\documentclass[reqno,12pt]{article}

\usepackage{amsmath,amssymb,amsthm,mathtools,bm,mathrsfs}
\usepackage{geometry}
\usepackage{microtype}
\usepackage{enumitem}
\usepackage{booktabs,array,aliascnt,fancyhdr}
\usepackage[numbers,sort&compress]{natbib}
\usepackage[hidelinks]{hyperref}
\usepackage[nameinlink,noabbrev]{cleveref}
\allowdisplaybreaks[3]
\numberwithin{equation}{section}
\setlist{itemsep=0.2em,topsep=0.35em}
\makeatletter
\renewcommand\section{\@startsection{section}{1}{\z@}%
  {-2.2ex \@plus -0.5ex \@minus -0.2ex}%
  {1.1ex \@plus 0.2ex}%
  {\normalfont\large\bfseries\raggedright}}
\makeatother

\newcommand{\RunningTitle}{Planar Viscous Contact Waves in the Whole Space}
\theoremstyle{plain}
\newtheorem{theorem}{Theorem}[section]
\newtheorem{lemma}{Lemma}[section]
\newaliascnt{proposition}{lemma}

\aliascntresetthe{proposition}
\newaliascnt{corollary}{lemma}
\newtheorem{corollary}[corollary]{Corollary}
\aliascntresetthe{corollary}
\theoremstyle{remark}
\newaliascnt{remark}{lemma}
\newtheorem{remark}[remark]{Remark}
\aliascntresetthe{remark}
\crefname{theorem}{Theorem}{Theorems}
\Crefname{theorem}{Theorem}{Theorems}
\crefname{lemma}{Lemma}{Lemmas}
\Crefname{lemma}{Lemma}{Lemmas}
\crefname{proposition}{Proposition}{Propositions}
\Crefname{proposition}{Proposition}{Propositions}
\crefname{corollary}{Corollary}{Corollaries}
\Crefname{corollary}{Corollary}{Corollaries}
\crefname{remark}{Remark}{Remarks}
\Crefname{remark}{Remark}{Remarks}
\crefname{equation}{equation}{equations}
\Crefname{equation}{Equation}{Equations}
\crefname{section}{Section}{Sections}
\Crefname{section}{Section}{Sections}
\crefname{appendix}{Appendix}{Appendices}
\Crefname{appendix}{Appendix}{Appendices}

\newcommand{\dd}{\,\mathrm d}
\newcommand{\R}{\mathbb R}
\newcommand{\pa}{\partial}
\newcommand{\Div}{\operatorname{div}}

\title{\bfseries Asymptotic Stability of Planar Viscous 
Contact Waves  
for the Three-Dimensional Full Compressible 
Navier--Stokes Equations in the Whole Space\thanks{This work is partially supported by the National Natural Science Foundation of
China (12471227, 12494541).}}
\author{Jing Li$^a$; Ziyin Liu$^b$
\thanks{a. Institute of Applied Mathematics \& State Key Laboratory of Mathematical Sciences, AMSS,\\ \& Hua Loo-Keng Key Laboratory of Mathematics,
Chinese Academy of Sciences,\\ Beijing, 100190 People's Republic of China, \\ and School of Mathematics and Computer Sciences, \\ Institute of Mathematics and Interdisciplinary Sciences,\\ Nanchang
University, Nanchang, 330031, People's Republic of China (ajingli@gmail.com).\\
b. Institute of Applied Mathematics, Academy of Mathematics and Systems Science,\\
Chinese Academy of Sciences, Beijing 100190, People's Republic of China\\
(a1909462594@gmail.com).}}
\date{}
\hypersetup{
  pdftitle={Asymptotic Stability of Planar Viscous Contact Waves  for the Three-Dimensional Full Compressible Navier-Stokes Equations in the Whole Space},
  pdfauthor={Jing Li; Ziyin Liu},
  pdfkeywords={compressible Navier-Stokes equations; viscous contact wave; planar wave; asymptotic stability; whole space}
}

\begin{document}
\maketitle
\begin{abstract}
The asymptotic stability of planar viscous contact waves is established
for the three-dimensional full compressible Navier--Stokes equations in
the whole space. The exact evolving planar background is constructed
in a companion one-dimensional paper and may arise from a large
integrable perturbation. The three-dimensional remainder has small
energy, prescribed second-derivative bounds for density and velocity,
and a prescribed first-derivative bound for temperature. Its initial
temperature lies in a fixed positive interval; its second derivatives
need only be finite. A unique global positive strong solution and
algebraic decay in the uniform norm are obtained. The wave-strength
threshold depends on lower-order planar bounds, while the allowed
remainder size may depend on the complete planar datum.
\end{abstract}



\section{Introduction}\label{sec:intro}

We study the stability of planar viscous contact waves for the full
compressible Navier--Stokes equations in $\R^3$:
\begin{equation}\label{eq:NSF}
\begin{cases}
 \rho_t+\Div(\rho u)=0,\\
 (\rho u)_t+\Div(\rho u\otimes u)+\nabla p=\Div\mathbb S,\\
 (\rho E)_t+\Div(\rho Eu+pu)=\kappa\Delta\theta+\Div(\mathbb Su),
\end{cases}
\end{equation}
where $t>0$ and $x=(x_1,x_2,x_3)\in\R^3$. Here $\rho$,
$u=(u_1,u_2,u_3)^\top$, $p$ and $\theta$ are the density, velocity,
pressure and absolute temperature, $E=c_v\theta+\frac12|u|^2$ is the
total energy per unit mass, and the viscous stress tensor is
\begin{equation}\label{eq:stress}
 \mathbb S=2\mu\mathbb D(u)+\lambda(\Div u)\mathbb I,\qquad
 \mathbb D(u)=\frac{\nabla u+(\nabla u)^\top}2,
\end{equation}
with $\mathbb I$ the $3\times3$ identity matrix. The viscosity
coefficients $\mu,\lambda$ and the heat conductivity $\kappa$ are
constants with
\begin{equation}\label{eq:cf-hyp}
 \mu>0,\qquad 2\mu+3\lambda\ge0,\qquad \kappa>0,
\end{equation}
and the gas is polytropic and ideal,
\begin{equation}\label{eq:state}
 p=R\rho\theta,\qquad c_v=\frac R{\gamma-1},
\end{equation}
where $R>0$ is the gas constant and $\gamma>1$ is the adiabatic exponent.
We write $c_p=\gamma c_v$, $\nu=2\mu+\lambda$ and $c=\frac\kappa{c_p}$;
by \eqref{eq:cf-hyp}, $\nu\ge\frac{4\mu}3>0$.

We study the Cauchy problem for \eqref{eq:NSF} with the initial data
\begin{equation}\label{eq:data}
 (\rho,u,\theta)(x,0)=(\rho_0,u_0,\theta_0)(x),\qquad x\in\R^3,
\end{equation}
and the far-field condition
\begin{equation}\label{eq:far}
 (\rho,u,\theta)(x,t)\to(\rho_\pm,u_\pm,\theta_\pm)\quad
 \text{as }x_1\to\pm\infty,
\end{equation}
uniformly in $x'=(x_2,x_3)\in\R^2$ for each $t\ge0$, where
$\rho_\pm,\theta_\pm>0$ and $u_\pm=u_{1\pm}e_1$, $e_1=(1,0,0)^\top$, are
constant states.

We assume that
\begin{equation}\label{eq:ct-st}
 u_{1-}=u_{1+}=:u_*,\qquad R\rho_-\theta_-=R\rho_+\theta_+,\qquad
 \theta_-\ne\theta_+ .
\end{equation}
The corresponding Euler contact discontinuity travels at speed $u_*$.
We put
\begin{equation}\label{e:end}
 p_*=R\rho_\pm\theta_\pm,\qquad \delta=|\theta_+-\theta_-|>0,\qquad \tau=1+t,
\end{equation}
and call $(R,\gamma,\mu,\lambda,\kappa,\rho_-,\theta_-)$ the fixed data.

One-dimensional contact-wave stability was established
in~\cite{HuangMatsumuraXin2006}. Huang and
Wang~\cite{HuangWang2016} treated large perturbations, while Liu, Wang
and Xu~\cite{LiuWangXu2026} obtained optimal decay for small Sobolev
perturbations with small corrected $L^2$ primitives.

In three dimensions, Matsumura and Nishida~\cite{MatsumuraNishida1980}
proved global existence and convergence for small $H^3$ perturbations
of a constant state. Huang and Li~\cite{HuangLi2018} allowed large
higher Sobolev norms under small initial energy. On
$\R\times\mathbb T^2$, planar rarefaction waves for the full system
and viscous shocks for the isentropic system are stable under small
$H^2$ perturbations~\cite{LiWangWang2018,WangWang2026}.
Duan, Huang, Li and Xu~\cite{DuanHuangLiXu2025} proved Eulerian
stability of planar viscous entropy waves on the same domain, with
optimal uniform decay. Their longitudinal means are small as well;
the precise assumptions are compared below.

The companion paper~\cite{LiLiuContact1D2026} supplies the exact
one-dimensional backgrounds, their optimal decay at a fixed Eulerian
position, and the coefficient and pressure estimates used here.
We prove their stability under three-dimensional perturbations in
$\R^3$, allowing a large planar perturbation and a continuum of
transverse frequencies approaching zero.

\subsection{The planar background}

The planar background used below is supplied by the companion
paper~\cite{LiLiuContact1D2026}. The Eulerian coordinate $x$, derivative
$\partial_x$ and material derivative $D_t$ in that paper are denoted
here by $x_1$, $\partial_1$ and $D_b$, respectively. Its Eulerian solution
$(\rho,u,\theta)$ and initial data $(\rho_0,u_0,\theta_0)$ are denoted
here by $(\rho_b,u_b,\theta_b)$ and $(\rho_0^b,u_0^b,\theta_0^b)$,
respectively, to distinguish them from the three-dimensional solution
and its initial data. We recall the normalization and data class.
Let $\widehat\Theta$ be the positive monotone self-similar profile with
end states $\theta_\pm$ for
$\Theta_t=\frac{\kappa p_*}{c_pR}\partial_y(\Theta_y/\Theta)$, and put

\begin{equation}\label{eq:prof}
 \Theta_c(y,t)=\widehat\Theta\Bigl(\frac y{\sqrt{1+t}}\Bigr),\qquad
 V_c=\frac{R\Theta_c}{p_*},\qquad
 U_c=u_*+\frac{\kappa}{c_p}\frac{\pa_y\Theta_c}{\Theta_c},
\end{equation}

Write $H_c=U_c-u_*$ and fix the Eulerian phase by

\begin{equation}\label{e:Xc}
 \bar X(y,t)=\int_0^yV_c(z,t)\dd z+\int_0^tU_c(0,r)\dd r,
 \qquad \bar Y(\cdot,t)=\bar X(\cdot,t)^{-1};
\end{equation}

The Eulerian contact wave is

\begin{equation}\label{eq:ct-w}
 \bar Q=(\bar\rho,\bar u,\bar\theta)
 =\left(\frac{p_*}{R\bar\theta},\,
 u_*+\frac{\kappa R}{c_pp_*}\partial_1\bar\theta,\,\bar\theta\right).
\end{equation}
Here $\bar\theta=\Theta_c\circ\bar Y$; equivalently,
$\bar Q=(V_c^{-1},U_c,\Theta_c)\circ\bar Y$. It satisfies

\begin{equation}\label{eq:ct-euler}
\begin{aligned}
 (\partial_t+u_*\partial_1)\bar\theta
 &=\frac{\kappa R}{c_pp_*}
   \bigl(\bar\theta\partial_1^2\bar\theta-(\partial_1\bar\theta)^2\bigr),\\
 \bar\theta(\pm\infty,t)&=\theta_\pm .
\end{aligned}
\end{equation}

For $e_0=(e_\rho,e_u,e_\theta)$, put
\begin{equation}\label{eq:bg-dat}
 W_0=(\rho_0^b,u_0^b,\theta_0^b)=\bar Q(\cdot,0)+e_0.
\end{equation}
The background $(\rho_b,u_b,\theta_b)$ is the solution of
\begin{equation}\label{eq:bg-eqs}
\begin{cases}
 (\rho_b)_t+u_b\partial_1\rho_b+\rho_b\partial_1u_b=0,\\
 \rho_b\bigl((u_b)_t+u_b\partial_1u_b\bigr)+\partial_1p_b
       =\nu\partial_1^2u_b,\\
 c_v\rho_b\bigl((\theta_b)_t+u_b\partial_1\theta_b\bigr)
       +p_b\partial_1u_b
       =\kappa\partial_1^2\theta_b+\nu(\partial_1u_b)^2,
\end{cases}
\qquad p_b=R\rho_b\theta_b,
\end{equation}
with the initial value $W_0$, where $\partial_1=\partial_{x_1}$.

Thus $(\rho_b,u_be_1,\theta_b)$ is a planar solution of \eqref{eq:NSF}.
For the initial mass labels, put

\begin{equation}\label{e:Y0}
 Y_0(x_1)=\int_0^{x_1}\rho_0^b(z)\dd z,\qquad X_0=Y_0^{-1},
\end{equation}

The initial conserved excess is defined directly from these data:
\begin{equation}\label{e:mass-initial}
\begin{gathered}
 \mathcal U_0=\Bigl(\frac1{\rho_0^b},u_0^b-u_*,
     c_v\theta_0^b+\frac{(u_0^b-u_*)^2}2\Bigr)\circ X_0,\\
 \widehat{\mathcal U}=(V_c,H_c,c_v\Theta_c),\qquad
 Q=(Q_v,Q_{\rm mom},Q_E)=\int_\R
       [\mathcal U_0-\widehat{\mathcal U}(y,0)]\dd y .
\end{gathered}
\end{equation}
For $A\ge1$, assume $e_\rho\in H^3(\R)$,
$(e_u,e_\theta)\in H^2(\R)$ and $e_0\in L^1(\R)$, and put
$\Phi(z)=z-\log z-1$. The quantitative data condition is

\begin{equation}\label{e:dat}
\begin{gathered}
 \Vert e_\rho\Vert_{H^1(\R)}+\Vert e_u\Vert_{L^2(\R)}^2+|Q|
 +\int_\R\Phi\left(\frac{\theta_0^b(x_1)}
                             {\bar\theta(x_1,0)}\right)\dd x_1\le A,\\
 \rho_0^b\ge A^{-1},\qquad \theta_0^b\ge A^{-1}.
\end{gathered}
\end{equation}

\begin{lemma}[Imported planar background]\label{l:background}
For the planar data just specified, under \eqref{eq:cf-hyp},
\eqref{eq:ct-st} and \eqref{e:dat}, the
one-dimensional stability theorem~\cite[Theorem~1.1]{LiLiuContact1D2026}
applies with viscosity $\nu=2\mu+\lambda\ge4\mu/3>0$ and heat
conductivity $\kappa>0$. It gives a threshold $\delta_0(A)>0$,
depending only on $A$ and the fixed data and
independent of $u_*$, such that $0<\delta\le\delta_0(A)$ yields a unique
global background solving \eqref{eq:bg-eqs} with initial value
\eqref{eq:bg-dat}. For every $S<\infty$ it satisfies

\begin{equation}\label{e:brg}
\begin{aligned}
 &\rho_b-\bar\rho\in C([0,S];H^3(\R)),\qquad
 (u_b-\bar u,\theta_b-\bar\theta)\in C([0,S];H^2(\R)),\\
 &\partial_1(u_b,\theta_b)\in L^2(0,S;H^2(\R)),\\
 &(\rho_b)_t\in C([0,S];H^1(\R))\cap L^2(0,S;H^2(\R)),\\
 &((u_b)_t,(\theta_b)_t)\in C([0,S];L^2(\R))\cap L^2(0,S;H^1(\R)),
\end{aligned}
\end{equation}
and $\rho_b,\theta_b$ lie between two positive constants on
$\R\times[0,\infty)$.
Moreover, for $t\ge0$, with $\bar Q$ of \eqref{eq:ct-w},
\begin{equation}\label{e:r1}
 \Vert(\rho_b,u_b,\theta_b)(t)-\bar Q(t)\Vert_{L^\infty(\R)}
 +\Vert p_b(t)-p_*\Vert_{L^\infty(\R)}
 \le C_b(1+t)^{-1/2}.
\end{equation}

The finite higher Sobolev norms and the full
$L^1$ size of $e_0$ may affect these background constants and entrance
times, but impose no further restriction on $\delta_0(A)$.
\end{lemma}

\subsection{The perturbation}

For a solution of \eqref{eq:NSF}, define its perturbation of the
exact planar background by
\begin{equation}\label{eq:pert}
 q=(\phi,\psi,\zeta)=(\rho-\rho_b,\,u-u_be_1,\,\theta-\theta_b),
 \qquad \Pi=p-p_b.
\end{equation} The initial data are
\begin{equation}\label{e:3dat}
 (\rho,u,\theta)(x,0)=(\rho_0^b,u_0^be_1,\theta_0^b)(x_1)+q_0(x),
 \quad q_0\in H^2(\R^3),\quad
 \epsilon=\Vert q_0\Vert_{H^2(\R^3)}.
\end{equation}

\subsection{Notation and constants}

Sobolev spaces are inhomogeneous, norms of vector-valued functions use
the Euclidean norm of the vector, and $\Vert\cdot\Vert$ is the $L^2$
norm. The domains are $\R_y$ for functions of the mass coordinate,
$\R_{x_1}$ for Eulerian backgrounds and $\R^3$ for three-dimensional
perturbations; a norm without a domain in
Sections~\ref{sec:pre}--\ref{sec:proof} is over $\R^3$. We use
\begin{equation}\label{e:der}
 \partial_1=\partial_{x_1},\qquad
 D_b=\partial_t+u_b\partial_1,\qquad
 D_u=\partial_t+u\cdot\nabla,
\end{equation}
the transverse operators $\nabla'=(\partial_2,\partial_3)$,
$\Delta'=\partial_2^2+\partial_3^2$ and $\Div'$, and $\tau_s=1+s$ for a
time $s$. 

Constants denoted by $C_M$ depend only on an auxiliary bound $M$
for the planar datum and the fixed data,
together with fixed parameters displayed in the statement; constants
denoted by $C_b$ may also depend on the background, and $C_{b,S}$ also on
a fixed starting time $S$. They are independent of the final time and of
the three-dimensional datum $q_0$. The notation $\omega_b(s)$ denotes a
nonnegative quantity depending on the background with
$\omega_b(s)\to0$ as $s\to\infty$; it may change between estimates.
Additional dependence on a multiplier bound
or a fixed coefficient is displayed. Gaussian exponents are written $c_G$ or $c_j$. The order of
choosing constants is given in the proof of Theorem~\ref{thm:main} in
Section~\ref{sec:proof}.

\subsection{Main result}

For an interval $I\subset[0,\infty)$, a solution of \eqref{eq:NSF} on $\R^3\times I$
is called strong if
\begin{equation}\label{eq:class}
 q\in C(I;H^2(\R^3)),\qquad
 (\psi,\zeta)\in L^2(I;H^3(\R^3)),\qquad
 q_t\in L^2(I;H^1(\R^3)),
\end{equation}
and $\rho,\theta$ are bounded below by positive constants on $\R^3\times I$.
Local existence and continuation are given in Lemma~\ref{l:local}
of Appendix~\ref{app:local}. Our main result is as follows.

\begin{theorem}[Three-dimensional stability]\label{thm:main}
Let $A\ge1$ and let the planar datum satisfy the assumptions of
Lemma~\ref{l:background}. Let $(\rho_b,u_b,\theta_b)$ be the background
supplied by that lemma.
After decreasing $\delta_0(A)$ using only $A$ and the fixed data,
the following holds for $0<\delta\le\delta_0(A)$.
Fix $M_\rho,M_u,M_\theta\ge1$ and
$0<m_\theta\le\Theta_\theta<\infty$. There are $\epsilon_0>0$ and
$C$, depending only on $A$, the fixed data, $\theta_+$,
$e_0$ and these five bounds, such that for every $q_0=(\phi_0,\psi_0,\zeta_0)$
with \eqref{e:3dat} and
\begin{equation}\label{e:3small}
\begin{gathered}
 \Vert\nabla^2\phi_0\Vert_{L^2}\le M_\rho,\qquad
 \Vert\nabla^2\psi_0\Vert_{L^2}\le M_u,\qquad
 \Vert\nabla\zeta_0\Vert_{L^2}\le M_\theta,\\
 m_\theta\le\theta_0^b(x_1)+\zeta_0(x)\le\Theta_\theta,
 \qquad \Vert q_0\Vert_{L^2}\le\epsilon_0,
\end{gathered}
\end{equation}
\eqref{eq:NSF} has a unique global strong solution with the initial
data \eqref{e:3dat}.
For every $T_e<\infty$, $\phi_t\in C([0,T_e];H^1(\R^3))$; at each
$t\ge0$, \eqref{eq:far} holds uniformly in $x'$; and, for $t\ge0$, with
$\bar Q$ of \eqref{eq:ct-w},
\begin{equation}\label{e:r3}
 \Vert(\phi,\zeta)(t)\Vert_{L^\infty(\R^3)}
 +\Vert\psi(t)\Vert_{L^\infty(\R^3)}
 +\Vert\Pi(t)\Vert_{L^\infty(\R^3)}
 \le C\epsilon(1+t)^{-1/20},
\end{equation}
\begin{equation}\label{e:rt}
\begin{aligned}
 &\Vert(\rho-\bar\rho,\theta-\bar\theta)(t)\Vert_{L^\infty(\R^3)}
 +\Vert u(t)-\bar u(t)e_1\Vert_{L^\infty(\R^3)}
 +\Vert p(t)-p_*\Vert_{L^\infty(\R^3)}\\
 &\hspace{15mm}\le C\epsilon(1+t)^{-1/20}+C_b(1+t)^{-1/2}.
\end{aligned}
\end{equation}
The three derivative bounds need not be small. Neither $\epsilon_0$
nor $C$ depends on $\Vert\nabla^2\zeta_0\Vert_{L^2}$; this norm is
required only to be finite. The density is initially positive by
\eqref{e:3small} and interpolation. In \eqref{e:r3} and \eqref{e:rt},
$\epsilon=\Vert q_0\Vert_{H^2}$ is the finite amplitude of
\eqref{e:3dat}, not the norm on which smallness is imposed.
\end{theorem}

\begin{remark}[Comparison with periodic planar-wave results]\label{r:whole-space}
The planar rarefaction, shock and entropy-wave results
in~\cite{LiWangWang2018,WangWang2026,DuanHuangLiXu2025} concern
$\R\times\mathbb T^2$, with distinct longitudinal end states and
periodic transverse directions.
On the unit torus, the transverse frequencies satisfy
$k\in2\pi\mathbb Z^2$, so $|k|\ge2\pi$ for every nonzero mode. Hence
\[
 \Vert f-\langle f\rangle\Vert_{L^2(\R\times\mathbb T^2)}
 \le C\Vert\nabla'f\Vert_{L^2(\R\times\mathbb T^2)},
 \qquad \langle f\rangle(x_1)=\int_{\mathbb T^2}f(x_1,x')\dd x'.
\]
Thus the one-dimensional mean $k=0$ is the only mode not controlled
by transverse dissipation; the mean and the other modes remain
nonlinearly coupled. In $\R^3$, $k\in\R^2$ can approach zero
continuously, so the transverse dissipation $|k|^2$ has no positive
lower bound for $k\ne0$. 
The whole-space problem therefore requires control of a continuum
of slowly damped transverse modes and their interaction with the
planar background; see Section~\ref{s:freq}.
\end{remark}

\begin{remark}[Perturbation size and decay]\label{r:DHLX}
The nonzero-mass theorem of Duan, Huang, Li and
Xu~\cite[Theorem~2.1]{DuanHuangLiXu2025} assumes small $H^3$
perturbations, small $L^1$ conserved perturbations, and small $H^3(\R)$
primitives of their transverse averages after the excess-mass
correction. Hence the longitudinal mean is small too. Their uniform
decay rate is $(1+t)^{-1/2}$, and the nonzero transverse modes decay
exponentially. Here the planar background may be large, and $q_0$
requires no $L^1$ or primitive condition. The planar rate remains
$(1+t)^{-1/2}$ by \eqref{e:r1}; the whole-space remainder has the
separate rate \eqref{e:r3}.
\end{remark}

\begin{remark}[Large evolving planar backgrounds]\label{r:thresholds}
He, Huang and Wang~\cite[Theorem~1.2]{HeHuangWang2020} establish
conditional stability of a given global smooth reference solution
with a single constant far-field state, assuming uniform regularity
bounds and initial deviations in $L^1\cap H^s$, $s>3/2$.
Here the reference solution is constructed from large planar data
connecting distinct end states. It lies outside that integrable
whole-space setting: its initial deviation from the contact profile
satisfies
\[
 e_0\not\equiv0\quad\Longrightarrow\quad
 \int_{\R^3}|e_0(x_1)|^2\dd x=\infty.
\]
Theorem~\ref{thm:main} proves its stability under three-dimensional
perturbations satisfying \eqref{e:3small}; smallness is imposed on
the remainder relative to this exact background.

For fixed $A$ in \eqref{e:dat}, the higher planar Sobolev norms and
the full $L^1$ norm may be arbitrarily large, while $\delta_0(A)$
remains unchanged. The contact profile has a fixed initial width,
and $\delta_0(A)$ is independent of $u_*$. The allowed remainder
size $\epsilon_0$ may depend on the complete planar datum and the
bounds in \eqref{e:3small}.
\end{remark}

\begin{remark}[Small energy and componentwise bounds]\label{r:small-energy}
On fixed positive density and temperature intervals, the initial
relative energy defined in \eqref{e:rel3} is equivalent to
$\Vert q_0\Vert_{L^2}^2$. The small-energy theory of Huang and
Li~\cite{HuangLi2018} concerns perturbations of a constant state;
here the reference state is a positive, nonconstant planar solution.
Condition~\eqref{e:3small} retains quantitative second-derivative
bounds for the density and velocity perturbations, but only a
first-derivative bound for the temperature perturbation. For suitable
fixed $M_\theta,m_\theta,\Theta_\theta$, it allows the family
\[
 \phi_0^{(N)}=\psi_0^{(N)}=0,\qquad
 \zeta_0^{(N)}=N^{-1}\chi(x)\sin(Nx_1),
 \qquad 0\ne\chi\in C_c^\infty(\R^3),
\]
for which, as $N\to\infty$,
\[
 \Vert q_0^{(N)}\Vert_{L^2}=O(N^{-1}),\qquad
 \Vert\nabla\zeta_0^{(N)}\Vert_{L^2}=O(1),\qquad
 \Vert\nabla^2\zeta_0^{(N)}\Vert_{L^2}\asymp N.
\]
Choose the positive temperature interval to contain the range of
$\theta_0^b$ in its interior. Then all sufficiently large $N$ satisfy
\eqref{e:3small} with the same $L^2$ threshold.
Lemma~\ref{l:thermal-layer} gives a uniform $H^2$ bound at time one,
independently of the growing initial temperature Hessian.
\end{remark}

\subsection{Difficulties and method}

The background derivatives decay too slowly for direct energy estimates.
The temperature equation \eqref{eq:phzeta} contains
$-c_v\rho\psi_1\partial_1\theta_b$, while \eqref{e:E1} gives only
\[
 \int_S^T\Vert\psi_1\partial_1\theta_b\Vert^2\dd t
 \le\bigl(C_0\delta^2+\omega_b(S)\bigr)
 \sup_{S\le t\le T}\Vert\psi_1(t)\Vert^2
 \log\frac{1+T}{1+S}.
\]
Thus a direct estimate of $\Vert\zeta\Vert^2$ loses a logarithm.
The acoustic masses need not be small either, so choosing small
$\delta$ does not absorb their strain.

\paragraph{The weighted perturbation energy.}
The energy \eqref{eq:E} weights the undifferentiated temperature by
$\tau^{-1}$. The time weight $((1+s)/(1+t))^{1/128}$ supplies damping
of order $\tau^{-1}$ to absorb $C\delta\tau^{-1}$ coefficients;
\eqref{e:Ns} combines this energy and its dissipation.

The density is recovered through the entropy decomposition of
Lemma~\ref{l:ent}:
\[
 \chi=\frac{\zeta-\Pi/(c_p\rho)}{\theta_b}=\eta+f,\qquad
 \phi=\frac{\rho_b\Pi}{\gamma p_b}-\rho(\eta+f).
\]
Here $\eta$ carries the initial entropy and $f$ its forced response.
The velocity $V=\psi-(u_b-u_*)\eta e_1$ removes the acceleration term
containing $\eta$. Lemma~\ref{l:mr} constructs a positive energy
that also controls this entropy history. Its remaining terms are the
strain, thermal and contact works in \eqref{e:works}.

The acoustic strain has size $C_M\tau^{-1}$ without a small factor.
We use the correction of Lemma~\ref{l:cor} at low transverse
frequencies, a decreasing Fourier multiplier at intermediate
frequencies, and weighted crossing identities at high frequencies.
Lemmas~\ref{l:low}--\ref{l:tf} control these regions and their overlaps;
late starting times make their time-boundary terms small.
Lemma~\ref{l:tw} treats the thermal and nonlinear works.

\paragraph{The contact work and the pressure decomposition.}
The contact work contains the background acceleration, and
$\rho_bD_bu_b=\nu\partial_1^2u_b-\partial_1p_b$.
For the decomposition $P_b-p_*=g_1+g_2$ supplied by
\cite{LiLiuContact1D2026} and recalled in Lemma~\ref{l:prs},
Lemma~\ref{l:cw} absorbs the finite weighted budget of $g_1$ after a
late start. Integration by parts makes the Gaussian observation of
$g_2$ multiply the accumulated energy; its logarithmic budget has a
small slope and is treated by Gronwall's inequality. Choosing $\delta$
first and the starting time afterwards gives the small polynomial
growth in Lemma~\ref{l:cl}, with a wave-strength threshold independent
of the higher planar norms for fixed $A$.

\paragraph{Decay and continuation.}
Despite this small energy growth, Lemma~\ref{l:age} controls
$\nabla(\Pi,\psi)$ with a time weight. Adjoint smoothing
(Lemma~\ref{l:adj}) controls the low-frequency temperature, and local
damping (Lemma~\ref{l:damp}) recovers the uniformly local and uniform
norms. Lemma~\ref{l:dec} then gives $(1+t)^{-1/20}$ decay and small
$\Vert(\Pi,\psi)\Vert_{L^3}$, closing continuation.

For the componentwise bounds \eqref{e:3small}, relative energy about
the exact planar background keeps the $L^2$ norm small.
Lemma~\ref{l:thermal-layer} controls density and velocity in $H^2$
and temperature in $H^1$; a time-weighted thermal estimate gives a
uniform full $H^2$ bound at time one, independently of the initial
temperature Hessian. Lemma~\ref{l:low-finite} propagates this bound
on fixed finite intervals. The high-frequency damping of
Lemma~\ref{l:low-damp} then reaches the polynomial late-time
threshold in Lemma~\ref{l:late-restart}.

Section~\ref{s:bg-interface} records the background estimates
from \cite{LiLiuContact1D2026}. Sections~\ref{sec:pre}
and~\ref{sec:en} give the perturbation estimates;
Sections~\ref{s:dec} and~\ref{sec:proof} prove decay and complete
the proof of Theorem~\ref{thm:main}. Appendix~\ref{app:local}
contains the local theory and the estimates for small energy
with prescribed derivative bounds.
\section{Background estimates from the companion paper}\label{s:bg-interface}
We restate the background estimates of~\cite{LiLiuContact1D2026},
including the coefficient and pressure bounds from its Appendices~B
and~D. Norms are over $\R_y$ in mass coordinates and $\R_{x_1}$ in
Eulerian coordinates. Constants and tails follow the conventions above.
For each fixed datum choose a finite auxiliary $M\ge A$ such that

\begin{equation}\label{e:dat-aux}
\begin{gathered}
 \Vert e_\rho\Vert_{H^3(\R)}
 +\Vert(e_u,e_\theta)\Vert_{H^2(\R)}
 +\Vert e_0\Vert_{L^1(\R)}\le M,\\
 \rho_0^b\ge M^{-1},\qquad \theta_0^b\ge M^{-1}.
\end{gathered}
\end{equation}

\paragraph{Uniformity of the imported estimates.}
Theorem~1.1 and Proposition~D.20 of~\cite{LiLiuContact1D2026} show that,
after decreasing $\delta_0(A)$ using only $A$ and the fixed data, all
estimates below hold for $0<\delta\le\delta_0(A)$ after a sufficiently
late background-dependent entrance time $S_0$. The acoustic masses and
Cole--Hopf factors are controlled by $A$. The persistent contact and
absorption coefficients are controlled by $A$ and the fixed data; the
constant $C_P$ in the pressure estimate depends only on the fixed data.
Higher background derivatives, finite source budgets and errors that
vanish with the entrance time retain their $C_b$ dependence. In particular,
this assertion does not replace every $C_M$ or $C_b$ by $C_A$.
The Gaussian exponents are fixed before the datum; the threshold may
also depend on those fixed exponents. The pressure constant $C_P$ is
independent of them, the acoustic masses and the background.

\subsection{Mass coordinates and the diffusion waves}\label{s:bg-coordinates}
Let $(V_b,U_b,\Theta_b)$ be the mass-coordinate representation of the
background of Lemma~\ref{l:background}. Its initial value is
$(1/\rho_0^b,u_0^b,\theta_0^b)\circ X_0$. Put
$P_b=R\Theta_b/V_b$. Then

\begin{equation}\label{eq:m-map}
 X(y,t)=\int_0^yV_b(z,t)\dd z+\int_0^tU_b(0,r)\dd r,
 \qquad Y(\cdot,t)=X(\cdot,t)^{-1}
\end{equation}
satisfies $X_y=V_b$, $X_t=U_b$ and $X(\cdot,0)=X_0$, and
\begin{equation}\label{eq:bg-eul}
 Q_b=(\rho_b,u_b,\theta_b)=\left(\frac1{V_b},U_b,\Theta_b\right)\circ Y,
 \qquad p_b=P_b\circ Y.
\end{equation}

We use $\nabla_m=(\partial_y,\partial_2,\partial_3)$. For a function $F(x_1,x',t)$ we write
$F\circ X=F(X(y,t),x',t)$; in particular,
\begin{equation}\label{e:chain}
 \partial_y(F\circ X)=V_b(\partial_1F)\circ X,\qquad
 \partial_t(F\circ X)=(D_bF)\circ X,\qquad
 \dd x_1=V_b\dd y.
\end{equation}
Mass-coordinate functions are composed with $Y$ when used in Eulerian
coordinates. The maps $X$ and $\bar X$ retain their initial
normalization at every later starting time.
For three-dimensional functions, this is a pullback by the
one-dimensional background map.

Put
\begin{equation}\label{e:Um}
 b_m=U_b-u_*,\qquad
 Z_m=(\varphi_m,\upsilon_m,\vartheta_m)
    =(V_b-V_c,\,b_m-H_c,\,\Theta_b-\Theta_c).
\end{equation}

The contact profile satisfies the following Gaussian estimates from
\cite[Section~2.1]{LiLiuContact1D2026}:

If
$\delta$ is sufficiently small, then for $\pm y\ge0$ and $t\ge0$ (here
and below, either all upper or all lower signs are taken)
\begin{equation}\label{eq:ct-gau}
 (1+t)|\pa_y^2\Theta_c|+(1+t)^{1/2}|\pa_y\Theta_c|+|\Theta_c-\theta_\pm|
 \le C\delta e^{-\frac{c_G y^2}{1+t}},
\end{equation}
where $C,c_G>0$ are independent of $y$, $t$ and $\delta$.

\begin{equation}\label{eq:pr-gau}
\begin{aligned}
 |\pa_y^k(V_c,\Theta_c)|
 &\le C\delta(1+t)^{-\frac k2}e^{-\frac{c_G y^2}{1+t}} &&(k\ge1),\\
 |\pa_y^k(U_c-u_*)|
 &\le C\delta(1+t)^{-\frac{k+1}2}e^{-\frac{c_G y^2}{1+t}} &&(k\ge0),
\end{aligned}
\end{equation}
where $C$ depends on $k$, and $C,c_G$ are independent of $y$, $t$ and
$\delta$.

The mass-coordinate regularity class is, for every $0<T_e<\infty$
and $t_0\in(0,T_e)$,

\begin{equation}\label{e:bgcl}
\begin{aligned}
 &\varphi_m\in C([0,T_e];H^3),\qquad
 (\upsilon_m,\vartheta_m)\in C([0,T_e];H^2)\cap C([t_0,T_e];H^3),\\
 &\partial_y(\upsilon_m,\vartheta_m)
       \in L^2(0,T_e;H^2)\cap L^2(t_0,T_e;H^3),\\
 &(V_b)_t\in C([0,T_e];H^1)\cap L^2(0,T_e;H^2),\qquad
 ((U_b)_t,(\Theta_b)_t)\in C([0,T_e];L^2)\cap L^2(0,T_e;H^1).
\end{aligned}
\end{equation}

\begin{lemma}[Mass-coordinate background regularity]\label{l:bg}
The background of Lemma~\ref{l:background} satisfies \eqref{e:bgcl} on
every finite time interval~\cite[Lemma~2.2]{LiLiuContact1D2026}.
There are $c_b,C_b>0$ such that
$c_b\le V_b,\Theta_b\le C_b$ on $\R\times[0,\infty)$.
The maps \eqref{eq:m-map} satisfy $X_y=V_b$, $X_t=U_b$ and
$X(\cdot,0)=X_0$; \eqref{eq:bg-eul} agrees with the Eulerian background
in Lemma~\ref{l:background}.
\end{lemma}

For the centered energy and the conserved phase, put
\[
 E_m=c_v\Theta_b+\frac{b_m^2}2,\qquad \mathcal U=(V_b,b_m,E_m).
\]
The mass-coordinate conserved excess equals the $Q$ in
\eqref{e:mass-initial}; define

\begin{equation}\label{e:mass}
\begin{gathered}
 v_\pm=\frac1{\rho_\pm},\qquad \Delta_v=v_+-v_-,\\
 c_\pm=\sqrt{\frac{\gamma p_*}{v_\pm}},\qquad s_\pm=\pm c_\pm,\qquad
 r_0=\Bigl(1,0,\frac{p_*}{\gamma-1}\Bigr),\qquad r_\sigma=(-1,s_\sigma,p_*),\\
 m_0=\frac{\gamma-1}\gamma\Bigl(Q_v+\frac{Q_E}{p_*}\Bigr),\qquad h_*=\frac{m_0}{\Delta_v},\\
 m_-=\frac{c_+(m_0-Q_v)-Q_{\rm mom}}{c_-+c_+},\qquad
 m_+=\frac{c_-(m_0-Q_v)+Q_{\rm mom}}{c_-+c_+}.
\end{gathered}
\end{equation}

The companion conservation result gives
$Q=h_*\Delta_vr_0+m_-r_-+m_+r_+$ and
$|m_-|+|m_+|+\delta|h_*|\le C_A$.
For the acoustic diffusion coefficients put

\begin{equation}\label{e:ac}
 d(n)=\frac{\nu+(\gamma-1)c}{2n},\qquad
 s_\sigma(n)=\sigma\sqrt{\frac{\gamma p_*}n},\qquad
 k(n)=\frac{\gamma+1}{4n},\qquad
 \alpha=\frac{\gamma+1}{2(\nu+(\gamma-1)c)},
\end{equation}
so that $\alpha d(n)=k(n)$, and let $d_\sigma=d(v_\sigma)$, $k_\sigma=k(v_\sigma)$;

For $\sigma=\pm1$, the diffusion waves and the remainder are

\begin{equation}\label{e:CH}
\begin{gathered}
 F_d(\xi)=\frac1{\sqrt{4\pi d}}\int_{-\infty}^\xi e^{-\frac{z^2}{4d}}\dd z,\qquad
 p_\sigma=1+\bigl(e^{-\alpha s_\sigma m_\sigma}-1\bigr)
   F_{d_\sigma}\Bigl(\frac{y-s_\sigma t}{\sqrt\tau}\Bigr),\\
 A_\sigma=-\frac1\alpha\log p_\sigma,\qquad a_\sigma=\partial_yA_\sigma,\qquad
 \vartheta_\sigma=\frac{a_\sigma}{s_\sigma},\qquad
 R_U=\mathcal U-\widehat{\mathcal U}(y+h_*,t)
   -\sum_{\sigma=\pm1}\vartheta_\sigma r_\sigma .
\end{gathered}
\end{equation}

The next two results are~\cite[Lemmas~3.2 and~4.12]{LiLiuContact1D2026}.
\begin{lemma}[Diffusion-wave bounds]\label{l:pr}
For the masses in \eqref{e:mass}, $p_\sigma$ takes values between $1$
and $e^{-\alpha s_\sigma m_\sigma}$ and $A_\sigma$ between $0$ and
$s_\sigma m_\sigma$. Moreover,
\[
\begin{gathered}
 A_\sigma(-\infty,t)=0,\quad A_\sigma(+\infty,t)=s_\sigma m_\sigma,
 \qquad \int_\R\vartheta_\sigma\dd y=m_\sigma,\\
 (\vartheta_\sigma)_t+s_\sigma(\vartheta_\sigma)_y
   +2k_\sigma s_\sigma\vartheta_\sigma(\vartheta_\sigma)_y
   =d_\sigma(\vartheta_\sigma)_{yy},\\
 |\partial_y^ja_\sigma|+|\partial_y^j\vartheta_\sigma|
 \le C_{A,j}\tau^{-\frac{j+1}2}
      e^{-\frac{(y-s_\sigma t)^2}{8d_\sigma\tau}}
 \qquad(j\ge0).
\end{gathered}
\]
Here $C_{A,j}$ depends only on $A$, $j$ and the fixed data.
\end{lemma}
\begin{lemma}[Diffusion-wave remainder]\label{l:asy}
\[
 \lim_{t\to\infty}\sqrt\tau\,\Vert R_U(t)\Vert_{L^\infty}=0,
 \qquad \sup_{t\ge0}\sqrt\tau\,\Vert R_U(t)\Vert_{L^\infty}\le C_b.
\]
\end{lemma}

\subsection{Characteristic budgets}
The entropy history uses the characteristic budgets
of~\cite[Lemma~2.5]{LiLiuContact1D2026}, with the quantities

\begin{equation}\label{e:char}
\begin{gathered}
 m=V_b\Bigl(\frac{P_b}{p_*}\Bigr)^{1/\gamma},\qquad e_m=m-V_c,\qquad
 H_T=c\partial_y\log\Theta_b,\\
 J_a=\frac{2\sqrt{\gamma p_*}\sqrt m}{\gamma-1}
 \Bigl[\Bigl(\frac{P_b}{p_*}\Bigr)^{\frac{\gamma-1}{2\gamma}}-1\Bigr],\qquad
 q_\pm=m^{-1/4}(\upsilon_m\pm J_a),\\
 \Gamma=\bigl((b_m)_y,(P_b)_y,H_T-H_c,(e_m)_y,(q_+)_y,(q_-)_y,
          m_yq_+,m_yq_-,H_Tq_+,H_Tq_-\bigr),\\
 c_\sigma(z,m)=\sqrt{\gamma p_*}\,m^{-1/2}
    \Bigl(1+\frac{\sigma(\gamma-1)m^{-1/4}z}{4\sqrt{\gamma p_*}}\Bigr)^{\frac{\gamma+1}{\gamma-1}},
 \qquad
 E_\sigma(q,m)=\int_0^q\frac{z\dd z}{c_\sigma(z,m)}.
\end{gathered}
\end{equation}

\begin{lemma}[Imported characteristic budgets]\label{l:dep}
For some $S_0\ge1$ there are $J_\pm,S_\pm$ on
$\R\times[S_0,\infty)$ such that

\begin{equation}\label{e:hcur}
\begin{gathered}
 \partial_t[E_\sigma(q_\sigma,m)]
       +\sigma\partial_y\frac{q_\sigma^2}{2}=(J_\sigma)_y+S_\sigma,\qquad
 C_*^{-1}q_\sigma^2\le E_\sigma(q_\sigma,m)\le C_*q_\sigma^2,\\
 C_*^{-1}\bigl((P_b-p_*)^2+\upsilon_m^2\bigr)\le q_+^2+q_-^2
 \le C_*\bigl((P_b-p_*)^2+\upsilon_m^2\bigr),\\
 |J_\sigma|+|\partial_y[E_\sigma(q_\sigma,m)]|\le C_*|q_\sigma||\Gamma|,
 \qquad |(e_m)_y|\le C_*|\Gamma|,
\end{gathered}
\end{equation}
Here $C_*$ depends only on the fixed data. The identities hold in distributions;
$E_\sigma(q_\sigma,m)\in C([S_0,T];L^1)$ and
$J_\sigma,S_\sigma\in L^1(\R\times(S_0,T))$ for finite $T>S_0$. The tail
\[
 \sup_{t\ge s}\Vert(e_m,q_+,q_-)(t)\Vert_{H^1}^2
 +\int_s^\infty\bigl(\Vert\Gamma\Vert_{L^\infty}^2
       +\Vert(S_+,S_-)\Vert_{L^1}\bigr)\dd t
\]
is finite for $s=S_0$ and tends to zero as $s\to\infty$.
\end{lemma}

\subsection{Coefficient estimates}
We retain the mass-coordinate notation needed for the acoustic and
contact works. The next two statements are the coefficient results of
\cite[Lemma~B.1 and Corollary~B.2]{LiLiuContact1D2026}.

\begin{lemma}\label{l:mcoef}
Put $Z=(V_b,b_m,\Theta_b)$,
$(V_h,H_h,\Theta_h)=(V_c,H_c,\Theta_c)(y+h_*,t)$ and
\begin{equation}\label{e:cmp}
\begin{gathered}
 (v_a,b_a,E_a)=\widehat{\mathcal U}(y+h_*,t)
   +\sum_{\sigma=\pm1}\vartheta_\sigma r_\sigma,
 \qquad T_a=\frac1{c_v}\Bigl(E_a-\frac{b_a^2}2\Bigr),\qquad
 p_a=\frac{RT_a}{v_a}.
\end{gathered}
\end{equation}
For $t\ge S_0$, $Z$ and $(v_a,b_a,T_a)$ take values in
$\mathcal K=[\frac1{2\rho_-},\frac2{\rho_-}]\times[-1,1]
\times[\frac{\theta_-}2,2\theta_-]$.

The following derivative bounds hold:

\begin{equation}\label{e:munw}
\begin{gathered}
 \lim_{t\to\infty}\Vert Z_y\Vert_{H^2}=0,\qquad
 \sup_{t\ge S_0}\sqrt\tau\bigl(\Vert b_m\Vert_{L^\infty}
       +\Vert Z_y\Vert_{L^\infty}\bigr)<\infty,\\
 \int_{S_0}^\infty\bigl(\Vert Z_{yy}\Vert_{H^1}^2
   +\Vert\partial_y^4(b_m,\Theta_b)\Vert_{L^2}^2+\Vert(b_m)_y\Vert_{L^2}^2
   +\Vert(P_b)_y\Vert_{H^2}^2+\Vert Z_y\Vert_{L^4}^4\bigr)\dd t<\infty.
\end{gathered}
\end{equation}

Finally, fix $0<c_0<\frac18\min\{c_-,c_+\}$, let $\sigma=\pm1$, and put
$\Omega_\sigma(t)=\{y:\sigma(y+h_*)\ge c_0\tau\}$,
$Z_{\sigma,*}=(v_\sigma,0,\theta_\sigma)$,
$Z_{\sigma,a}=Z_{\sigma,*}+\vartheta_\sigma(-1,s_\sigma,\frac{p_*}{c_v})
-(0,0,\frac{s_\sigma^2\vartheta_\sigma^2}{2c_v})$ and
\[
 \mathsf B=\begin{pmatrix}0&-1&0\\-\frac{P_b}{V_b}&0&\frac R{V_b}\\
   0&\frac{P_b}{c_v}&0\end{pmatrix},\quad
 \mathsf V=\operatorname{diag}\Bigl(0,\frac\nu{V_b},\frac\kappa{c_vV_b}\Bigr),\quad
 \mathsf Q=\Bigl(0,\ \nu(b_m)_y\partial_y\frac1{V_b},\ 
   \frac{\kappa(\Theta_b)_y}{c_v}\partial_y\frac1{V_b}
   +\frac{\nu(b_m)_y^2}{c_vV_b}\Bigr)^{\!\top}.
\]
For $\varsigma\in\R$ and smooth functions $a,f$ on $\mathcal K$, with $f_Z$
the row of partial derivatives of $f$, put
$L=[af_Z(\varsigma I-\mathsf B)](Z_{\sigma,*})$ and
\[
 F=L(Z-Z_{\sigma,*})+af_Z\mathsf VZ_y,\qquad
 H=\bigl(af_Z(\varsigma I-\mathsf B)-L\bigr)Z_y+af_Z\mathsf Q
   -(af_Z\mathsf V)_yZ_y,
\]
where $a,f_Z,\mathsf B,\mathsf V$ are evaluated at $Z$. Then, for $s\ge S_0$,
\begin{equation}\label{e:E10}
\begin{gathered}
 \sup_{t\ge s}\sqrt\tau\bigl(\Vert Z-Z_{\sigma,a}\Vert_{L^\infty(\Omega_\sigma(t))}
   +\Vert Z_y\Vert_{L^\infty(\Omega_\sigma(t))}\bigr)\le\omega_b(s),\qquad
 a(Z)(\partial_t+\varsigma\partial_y)f(Z)=\partial_yF+H,\\
 \sup_{t\ge S_0}\sqrt\tau\Vert F\Vert_{L^\infty(\Omega_\sigma(t))}\le C_b,\qquad
 \sup_{t\ge s}\tau\Vert H\Vert_{L^\infty(\Omega_\sigma(t))}\le\omega_b(s),
\end{gathered}
\end{equation}
with $C_b,\omega_b$ depending also on $a,f,\varsigma$; if $\varsigma=s_\sigma$,
then $\sup_{t\ge s}\sqrt\tau\Vert F\Vert_{L^\infty(\Omega_\sigma(t))}\le\omega_b(s)$.
\end{lemma}

\begin{corollary}\label{l:coef}
Put
\begin{equation}\label{e:bc}
 h_b=\partial_1u_b,\qquad A_b=D_bu_b,\qquad B_b=D_b\theta_b,\qquad
 g_T=\frac{\partial_1\theta_b}{\theta_b}.
\end{equation}
Then
\begin{equation}\label{e:bal}
 D_b\rho_b=-\rho_bh_b,\qquad \rho_bA_b+\partial_1p_b=\nu\partial_1h_b,\qquad
 c_v\rho_bB_b=\kappa\partial_1^2\theta_b+\nu h_b^2-p_bh_b,
\end{equation}
and there is $C_0$, depending only on the fixed data, such that, for
$s\ge S_0$, with $V_h,\Theta_h,Z,\mathcal K$ as in Lemma~\ref{l:mcoef},
\begin{gather}
 \sup_{t\ge s}\Vert(\rho_b-\rho_-,u_b-u_*,\theta_b-\theta_-)\Vert_{L^\infty}
   \le C_0\delta+\omega_b(s),\label{e:E0}\\
 \sup_{t\ge s}\sqrt\tau\Vert u_b-u_*\Vert_{L^\infty}\le C_A+\omega_b(s),
   \label{e:cbg}
\end{gather}
\begin{equation}\label{e:E1}
\begin{gathered}
 \sup_{t\ge s}\tau\Vert\partial_1(\rho_b,\theta_b)\Vert_{L^\infty}^2
   \le C_0\delta^2+\omega_b(s),\qquad
 \partial_1(\rho_b,\theta_b)=\Bigl(-\frac{(V_h)_y}{V_b^3},
   \frac{(\Theta_h)_y}{V_b}\Bigr)\circ Y+k_1,\\
 \Bigl|\Bigl(-\frac{(V_h)_y}{V_b^3},\frac{(\Theta_h)_y}{V_b}\Bigr)\circ Y\Bigr|
   \le C_0\delta\tau^{-1/2},\qquad
 \sup_{t\ge s}\tau\Vert k_1\Vert_{L^\infty}^2
   +\int_s^\infty\bigl(\Vert k_1\Vert_{L^\infty}^2+\Vert k_1\Vert_{L^2}^2\bigr)\dd t
   \le\omega_b(s),
\end{gathered}
\end{equation}
\begin{gather}
 \sup_{t\ge s}\tau\Vert(h_b,\partial_1p_b)\Vert_{L^\infty}^2\le\omega_b(s),
   \label{e:E2}\\
 \int_s^\infty\Vert(h_b,\partial_1p_b,A_b,B_b,\partial_1h_b,\partial_1^2p_b,
   \partial_1^2\theta_b,\partial_1A_b,\partial_1B_b)\Vert_{L^\infty}^2\dd t
   \le\omega_b(s),\label{e:E3}
\end{gather}
\begin{equation}\label{e:curv}
\begin{gathered}
 \partial_1^2\rho_b=c_\rho+e_{\rho,2},\qquad
 \partial_1^2\theta_b=c_\theta+e_{\theta,2},\\
 c_\rho=\Bigl(-\frac{(V_h)_{yy}}{V_b^4}+\frac{3(V_h)_y^2}{V_b^5}\Bigr)\circ Y,
 \qquad
 c_\theta=\Bigl(\frac{(\Theta_h)_{yy}}{V_b^2}
   -\frac{(V_h)_y(\Theta_h)_y}{V_b^3}\Bigr)\circ Y,\\
 |c_\rho|+|c_\theta|\le\frac{C_0\delta}\tau,\qquad
 \int_s^\infty\bigl(\sqrt\tau\Vert(e_{\rho,2},e_{\theta,2})\Vert_{L^2}^2
   +\tau\Vert(e_{\rho,2},e_{\theta,2})\Vert_{L^\infty}^2\bigr)\dd t
   \le\omega_b(s),
\end{gathered}
\end{equation}
\begin{gather}
 \int_s^\infty\bigl(\Vert\partial_1^3\rho_b\Vert_{L^2}^2
   +\Vert(h_b,B_b)\Vert_{L^2}^2
   +\sqrt\tau\Vert(\partial_1h_b,\partial_1^2h_b,\partial_1A_b,\partial_1B_b)
     \Vert_{L^2}^2\bigr)\dd t\le\omega_b(s),\label{e:E5}\\
 \sup_{t\ge s}\bigl(\Vert\partial_1(\rho_b,u_b,\theta_b)\Vert_{H^2}
   +\Vert(A_b,B_b)\Vert_{H^1}\bigr)\le\omega_b(s),\label{e:E6}
\end{gather}
\begin{equation}\label{e:unwE}
\begin{gathered}
 \sup_{t\ge S_0}\sqrt\tau\bigl(\Vert u_b-u_*\Vert_{L^\infty}
   +\Vert\partial_1Q_b\Vert_{L^\infty}\bigr)<\infty,\\
 \int_s^\infty\bigl(\Vert\partial_1^2Q_b\Vert_{H^1}^2
   +\Vert\partial_1^4(u_b,\theta_b)\Vert_{L^2}^2+\Vert\partial_1u_b\Vert_{L^2}^2
   +\Vert\partial_1p_b\Vert_{H^1}^2\bigr)\dd t\le\omega_b(s),
\end{gathered}
\end{equation}
\begin{equation}\label{e:E8}
\begin{gathered}
 \sup_{t\ge S_0}\tau\Vert(u_b-u_*,\partial_1\rho_b,\partial_1\theta_b,
   \partial_1p_b,h_b)\Vert_{L^\infty}^2<\infty,\qquad
 \int_{S_0}^\infty\tau^a\Vert(A_b,B_b,h_b,\partial_1p_b)\Vert_{L^\infty}^2\dd t
   <\infty\quad(a<\tfrac23),\\
 \int_{S_0}^\infty\tau\Bigl(\Vert(\partial_1h_b,\partial_1^2p_b)\Vert_{L^\infty}^2
   +\Bigl\Vert\Bigl(\partial_1^2p_b,\partial_1A_b,\partial_1B_b,
     \partial_1\frac{\partial_1^2\theta_b}{\rho_b\theta_b}\Bigr)\Bigr\Vert_{L^2}^2
   \Bigr)\dd t<\infty,\\
 \int_{S_0}^\infty\tau^{1/4}\Bigl(\Bigl\Vert\Bigl(\partial_1^2\log\rho_b,
   \partial_1^2\log\theta_b,\frac{\partial_1^2\theta_b}{\rho_b\theta_b}\Bigr)
   \Bigr\Vert_{L^\infty}^2+\Vert\partial_1^3\log\rho_b\Vert_{L^2}^2\Bigr)\dd t
   <\infty,
\end{gathered}
\end{equation}
and, for smooth functions $f,g$ on $\mathcal K$, integers $n\ge0$,
$c_1>0$ and $\xi_\sigma=\frac{Y-s_\sigma t}{\sqrt\tau}$,
\begin{equation}\label{e:floc}
 \int_s^\infty\sqrt\tau\,\Bigl\Vert\bigl(g(Z)\circ Y\bigr)|\xi_\sigma|^n
   e^{-c_1\xi_\sigma^2}\bigl(|\partial_1(f(Z)\circ Y)|
   +|D_b(f(Z)\circ Y)|\bigr)\Bigr\Vert_{L^\infty}^2\dd t\le\omega_b(s),
\end{equation}
with $\omega_b$ depending also on $f,g,n,c_1$.
\end{corollary}

For the reference pressure used in the contact work, define
\begin{equation}\label{e:Uhat}
 \widehat\theta=\Theta_c-\frac{H_c^2}{2c_v},\qquad
 \widehat p=\frac{R\widehat\theta}{V_c}
      =p_*-\frac{(\gamma-1)H_c^2}{2V_c}.
\end{equation}

\subsection{Gaussian observations and pressure decomposition}

For $a>0$, $0\le r\le T\le\infty$ and a function $f$ on $\R\times(r,T)$ put
\begin{equation}\label{e:nobs}
\begin{aligned}
 A_f^a(t)&=\sqrt\tau\int_\R e^{-\frac{a(y+h_*)^2}{2\tau}}|f(y,t)|^2\dd y,
 \qquad
 \mathscr O^a_{r,T}(f)=\int_r^TA_f^a(t)\dd t,\\
 \mathscr N^a_{r,T}(f)&=\int_r^T\Bigl[\tau^{\frac14}
   \bigl\Vert e^{-\frac{a(y+h_*)^2}{4\tau}}f_y\bigr\Vert_{L^2}
   +\tau^{-\frac14}\bigl\Vert e^{-\frac{a(y+h_*)^2}{4\tau}}f\bigr\Vert_{L^2}
   \Bigr]\dd t .
\end{aligned}
\end{equation}

Fix $0<a_s\le c_G/4$ and an auxiliary $0<b_0<a_s/64$, and write
$A_f=A_f^{a_s}$, $\mathscr O=\mathscr O^{a_s}$ and
$\mathscr N=\mathscr N^{a_s}$. The exponent choices are made before
fixing the datum; all entrance times may subsequently be enlarged.

We also use the Gaussian estimate of~\cite[Lemma~D.5]{LiLiuContact1D2026}.
\begin{lemma}\label{l:bulk}
Let $0<a\le c_G$. For every sufficiently late fixed $r$,
\[
 \mathscr O^a_{r,\infty}\bigl(\pa_y(P_b-p_a)\bigr)
 +\mathscr O^a_{r,\infty}\bigl(\pa_y(b_m-b_a)\bigr)<\infty .
\]
\end{lemma}

\begin{lemma}[Imported pressure decomposition]\label{l:prs}
For $0<\delta\le\delta_0(A)$ and the fixed Gaussian exponents above,
there are a background-dependent $r\ge S_0$ and functions $g_1,g_2$ on
$\R\times[r,\infty)$, independent of the final time, such that
\begin{equation}\label{e:pspl}
 P_b-p_*=g_1+g_2,\qquad \mathscr N_{r,\infty}(g_1)<\infty.
\end{equation}
For every $\xi>0$ there is a finite $C_{b,r,\xi}$ such that
\begin{equation}\label{e:prlog}
 \mathscr O_{r,T}(g_2)\le C_{b,r,\xi}
       +(4C_P\delta^2+\xi)\log\frac{2+T}{2+r}\qquad(T\ge r).
\end{equation}
Here $C_P$ depends only on the fixed data and is independent of $A$,
$M$, $a_s$, $b_0$, $r$ and the background. The threshold uniformity is
the one stated above, as proved in
\cite[Proposition~D.20 and Lemma~D.19]{LiLiuContact1D2026}.
\end{lemma}
\section{Perturbation equations and basic estimates}\label{sec:pre}

By Galilean invariance, the three-dimensional estimates are proved in
the frame $(x_1-u_*t,x')$, with the velocity replaced by $u-u_*e_1$.
Spatial Sobolev and uniform norms of the perturbations and their
initial values are unchanged, as are all physical coefficients.
Thus the constants and smallness
thresholds in these estimates can be chosen independently of $u_*$.
Returning to the original frame preserves all stated regularity classes:
the time derivative only acquires $-u_*\partial_1$ of the translated
function, which has the required spatial regularity for each fixed
$u_*$. No uniform bound on those time-derivative norms is asserted.

Except in Lemmas~\ref{l:uniq}--\ref{l:fp},
$(\rho_b,u_b,\theta_b)$ is the background of Lemma~\ref{l:bg} with
$0<\delta\le\delta_1(M)$, the bounds of Corollary~\ref{l:coef} hold for
$s\ge S$, $\delta_1(M)$ may be decreased and $S$ increased in each lemma,
and $(\rho,u,\theta)$ solves \eqref{eq:NSF} on $\R^3\times[s,T_e]$,
$0\le s<T_e<\infty$, with $q$ in the class \eqref{eq:class} on $[s,T_e]$
and $\rho,\theta$ bounded below by positive constants. By \eqref{e:bgcl},
\eqref{e:chain} and the coordinate-conversion argument
in~\cite[Lemma~2.2]{LiLiuContact1D2026}, applied on $[t_0,T]$,
\begin{equation}\label{e:lstr}
 Q_b-\bar Q\in C([t_0,T];H^3(\R)),\qquad
 (u_b-\bar u,\theta_b-\bar\theta)\in L^2(t_0,T;H^4(\R))
 \qquad(0<t_0<T<\infty).
\end{equation}

\subsection{Norms}

For $s\le t$, set
\begin{equation}\label{e:w}
\begin{gathered}
 q_a=\frac1{128},\\
 w_s(t)=\left(\frac{\tau_s}{1+t}\right)^{q_a},\qquad
 I_s=\Vert q(s)\Vert_{H^2}^2.
\end{gathered}
\end{equation}
and
\begin{align}
 E_*(t)&=\Vert(\Pi,\psi)(t)\Vert_{H^1}^2+\Vert\nabla\zeta(t)\Vert_{H^1}^2
   +\tau^{-1}\Vert\zeta(t)\Vert_{L^2}^2,\label{eq:E}\\
 D_*(t)&=\Vert\nabla\Pi(t)\Vert_{L^2}^2+\Vert\nabla\psi(t)\Vert_{H^1}^2
   +\Vert\nabla^2\zeta(t)\Vert_{H^1}^2
   +\tau^{-1}\Vert\nabla\zeta(t)\Vert_{L^2}^2
   +\tau^{-2}\Vert\zeta(t)\Vert_{L^2}^2,\label{eq:D}\\
 N_o^2(s,t)&=\sup_{s\le r\le t}w_s(r)E_*(r)
   +\int_s^tw_s(r)\left(D_*(r)+\frac{q_aE_*(r)}{1+r}\right)\dd r.
   \label{e:Ns}
\end{align}
For a triple $f=(P_f,v_f,z_f)$, $E[f]$, $D[f]$ and $N[f](s,t)$ are the
expressions \eqref{eq:E}--\eqref{e:Ns} with $(\Pi,\psi,\zeta)$ replaced
by $(P_f,v_f,z_f)$; thus
\begin{equation}\label{e:qo}
 q_o=(\Pi,\psi,\zeta),\qquad N_o(s,t)=N[q_o](s,t).
\end{equation}
For a scalar or vector function $Z$ on $\R^3\times[s,t]$,
\begin{equation}\label{e:N0}
 N_0(Z;s,t)^2=\mathop{\rm ess\,sup}_{s\le r\le t}w_s(r)\Vert Z(r)\Vert_{L^2}^2
 +\int_s^tw_s(r)\left(\Vert\nabla Z(r)\Vert_{L^2}^2
       +\frac{q_a\Vert Z(r)\Vert_{L^2}^2}{1+r}\right)\dd r .
\end{equation}
For $k=1,2$,
\begin{equation}\label{e:ul}
 \Vert f\Vert_{H^k_{\rm uloc}}=\sup_{z\in\R^3}\Vert f\Vert_{H^k(z+(-1,1)^3)},
\end{equation}
and on $[s,T_e]$
\begin{equation}\label{e:mu}
 \mu_*=\mathop{\rm ess\,sup}_{s\le t\le T_e}
 \bigl(\Vert q(t)\Vert_{H^2_{\rm uloc}}+\Vert\Pi(t)\Vert_{L^3}
       +\Vert\psi(t)\Vert_{L^3}\bigr).
\end{equation}

\subsection{Equations}

Set
\begin{equation}\label{e:Lam}
 \mathcal L=\mu\Delta+(\mu+\lambda)\nabla\Div,\qquad
 Q_\mu(\nabla u)=2\mu|\mathbb D(u)|^2+\lambda(\Div u)^2,\qquad
 H[\psi]=2h_b(2\mu\partial_1\psi_1+\lambda\Div\psi).
\end{equation}
Then
\begin{align}
 D_u\phi+\rho\Div\psi
 &=-\psi_1\partial_1\rho_b-\phi\partial_1u_b,\label{eq:ph-phi}\\
 \rho D_u\psi+\nabla\Pi-\mathcal L\psi
 &=-\phi D_bu_b\,e_1-\rho\psi_1\partial_1u_b\,e_1,\label{eq:ph-psi}\\
 c_v\rho D_u\zeta+p\Div\psi-\kappa\Delta\zeta
 &=-c_v\phi D_b\theta_b-c_v\rho\psi_1\partial_1\theta_b
   -\Pi\partial_1u_b+H[\psi]+Q_\mu(\nabla\psi),\label{eq:phzeta}
\end{align}
and, with
\begin{equation}\label{e:ps3}
\begin{aligned}
 f_p&=-\partial_1p_b\,\psi_1-\gamma h_b\Pi
                  +(\gamma-1)(H[\psi]+Q_\mu(\nabla\psi)),\qquad
 f_m=-\rho h_b\psi_1e_1-\phi A_be_1,\\
 f_z&=-\partial_1\theta_b\,\psi_1-\frac\phi\rho B_b
       -\frac{h_b\Pi}{c_v\rho}
       +\frac{H[\psi]+Q_\mu(\nabla\psi)}{c_v\rho},
\end{aligned}
\end{equation}
\begin{equation}\label{e:peq3}
\begin{cases}
 D_u\Pi+\gamma p\Div\psi-(\gamma-1)\kappa\Delta\zeta=f_p,\\
 \rho D_u\psi+\nabla\Pi-\mathcal L\psi=f_m,\\
 D_u\zeta-\frac{\kappa}{c_v\rho}\Delta\zeta
                       +(\gamma-1)\theta\Div\psi=f_z.
\end{cases}
\end{equation}
The state law and the total-energy equation give
\[
 c_v\rho D_u\theta+p\Div u=\kappa\Delta\theta+Q_\mu(\nabla u),\qquad
 D_up+\gamma p\Div u=(\gamma-1)\bigl(\kappa\Delta\theta+Q_\mu(\nabla u)\bigr).
\]
Subtract the background equations, using
\[
 D_uf=D_bf+\psi_1\partial_1f\quad(f=f(x_1,t)),\qquad
 \Div u=h_b+\Div\psi,
\]
\[
 Q_\mu(\nabla u)=\nu h_b^2+H[\psi]+Q_\mu(\nabla\psi).
\]
This gives \eqref{eq:ph-phi}--\eqref{eq:phzeta} and \eqref{e:peq3}.
For $q$ in \eqref{eq:class}, the products are locally integrable by
\eqref{e:brg}, so the identities also hold in distributions.

With $\nabla\Pi=R\theta\nabla\phi+R\rho\nabla\zeta
+R(\zeta\partial_1\rho_b+\phi\partial_1\theta_b)e_1$, \eqref{eq:ph-phi},
\eqref{eq:ph-psi} divided by $\rho$ and the third row of \eqref{e:peq3}
become
\begin{equation}\label{e:cpeq}
\begin{cases}
 (\partial_t+u_*\partial_1)\phi+\rho_-\Div\psi=N_\phi,\\
 (\partial_t+u_*\partial_1)\psi-\frac1{\rho_-}\mathcal L\psi
            +\frac{R\theta_-}{\rho_-}\nabla\phi+R\nabla\zeta=N_\psi,\\
 (\partial_t+u_*\partial_1)\zeta-\frac{\kappa}{c_v\rho_-}\Delta\zeta
                              +\frac{R\theta_-}{c_v}\Div\psi=N_\zeta ,
\end{cases}
\end{equation}
where
\begin{equation}\label{e:nphys}
\begin{aligned}
 N_\phi={}&-((u_b-u_*)e_1+\psi)\cdot\nabla\phi
       -(\rho-\rho_-)\Div\psi-(\partial_1\rho_b)\psi_1-h_b\phi,\\
 N_\psi={}&-((u_b-u_*)e_1+\psi)\cdot\nabla\psi
       +\left(\frac1\rho-\frac1{\rho_-}\right)\mathcal L\psi
       -\left(\frac{R\theta}{\rho}-\frac{R\theta_-}{\rho_-}\right)\nabla\phi\\
 &-\frac{\phi}{\rho}A_be_1-h_b\psi_1e_1
       -\frac R\rho(\phi\partial_1\theta_b+\zeta\partial_1\rho_b)e_1,\\
 N_\zeta={}&-((u_b-u_*)e_1+\psi)\cdot\nabla\zeta
       +\frac\kappa{c_v}\left(\frac1\rho-\frac1{\rho_-}\right)\Delta\zeta
       -\frac R{c_v}(\theta-\theta_-)\Div\psi+f_z .
\end{aligned}
\end{equation}

\subsection{Local theory}

Let $0<T_f<\infty$, and let $(\rho_b,u_b,\theta_b)$ be a solution of
\eqref{eq:bg-eqs} on $\R\times[0,T_f]$ with
\begin{equation}\label{e:locb}
\begin{aligned}
 &(\rho_b,u_b,\theta_b)\in C([0,T_f];L^\infty(\R)),\qquad
 \inf_{\R\times[0,T_f]}\rho_b>0,\qquad
 \inf_{\R\times[0,T_f]}\theta_b>0,\\
 &\partial_1\rho_b\in C([0,T_f];H^2(\R)),\qquad
 \partial_1(u_b,\theta_b)
     \in C([0,T_f];H^1(\R))\cap L^2(0,T_f;H^2(\R)).
\end{aligned}
\end{equation}
By
Lemma~\ref{l:bg}, \eqref{e:brg} and \eqref{eq:pr-gau}, the background
of Lemma~\ref{l:bg} satisfies \eqref{e:locb} for every $T_f$. Besides
Lemma~\ref{l:local}, the following three results are used; they are
proved with it in Appendix~\ref{app:local}.

\begin{lemma}\label{l:uniq}
Let $0\le t_1<t_2\le T_f$. Two solutions of \eqref{eq:NSF} on
$\R^3\times[t_1,t_2]$ with $q$ in the class \eqref{eq:class} on
$[t_1,t_2]$, $\rho,\theta$ bounded below by positive constants and the
same value at $t_1$ coincide. Every such solution satisfies
$\phi_t\in C([t_1,t_2];H^1)$, and, if $(\rho_b,u_b,\theta_b)$ tends to
$(\rho_\pm,u_{1\pm},\theta_\pm)$ as $x_1\to\pm\infty$, also \eqref{eq:far}
at each $t\in[t_1,t_2]$, uniformly in $x'$.
\end{lemma}

\begin{lemma}\label{l:cont}
Let $0<T_{\rm fin}<T_f$, and let $(\rho,u,\theta)$ solve \eqref{eq:NSF}
on $\R^3\times[0,T_{\rm fin})$ with $q$ in \eqref{eq:class} on every
$[0,T]$, $T<T_{\rm fin}$, and $\rho,\theta$ bounded below by positive
constants on each such $[0,T]$. If
\begin{equation}\label{e:cont}
 \sup_{0\le t<T_{\rm fin}}\Vert q(t)\Vert_{H^2}<\infty,
 \qquad
 \inf_{\R^3\times[0,T_{\rm fin})}\min\{\rho,\theta\}>0,
\end{equation}
then $q(t)$ converges in $H^2$ as $t\to T_{\rm fin}$, and for some
$\ell>0$ the solution extends to $\R^3\times[0,T_{\rm fin}+\ell]$, with
$q$ in \eqref{eq:class} on $[0,T_{\rm fin}+\ell]$ and $\rho,\theta$
bounded below by a positive constant.
\end{lemma}

\begin{lemma}\label{l:fp}
For $0<T\le T_f$ there are $\epsilon_{b,T}>0$ and $K_{b,T}<\infty$,
depending only on $T$, $(R,\gamma,\mu,\lambda,\kappa)$, the lower
bounds in \eqref{e:locb} and the norms in \eqref{e:locb} on $[0,T]$,
such that for every $q_0\in H^2(\R^3)$ with
$\Vert q_0\Vert_{H^2}\le\epsilon_{b,T}$, \eqref{eq:NSF} with $q(0)=q_0$
has a solution on $\R^3\times[0,T]$ with $q$ in \eqref{eq:class} on
$[0,T]$,
\begin{equation}\label{e:fp}
 \sup_{0\le r\le T}\Vert q(r)\Vert_{H^2}^2
 +\int_0^{T}\bigl(\Vert\nabla(\psi,\zeta)(r)\Vert_{H^2}^2
                       +\Vert q_t(r)\Vert_{H^1}^2\bigr)\dd r
 \le K_{b,T}^2\Vert q_0\Vert_{H^2}^2,
\end{equation}
$\rho\ge\frac12\inf_{\R\times[0,T]}\rho_b$ and
$\theta\ge\frac12\inf_{\R\times[0,T]}\theta_b$.
\end{lemma}

\subsection{Basic estimates}

For $f,\partial_1f\in L^2(\R^3)$,
\begin{equation}\label{eq:planar}
 \sup_{x_1\in\R}\Vert f(x_1,\cdot)\Vert_{L^2(\R^2)}^2
 \le2\Vert f\Vert_{L^2(\R^3)}\Vert\partial_1f\Vert_{L^2(\R^3)},
\end{equation}
since $\Vert f(x_1,\cdot)\Vert_{L^2(\R^2)}^2
=2\int_{-\infty}^{x_1}\int_{\R^2}(f\partial_1f)(y_1,x')\dd x'\dd y_1$
for $f\in C_c^\infty(\R^3)$. For $a=a(x_1)$ and $j=1,2$, the Leibniz
formula, with \eqref{eq:planar} for the term containing $\partial_1^ja$,
gives
\begin{equation}\label{e:pmul}
 \Vert aF\Vert_{H^j}
 \le C\bigl(\Vert a\Vert_{L^\infty(\R)}
                 +\Vert\partial_1a\Vert_{H^{j-1}(\R)}\bigr)\Vert F\Vert_{H^j}.
\end{equation}
In particular, $\Pi=R(\theta_b\phi+\rho_b\zeta+\phi\zeta)$ and
\eqref{e:brg} give
$\Vert\Pi\Vert_{H^2}\le C_{b,T}(1+\Vert q\Vert_{H^2})\Vert q\Vert_{H^2}$
on $[0,T]$.

Let $Q_j=j+(0,1)^3$, $j\in\mathbb Z^3$. For $k=1,2$,
$\sum_j\Vert f\Vert_{H^k(Q_j)}^2=\Vert f\Vert_{H^k}^2$ and
$\Vert f\Vert_{H^k(Q_j)}\le\Vert f\Vert_{H^k_{\rm uloc}}$; moreover
$\Vert f\Vert_{L^6(Q_j)}\le C\Vert f\Vert_{H^1(Q_j)}$ and
$\Vert f\Vert_{L^\infty(Q_j)}\le C\Vert f\Vert_{H^2(Q_j)}$ with $C$
independent of $j$.

\begin{lemma}\label{l:st}
There is $\mu_0>0$, depending on the background, such that
$\mu_*\le\mu_0$ implies, on $\R^3\times[s,T_e]$,
\begin{equation}\label{e:stb}
 C_b^{-1}\le\rho,\theta,p\le C_b,\qquad
 \Vert(q,\Pi)\Vert_{L^\infty}
 +\Vert(q,\Pi)\Vert_{H^2_{\rm uloc}}\le C_b\mu_*,
\end{equation}
\begin{equation}\label{e:state}
 \phi=\frac{\Pi}{R\theta}-\frac{\rho_b\zeta}{\theta},\qquad
 \nabla\phi=\frac{\nabla\Pi}{R\theta}
 -\frac\rho\theta\nabla\zeta
 -\frac{\zeta\partial_1\rho_b+\phi\partial_1\theta_b}{\theta}e_1,\qquad
 \Vert\phi\Vert_{H^1}\le C_b\Vert(\Pi,\zeta)\Vert_{H^1}.
\end{equation}
If $s\ge S$, then \eqref{e:stb} and \eqref{e:state} hold with $C_M$ in
place of $C_b$ for some $\mu_0$ depending only on $M$ and the fixed data,
and, for $s\le t\le T_e$,
\begin{equation}\label{e:trace}
 \Vert\nabla\phi(t)\Vert_{L^2}^2\le C_ME_*(t),\qquad
 \sup_{x_1\in\R}\Vert(\phi,\zeta)(x_1,\cdot,t)\Vert_{L^2(\R^2)}^2
 \le C_M\sqrt\tau E_*(t).
\end{equation}
For $G\in H^1$ and $F\in H^1_{\rm uloc}$, resp. $F\in H^2_{\rm uloc}$,
\begin{equation}\label{e:prod}
 \Vert FG\Vert_{L^2}\le C\Vert F\Vert_{H^1_{\rm uloc}}\Vert G\Vert_{H^1},
 \qquad
 \Vert FG\Vert_{H^1}\le C\Vert F\Vert_{H^2_{\rm uloc}}\Vert G\Vert_{H^1}.
\end{equation}
For scalar $a\in H^2$, resp. $a,d\in\{\phi,\psi_1,\psi_2,\psi_3,\zeta\}$,
\begin{equation}\label{e:gn4}
 \Vert\nabla a\Vert_{L^4}^2\le C\Vert a\Vert_{L^\infty}\Vert\nabla^2a\Vert_{L^2},
 \qquad
 \Vert\nabla a\,\nabla d\Vert_{L^2}
 \le C\mu_*\bigl(\Vert\nabla^2a\Vert_{L^2}+\Vert\nabla^2d\Vert_{L^2}\bigr).
\end{equation}
\end{lemma}
\begin{proof}
The embedding $H^2(Q_j)\subset L^\infty(Q_j)$ gives
$\Vert q\Vert_{L^\infty}\le C\mu_*$. By Lemma~\ref{l:bg},
$C_b^{-1}\le\rho_b,\theta_b\le C_b$; choosing $C\mu_0$ below half the
lower bound gives the bounds for $\rho$, $\theta$ and $p=R\rho\theta$. The identity
\[
 \Pi=R(\theta_b\phi+\rho_b\zeta+\phi\zeta)=R(\theta\phi+\rho_b\zeta)
\]
and the algebra property of $H^2(Q_j)$ give
$\Vert\Pi\Vert_{H^2_{\rm uloc}}\le C_b(\mu_*+\mu_*^2)\le C_b\mu_*$,
with $\mu_0\le1$. Here
$\Vert(\rho_b,\theta_b)\Vert_{L^\infty}
+\Vert\partial_1(\rho_b,\theta_b)\Vert_{H^1(\R)}\le C_b$
by \eqref{e:brg}, \eqref{eq:pr-gau} and \eqref{e:E6}.
The local embedding also gives the $L^\infty$ bound for $\Pi$.
Solving the second form of $\Pi$ for $\phi$ and differentiating, with
$\nabla\theta=(\partial_1\theta_b)e_1+\nabla\zeta$, gives the identities
in \eqref{e:state}; they and
$\Vert\partial_1(\rho_b,\theta_b)\Vert_{L^\infty}\le C_b$ give the $H^1$
bound. For $s\ge S$, \eqref{e:E0} and \eqref{e:E6} give
$|\rho_b-\rho_-|\le\frac{\rho_-}2$, $|\theta_b-\theta_-|\le\frac{\theta_-}2$
and $\Vert\partial_1(\rho_b,\theta_b)\Vert_{H^1(\R)}\le1$ after
decreasing $\delta_1(M)$ and increasing $S$, so the preceding constants
and $\mu_0$ depend only on $M$ and the fixed data.

Let $s\ge S$. By \eqref{e:state} and \eqref{e:E1},
\[
 \Vert\phi\Vert_{L^2}\le C_M\Vert(\Pi,\zeta)\Vert_{L^2},\qquad
 \Vert\nabla\phi\Vert_{L^2}
 \le C_M\bigl(\Vert\nabla\Pi\Vert_{L^2}+\Vert\nabla\zeta\Vert_{L^2}
       +\tau^{-\frac12}\Vert(\Pi,\zeta)\Vert_{L^2}\bigr),
\]
and the square of the last bound is at most $C_ME_*$. By
\eqref{eq:planar},
$\sup_{x_1}\Vert\Pi(x_1,\cdot)\Vert_{L^2(\R^2)}^2
\le2\Vert\Pi\Vert_{L^2}\Vert\partial_1\Pi\Vert_{L^2}\le E_*$ and
\[
 \sup_{x_1}\Vert\zeta(x_1,\cdot)\Vert_{L^2(\R^2)}^2
 \le2\tau^{\frac12}\bigl(\tau^{-\frac12}\Vert\zeta\Vert_{L^2}\bigr)
       \Vert\partial_1\zeta\Vert_{L^2}
 \le\sqrt\tau E_*;
\]
the first identity of \eqref{e:state} on each slice $\{x_1\}\times\R^2$
gives \eqref{e:trace}.

H\"older's inequality and the embeddings on $Q_j$ give
\[
\begin{aligned}
 \Vert FG\Vert_{L^2(Q_j)}
 &\le\Vert F\Vert_{L^6(Q_j)}\Vert G\Vert_{L^3(Q_j)}
 \le C\Vert F\Vert_{H^1_{\rm uloc}}\Vert G\Vert_{H^1(Q_j)},\\
 \Vert\nabla(FG)\Vert_{L^2(Q_j)}
 &\le\Vert F\Vert_{L^\infty(Q_j)}\Vert\nabla G\Vert_{L^2(Q_j)}
       +\Vert\nabla F\Vert_{L^6(Q_j)}\Vert G\Vert_{L^3(Q_j)}
 \le C\Vert F\Vert_{H^2_{\rm uloc}}\Vert G\Vert_{H^1(Q_j)};
\end{aligned}
\]
squaring and summing over $j$ gives \eqref{e:prod}.
The first bound in \eqref{e:gn4} is the Gagliardo--Nirenberg
inequality~\cite[Lecture~II]{Nirenberg1959}.
Its application to $a,d$, followed by H\"older's inequality and
$\Vert(a,d)\Vert_{L^\infty}\le C\mu_*$, gives
\[
 \Vert\nabla a\,\nabla d\Vert_{L^2}
 \le C\mu_*(\Vert\nabla^2a\Vert_{L^2}\Vert\nabla^2d\Vert_{L^2})^{\frac12}
 \le C\mu_*(\Vert\nabla^2a\Vert_{L^2}+\Vert\nabla^2d\Vert_{L^2}).
\]
\end{proof}

Let $s\ge S$, and let $Z$ be a scalar or vector function with
$N_0(Z;s,t)<\infty$. Since
$\int_s^t\frac{w_s}{1+r}\Vert Z\Vert_{L^2}^2\dd r\le q_a^{-1}N_0(Z;s,t)^2$
by \eqref{e:N0}, \eqref{e:E1} and \eqref{e:E3} give
\begin{equation}\label{e:E1Z}
\begin{aligned}
 \int_s^tw_s\Vert(|\partial_1\rho_b|+|\partial_1\theta_b|)Z\Vert_{L^2}^2\dd r
 &\le\frac{2(C_0\delta^2+\omega_b(s))}{q_a}N_0(Z;s,t)^2,\\
 \int_s^tw_s\Vert aZ\Vert_{L^2}^2\dd r
 &\le\omega_b(s)\mathop{\rm ess\,sup}_{s\le r\le t}w_s\Vert Z\Vert_{L^2}^2
 \le\omega_b(s)N_0(Z;s,t)^2
\end{aligned}
\end{equation}
for $a\in\{h_b,\partial_1p_b,A_b,B_b,\partial_1h_b,\partial_1^2p_b,
\partial_1^2\theta_b,\partial_1A_b,\partial_1B_b\}$.

\subsection{The second-order estimate}

\begin{lemma}\label{l:ha}
There is $\mu_0>0$, depending only on $M$ and the fixed data, such that
$s\ge S$ and $\mu_*\le\mu_0$ imply, for $s<t<T_e$,
\begin{equation}\label{e:ha}
\begin{aligned}
 &\sup_{s\le r\le t}w_s(r)
       \Vert\nabla^2(\phi,\psi)(r)\Vert_{L^2}^2
 +\int_s^t w_s(r)\left[
       \Vert\nabla^2\phi\Vert_{L^2}^2+\Vert\nabla^3\psi\Vert_{L^2}^2
       +\frac{q_a}{1+r}\Vert\nabla^2(\phi,\psi)\Vert_{L^2}^2
                         \right]\dd r\\
 &\qquad\le C_M\bigl(I_s+N_o^2(s,t)\bigr).
\end{aligned}
\end{equation}
\end{lemma}
\begin{proof}
Put
\[
 X=N_o^2(s,t)+\int_s^tw_s\bigl(\Vert\nabla^2\phi\Vert_{L^2}^2
       +\Vert\nabla^3\psi\Vert_{L^2}^2\bigr)\dd r .
\]
For $a_0=a_0(x_1,t)$ and $z_0\in\{\phi,\psi,\zeta\}$,
$\Vert a_0z_0\Vert_{L^2}^2\le\Vert a_0\Vert_{L^2(\R)}^2
\sup_{x_1}\Vert z_0(x_1,\cdot)\Vert_{L^2(\R^2)}^2$, and the last factor is
at most $C_M\sqrt\tau E_*$ by \eqref{e:trace} and \eqref{eq:planar}, and
at most $E_*$ if $z_0=\psi$. Hence
\[
 \int_s^t w_s\Vert a_0z_0\Vert_{L^2}^2\dd r
 \le C_MN_o^2(s,t)\int_s^t\sqrt{1+r}\,
                      \Vert a_0(r)\Vert_{L^2(\R)}^2\dd r,
\]
where $\sqrt{1+r}$ may be omitted if $z_0=\psi$. By \eqref{e:curv} and
\eqref{e:E5},
\begin{equation}\label{e:hpl}
 \int_s^t w_s\Vert a_0z_0\Vert_{L^2}^2\dd r\le\omega_b(s)N_o^2(s,t)
\end{equation}
for $a_0\in\{e_{\rho,2},e_{\theta,2},\partial_1h_b,\partial_1^2h_b,
\partial_1A_b,\partial_1B_b\}$, $z_0\in\{\phi,\psi,\zeta\}$, and for
$a_0=\partial_1^3\rho_b$, $z_0=\psi$. By $|c_\rho|+|c_\theta|\le\frac{C_M\delta}\tau$
in \eqref{e:curv}, $\Vert\phi\Vert_{L^2}\le C_M\Vert(\Pi,\zeta)\Vert_{L^2}$,
$\tau^{-2}\Vert(\Pi,\psi)\Vert_{L^2}^2\le\tau^{-1}E_*$ and
$\tau^{-2}\Vert\zeta\Vert_{L^2}^2\le D_*$,
\[
 \int_s^t w_s\Vert(c_\rho,c_\theta)(\phi,\psi,\zeta)\Vert_{L^2}^2\dd r
 \le C_M\delta^2\int_s^t\frac{w_s}{(1+r)^2}
                            \Vert(\Pi,\psi,\zeta)\Vert_{L^2}^2\dd r
 \le C_M\delta^2N_o^2(s,t),
\]
so \eqref{e:hpl} holds for $a_0\in\{\partial_1^2\rho_b,\partial_1^2\theta_b\}$
with $C_M\delta^2+\omega_b(s)$ in place of $\omega_b(s)$.

We estimate $\nabla(N_\phi,N_\psi,N_\zeta)$ using \eqref{e:nphys}.
By \eqref{e:stb}, its coefficients are bounded; the terms have six forms.
First, $a\,Y$ with $|a|\le|u_b-u_*|+|\rho-\rho_-|+|\theta-\theta_-|+|\psi|$
and $Y\in\{\nabla^2\phi,\nabla^2\psi,\nabla^3\psi,\nabla^2\zeta,
\nabla^3\zeta\}$; by \eqref{e:E0} and \eqref{e:stb},
$\Vert a\Vert_{L^\infty}\le C_0\delta+\omega_b(s)+C\mu_*$.

Second,
products $\nabla a\,\nabla d$, $\nabla a\,\nabla^2d$ and
$|\nabla\psi|^2\nabla\phi$ with $a,d\in\{\phi,\psi_k,\zeta\}$ and
$d\neq\phi$ in $\nabla^2d$; by \eqref{e:gn4}, \eqref{e:prod} with
$\Vert\nabla a\Vert_{H^1_{\rm uloc}}\le\mu_*$, and
$\Vert|\nabla\psi|^2\nabla\phi\Vert_{L^2(Q_j)}
\le\Vert\nabla\psi\Vert_{L^6(Q_j)}^2\Vert\nabla\phi\Vert_{L^6(Q_j)}
\le C\mu_*^2\Vert\nabla\psi\Vert_{H^1(Q_j)}$, their norms are at most
$C\mu_*(\Vert\nabla\psi\Vert_{H^2}+\Vert\nabla^2(\phi,\zeta)\Vert_{L^2}
+\Vert\nabla^3\zeta\Vert_{L^2})$.

Third,
$(\partial_1\rho_b,\partial_1\theta_b)$ times
$\nabla\phi,\nabla\psi,\nabla\zeta,\nabla^2\psi,\nabla^2\zeta$; we use
\eqref{e:E1Z} with $N_0(\nabla\phi;s,t)^2\le C_MX$ by \eqref{e:trace},
$N_0(\nabla(\psi,\zeta);s,t)\le N_o(s,t)$, and
$\Vert\partial_1(\rho_b,\theta_b)\Vert_{L^\infty}^2\le C_0\delta^2+\omega_b(s)$
for the second derivatives.

Fourth, $h_b\nabla(\phi,\psi,\zeta,\Pi)$,
$h_b\nabla^2\psi$, $(A_b,B_b)\nabla\phi$, $\partial_1h_b(\psi,\Pi,\nabla\psi)$
and $(\partial_1^2\theta_b)\psi$; we use \eqref{e:E1Z} with
$\sup_{s\le r\le t}w_s\Vert(\psi,\Pi,\nabla\phi,\nabla\psi,\nabla\zeta,
\nabla\Pi)\Vert_{L^2}^2\le C_MN_o^2(s,t)$, and
$\Vert h_b\Vert_{L^\infty}\le\omega_b(s)$ by \eqref{e:E6} for
$h_b\nabla^2\psi$.

Fifth, $(A_b,B_b,h_b)(\partial_1\rho_b)(\phi,\Pi)$ and
$(\partial_1\rho_b,\partial_1\theta_b)^2(\phi,\zeta)$; by \eqref{e:E1} and
\eqref{e:state}, their squared norms are at most
$C_M\Vert(A_b,B_b,h_b)\Vert_{L^\infty}^2E_*$ and
$C_M(C_0\delta^2+\omega_b(s))^2\tau^{-1}E_*$, and we use \eqref{e:E3} and
\eqref{e:Ns}.

Sixth, $(\partial_1^2\rho_b,\partial_1^2\theta_b)
(\phi,\psi,\zeta)$ and $(\partial_1h_b,\partial_1A_b,\partial_1B_b)\phi$,
bounded by \eqref{e:hpl}. Altogether,
\begin{equation}\label{e:hN}
 \int_s^t w_s\Vert\nabla(N_\phi,N_\psi,N_\zeta)\Vert_{L^2}^2\dd r
 \le\bigl[C_M(\delta+\mu_*)+\omega_b(s)\bigr]^2X .
\end{equation}

Apply $\partial_{ij}=\partial_i\partial_j$ to \eqref{e:cpeq}, multiply
the three rows by $\frac{R\theta_-}{\rho_-}\partial_{ij}\phi$,
$\rho_-\partial_{ij}\psi$ and $\frac{c_v\rho_-}{\theta_-}\partial_{ij}\zeta$,
integrate over $\R^3$ and sum over $i,j$. The terms
$\pm R\theta_-\int_{\R^3}\partial_{ij}\phi\Div\partial_{ij}\psi\dd x$ and
$\pm R\rho_-\int_{\R^3}\partial_{ij}\zeta\Div\partial_{ij}\psi\dd x$ cancel.
In the first row write, by \eqref{e:nphys},
\[
 \partial_{ij}N_\phi=-w\cdot\nabla\partial_{ij}\phi+G_{ij},\qquad
 G_{ij}=-[\partial_{ij},w\cdot\nabla]\phi
 -\partial_{ij}\bigl((\rho-\rho_-)\Div\psi
                      +(\partial_1\rho_b)\psi_1+h_b\phi\bigr),
\]
with $w=(u_b-u_*)e_1+\psi$ and $\Div w=h_b+\Div\psi$. With
\[
 E_2=\frac12\int_{\R^3}\left[
       \frac{R\theta_-}{\rho_-}|\nabla^2\phi|^2
       +\rho_-|\nabla^2\psi|^2
       +\frac{c_v\rho_-}{\theta_-}|\nabla^2\zeta|^2\right]\dd x,
\]
we obtain
\begin{equation}\label{e:h2id}
\begin{aligned}
 &\frac{\dd}{\dd t}E_2+\sum_{i,j}\int_{\R^3}\Bigl[
       \mu|\nabla\partial_{ij}\psi|^2+(\mu+\lambda)(\Div\partial_{ij}\psi)^2
       +\frac{\kappa}{\theta_-}|\nabla\partial_{ij}\zeta|^2\Bigr]\dd x\\
 &\quad=\sum_{i,j}\int_{\R^3}\Bigl[\frac{R\theta_-}{\rho_-}
       \Bigl(\frac{\Div w}{2}|\partial_{ij}\phi|^2+\partial_{ij}\phi\,G_{ij}\Bigr)
       -\rho_-\partial_i\partial_{ij}\psi\cdot\partial_jN_\psi
       -\frac{c_v\rho_-}{\theta_-}\partial_i\partial_{ij}\zeta\,
                  \partial_jN_\zeta\Bigr]\dd x .
\end{aligned}
\end{equation}
By \eqref{e:nphys}, $G_{ij}$ consists of $(\rho-\rho_-)\Div\partial_{ij}\psi$,
of products $\nabla\psi\,\nabla^2\phi$ and $\nabla\phi\,\nabla^2\psi$, and of
bounded multiples of
\[
 h_b\nabla^2\phi,\quad (\partial_1h_b)\nabla\phi,\quad
 (\partial_1\rho_b)\nabla^2\psi,\quad (\partial_1^2\rho_b)\nabla\psi,\quad
 (\partial_1^3\rho_b)\psi,\quad (\partial_1^2h_b)\phi .
\]
On each $Q_j$, $\Vert\nabla\psi\Vert_{L^\infty(Q_j)}\le
C\Vert\nabla\psi\Vert_{H^2(Q_j)}$,
$\Vert\nabla^2\psi\Vert_{L^3(Q_j)}\le C\Vert\nabla^2\psi\Vert_{H^1(Q_j)}$,
$\Vert\nabla^2\phi\Vert_{L^2(Q_j)}\le\mu_*$ and
$\Vert\nabla\phi\Vert_{L^6(Q_j)}\le C\mu_*$; summing over $j$ with the
Cauchy--Schwarz inequality,
\[
 \int_{\R^3}\bigl(|\nabla\psi||\nabla^2\phi|^2
       +|\nabla^2\psi||\nabla\phi||\nabla^2\phi|\bigr)\dd x
 \le C\mu_*\Vert\nabla^2\phi\Vert_{L^2}\Vert\nabla\psi\Vert_{H^2}.
\]
By this inequality, $\Vert\rho-\rho_-\Vert_{L^\infty}\le C_0\delta+\omega_b(s)
+C\mu_*$, $\Vert(h_b,\partial_1\rho_b,\partial_1^2\rho_b)\Vert_{L^\infty}
\le\omega_b(s)$ by \eqref{e:E6}, \eqref{e:E1Z} for $(\partial_1h_b)\nabla\phi$,
\eqref{e:hpl} for the last two terms, and, for the last two integrals in
\eqref{e:h2id}, \eqref{e:hN} and $\int_s^tw_s\Vert\nabla^3\zeta\Vert_{L^2}^2
\dd r\le N_o^2(s,t)$, the Cauchy--Schwarz inequality in
$\int_s^tw_s\,\cdot\,\dd r$ gives
\begin{equation}\label{e:h2R}
 \int_s^tw_s|\mathcal R_2|\dd r\le\bigl[C_M(\delta+\mu_*)+\omega_b(s)\bigr]X,
\end{equation}
where $\mathcal R_2$ is the right-hand side of \eqref{e:h2id}.

Differentiating
$J_1=\sum_i\int_{\R^3}\partial_i\psi\cdot\nabla\partial_i\phi\dd x$ by
\eqref{e:cpeq} and integrating by parts in the term containing
$\partial_t\partial_i\phi$ give
\begin{equation}\label{e:J1id}
\begin{aligned}
 \frac{\dd}{\dd t}J_1+\frac{R\theta_-}{\rho_-}\Vert\nabla^2\phi\Vert_{L^2}^2
 ={}&\sum_i\int_{\R^3}\Bigl[\frac1{\rho_-}\mathcal L\partial_i\psi\cdot
       \nabla\partial_i\phi
       -R\nabla\partial_i\zeta\cdot\nabla\partial_i\phi
       +\rho_-(\Div\partial_i\psi)^2\\
 &+\partial_iN_\psi\cdot\nabla\partial_i\phi
       -\Div\partial_i\psi\,\partial_iN_\phi\Bigr]\dd x,
\end{aligned}
\end{equation}
Young's inequality gives, with $C_2$ depending only on the fixed data,
\[
 \frac{\dd}{\dd t}J_1+\frac{R\theta_-}{2\rho_-}\Vert\nabla^2\phi\Vert_{L^2}^2
 \le C_2\bigl(\Vert\nabla^3\psi\Vert_{L^2}^2+\Vert\nabla^2(\psi,\zeta)\Vert_{L^2}^2
       +\Vert\nabla(N_\phi,N_\psi)\Vert_{L^2}^2\bigr).
\]
The $\phi$-terms of \eqref{e:h2id} hold in integrated form by
\eqref{e:phid} in Appendix~\ref{app:local}, applied to $\partial_{ij}\phi$; the other terms of
\eqref{e:h2id} and \eqref{e:J1id} follow from
$\partial_{ij}(\psi,\zeta)\in L^2(s,t;H^1)$,
$\partial_t\partial_{ij}(\psi,\zeta)\in L^2(s,t;H^{-1})$ and
$\phi_t\in C([s,t];H^1)$ (Lemma~\ref{l:uniq}).

Let $\varepsilon_1=\min\{1,\frac{\mu}{2C_2}\}$. Add $\varepsilon_1$ times
the last inequality to \eqref{e:h2id}, use $\mu+\lambda>0$, multiply by
$w_s$, use $\frac{\dd}{\dd t}w_s=-\frac{q_aw_s}{1+t}$ and integrate over
$(s,t')$, $t'\le t$. Since
$\varepsilon_1|J_1|\le\Vert\nabla\psi\Vert_{L^2}\Vert\nabla^2\phi\Vert_{L^2}
\le\frac12E_2+C\Vert\nabla\psi\Vert_{L^2}^2$, $E_2(s)+|J_1(s)|\le CI_s$
and $\sup_{s\le r\le t}w_s\Vert\nabla\psi\Vert_{L^2}^2
+\int_s^tw_s\bigl(\Vert\nabla^2(\psi,\zeta)\Vert_{L^2}^2
+\frac{q_a}{1+r}\Vert\nabla\psi\Vert_{L^2}^2\bigr)\dd r\le N_o^2(s,t)$,
\eqref{e:hN} and \eqref{e:h2R} give
\[
\begin{aligned}
 &\sup_{s\le r\le t}w_sE_2
 +\int_s^tw_s\left[\frac{q_aE_2}{1+r}
       +\frac1{C}\bigl(\Vert\nabla^2\phi\Vert_{L^2}^2
                        +\Vert\nabla^3\psi\Vert_{L^2}^2\bigr)\right]\dd r\\
 &\qquad\le C\bigl(I_s+N_o^2(s,t)\bigr)
     +C\bigl[C_M(\delta+\mu_*)+\omega_b(s)\bigr]X .
\end{aligned}
\]
Decrease $\delta_1(M)$ and $\mu_0$, depending only on $M$ and the fixed
data, and then increase $S$, so that
$C^2[C_M(\delta+\mu_*)+\omega_b(s)]\le\frac12$. The part of $X$ other than
$N_o^2$ is absorbed, and $E_2\ge C^{-1}\Vert\nabla^2(\phi,\psi)\Vert_{L^2}^2$
gives \eqref{e:ha}.
\end{proof}

\section{Weighted energy estimates}\label{sec:en}

Throughout this section \eqref{e:dat-aux} holds, $(\rho_b,u_b,\theta_b)$ is the
background of Lemma~\ref{l:bg} with $0<\delta\le\delta_1(M)$, $X,Y$ are the
maps \eqref{eq:m-map}, and $(\rho,u,\theta)$ solves \eqref{eq:NSF} on
$\R^3\times[s,T_e]$ with $q$ in the class \eqref{eq:class},
$S\le s<T_e<\infty$ and $\mu_*\le\mu_0$; thus Lemmas~\ref{l:st}
and~\ref{l:ha} apply. In each lemma $\delta_1(M)$ and $\mu_0$ may be
decreased, depending only on $M$ and the fixed data, and $S$ may be
increased, depending also on the background. Every function defined below
by an evolution equation on $[s,T_e]$ takes its initial value at $s$. All
estimates hold on every interval $[s,t]$ with $s<t<T_e$, with constants
independent of $t$ and $T_e$, and we write
\[
 \Vert a\Vert_{L^2_w}^2=\int_s^tw_s(r)\Vert a(r)\Vert_{L^2}^2\dd r,
 \qquad N_o=N_o(s,t).
\]
The identities are written for smooth solutions decaying at infinity; for
the solutions above they follow by spatial regularization, with the
$H^{-1}$--$H^1$ pairing for the inertia terms and with \eqref{e:ha}.

\subsection{Energy identities for a carrier}

We choose a mixed velocity--temperature energy whose dissipation controls
$\nabla v$ and $\nabla^2z$; a screened test supplies the pressure dissipation.
A carrier is a triple $(\rho_c,\theta_c,u_c)$ of functions on
$\R^3\times[s,T_e]$, with $u_c$ vector valued, $\rho_c,\theta_c$ between two
positive constants and $\partial_t\rho_c+\Div(\rho_cu_c)=0$. We put
$p_c=R\rho_c\theta_c$, $D_c=\partial_t+u_c\cdot\nabla$, $J=(1-\Delta)^{-1}$
and
\begin{equation}\label{e:FC}
\begin{gathered}
 d_c=\frac{(\gamma-1)\kappa}{\gamma p_c},\qquad
 j_c=\rho_cd_c,\qquad
 k_c=\frac{d_c}{(\gamma-1)\theta_c}\Bigl(\nu+\frac\kappa{c_v}\Bigr),\\
 e_c(P,v,z)=\frac{P^2}{2\gamma p_c}+\frac{\rho_c|v|^2}2
   -j_cv\cdot\nabla z+\frac{k_c|\nabla z|^2}2 .
\end{gathered}
\end{equation}
The subscripts $N$ and $b$ in place of $c$ refer to the carriers
$(\rho,\theta,u)$ and $(\rho_b,\theta_b,u_be_1)$. Thus $p_N=p$, $D_N=D_u$,
$D_b$ is as in \eqref{e:der}, and, since $\gamma R=(\gamma-1)c_p$,
\[
 d_b=\frac c{\rho_b\theta_b},\qquad
 j_N=\frac c\theta,\qquad
 k_N=\frac{c(c_v\nu+\kappa)}{(\gamma-1)c_v\rho\theta^2}.
\]
For complex-valued functions, a product $ab$ of two factors in the quadratic
and bilinear expressions below stands for $\operatorname{Re}(a\bar b)$, and
$a^2$ for $|a|^2$.

\begin{lemma}\label{l:FC}
Let $(P,v,z)$ solve
\begin{equation}\label{e:fcsy}
\begin{cases}
 D_cP+\gamma p_c\Div v-(\gamma-1)\kappa\Delta z=f_p,\\
 \rho_cD_cv+\nabla P-\mathcal Lv=f_m,\\
 D_cz-\frac{\kappa}{c_v\rho_c}\Delta z+(\gamma-1)\theta_c\Div v=f_z,
\end{cases}
\end{equation}
put $\mathbf f=(f_p,f_m,f_z)$ and $\mathbb S(v)=2\mu\mathbb D(v)+\lambda(\Div v)\mathbb I$,
and let $B=B(x_1,t)$ be real with $B,\partial_tB,\partial_1B$ bounded. With
\begin{equation}\label{e:fcd}
 \mathcal D_c^B(v,z)=\int_{\R^3}B\Bigl[\mu|\nabla v|^2
   +\bigl(\mu+\lambda+(\gamma-1)j_c\theta_c\bigr)(\Div v)^2
   +\frac{\kappa k_c}{c_v\rho_c}(\Delta z)^2
   -2(\gamma-1)k_c\theta_c\Div v\,\Delta z\Bigr]\dd x,
\end{equation}
\begin{equation}\label{e:fcs}
 \mathscr S_c^B(P,v,z;\mathbf f)=\int_{\R^3}B\Bigl[(v-d_c\nabla z)\cdot f_m
   +\frac{Pf_p}{\gamma p_c}+\{\Div(j_cv)-\Div(k_c\nabla z)\}f_z\Bigr]\dd x,
\end{equation}
\begin{equation}\label{e:fcg}
\begin{aligned}
 \mathscr G_c^B(P,v,z)=\int_{\R^3}B\bigg[&
   \frac{\Div u_c-D_c\log p_c}{2\gamma p_c}P^2-P\nabla d_c\cdot\nabla z
   +\nu(\nabla d_c\cdot\nabla z)\Div v\\
 &+\mu\sum_{i,l}(\partial_ld_c)(\partial_lv_i-\partial_iv_l)\partial_iz
   +(D_c\log p_c)\bigl(j_cv\cdot\nabla z-k_c|\nabla z|^2\bigr)\\
 &+(\nabla j_c\cdot v-\nabla k_c\cdot\nabla z)
     \Bigl(\frac{\kappa}{c_v\rho_c}\Delta z-(\gamma-1)\theta_c\Div v\Bigr)\\
 &+\sum_{i,l}(\partial_iu_{c,l})(j_cv_i-k_c\partial_iz)\partial_lz
   \bigg]\dd x
\end{aligned}
\end{equation}
and
\[
\begin{aligned}
 n_c={}&P(v_1-d_c\partial_1z)-\bigl(\mathbb S(v)(v-d_c\nabla z)\bigr)_1
   +\mu\sum_{l=1}^3(v_l\partial_lv_1-v_1\partial_lv_l)
   -2\mu d_c\sum_{l=1}^3(\partial_lv_1\partial_lz-\partial_lv_l\partial_1z)\\
 &+(j_cv_1-k_c\partial_1z)
   \Bigl(\frac{\kappa}{c_v\rho_c}\Delta z-(\gamma-1)\theta_c\Div v+f_z\Bigr),
\end{aligned}
\]
one has
\begin{equation}\label{e:wFC}
 \frac{\dd}{\dd t}\int_{\R^3}Be_c(P,v,z)\dd x+\mathcal D_c^B(v,z)
 =\mathscr S_c^B(P,v,z;\mathbf f)+\mathscr G_c^B(P,v,z)
 +\int_{\R^3}\bigl[(D_cB)e_c(P,v,z)+(\partial_1B)n_c\bigr]\dd x .
\end{equation}
We write $\mathcal D_c$, $\mathscr S_c$, $\mathscr G_c$ for $B=1$; then
\begin{equation}\label{e:fcid}
 \frac{\dd}{\dd t}\int_{\R^3}e_c(P,v,z)\dd x+\mathcal D_c(v,z)
 =\mathscr S_c(P,v,z;\mathbf f)+\mathscr G_c(P,v,z).
\end{equation}
Moreover
\[
\begin{gathered}
 e_c(P,v,z)=\frac12\Bigl[\frac{P^2}{\gamma p_c}+\rho_c|v-d_c\nabla z|^2
   +(k_c-\rho_cd_c^2)|\nabla z|^2\Bigr],\qquad
 k_c-\rho_cd_c^2=\frac{d_c}{(\gamma-1)\theta_c}\Bigl(\nu+\frac\kappa{c_p}\Bigr),\\
 \mathcal D_c(v,z)=\int_{\R^3}\Bigl[\mu|\nabla\times v|^2+\frac\nu\gamma(\Div v)^2
   +\frac{\kappa k_c}{c_v\rho_c}\Bigl(\Delta z
     -\frac{(\gamma-1)c_v\rho_c\theta_c}{\kappa}\Div v\Bigr)^2\Bigr]\dd x,
\end{gathered}
\]
so that, with $C$ depending only on the fixed data and the bounds of
$\rho_c,\theta_c$,
\[
\begin{gathered}
 \frac{P^2+|v|^2+|\nabla z|^2}C\le e_c(P,v,z)\le C(P^2+|v|^2+|\nabla z|^2),\\
 \frac{\Vert\nabla v\Vert_{L^2}^2+\Vert\nabla^2z\Vert_{L^2}^2}C
 \le\mathcal D_c(v,z)\le C\bigl(\Vert\nabla v\Vert_{L^2}^2+\Vert\nabla^2z\Vert_{L^2}^2\bigr).
\end{gathered}
\]
With $[D_c,\nabla J]P=u_c\cdot\nabla(\nabla JP)-\nabla J(u_c\cdot\nabla P)$ and
\begin{equation}\label{e:fck}
\begin{aligned}
 \mathscr K_c(P,v;f_p,f_m)&=\int_{\R^3}\bigl(f_m\cdot\nabla JP
   +\rho_cv\cdot\nabla Jf_p\bigr)\dd x,\\
 \mathscr P_c(P,v,z)&=\int_{\R^3}\Bigl(\mathcal Lv\cdot\nabla JP
   +\rho_cv\cdot\nabla J\bigl[(\gamma-1)\kappa\Delta z-\gamma p_c\Div v\bigr]
   +\rho_cv\cdot[D_c,\nabla J]P\Bigr)\dd x,
\end{aligned}
\end{equation}
one has
\begin{equation}\label{e:fcsc}
 \frac{\dd}{\dd t}\int_{\R^3}\rho_cv\cdot\nabla JP\dd x
 +\Vert\nabla J^{1/2}P\Vert_{L^2}^2
 =\mathscr P_c(P,v,z)+\mathscr K_c(P,v;f_p,f_m),
\end{equation}
\begin{equation}\label{e:fcwv}
 \frac{\dd}{\dd t}\Bigl(\frac1{2\tau}\int_{\R^3}\rho_cz^2\dd x\Bigr)
 +\frac1{2\tau^2}\int_{\R^3}\rho_cz^2\dd x
 +\frac{\kappa}{c_v\tau}\Vert\nabla z\Vert_{L^2}^2
 =\frac1\tau\int_{\R^3}\rho_cz\bigl(f_z-(\gamma-1)\theta_c\Div v\bigr)\dd x .
\end{equation}
For a spatial multi-index $\alpha$, $(\partial^\alpha P,\partial^\alpha v,\partial^\alpha z)$
solves \eqref{e:fcsy} with the sources
\begin{equation}\label{e:fcal}
\begin{aligned}
 f_p^\alpha&=\partial^\alpha f_p-[\partial^\alpha,u_c\cdot\nabla]P
   -\gamma[\partial^\alpha,p_c]\Div v,\\
 f_m^\alpha&=\partial^\alpha f_m-[\partial^\alpha,\rho_c]D_cv
   -\rho_c[\partial^\alpha,u_c\cdot\nabla]v,\\
 f_z^\alpha&=\partial^\alpha f_z-[\partial^\alpha,u_c\cdot\nabla]z
   +\Bigl[\partial^\alpha,\frac{\kappa}{c_v\rho_c}\Bigr]\Delta z
   -(\gamma-1)[\partial^\alpha,\theta_c]\Div v,
\end{aligned}
\end{equation}
where $[A,G]h=A(Gh)-G(Ah)$ and a function denotes the multiplication by it.
For $c=b$ all these identities also hold for the Fourier transforms in $x'$,
with $\nabla$, $\Delta$ and $J$ replaced by $(\partial_1,i\xi_2,i\xi_3)$,
$\partial_1^2-|\xi'|^2$ and $(1-\partial_1^2+|\xi'|^2)^{-1}$; the constants
of the two equivalences do not depend on $\xi'$.
\end{lemma}

\begin{proof}
Put $\mathcal Ta=D_ca+(\Div u_c)a$. The equation of $\rho_c$ gives
$\mathcal T(\rho_ca)=\rho_cD_ca$ and
\[
 \frac{\dd}{\dd t}\int_{\R^3}Ba\dd x
 =\int_{\R^3}\bigl[(D_cB)a+B\mathcal Ta\bigr]\dd x .
\]
By \eqref{e:FC}, $d_c$ and $\frac{k_c}{\rho_c}$ are constant multiples of
$p_c^{-1}$ and $p_c^{-2}$; hence $D_cd_c=-d_cD_c\log p_c$ and
$D_c\frac{k_c}{\rho_c}=-\frac{2k_c}{\rho_c}D_c\log p_c$. Moreover
$D_c\partial_lz=\partial_lD_cz-\sum_i(\partial_lu_{c,i})\partial_iz$. Writing
each term of $e_c$ as $\rho_c$ times a function and using the first two rows
of \eqref{e:fcsy}, we obtain
\[
\begin{aligned}
 \mathcal T\frac{P^2}{2\gamma p_c}
 ={}&\frac{Pf_p}{\gamma p_c}-P\Div v+d_cP\Delta z
   +\frac{\Div u_c-D_c\log p_c}{2\gamma p_c}P^2,\\
 \mathcal T\frac{\rho_c|v|^2}2={}&v\cdot(f_m-\nabla P+\mathcal Lv),\\
 \mathcal T(-j_cv\cdot\nabla z)
 ={}&-d_c\nabla z\cdot(f_m-\nabla P+\mathcal Lv)-j_cv\cdot\nabla D_cz
   +j_c(D_c\log p_c)v\cdot\nabla z
   +j_c\sum_{i,l}v_i(\partial_iu_{c,l})\partial_lz,\\
 \mathcal T\frac{k_c|\nabla z|^2}2
 ={}&k_c\nabla z\cdot\nabla D_cz-k_c(D_c\log p_c)|\nabla z|^2
   -k_c\sum_{i,l}\partial_iz(\partial_iu_{c,l})\partial_lz .
\end{aligned}
\]
In the sum we use
\[
\begin{gathered}
 -P\Div v-v\cdot\nabla P+d_cP\Delta z+d_c\nabla z\cdot\nabla P
 =-\Div\bigl(P(v-d_c\nabla z)\bigr)-P\nabla d_c\cdot\nabla z,\\
 (k_c\nabla z-j_cv)\cdot\nabla D_cz
 =\Div\bigl((k_c\nabla z-j_cv)D_cz\bigr)
   +\{\Div(j_cv)-\Div(k_c\nabla z)\}D_cz,
\end{gathered}
\]
and, since $\Div\mathbb S(v)=\mathcal Lv$ and $\mathbb S(v)$ is symmetric,
\[
\begin{aligned}
 (v-d_c\nabla z)\cdot\mathcal Lv
 ={}&\Div\bigl(\mathbb S(v)(v-d_c\nabla z)\bigr)
   -\mathbb S(v):\nabla v+\mathbb S(v):\nabla(d_c\nabla z),\\
 \mathbb S(v):\nabla v={}&\mu|\nabla v|^2+(\mu+\lambda)(\Div v)^2+\mu\Div\Lambda,\\
 \mathbb S(v):\nabla(d_c\nabla z)
 ={}&\nu d_c\Div v\,\Delta z+\nu(\nabla d_c\cdot\nabla z)\Div v
   +\mu\sum_{i,l}(\partial_ld_c)(\partial_lv_i-\partial_iv_l)\partial_iz
   +2\mu\Div(d_c\Gamma),
\end{aligned}
\]
where $\Lambda_i=\sum_l(v_l\partial_lv_i-v_i\partial_lv_l)$ and
$\Gamma_i=\sum_l(\partial_lv_i\partial_lz-\partial_lv_l\partial_iz)$; the
last two lines follow from
$\Div\Lambda=\sum_{i,l}\partial_iv_l\partial_lv_i-(\Div v)^2$ and
$\Div\Gamma=\sum_{i,l}\partial_lv_i\partial_{il}z-\Div v\,\Delta z$.
Finally substitute
$D_cz=f_z+\frac{\kappa}{c_v\rho_c}\Delta z-(\gamma-1)\theta_c\Div v$ and
$\Div(j_cv)-\Div(k_c\nabla z)=j_c\Div v-k_c\Delta z+\nabla j_c\cdot v-\nabla k_c\cdot\nabla z$.
The products $\Div v\,\Delta z$ then carry the coefficient
\[
 \nu d_c+\frac{\kappa j_c}{c_v\rho_c}+(\gamma-1)\theta_ck_c
 =d_c\Bigl(\nu+\frac\kappa{c_v}\Bigr)+(\gamma-1)\theta_ck_c
 =2(\gamma-1)\theta_ck_c
\]
by the definition of $k_c$, and $(\Div v)^2$, $(\Delta z)^2$ receive
$-(\gamma-1)\theta_cj_c$ and $-\frac{\kappa k_c}{c_v\rho_c}$. Collecting all
terms, $\mathcal Te_c(P,v,z)$ equals the integrand of
$-\mathcal D_c+\mathscr S_c+\mathscr G_c$ plus
\[
 \Div\Bigl[-P(v-d_c\nabla z)+\mathbb S(v)(v-d_c\nabla z)-\mu\Lambda
   +2\mu d_c\Gamma+(k_c\nabla z-j_cv)D_cz\Bigr].
\]
The first component of this flux is $-n_c$; multiplication by $B$ and integration give
\eqref{e:wFC}.

Completing the square in $e_c$ gives its stated formula, since
$R=(\gamma-1)c_v$ and $c_p=\gamma c_v$. For the dissipation,
\[
 \nu+(\gamma-1)j_c\theta_c-\frac{(\gamma-1)^2c_v\rho_c\theta_c^2k_c}{\kappa}
 =\nu+\frac{\gamma-1}\gamma\frac\kappa{c_v}
 -\frac{\gamma-1}\gamma\Bigl(\nu+\frac\kappa{c_v}\Bigr)=\frac\nu\gamma>0.
\]
Together with
$\int|\nabla v|^2=\int(|\nabla\times v|^2+(\Div v)^2)$, this gives the
square formula for $\mathcal D_c$. The two equivalences follow from these
uniformly positive quadratic forms and
$\Vert\Delta z\Vert_{L^2}=\Vert\nabla^2z\Vert_{L^2}$.

For \eqref{e:fcsc}, the equation of $\rho_c$ and $D_c\nabla J=\nabla JD_c+[D_c,\nabla J]$ give
\[
\begin{aligned}
 \frac{\dd}{\dd t}\int_{\R^3}\rho_cv\cdot\nabla JP\dd x
 ={}&\int_{\R^3}\Bigl[(f_m-\nabla P+\mathcal Lv)\cdot\nabla JP
   +\rho_cv\cdot[D_c,\nabla J]P\Bigr]\dd x\\
 &+\int_{\R^3}\rho_cv\cdot\nabla J\bigl(f_p-\gamma p_c\Div v
   +(\gamma-1)\kappa\Delta z\bigr)\dd x,
\end{aligned}
\]
and $\int_{\R^3}\nabla P\cdot\nabla JP\dd x=\Vert\nabla J^{1/2}P\Vert_{L^2}^2$
by Plancherel's identity. Testing the third row of \eqref{e:fcsy} by
$\frac{\rho_cz}\tau$ gives \eqref{e:fcwv}, since $\frac{\dd\tau}{\dd t}=1$.
Applying $\partial^\alpha$ to \eqref{e:fcsy} gives \eqref{e:fcal}. For $c=b$
the coefficients and $B$ depend only on $(x_1,t)$, so the Fourier transform
in $x'$ commutes with them; the same computation with
$\partial_j=i\xi_j$, $j=2,3$, and real parts gives the identities on each
fiber, and the two squares are pointwise in $\xi'$.
\end{proof}

For $c=b$, \eqref{e:fcsy} reads
\begin{equation}\label{e:forced}
\begin{cases}
 D_bP+\gamma p_b\Div v-(\gamma-1)\kappa\Delta z=f_p,\\
 \rho_bD_bv+\nabla P-\mathcal Lv=f_m,\\
 D_bz-\frac{\kappa}{c_v\rho_b}\Delta z+(\gamma-1)\theta_b\Div v=f_z,
\end{cases}
\end{equation}
and $\Div(u_be_1)=h_b$, $D_b\log p_b=\frac{B_b}{\theta_b}-h_b$ by \eqref{e:bc}
and \eqref{e:bal}. Let $b=u_b-u_*$.

\begin{lemma}\label{l:bessel}
The operator $J$ is the convolution with
$K_J(x)=\frac{e^{-|x|}}{4\pi|x|}$; $\Vert K_J\Vert_{L^1}=1$,
$\nabla K_J\in L^1$, and
$|\nabla^2K_J(x)|+|x||\nabla^3K_J(x)|\le C|x|^{-3}e^{-\frac{|x|}2}$ for $x\ne0$.
For $h,g\in L^2$ and $1\le i\le3$,
\begin{equation}\label{e:Jsym}
\begin{gathered}
 \Vert\nabla Jh\Vert_{L^2}\le\min\Bigl\{\frac12\Vert h\Vert_{L^2},
   \Vert\nabla h\Vert_{L^2},\Vert h\Vert_{H^{-1}}\Bigr\},\qquad
 \Vert\nabla^2Jh\Vert_{L^2}\le\Vert\nabla J^{1/2}h\Vert_{L^2}
   \le\min\{\Vert h\Vert_{L^2},\Vert\nabla h\Vert_{L^2}\},\\
 \Vert\nabla J\partial_ih\Vert_{L^2}\le\Vert h\Vert_{L^2},\qquad
 \Vert\nabla J\Div g\Vert_{L^2}\le\Vert g\Vert_{L^2}.
\end{gathered}
\end{equation}
For $a=a(x_1,t)$ with $\partial_1a$ bounded,
\begin{equation}\label{e:Jcom}
 \Vert a\partial_1\nabla Jh-\nabla J(a\partial_1h)\Vert_{L^2}
 \le C\Vert\partial_1a\Vert_{L^\infty}
   \min\{\Vert h\Vert_{L^2},\Vert\partial_1h\Vert_{L^2}\}.
\end{equation}
For $1\le i,j\le3$, with $\Vert\cdot\Vert_{H^4_{\rm uloc}}$ defined as in
\eqref{e:ul},
\begin{equation}\label{e:Jloc}
 \Vert Jh\Vert_{L^3}+\Vert\nabla Jh\Vert_{L^3}+\Vert\partial_{ij}Jh\Vert_{L^3}
 \le C\Vert h\Vert_{L^3},
 \qquad
 \Vert Jh\Vert_{H^4_{\rm uloc}}+\sum_{l=0}^2\Vert\nabla^lJh\Vert_{L^\infty}
 \le C\Vert h\Vert_{H^2_{\rm uloc}} .
\end{equation}
For $k\in\R^2$, the operator $(1-\partial_1^2+|k|^2)^{-1}$ on $L^2(\R)$ has
the kernel $\frac{e^{-m_k|x_1-y_1|}}{2m_k}$, $m_k=\sqrt{1+|k|^2}$; for
measurable $A_1,A_2\subset\R$ at distance $d\ge0$, the operator
$h\mapsto\mathbf 1_{A_1}(1-\partial_1^2+|k|^2)^{-1}(\mathbf 1_{A_2}h)$ has
norm at most $e^{-d}$ on $L^2(\R)$.
\end{lemma}

\begin{proof}
$K_J$ is the radial solution of $(1-\Delta)K_J=\delta_0$ decaying at
infinity; $\int_{\R^3}K_J\dd x=\int_0^\infty re^{-r}\dd r=1$, and
$|\partial^\beta K_J(x)|\le C_\beta(|x|^{-1-|\beta|}+|x|^{-1})e^{-|x|}$ for
$x\ne0$ gives the kernel bounds. The symbols of $\nabla J$, $\nabla^2J$,
$\nabla J^{1/2}$, $\nabla J\partial_i$ and $\nabla J\Div$ are bounded in
modulus by $\frac{|\xi|}{1+|\xi|^2}$, $\frac{|\xi|^2}{1+|\xi|^2}$,
$\frac{|\xi|}{(1+|\xi|^2)^{1/2}}$ and $\frac{|\xi|^2}{1+|\xi|^2}$, and
$\frac{|\xi|}{1+|\xi|^2}\le\min\{\frac12,|\xi|,(1+|\xi|^2)^{-1/2}\}$,
$\frac{|\xi|^2}{1+|\xi|^2}\le\frac{|\xi|}{(1+|\xi|^2)^{1/2}}\le\min\{1,|\xi|\}$;
Plancherel's identity gives \eqref{e:Jsym}.

Since $\partial_1$ commutes with $\nabla J$, the operator in \eqref{e:Jcom}
is $[a,\nabla J](\partial_1h)$, whose kernel $(a(x_1)-a(y_1))\nabla K_J(x-y)$
is bounded by $\Vert\partial_1a\Vert_{L^\infty}|x-y||\nabla K_J(x-y)|$; this
is an integrable function of $x-y$, and Young's inequality gives the bound
by $\Vert\partial_1h\Vert_{L^2}$. Moreover
$[a,\nabla J]\partial_1h=[a,\nabla J\partial_1]h+\nabla J((\partial_1a)h)$,
and an integration by parts in $y_1$ shows that $[a,\nabla J\partial_1]$ has
the kernel $(a(x_1)-a(y_1))\nabla\partial_1K_J(x-y)$, $x\ne y$, bounded by
$C\Vert\partial_1a\Vert_{L^\infty}|x-y|^{-2}e^{-\frac{|x-y|}2}$. Young's inequality
and \eqref{e:Jsym} give the bound by $\Vert h\Vert_{L^2}$.

Young's inequality gives the $L^3$ bounds for $J$ and $\nabla J$.
As an operator, $\partial_{ij}J$ is a constant multiple of $I$ plus
$\operatorname{p.v.}(\partial_{ij}K_J)*$.
The $L^2$ bound \eqref{e:Jsym} and the kernel estimates
$|\partial_{ij}K_J(x)|\le C|x|^{-3}$, $|\nabla\partial_{ij}K_J(x)|\le C|x|^{-4}$
give its $L^3$ boundedness by the Calder\'on--Zygmund theorem
\cite{Stein1970}. For the second bound in
\eqref{e:Jloc}, let $\{\varphi_n\}_{n\in\mathbb Z^3}$ be a smooth partition
of unity with $\operatorname{supp}\varphi_n\subset n+(-1,1)^3$ and bounded
derivatives, and let $Q=x_0+(-1,1)^3$. For $|n-x_0|_\infty\le4$,
$\Vert J(\varphi_nh)\Vert_{H^4}\le\Vert\varphi_nh\Vert_{H^2}
\le C\Vert h\Vert_{H^2_{\rm uloc}}$; for $|n-x_0|_\infty>4$ the distance of
$Q$ to $\operatorname{supp}\varphi_n$ is at least $\frac12|n-x_0|_\infty$,
and $|\partial^\beta K_J(z)|\le C_\beta e^{-|z|}$ for $|z|\ge1$ gives
$\Vert J(\varphi_nh)\Vert_{W^{4,\infty}(Q)}
\le Ce^{-\frac{|n-x_0|_\infty}2}\Vert h\Vert_{H^2_{\rm uloc}}$. Summation over $n$
and the embedding $H^4(Q)\subset W^{2,\infty}(Q)$ prove \eqref{e:Jloc}.

The kernel of $(1-\partial_1^2+|k|^2)^{-1}$ is the decaying solution of
$(m_k^2-\partial_1^2)G=\delta_{y_1}$. For $x_1\in A_1$,
$\int_{A_2}\frac{e^{-m_k|x_1-y_1|}}{2m_k}\dd y_1
\le\int_{|w|\ge d}\frac{e^{-m_k|w|}}{2m_k}\dd w=\frac{e^{-m_kd}}{m_k^2}\le e^{-d}$,
and the same holds with $A_1,A_2$ exchanged; Schur's test gives the last
assertion.
\end{proof}

\subsection{The entropy history and the contact part}

With $m$ of \eqref{e:char}, $h_*$ of \eqref{e:mass} and $g_T$ of
\eqref{e:bc}, put
\begin{equation}\label{e:ent}
 \chi=\frac1{\theta_b}\Bigl(\zeta-\frac{\Pi}{c_p\rho}\Bigr),\qquad
 f_h=\Bigl(\log\frac{m}{V_c(y+h_*,t)}\Bigr)\circ Y,
\end{equation}
and let $\eta$ and $f$ solve
\begin{equation}\label{e:hist}
\begin{aligned}
 \rho D_u\eta-\frac c{\theta_b}\Delta(\theta_b\eta)&=S_\eta,
 &\eta(s)&=\chi(s),\\
 \rho D_uf-\frac c{\theta_b}\Delta(\theta_bf)
 &=-\rho_bG_c\psi_1,&f(s)&=0,
\end{aligned}
\end{equation}
where $G_c=\partial_1\{\log V_c(Y+h_*,t)\}=\rho_b(\partial_y\log V_c)(Y+h_*,t)$, with
$\log V_c(Y+h_*,t)$ the function $x\mapsto\log V_c(Y(x_1,t)+h_*,t)$, and
\begin{equation}\label{e:Fn}
\begin{gathered}
 S_\eta=\frac c{\theta_b}\Delta\Bigl(\frac{\Pi}{c_p\rho}\Bigr)
 -\frac{\kappa\partial_1^2\theta_b+\nu h_b^2}{c_p\theta_b}\frac\Pi{p_b}
 -\rho_b\psi_1\partial_1f_h+\frac{H[\psi]}{c_p\theta_b}+F_n,\\
 F_n=\frac{Q_\mu(\nabla\psi)-\Pi\Div\psi}{c_p\theta_b}
 -\frac{\rho_b}{\gamma}\frac\Pi{p_b}g_T\psi_1 .
\end{gathered}
\end{equation}
Let $\omega_c$ be a smooth even function with $0\le\omega_c\le1$,
$\omega_c=0$ on $[-2,2]$ and $\omega_c=1$ outside $[-3,3]$, fix
$0<c_0<\frac18\min\{c_-,c_+\}$ with $c_\pm$ of \eqref{e:mass}, and put
\begin{equation}\label{e:cut}
 \chi_c=1-\omega_c\Bigl(\frac{Y+h_*}{c_0\tau}\Bigr),\quad
 F=\chi_cf,\quad f_{\rm out}=(1-\chi_c)f,\quad
 V=\psi-b\eta e_1,\quad \Phi(s,t)=N_0(\eta;s,t)^2 .
\end{equation}
For $\epsilon_K,\theta_W\in(0,1]$, which are fixed in Lemma~\ref{l:mr},
put
\begin{equation}\label{e:test}
\begin{aligned}
 W^0&=V-d_N\nabla\zeta+\epsilon_K\nabla J\Pi,&
 Z^0&=\Div(j_NV)-\Div(k_N\nabla\zeta),\\
 W^i&=\partial_i\psi-d_N\nabla\partial_i\zeta+\epsilon_K\nabla J\partial_i\Pi,&
 Z^i&=\Div(j_N\partial_i\psi)-\Div(k_N\nabla\partial_i\zeta)
 \quad(1\le i\le3),
\end{aligned}
\end{equation}
and put
\begin{equation}\label{e:zN}
 z_N=\frac1{c_p\theta_b}\Bigl(\frac\Pi\rho-\frac{\Pi}{\rho_b}\Bigr)
 =-\frac{(\gamma-1)\phi\Pi}{\gamma\rho p_b},
\end{equation}
where the last identity follows from
$\frac1\rho-\frac1{\rho_b}=-\frac\phi{\rho\rho_b}$, $p_b=R\rho_b\theta_b$
and $R/c_p=(\gamma-1)/\gamma$.

Put
\begin{equation}\label{e:Geta}
 G_\eta=\mathcal L(b\eta e_1)-\rho bD_u\eta\,e_1-\rho h_b\psi_1\eta\,e_1.
\end{equation}
Since $\rho=\frac p{R\theta}$,
\begin{equation}\label{e:rhoD}
 \nabla\rho=\frac{\nabla\Pi+(\partial_1p_b)e_1}{R\theta}
 -\frac\rho\theta\bigl(\nabla\zeta+(\partial_1\theta_b)e_1\bigr),\qquad
 \nabla\Bigl(\frac{\Pi}{c_p\rho}\Bigr)
 =\frac{\nabla\Pi}{c_p\rho}-\frac{\Pi}{c_p\rho^2}\nabla\rho,
\end{equation}
and \eqref{e:hist}, \eqref{e:Lam} and $D_u=D_b+\psi\cdot\nabla$ give
\begin{equation}\label{e:etaD}
 D_b\eta=\frac{c}{\rho\theta_b}\Delta\Bigl(\theta_b\eta+\frac{\Pi}{c_p\rho}\Bigr)
 -\frac{\kappa\partial_1^2\theta_b+\nu h_b^2}{c_p\rho\theta_b}\frac\Pi{p_b}
 -\frac{\rho_b}\rho\psi_1\partial_1f_h
 +\frac{H[\psi]}{c_p\rho\theta_b}+\frac{F_n}\rho-\psi\cdot\nabla\eta .
\end{equation}
By \eqref{e:Um} and \eqref{e:CH}, $b\circ X=b_m=H_c(y+h_*,t)
+\sum_\sigma s_\sigma\vartheta_\sigma+(R_U)_2$. Lemma~\ref{l:asy},
Lemma~\ref{l:pr} with $|m_\pm|\le C_M$ and \eqref{eq:pr-gau} give, for
$t\ge S$,
\begin{equation}\label{e:bsup}
 \sqrt\tau\,\Vert b\Vert_{L^\infty(\R)}\le C_M .
\end{equation}
Below, $\langle\cdot,\cdot\rangle$ is the pairing of $H^{-1}$ and $H^1$.

\begin{lemma}\label{l:ent}
The equations \eqref{e:hist} have unique solutions
$\eta,f\in C([s,T_e];L^2)\cap L^2(s,T_e;H^1)$, and
\begin{equation}\label{e:espl}
 \chi=\eta+f,\qquad
 \rho(1+\eta+f)=\rho_b\Bigl(1+\frac\Pi{\gamma p_b}\Bigr),\qquad
 \phi=\frac{\rho_b\Pi}{\gamma p_b}-\rho(\eta+f),
\end{equation}
\begin{equation}\label{e:etaE}
 \frac12\frac{\dd}{\dd t}\int_{\R^3}\rho\eta^2\dd x
 +c\Vert\nabla\eta\Vert_{L^2}^2
 =c\int_{\R^3}g_T^2\eta^2\dd x+\langle S_\eta,\eta\rangle .
\end{equation}
With $\epsilon_f(s)=\sup_{r\ge s}\sqrt{1+r}\,\Vert f_h(r)\Vert_{L^\infty(\R)}$,
$\lim_{s\to\infty}\epsilon_f(s)=0$. Moreover,
\begin{equation}\label{e:eta}
\begin{aligned}
 \Phi(s,t)&\le C_M\bigl(I_s+N_o^2\bigr),\\
 \Vert b\eta\Vert_{L^2_w}^2
 &\le\bigl(C_M\delta^2+\omega_b(s)\bigr)\bigl(I_s+N_o^2\bigr),\\
 N_0(b\eta;s,t)^2&\le\omega_b(s)\Phi(s,t),
\end{aligned}
\end{equation}
\begin{equation}\label{e:Seta}
 \int_s^tw_s|\langle S_\eta,\eta\rangle|\dd r
 \le\frac c4\Vert\nabla\eta\Vert_{L^2_w}^2
 +C_M\Vert\nabla\Pi\Vert_{L^2_w}^2
 +\bigl[C_M(\delta+\mu_*)+\omega_b(s)\bigr]\bigl(N_o^2+\Phi\bigr).
\end{equation}
There are $\lambda_0,k_0>0$, depending only on $M$ and the fixed data,
such that, with $S$ independent of $L_0$, for every $L_0\ge0$,
$\sigma=\pm1$ and $\varpi_\sigma=\exp\{\lambda_0(\sigma(Y+h_*)-c_0\tau)\}$,
\begin{equation}\label{e:ext}
\begin{aligned}
 &\sup_{s\le r\le t}w_s\tau^{L_0}\int_{\R^3}\varpi_\sigma f^2\dd x
 +\int_s^tw_s\tau^{L_0}\int_{\R^3}\varpi_\sigma\bigl(f^2+|\nabla f|^2\bigr)\dd x\dd r
 \le C_{M,L_0}\delta^2e^{-k_0s}N_0(\psi;s,t)^2,\\
 &\sup_{s\le r\le t}w_s\tau^{L_0}\Vert f_{\rm out}\Vert_{L^2}^2
 +\int_s^tw_s\tau^{L_0}\Vert f_{\rm out}\Vert_{H^1}^2\dd r
 \le C_{M,L_0}\delta^2e^{-k_0s}N_0(\psi;s,t)^2 .
\end{aligned}
\end{equation}
Further,
\begin{equation}\label{e:ehp}
\begin{gathered}
 |z_N|\le C_M\mu_*|\Pi|,\qquad
 \Vert\nabla z_N\Vert_{L^2_w}\le C_M\mu_*N_o,\qquad
 \Bigl\Vert\nabla\Bigl(\frac{\Pi}{c_p\rho}\Bigr)\Bigr\Vert_{L^2_w}\le C_MN_o,\\
 \Bigl\Vert\nabla\Bigl(\frac{\Pi}{c_p\rho}\Bigr)-\frac{\nabla\Pi}{c_p\rho}\Bigr\Vert_{L^2_w}
 \le\bigl[C_M(\delta+\mu_*)+\omega_b(s)\bigr]N_o,\qquad
 \int_s^tw_s\Vert F_n\Vert_{L^{6/5}}^2\dd r\le C_M\mu_*^2N_o^2,
\end{gathered}
\end{equation}
\begin{equation}\label{e:etst}
\begin{gathered}
 N_0(W^0;s,t)^2+\Vert Z^0\Vert_{L^2_w}^2+\Vert\Div(\rho V)\Vert_{L^2_w}^2
 \le C_M\bigl(I_s+N_o^2+\Phi\bigr),\\
 \sum_{i=1}^3\bigl(\Vert W^i\Vert_{L^2_w}^2+\Vert\nabla W^i\Vert_{L^2_w}^2
 +\Vert Z^i\Vert_{L^2_w}^2\bigr)\le C_MN_o^2,\\
 \Vert(\nabla\rho)\cdot V\Vert_{L^2_w}
 \le\bigl[C_M(\delta+\mu_*)+\omega_b(s)\bigr]\bigl(N_o+\Phi^{\frac12}\bigr),
\end{gathered}
\end{equation}
and, for every vector function $T$ with $N_0(T;s,t)<\infty$,
\begin{equation}\label{e:Getb}
 \int_s^tw_s|\langle G_\eta,T\rangle|\dd r
 \le\omega_b(s)\bigl(N_o^2+\Phi\bigr)^{\frac12}N_0(T;s,t).
\end{equation}
With
\begin{equation}\label{e:Reta}
\begin{aligned}
 \mathcal R_\eta={}&\langle G_\eta,W^0\rangle
 +\int_{\R^3}\Bigl[(\partial_1\Pi)b\eta
 +Z^0\bigl\{B_b\eta-(\gamma-1)\theta\partial_1(b\eta)\bigr\}
 +\frac{\theta_W}\tau\rho\zeta B_b\eta\Bigr]\dd x\\
 &-\epsilon_K\int_{\R^3}\rho V\cdot\nabla J\{\gamma p\partial_1(b\eta)\}\dd x
 +\sum_{i=1}^3\int_{\R^3}\bigl[\rho(W^i)_1\partial_i(A_b\eta)
 +Z^i\partial_i(B_b\eta)\bigr]\dd x,\\
 \mathcal R_{\rm out}={}&\int_{\R^3}\Bigl[\rho A_bf_{\rm out}(W^0)_1
 +B_bf_{\rm out}Z^0+\frac{\theta_W}\tau\rho\zeta B_bf_{\rm out}\Bigr]\dd x\\
 &+\sum_{i=1}^3\int_{\R^3}\bigl[\rho(W^i)_1\partial_i(A_bf_{\rm out})
 +Z^i\partial_i(B_bf_{\rm out})\bigr]\dd x,
\end{aligned}
\end{equation}
one has
\begin{equation}\label{e:eret}
 \int_s^tw_s\bigl(|\mathcal R_\eta|+|\mathcal R_{\rm out}|\bigr)\dd r
 \le\bigl(C_M\delta+\omega_b(s)\bigr)\bigl(I_s+N_o^2+\Phi\bigr).
\end{equation}
\end{lemma}

\begin{proof}
For $v>0$ and $E_0>\frac{b_0^2}2$ put
$\mathfrak n(v,b_0,E_0)=\frac{\gamma-1}\gamma\log v
+\frac1\gamma\log\frac{(\gamma-1)(E_0-\frac{b_0^2}2)}{p_*}$. By \eqref{e:char}
and $P_b=\frac{(\gamma-1)(E_m-\frac{b_m^2}2)}{V_b}$, $\log m=\mathfrak n(\mathcal U)$;
by \eqref{eq:prof}, $(\gamma-1)c_v\Theta_c=p_*V_c$, so
$\mathfrak n(\widehat{\mathcal U})=\log V_c+\frac1\gamma
\log\bigl(1-\frac{(\gamma-1)H_c^2}{2p_*V_c}\bigr)$. With
$\widehat{\mathcal U}^h=\widehat{\mathcal U}(y+h_*,t)$, \eqref{e:CH} gives
\[
 f_h\circ X=\mathfrak n\Bigl(\widehat{\mathcal U}^h
 +\sum_\sigma\vartheta_\sigma r_\sigma+R_U\Bigr)
 -\mathfrak n(\widehat{\mathcal U}^h)
 +\frac1\gamma\log\Bigl(1-\frac{(\gamma-1)H_c^2}{2p_*V_c}\Bigr)(y+h_*,t).
\]
At $\mathcal U_{\sigma,*}=(v_\sigma,0,\frac{p_*v_\sigma}{\gamma-1})$,
$D\mathfrak n(\mathcal U_{\sigma,*})r_\sigma
=-\frac{\gamma-1}{\gamma v_\sigma}+\frac{\gamma-1}{\gamma v_\sigma}=0$. For $t\ge S$,
Lemmas~\ref{l:asy} and~\ref{l:pr} and \eqref{eq:ct-gau} place the segment
between $\widehat{\mathcal U}^h$ and $\mathcal U$ in the convex set
$\{\frac12\min v_\pm\le v\le2\max v_\pm,\ |b_0|\le1,\
\frac{c_v}2\min\theta_\pm\le E_0-\frac{b_0^2}2,\ E_0\le C\}$, on which
$\mathfrak n$ and its first two derivatives are bounded by constants of the
fixed data. Taylor's formula, $|D\mathfrak n(\widehat{\mathcal U}^h)r_\sigma|
\le C|\widehat{\mathcal U}^h-\mathcal U_{\sigma,*}|$ and $|H_c|\le C\delta\tau^{-\frac12}$
give
\[
 |f_h\circ X|\le C\sum_\sigma|\vartheta_\sigma|
 |\widehat{\mathcal U}^h-\mathcal U_{\sigma,*}|
 +C|R_U|+C_M\tau^{-1}.
\]
By \eqref{eq:ct-gau}, $|\widehat{\mathcal U}^h-\mathcal U_{\sigma,*}|\le
C\delta e^{-\frac{c_G(y+h_*)^2}\tau}$ if $\sigma(y+h_*)\ge0$, and $\le C\delta$
otherwise; by Lemma~\ref{l:pr},
$|\vartheta_\sigma|\le C_M\tau^{-\frac12}e^{-\frac{(y-s_\sigma t)^2}{8d_\sigma\tau}}$.
If $|y-s_\sigma t|\le\frac{|s_\sigma|t}2$, then
$\sigma(y+h_*)\ge\frac{|s_\sigma|t}2-|h_*|$. Hence
$|\vartheta_\sigma||\widehat{\mathcal U}^h-\mathcal U_{\sigma,*}|\le C_be^{-c_bt}$,
\[
 \Vert f_h(t)\Vert_{L^\infty}\le C_M\Vert R_U(t)\Vert_{L^\infty}
 +C_M\tau^{-1}+C_be^{-c_bt},
\]
and Lemma~\ref{l:asy} gives $\epsilon_f(s)\to0$.

Subtracting from $c_p\rho D_u\theta-D_up=\kappa\Delta\theta+Q_\mu(\nabla u)$,
which follows from the four equations displayed after \eqref{e:peq3}, its
background form and using $Q_\mu(\nabla u)=\nu h_b^2+H[\psi]+Q_\mu(\nabla\psi)$,
we obtain
\[
 \rho D_u\zeta-\frac1{c_p}D_u\Pi
 =c\Delta\zeta+\frac{H[\psi]+Q_\mu(\nabla\psi)}{c_p}
 -\phi B_b-\rho(\partial_1\theta_b)\psi_1+\frac{\partial_1p_b}{c_p}\psi_1 .
\]
We substitute $\theta_b\chi=\zeta-\Pi/(c_p\rho)$ and use $\rho D_u\rho^{-1}=h_b+\Div\psi$,
$\Pi=R(\theta_b\phi+\rho\zeta)$, the last identity of \eqref{e:bal},
$R/c_p=(\gamma-1)/\gamma$, $\gamma R=(\gamma-1)c_p$ and
\[
 \rho\theta-\frac\Pi{c_p}=\rho_b\theta_b\Bigl(1+\frac\Pi{\gamma p_b}\Bigr),
 \qquad
 \partial_1\log(m\circ Y)=g_T-\frac{\gamma-1}\gamma\frac{\partial_1p_b}{p_b}
 =\partial_1f_h+G_c;
\]
the $B_b$ terms give $-\frac{c_v\rho_bB_b}{c_p\theta_b}\frac\Pi{p_b}$, and
\[
 \rho D_u\chi-\frac c{\theta_b}\Delta(\theta_b\chi)
 =S_\eta-\rho_bG_c\psi_1.
\]
The equations \eqref{e:hist} are linear and uniformly parabolic, with
sources in $L^2(s,T_e;H^{-1})$. A Galerkin construction uses the regularized
coefficients
\[
 \rho_\varepsilon=\rho_b+\phi*\varphi_\varepsilon,\qquad
 \rho_\varepsilon u_\varepsilon=\rho_bu_be_1
 +(\rho u-\rho_bu_be_1)*\varphi_\varepsilon,
\]
which preserve the continuity equation. Testing the first equation by
$\eta$ and using
$\nabla(\theta_b\eta)\cdot\nabla\frac\eta{\theta_b}=|\nabla\eta|^2-g_T^2\eta^2$
gives \eqref{e:etaE} and uniform bounds for passage to the limit.
The same identity for differences gives uniqueness and $\chi=\eta+f$.
Finally, $\theta_b\rho(1+\chi)=\rho\theta-\Pi/c_p$ gives \eqref{e:espl}.

Let $\zeta=\zeta_\ell+\zeta_{\rm h}$, $\widehat{\zeta_\ell}=\chi_\ell\widehat\zeta$,
with a smooth $\chi_\ell$ equal to $1$ for $|\xi|\le1$ and to $0$ for
$|\xi|\ge2$. Then $\Vert\nabla\zeta_\ell\Vert_{L^\infty}\le C\Vert\zeta\Vert_{L^\infty}$
and $\Vert\nabla\zeta_{\rm h}\Vert_{H^1}\le C\Vert\nabla^2\zeta\Vert_{L^2}$.
By \eqref{e:stb}, \eqref{e:mu}, $\Vert g\Vert_{L^6}\le C\Vert\nabla g\Vert_{L^2}$,
\eqref{e:gn4}, H\"older's inequality on the cubes $Q_j$ with
$\Vert\nabla(\Pi,\psi)\Vert_{L^3(Q_j)}\le C\mu_*$, and the integration by
parts $\int|\nabla\psi|^{\frac{12}5}=-\int\psi\cdot\Div(|\nabla\psi|^{\frac25}\nabla\psi)$
with exponents $3,2,6$,
\begin{equation}\label{e:pbd}
\begin{gathered}
 \Vert\Pi\nabla\zeta\Vert_{L^2}+\Vert\psi\nabla\zeta\Vert_{L^2}
 +\Vert\nabla\zeta\Vert_{L^4}^2\le C\mu_*\Vert\nabla^2\zeta\Vert_{L^2},\\
 \Vert(|\nabla\Pi|+|\nabla\psi|)\nabla\zeta\Vert_{L^2}
 \le C\mu_*\bigl(\Vert\nabla\Pi\Vert_{L^2}+\Vert\nabla\psi\Vert_{L^2}
 +\Vert\nabla^2\zeta\Vert_{L^2}\bigr),\\
 \Vert Q_\mu(\nabla\psi)\Vert_{L^{6/5}}\le C\mu_*\Vert\nabla^2\psi\Vert_{L^2},
 \qquad
 \Vert\Pi\Div\psi\Vert_{L^{6/5}}\le C\mu_*\Vert\nabla\psi\Vert_{L^2}.
\end{gathered}
\end{equation}
By \eqref{e:state}, $|\Pi\nabla\phi|\le C_M(\mu_*|\nabla\Pi|+|\Pi\nabla\zeta|
+\mu_*(|\partial_1\rho_b|+|\partial_1\theta_b|)|\Pi|)$; this, \eqref{e:rhoD},
$\nabla\frac\phi\rho=\frac{\rho_b\nabla\phi-\phi(\partial_1\rho_b)e_1}{\rho^2}$,
$\Vert\partial_1p_b\Vert_{L^\infty}\le C_M\tau^{-\frac12}$, \eqref{e:pbd},
\eqref{e:E1Z} and $N_0(\Pi;s,t)+N_0(\psi;s,t)\le2N_o$ give \eqref{e:ehp}.

By \eqref{e:Fn} and integration by parts,
\[
\begin{aligned}
 \langle S_\eta,\eta\rangle={}&-c\int\nabla\Bigl(\frac{\Pi}{c_p\rho}\Bigr)\cdot\nabla\frac\eta{\theta_b}
 -\int\frac{\kappa\partial_1^2\theta_b+\nu h_b^2}{c_p\theta_bp_b}\Pi\eta
 +\int f_h\partial_1(\rho_b\psi_1\eta)
 +\int\frac{H[\psi]}{c_p\theta_b}\eta+\langle F_n,\eta\rangle .
\end{aligned}
\]
We bound the five terms after multiplication by $w_s$ and integration over
$[s,t]$. The first is at most $\frac c{16}\Vert\nabla\eta\Vert_{L^2_w}^2
+C_M\Vert\nabla\Pi\Vert_{L^2_w}^2+\Vert g_T\eta\Vert_{L^2_w}^2
+C_M\Bigl\Vert\nabla\Bigl(\frac{\Pi}{c_p\rho}\Bigr)-\frac{\nabla\Pi}{c_p\rho}\Bigr\Vert_{L^2_w}(\Vert\nabla\eta\Vert_{L^2_w}
+\Vert g_T\eta\Vert_{L^2_w})$. For the second, \eqref{e:curv} and \eqref{e:E2}
give $\Vert\kappa\partial_1^2\theta_b+\nu h_b^2\Vert_{L^\infty}\le C_M\delta\tau^{-1}
+Cr_2+\omega_b(s)\tau^{-1}$ with $\int_s^\infty\tau r_2^2\le\omega_b(s)$, and
$\int_s^t\frac{w_s}{1+r}\Vert Z\Vert_{L^2}^2\dd r\le q_a^{-1}N_0(Z;s,t)^2$ gives
the bound $(C_M\delta+\omega_b(s))N_o\Phi^{\frac12}$. The third is at most
$C\epsilon_f(s)\tau^{-\frac12}\bigl(\Vert\partial_1\rho_b\Vert_{L^\infty}\Vert\psi\Vert_{L^2}
\Vert\eta\Vert_{L^2}+\Vert\nabla\psi\Vert_{L^2}\Vert\eta\Vert_{L^2}
+\Vert\psi\Vert_{L^2}\Vert\nabla\eta\Vert_{L^2}\bigr)$ and, with
$\Vert\partial_1\rho_b\Vert_{L^\infty}\le C_M\tau^{-\frac12}$ by \eqref{e:E1}, has
the bound $C_M\epsilon_f(s)N_o\Phi^{\frac12}$. The fourth is at most
$C\Vert h_b\Vert_{L^2(s,\infty;L^\infty)}\sup_r(\sqrt{w_s}\Vert\eta\Vert_{L^2})
\Vert\nabla\psi\Vert_{L^2_w}\le\omega_b(s)N_o\Phi^{\frac12}$ by \eqref{e:E3}, and
the last is at most $C\Vert F_n\Vert_{L^2_w(L^{6/5})}\Vert\nabla\eta\Vert_{L^2_w}
\le C_M\mu_*N_o\Phi^{\frac12}$ by $\Vert\eta\Vert_{L^6}\le C\Vert\nabla\eta\Vert_{L^2}$.
Together with \eqref{e:ehp} and \eqref{e:E1Z} for $g_T\eta$, this proves
\eqref{e:Seta}. We multiply \eqref{e:etaE} by $w_s$; since
$\partial_tw_s=-\frac{q_aw_s}{1+t}$, integration gives
\[
 \frac{w_s(t)}2\int\rho\eta^2(t)+\int_s^tw_s\Bigl[\frac{q_a}{2(1+r)}\int\rho\eta^2
 +c\Vert\nabla\eta\Vert_{L^2}^2\Bigr]\dd r
 =\frac12\int\rho\chi^2(s)+\int_s^tw_s\Bigl[c\int g_T^2\eta^2
 +\langle S_\eta,\eta\rangle\Bigr]\dd r .
\]
Here $\Vert\chi(s)\Vert_{L^2}^2\le C_MI_s$, \eqref{e:E1Z} bounds the $g_T$ term
by $C(\delta^2+\omega_b(s))\Phi$, and \eqref{e:Seta} applies. After decreasing
$\delta_1(M)$, $\mu_0$ and increasing $S$, the $\Phi$ terms on the right are
absorbed, which proves the first line of \eqref{e:eta}.

On $\R\times[S_0,\infty)$, with $E_\pm=E_\pm(q_\pm,m)$, $J_\pm$, $S_\pm$ and
$\Gamma$ of Lemma~\ref{l:dep}, put
\[
 \mathcal A=\int_{-\infty}^yE_+\dd z+\int_y^\infty E_-\dd z,\qquad
 \mathcal Q=q_+^2+q_-^2,\qquad
 S_{\rm pot}=\int_{-\infty}^yS_+\dd z+\int_y^\infty S_-\dd z .
\]
By \eqref{e:hcur}, $\partial_tE_\sigma+\partial_y(\frac\sigma2q_\sigma^2-J_\sigma)
=S_\sigma$ with $E_\sigma\ge0$, and integration of \eqref{e:hcur} over the two half-lines gives
\[
 \partial_t\mathcal A=-\frac{\mathcal Q}2+J_+-J_-+S_{\rm pot},\qquad
 \partial_y\mathcal A=E_+-E_- .
\]
Let $a=\mathcal A\circ Y$. By \eqref{e:chain}, $D_ba=(\partial_t\mathcal A)\circ Y$
and $\partial_1a=\rho_b(E_+-E_-)\circ Y$; thus $a\ge0$,
$|\partial_1a|\le C_M\mathcal Q\circ Y$, and Lemma~\ref{l:dep} gives
\[
 \sup_{r\ge s}\bigl(\Vert\mathcal A\Vert_{L^\infty}+\Vert\mathcal Q\Vert_{L^\infty}\bigr)
 +\int_s^\infty\bigl(\Vert\Gamma\Vert_{L^\infty}^2
 +\Vert S_{\rm pot}\Vert_{L^\infty}\bigr)\dd r\le\omega_b(s).
\]
Since $b\circ X=\upsilon_m+H_c$ and $\upsilon_m=\frac{m^{1/4}}2(q_++q_-)$ by
\eqref{e:Um} and \eqref{e:char}, \eqref{eq:pr-gau} gives
$b^2\le C_M\mathcal Q\circ Y+C\delta^2\tau^{-1}$. Testing \eqref{e:hist} by
$a\eta$, with $D_ua=D_ba+\psi_1\partial_1a$, gives
\[
\begin{aligned}
 &\frac12\frac{\dd}{\dd t}\int\rho a\eta^2
 +\frac14\int\rho(\mathcal Q\circ Y)\eta^2+c\int a|\nabla\eta|^2\\
 &=\frac12\int\rho\eta^2\bigl[(J_+-J_-+S_{\rm pot})\circ Y+\psi_1\partial_1a\bigr]
 -c\int\eta\nabla a\cdot\nabla\eta-c\int g_T\eta^2\partial_1a
 +c\int ag_T^2\eta^2+\langle S_\eta,a\eta\rangle .
\end{aligned}
\]
We multiply by $w_s$ and integrate; $\partial_tw_s\le0$ and $a\ge0$. Write
$\mathcal Y=\int\rho(\mathcal Q\circ Y)\eta^2$. By $|J_\sigma|\le C_*|q_\sigma||\Gamma|$,
$|\psi_1|\le C\mu_*$, $\Vert g_T\Vert_{L^\infty}\le C(\delta+\omega_b(s))$ and
Young's inequality, the terms on the right other than
$\langle S_\eta,a\eta\rangle$ are at most
\[
 \frac4{64}\mathcal Y+C\bigl(\Vert\Gamma\Vert_{L^\infty}^2
 +\Vert S_{\rm pot}\Vert_{L^\infty}\bigr)\Vert\eta\Vert_{L^2}^2
 +C_M\Vert\mathcal Q\Vert_{L^\infty}\Vert\nabla\eta\Vert_{L^2}^2
 +C\Vert\mathcal A\Vert_{L^\infty}\Vert g_T\eta\Vert_{L^2}^2
\]
after decreasing $\mu_0,\delta_1(M)$ and increasing $S$. In
$\langle S_\eta,a\eta\rangle$ we use the five terms above with $\eta$ replaced by
$a\eta$. Since $|\nabla(a\eta)|\le a|\nabla\eta|+C_M(\mathcal Q\circ Y)|\eta|$,
\[
 \Vert a\eta\Vert_{L^6}\le C\Vert\mathcal A\Vert_{L^\infty}\Vert\nabla\eta\Vert_{L^2}
 +C_M\Vert\mathcal Q\Vert_{L^\infty}^{\frac12}\mathcal Y^{\frac12};
\]
the terms with a derivative of $a$, namely from $\nabla(\Pi/(c_p\rho))$, from
$f_h\rho_b\psi_1\eta\partial_1a$ and from $F_n$, are at most
$\frac3{64}\mathcal Y+C_M\Vert\mathcal Q\Vert_{L^\infty}\bigl(\Vert\nabla(\Pi/(c_p\rho))\Vert_{L^2}^2
+\epsilon_f(s)^2\tau^{-1}\Vert\psi\Vert_{L^2}^2+\Vert F_n\Vert_{L^{6/5}}^2\bigr)$, and
the others are the terms bounded above, multiplied by
$\Vert\mathcal A\Vert_{L^\infty}$. With $\Vert\chi(s)\Vert_{L^2}^2\le C_MI_s$ for
the initial value, \eqref{e:ehp}, \eqref{e:E1Z} and the first line of
\eqref{e:eta}, we obtain
\[
 \frac18\int_s^tw_s\int\rho(\mathcal Q\circ Y)\eta^2\dd x\dd r
 \le\omega_b(s)\bigl(I_s+N_o^2+\Phi\bigr).
\]
With $b^2\le C_M\mathcal Q\circ Y+C\delta^2\tau^{-1}$ and the first line, this
gives the second line of \eqref{e:eta}. The third follows from
$\nabla(b\eta)=b\nabla\eta+h_b\eta e_1$, \eqref{e:bsup} and \eqref{e:E1Z}.

Since $D_uY=\rho_b\psi_1$ and $\partial_1Y=\rho_b$,
\[
 D_u\log\varpi_\sigma=\lambda_0(\sigma\rho_b\psi_1-c_0),\qquad
 \nabla\log\varpi_\sigma=\lambda_0\sigma\rho_be_1 .
\]
Let $C_M^{-1}\le\rho,\rho_b\le C_M$ by \eqref{e:stb}, take $\mu_0$ with
$|\rho_b\psi_1|\le\frac{c_0}4$, and $\lambda_0=\min\{\frac{c_0}{4cC_M^3},4c_Gc_0\}$,
$k_0=\frac{\lambda_0c_0}4$. By \eqref{eq:pr-gau},
$\varpi_\sigma G_c^2\le C_M\delta^2\tau^{-1}
\exp\{(\frac{\lambda_0^2}{8c_G}-\lambda_0c_0)\tau\}\le C_M\delta^2\tau^{-1}e^{-2k_0\tau}$.
Testing the second equation of \eqref{e:hist} by $\varpi_\sigma f$ gives
\[
\begin{aligned}
 &\frac12\frac{\dd}{\dd t}\int\rho\varpi_\sigma f^2
 +\frac{\lambda_0}2\int\rho(c_0-\sigma\rho_b\psi_1)\varpi_\sigma f^2
 +c\int\varpi_\sigma|\nabla f|^2\\
 &=c\int\varpi_\sigma g_T^2f^2-c\lambda_0\sigma\int\varpi_\sigma\rho_bf(\partial_1f+g_Tf)
 -\int\rho_bG_c\psi_1\varpi_\sigma f .
\end{aligned}
\]
The choice of $\lambda_0$ gives $c\lambda_0^2\rho_b^2\le\frac{\lambda_0c_0\rho}4$;
after increasing $S$ so that $\Vert g_T\Vert_{L^\infty}$ is small,
\[
 \frac{\dd}{\dd t}\int\rho\varpi_\sigma f^2+c_1\int\varpi_\sigma(f^2+|\nabla f|^2)
 \le C_M\delta^2\tau^{-1}e^{-2k_0\tau}\Vert\psi\Vert_{L^2}^2,
\]
with $c_1>0$ depending only on $M$ and the fixed data. With
$\sigma(Y+h_*)$ replaced by $\min\{\sigma(Y+h_*),K\}$ in $\varpi_\sigma$ the same
inequality holds with the same constants, and $K\to\infty$ gives it for
$\varpi_\sigma$. Since $f(s)=0$, $w_s$ is decreasing and
$\frac{\tau_t}{\tau_r}\le1+t-r$, Gronwall's inequality gives
\[
 w_s(t)\tau_t^{L_0}\int\rho\varpi_\sigma f^2(t)
 \le C_{M,L_0}\delta^2\int_s^te^{-\frac{c_1(t-r)}{2C_M}}e^{-k_0\tau_r}
 \frac{w_s(r)}{\tau_r}\Vert\psi(r)\Vert_{L^2}^2\dd r
 \le C_{M,L_0}\delta^2e^{-k_0s}N_0(\psi;s,t)^2 .
\]
Multiplying the inequality by $w_s\tau^{L_0}$ and integrating, we bound the
term $L_0w_s\tau^{L_0-1}\int\rho\varpi_\sigma f^2$ by the preceding estimate with
$L_0+1$; this proves the first line of \eqref{e:ext}. On the support of
$f_{\rm out}$ one of $\varpi_\pm$ is at least $1$, and
$|\nabla\chi_c|\le C_M\tau^{-1}$; this proves the second line.

By $\nabla d_N=-d_N\frac{\nabla p}p$, $\nabla j_N=-j_N\frac{\nabla\theta}\theta$,
$\nabla k_N=-k_N(\frac{\nabla p}p+\frac{\nabla\theta}\theta)$,
$\nabla p=\nabla\Pi+(\partial_1p_b)e_1$, $\nabla\theta=\nabla\zeta+(\partial_1\theta_b)e_1$
and \eqref{e:rhoD}, the quantities in \eqref{e:etst} are sums of terms of
$D_*$, of $\nabla(b\eta)$, and of products. The products are bounded by
\eqref{e:pbd}, \eqref{e:prod}, \eqref{e:E1Z}, \eqref{e:Jsym}, the second line of
\eqref{e:eta}, and, for $g\in\{\Pi,\zeta\}$, by
\[
 \Vert b\eta\nabla g\Vert_{L^2}^2\le\Vert b\Vert_{L^\infty}^2
 \sum_j\Vert\eta\Vert_{L^6(Q_j)}^2\Vert\nabla g\Vert_{L^3(Q_j)}^2
 \le C_M\mu_*^2\tau^{-1}\Vert\eta\Vert_{H^1}^2,
\]
which uses \eqref{e:bsup}. This proves \eqref{e:etst}.

By \eqref{e:hist} and \eqref{e:Fn}, for $T\in H^1$,
\[
\begin{aligned}
 \langle G_\eta,T\rangle={}&-\int\mathbb S(b\eta e_1):\nabla T
 +\int\nabla\Bigl(\frac{cbT_1}{\theta_b}\Bigr)\cdot\nabla\Bigl(\theta_b\eta+\frac{\Pi}{c_p\rho}\Bigr)
 +\int\frac{\kappa\partial_1^2\theta_b+\nu h_b^2}{c_p\theta_bp_b}b\Pi T_1\\
 &-\int f_h\partial_1(\rho_bb\psi_1T_1)
 -\int\frac{bH[\psi]}{c_p\theta_b}T_1-\langle F_n,bT_1\rangle
 -\int\rho h_b\psi_1\eta T_1 .
\end{aligned}
\]
We bound the weighted time integrals of these seven terms, in this order, by
$\Vert\nabla(b\eta)\Vert_{L^2_w}\Vert\nabla T\Vert_{L^2_w}$;
$\Vert\nabla(\frac{bT_1}{\theta_b})\Vert_{L^2_w}(\Phi^{\frac12}+N_o)$; the second
term of $\langle S_\eta,\eta\rangle$ with the factor $\Vert b\Vert_{L^\infty}$; the
third term with the factors $\Vert b\Vert_{L^\infty}$ and $\Vert(b,h_b)\Vert_{L^\infty}$;
$\Vert b\Vert_{L^\infty}\Vert h_b\Vert_{L^2(s,\infty;L^\infty)}N_0(T)N_o$;
$\Vert F_n\Vert_{L^2_w(L^{6/5})}\Vert\nabla(bT_1)\Vert_{L^2_w}$; and
$C\mu_*\Vert h_b\Vert_{L^\infty}\Vert\eta\Vert_{L^6}\Vert T\Vert_{L^2}$ with
$\Vert h_b\Vert_{L^\infty}\le\omega_b(s)\tau^{-\frac12}$ by \eqref{e:E2}. With
\eqref{e:eta}, \eqref{e:ehp}, \eqref{e:bsup} and \eqref{e:E1Z}, this proves
\eqref{e:Getb}.

In \eqref{e:Reta}, \eqref{e:Getb} with $T=W^0$ and \eqref{e:etst} bound
$\langle G_\eta,W^0\rangle$. The other terms are at most
\[
\begin{gathered}
 \Vert\nabla\Pi\Vert_{L^2_w}\Vert b\eta\Vert_{L^2_w},\qquad
 \Vert Z^0\Vert_{L^2_w}\bigl(\Vert B_b\eta\Vert_{L^2_w}
 +C\Vert\nabla(b\eta)\Vert_{L^2_w}\bigr),\qquad
 C\Vert\tau^{-1}\zeta\Vert_{L^2_w}\Vert B_b\eta\Vert_{L^2_w},\\
 C\Vert\Div(\rho V)\Vert_{L^2_w}\Vert\nabla(b\eta)\Vert_{L^2_w},\qquad
 C\sum_{i=1}^3\bigl(\Vert W^i\Vert_{L^2_w}+\Vert Z^i\Vert_{L^2_w}\bigr)
 \bigl(\Vert\nabla(A_b\eta)\Vert_{L^2_w}+\Vert\nabla(B_b\eta)\Vert_{L^2_w}\bigr),
\end{gathered}
\]
where the screened term was written as
$\epsilon_K\int\Div(\rho V)J\{\gamma p\partial_1(b\eta)\}$ and \eqref{e:Jsym} was
used. By \eqref{e:E6} and \eqref{e:E1Z},
$\Vert B_b\eta\Vert_{L^2_w}+\Vert\nabla(A_b\eta)\Vert_{L^2_w}
+\Vert\nabla(B_b\eta)\Vert_{L^2_w}\le\omega_b(s)\Phi^{\frac12}$; with
\eqref{e:eta} and \eqref{e:etst} this bounds $\mathcal R_\eta$. In
$\mathcal R_{\rm out}$,
$\int_s^tw_s|\int\rho A_bf_{\rm out}(W^0)_1|\le C\sup_r\Vert A_b\Vert_{L^\infty}
\Vert\sqrt\tau f_{\rm out}\Vert_{L^2_w}\Vert\tau^{-\frac12}W^0\Vert_{L^2_w}$, and the
other terms are bounded as for $\mathcal R_\eta$ with $f_{\rm out}$ in place of
$\eta$; the second line of \eqref{e:ext} with $L_0=0,1$ gives the bound
$C_M\delta e^{-\frac{k_0s}2}(I_s+N_o^2+\Phi)$. This proves \eqref{e:eret}.
\end{proof}

By \eqref{e:espl}, \eqref{e:zN} and \eqref{e:cut},
\begin{equation}\label{e:chio}
 \chi^b_o=\frac1{\theta_b}\Bigl(\zeta-\frac\Pi{c_p\rho_b}\Bigr)
 \quad\text{satisfies}\quad \chi^b_o=\eta+f+z_N,
\end{equation}
and $\zeta=\theta_b(F+\eta+f_{\rm out})+\Pi/(c_p\rho)$ gives, for $1\le i\le3$,
\begin{equation}\label{e:tz7}
\begin{aligned}
 \partial_i\zeta={}&\theta_b\partial_iF+\mathbf1_{\{i=1\}}(\partial_1\theta_b)F
 +\theta_b\partial_i\eta+\mathbf1_{\{i=1\}}(\partial_1\theta_b)\eta\\
 &+\partial_i\Bigl(\frac{\Pi}{c_p\rho}\Bigr)+\theta_b\partial_if_{\rm out}
 +\mathbf1_{\{i=1\}}(\partial_1\theta_b)f_{\rm out}.
\end{aligned}
\end{equation}

\begin{lemma}\label{l:hes}
Fix $\alpha_G\in(0,c_G]$ with $16c_0\alpha_G\le\lambda_0$ and so small that
Lemma~\ref{l:bulk} applies with the exponent $\alpha_G$. Put
\begin{equation}\label{e:EC}
\begin{gathered}
 \mathcal W=\exp\Bigl\{\frac{\alpha_G(Y+h_*)^2}{\tau}\Bigr\},\qquad
 E_C=\frac1\tau\int_{\R^3}\mathcal WF^2\dd x
 +\int_{\R^3}\mathcal W|\nabla F|^2\dd x,\\
 D_C=\int_{\R^3}\mathcal W\Bigl(\frac{F^2}{\tau^2}+\frac{|\nabla F|^2}\tau
 +|\nabla^2F|^2\Bigr)\dd x,\qquad
 \mathcal C_C(s,t)=\sup_{s\le r\le t}w_sE_C+\int_s^tw_sD_C\dd r .
\end{gathered}
\end{equation}
Then
\begin{equation}\label{e:Ch}
 \mathcal C_C(s,t)\le C_M\delta^2N_0(\psi;s,t)^2,
\end{equation}
for almost every $r\in[s,t]$
\begin{equation}\label{e:ctr}
\begin{aligned}
 \sup_{y\in\R}e^{\frac{\alpha_G(y+h_*)^2}{2\tau}}
 \Vert(F\circ X)(y,\cdot)\Vert_{L^2(\R^2)}^2&\le C_M\sqrt\tau E_C,\\
 \sup_{y\in\R}e^{\frac{\alpha_G(y+h_*)^2}{2\tau}}
 \Vert((\nabla F)\circ X)(y,\cdot)\Vert_{L^2(\R^2)}^2&\le C_M\sqrt\tau D_C,
\end{aligned}
\end{equation}
and
\begin{equation}\label{e:ABF}
\begin{aligned}
 \Vert A_bF\Vert_{L^2_w}^2+\Vert\nabla(A_bF)\Vert_{L^2_w}^2
 +\Vert\nabla(B_bF)\Vert_{L^2_w}^2&\le\omega_b(s)\mathcal C_C(s,t),\\
 \Vert B_bF\Vert_{L^2_w}^2&\le\bigl(C_M\delta^2+\omega_b(s)\bigr)\mathcal C_C(s,t).
\end{aligned}
\end{equation}
With
\begin{equation}\label{e:Tin}
 \mathcal T_{\rm in}=\int_{\R^3}\Bigl[\epsilon_K\rho A_bF\partial_1J\Pi
 +B_bFZ^0+\frac{\theta_W}\tau\rho\zeta B_bF\Bigr]\dd x
 +\sum_{i=1}^3\int_{\R^3}\bigl[\rho\partial_i(A_bF)(W^i)_1
 +\partial_i(B_bF)Z^i\bigr]\dd x,
\end{equation}
one has
\begin{equation}\label{e:Tib}
 \int_s^tw_s|\mathcal T_{\rm in}|\dd r
 \le\bigl(C_M\delta^2+\omega_b(s)\bigr)\bigl(I_s+N_o^2+\Phi\bigr).
\end{equation}
\end{lemma}

\begin{proof}
Since $D_bY=0$,
\[
 D_b\mathcal W=-\alpha_G\frac{(Y+h_*)^2}{\tau^2}\mathcal W,\qquad
 \nabla\log\mathcal W=\frac{2\alpha_G\rho_b(Y+h_*)}\tau e_1.
\]
By the choice of $\lambda_0$ in the proof of
Lemma~\ref{l:ent}, $4c\alpha_G\rho_b^2\le\frac\rho{16}$. By \eqref{eq:pr-gau} and
$\alpha_G\le c_G$, $\mathcal WG_c^2\le C_M\delta^2\tau^{-1}$. The second equation
of \eqref{e:hist} gives
\[
 \rho D_uF-\frac c{\theta_b}\Delta(\theta_bF)=-\rho_b\chi_cG_c\psi_1+\rho K_c,
 \qquad
 K_c=\Bigl[D_u\chi_c-\frac c\rho\Delta\chi_c-\frac{2c}\rho g_T\partial_1\chi_c\Bigr]f
 -\frac{2c}\rho\nabla\chi_c\cdot\nabla f .
\]
On the support of $K_c$, $2c_0\tau\le|Y+h_*|\le3c_0\tau$,
$|D_u\chi_c|+|\nabla\chi_c|+|\Delta\chi_c|\le C_M\tau^{-1}$ and, for the $\sigma$
with $\sigma(Y+h_*)\ge2c_0\tau$, $\mathcal W\le e^{9\alpha_Gc_0^2\tau}\le\varpi_\sigma$.
Hence \eqref{e:ext} gives
$\int_s^tw_s\int\mathcal WK_c^2\le C_M\delta^2e^{-k_0s}N_0(\psi;s,t)^2$. Testing by
$\mathcal WF$ gives
\[
\begin{aligned}
 &\frac12\frac{\dd}{\dd t}\int\rho\mathcal WF^2+c\int\mathcal W|\nabla F|^2
 +\frac{\alpha_G}2\int\rho\mathcal W\frac{(Y+h_*)^2}{\tau^2}F^2\\
 &=\frac12\int\rho\mathcal W(\psi\cdot\nabla\log\mathcal W)F^2
 -c\int\mathcal WF\nabla\log\mathcal W\cdot\nabla F\\
 &\quad-c\int\mathcal Wg_T\partial_1\log\mathcal W\,F^2+c\int\mathcal Wg_T^2F^2\\
 &\quad
 -\int\rho_b\mathcal W\chi_cG_c\psi_1F+\int\rho\mathcal WK_cF .
\end{aligned}
\]
By $\Vert\sqrt{\mathcal W}F\Vert_{L^6}\le C(\Vert\sqrt{\mathcal W}\nabla F\Vert_{L^2}
+\Vert\sqrt{\mathcal W}F\nabla\log\mathcal W\Vert_{L^2})$, $\Vert\psi\Vert_{L^3}\le\mu_*$, the
bound $4c\alpha_G\rho_b^2\le\frac\rho{16}$ and \eqref{e:E1} in the form
$|g_T|\le C_M\delta\tau^{-\frac12}+Ck_1$ with $\int_s^\infty\Vert k_1\Vert_{L^\infty}^2
\le\omega_b(s)$, the first four terms on the right are at most
$\frac c2\int\mathcal W|\nabla F|^2+\frac{\alpha_G}4\int\rho\mathcal W\frac{(Y+h_*)^2}{\tau^2}F^2
+C_M(\delta^2\tau^{-1}+k_1^2)\int\mathcal WF^2$ for small $\mu_0$, and the last two
are at most $C_M\delta\tau^{-\frac12}\Vert\sqrt{\mathcal W}F\Vert_{L^2}\Vert\psi\Vert_{L^2}
+C\Vert\sqrt{\mathcal W}K_c\Vert_{L^2}\Vert\sqrt{\mathcal W}F\Vert_{L^2}$. We multiply by
$\frac{w_s}\tau$, whose derivative is $-(1+q_a)\frac{w_s}{\tau^2}$, use $F(s)=0$,
$\int_s^t\frac{w_s}\tau\Vert\psi\Vert_{L^2}^2\le q_a^{-1}N_0(\psi;s,t)^2$ and
absorb, for small $\delta_1(M)$ and large $S$:
\[
 \sup_{[s,t]}\frac{w_s}\tau\int\mathcal WF^2
 +\int_s^tw_s\int\mathcal W\Bigl(\frac{|\nabla F|^2}\tau+\frac{F^2}{\tau^2}
 +\frac{\alpha_G(Y+h_*)^2F^2}{\tau^3}\Bigr)\dd x\dd r\le C_M\delta^2N_0(\psi;s,t)^2 .
\]
With $S_F=\frac{2c}\rho g_T\partial_1F+\frac c\rho\frac{\partial_1^2\theta_b}{\theta_b}F
-\frac{\rho_b}\rho\chi_cG_c\psi_1+K_c$, the equation reads
$D_bF-\frac c\rho\Delta F=-\psi\cdot\nabla F+S_F$. Testing it by
$-\Div(\mathcal W\nabla F)$ and using $\nabla D_bF=D_b\nabla F+h_b(\partial_1F)e_1$
gives
\[
\begin{aligned}
 &\frac12\frac{\dd}{\dd t}\int\mathcal W|\nabla F|^2
 +\int\frac c\rho\mathcal W(\Delta F)^2
 +\frac{\alpha_G}2\int\mathcal W\frac{(Y+h_*)^2}{\tau^2}|\nabla F|^2\\
 &=\int\mathcal Wh_b\Bigl(\frac{|\nabla F|^2}2-|\partial_1F|^2\Bigr)
 -\int\frac c\rho\mathcal W\Delta F\,\nabla\log\mathcal W\cdot\nabla F
 -\int(-\psi\cdot\nabla F+S_F)\Div(\mathcal W\nabla F) .
\end{aligned}
\]
Two integrations by parts give
$\int\mathcal W|\nabla^2F|^2\le C\int\mathcal W(\Delta F)^2+C\int\mathcal W|\nabla\log\mathcal W|^2|\nabla F|^2$,
and the transport term is at most
$C_M\mu_*\int\mathcal W\{(\Delta F)^2+|\nabla\log\mathcal W|^2|\nabla F|^2\}$. In the first term
we write $h_b=\partial_1b$ and integrate by parts; on the support of $F\circ X$,
$|y-s_\sigma t|\ge\frac{|s_\sigma|t}2-3c_0-|h_*|$, so \eqref{e:CH}, Lemma~\ref{l:pr},
Lemma~\ref{l:asy} and \eqref{eq:pr-gau} give $\tau b^2\le C\delta^2+\omega_b(s)$, and
the term is at most
$\varepsilon\int\mathcal W(|\nabla^2F|^2+|\nabla\log\mathcal W|^2|\nabla F|^2)
+C_\varepsilon(C\delta^2+\omega_b(s))\tau^{-1}\int\mathcal W|\nabla F|^2$. By \eqref{e:curv},
$\int_s^tw_s\int\mathcal W|\partial_1^2\theta_bF|^2\le C\delta^2\int_s^t\frac{w_s}{\tau^2}
\int\mathcal WF^2+C\int_s^\infty\tau r_2^2\dd r\sup\frac{w_s}\tau\int\mathcal WF^2$, and the
other terms of $S_F$ are bounded as above. Multiplying by $w_s$, integrating
and absorbing gives \eqref{e:Ch}.

For $H=F$ and $H=\partial_iF$ let
$\mathcal U_H=e^{\frac{\alpha_G(y+h_*)^2}{4\tau}}H\circ X$. By
$\partial_y(H\circ X)=V_b(\partial_1H)\circ X$, $\dd x_1=V_b\dd y$ and
$\frac{z^2}\tau e^{\frac{\alpha_Gz^2}{2\tau}}\le C_{\alpha_G}e^{\frac{\alpha_Gz^2}\tau}$,
\[
\begin{gathered}
 \Vert\mathcal U_F\Vert_{L^2}^2\le C_M\tau E_C,\qquad
 \Vert\partial_y\mathcal U_F\Vert_{L^2}^2\le C_ME_C,\\
 \sum_i\Vert\mathcal U_{\partial_iF}\Vert_{L^2}^2\le C_M\tau D_C,\qquad
 \sum_i\Vert\partial_y\mathcal U_{\partial_iF}\Vert_{L^2}^2\le C_MD_C,
\end{gathered}
\]
and \eqref{eq:planar} in the variables $(y,x')$ gives \eqref{e:ctr}. For
$a=a(x_1,r)$, \eqref{e:ctr} gives
$\Vert aF\Vert_{L^2}^2\le C_ME_C\,A^{\alpha_G}_{a\circ X}$ with $A^{\alpha_G}$ of
\eqref{e:nobs}, and also
$\Vert aF\Vert_{L^2}^2\le C_M\sqrt\tau E_C\Vert a\Vert_{L^2(\R)}^2$. By
\eqref{e:chain}, $(\partial_1p_b)\circ X=\frac{\partial_yP_b}{V_b}$ and
$h_b\circ X=\frac{\partial_yb_m}{V_b}$. By \eqref{e:cmp}, with
$\widehat p=p_*-\frac{(\gamma-1)H_c^2}{2V_c}$ and Lemma~\ref{l:pr},
\[
 \bigl|\partial_y\bigl(p_a-\widehat p(y+h_*,t)\bigr)\bigr|
 +\Bigl|\partial_y\Bigl(b_a-H_c(y+h_*,t)\Bigr)\Bigr|
 \le C_M\sum_\sigma\bigl(|\partial_y\vartheta_\sigma|
 +\delta\tau^{-\frac12}|\vartheta_\sigma|\bigr)
 \le C_M\tau^{-1}\sum_\sigma e^{-\frac{(y-s_\sigma t)^2}{8d_\sigma\tau}},
\]
whose overlap with $e^{-\frac{\alpha_G(y+h_*)^2}{2\tau}}$ gives
$A^{\alpha_G}\le C_be^{-c_bt}$, while $|\partial_y\widehat p|\le C\delta^2\tau^{-\frac32}$
and $|\partial_yH_c|\le C\delta\tau^{-1}$ by \eqref{eq:pr-gau}. With
Lemma~\ref{l:bulk} for
$\partial_y(P_b-p_a)$ and $\partial_y(b_m-b_a)$, and $\int_s^tw_sD_C\ge
\int_s^t\frac{w_s}{\tau^2}\int\mathcal WF^2$ for the contact part of $h_b$,
\[
 \Vert(\partial_1p_b)F\Vert_{L^2_w}^2\le\omega_b(s)\mathcal C_C,\qquad
 \Vert h_bF\Vert_{L^2_w}^2+\Vert(\partial_1^2\theta_b)F\Vert_{L^2_w}^2
 \le(C\delta^2+\omega_b(s))\mathcal C_C,
\]
where the last term uses \eqref{e:curv} as above. By \eqref{e:E5} and the second
bound for $\Vert aF\Vert$, $\Vert(\partial_1h_b,\partial_1A_b,\partial_1B_b)F\Vert_{L^2_w}^2
\le\omega_b(s)\mathcal C_C$. Now \eqref{e:bal} gives $\rho_bA_b=\nu\partial_1h_b-\partial_1p_b$
and $|B_b|\le C(|\partial_1^2\theta_b|+|h_b|)$, hence the bounds for $A_bF$ and
$B_bF$; with $\nabla(A_bF)=A_b\nabla F+(\partial_1A_b)Fe_1$, the same for $B_b$,
and $\int_s^t w_s\Vert A_b\nabla F\Vert^2_{L^2}\le\sup_r(w_s\Vert\nabla F\Vert_{L^2}^2)
\int_s^\infty\Vert A_b\Vert_{L^\infty}^2$ by \eqref{e:E3}, this proves \eqref{e:ABF}.

The five terms of \eqref{e:Tin} are at most, by \eqref{e:Jsym}, \eqref{e:etst},
$\Vert\tau^{-1}\zeta\Vert_{L^2_w}\le N_o$ and \eqref{e:ABF},
\[
 \omega_b(s)\mathcal C_C^{\frac12}N_o,\quad
 \bigl(C_M\delta+\omega_b(s)\bigr)\mathcal C_C^{\frac12}(I_s+N_o^2+\Phi)^{\frac12},\quad
 \bigl(C_M\delta+\omega_b(s)\bigr)\mathcal C_C^{\frac12}N_o,\quad
 \omega_b(s)\mathcal C_C^{\frac12}N_o,\quad
 \omega_b(s)\mathcal C_C^{\frac12}N_o,
\]
and \eqref{e:Ch} with $N_0(\psi;s,t)\le N_o$ gives \eqref{e:Tib}.
\end{proof}

\begin{lemma}\label{l:mat}
The function $\eta$ of \eqref{e:hist} satisfies
\begin{equation}\label{e:emat}
 \sup_{s\le r\le t}w_s\Vert\nabla\eta\Vert_{L^2}^2
 +\int_s^tw_s\Bigl[\Vert\nabla^2\eta\Vert_{L^2}^2+\Vert D_u\eta\Vert_{L^2}^2
 +\Vert D_b\eta\Vert_{L^2}^2+\frac{q_a}{1+r}\Vert\nabla\eta\Vert_{L^2}^2\Bigr]\dd r
 \le C_b\bigl(I_s+N_o^2+\Phi\bigr).
\end{equation}
For $C_{\rm test}>0$ and every scalar function $T_{\rm th}$ with
$\Vert T_{\rm th}\Vert_{L^2_w}^2\le C_{\rm test}(I_s+N_o^2+\Phi)$,
\[
 \int_s^tw_s\Bigl|\int_{\R^3}\rho_bbD_b\eta\,T_{\rm th}\dd x\Bigr|\dd r
 \le\omega_{b,C_{\rm test}}(s)\bigl(I_s+N_o^2+\Phi\bigr).
\]
\end{lemma}

\begin{proof}
By $\Pi=R(\theta_b\phi+\rho\zeta)$,
$\Pi/(c_p\rho)=\frac R{c_p}(\zeta+\theta_b\frac\phi\rho)$.
Expanding $\Delta(\theta_b\frac\phi\rho)$, with \eqref{e:ha} for $\nabla^2\phi$,
\eqref{e:gn4} for $|\nabla\phi|^2$, and \eqref{e:trace}, \eqref{e:E1}, \eqref{e:curv}
for the products with $\partial_1\rho_b,\partial_1^2\rho_b,\partial_1\theta_b,
\partial_1^2\theta_b$ as in the proof of Lemma~\ref{l:ha}, gives
$\int_s^tw_s\Vert\Delta(\Pi/(c_p\rho))\Vert_{L^2}^2\le C_b(I_s+N_o^2)$. By \eqref{e:char},
$\partial_y(f_h\circ X)=\frac{(e_m)_y}m+\partial_yV_c\bigl(\frac1m-\frac1{V_c}\bigr)
+\partial_y\log V_c-(\partial_y\log V_c)(y+h_*,t)$. Here
$|(e_m)_y|\le C_*|\Gamma|$ by \eqref{e:hcur}, $|e_m|\le C_b\tau^{-\frac12}$ by
$\epsilon_f(S)<\infty$ and \eqref{eq:pr-gau}, and the last difference is at most
$C\delta|h_*|\tau^{-1}$; thus $\int_S^\infty\Vert\partial_1f_h\Vert_{L^\infty}^2\le C_b$.
With \eqref{e:E3}, $\Vert h_b\Vert_{L^\infty}\le C_b$,
$\Vert Q_\mu(\nabla\psi)\Vert_{L^2}\le C\mu_*\Vert\nabla^2\psi\Vert_{L^2}$ by
\eqref{e:gn4}, and \eqref{e:E1Z}, this gives
$\int_s^tw_s\Vert S_\eta\Vert_{L^2}^2\le C_b(I_s+N_o^2)$.

Dividing \eqref{e:hist} by $\rho$ gives
$D_u\eta-\frac c\rho\Delta\eta=\frac{2c}\rho g_T\partial_1\eta
+\frac c\rho\frac{\partial_1^2\theta_b}{\theta_b}\eta+\frac{S_\eta}\rho$. Testing by
$-\Delta\eta$ gives
\[
 \frac12\frac{\dd}{\dd t}\Vert\nabla\eta\Vert_{L^2}^2+\int\frac c\rho(\Delta\eta)^2
 =\frac12\int(\Div u)|\nabla\eta|^2-\sum_{i,j}\int(\partial_iu_j)\partial_i\eta\,\partial_j\eta
 -\int\Bigl(\frac{2c}\rho g_T\partial_1\eta
 +\frac c\rho\frac{\partial_1^2\theta_b}{\theta_b}\eta+\frac{S_\eta}\rho\Bigr)\Delta\eta .
\]
By H\"older's inequality on the cubes $Q_j$,
$\int|\nabla\psi||\nabla\eta|^2\le C\mu_*(\Vert\nabla\eta\Vert_{L^2}^2
+\Vert\nabla\eta\Vert_{L^2}\Vert\nabla^2\eta\Vert_{L^2})$, and
$\Vert\nabla^2\eta\Vert_{L^2}=\Vert\Delta\eta\Vert_{L^2}$; hence
\[
 \frac{\dd}{\dd t}\Vert\nabla\eta\Vert_{L^2}^2+c_1\Vert\nabla^2\eta\Vert_{L^2}^2
 \le C_b\Vert\nabla\eta\Vert_{L^2}^2
 +C_M\Vert\partial_1^2\theta_b\Vert_{L^\infty}^2\Vert\eta\Vert_{L^2}^2
 +C_M\Vert S_\eta\Vert_{L^2}^2 .
\]
We multiply by $w_s$ and integrate; $\Vert\nabla\chi(s)\Vert_{L^2}^2\le C_bI_s$ and
\eqref{e:E3} give the bounds for $\nabla\eta$ and $\nabla^2\eta$ in \eqref{e:emat}.
The equation gives the bound for $D_u\eta$, and
$\Vert D_b\eta-D_u\eta\Vert_{L^2}=\Vert\psi\cdot\nabla\eta\Vert_{L^2}
\le C\mu_*\Vert\nabla\eta\Vert_{L^2}$ that for $D_b\eta$. Finally,
\[
 \int_s^tw_s\Bigl|\int\rho_bbD_b\eta\,T_{\rm th}\Bigr|\dd r
 \le C_b\sup_{r\ge s}\Vert b(r)\Vert_{L^\infty}\Vert D_b\eta\Vert_{L^2_w}
 \Vert T_{\rm th}\Vert_{L^2_w}
 \le C_bC_{\rm test}^{\frac12}\sup_{r\ge s}\Vert b(r)\Vert_{L^\infty}
 \bigl(I_s+N_o^2+\Phi\bigr),
\]
and $\sup_{r\ge s}\Vert b(r)\Vert_{L^\infty}\to0$ by \eqref{e:bsup}.
\end{proof}

\subsection{The energy identity for the perturbation}

With $Z^i$ of \eqref{e:test} put
\[
 \mathcal W_{\rm heat}=\sum_{i=1}^3\int_{\R^3}\Bigl[(\gamma-1)\frac{\partial_i\Pi}{\gamma p}
 \partial_iQ_\mu(\nabla\psi)+Z^i\,\partial_i\Bigl(\frac{Q_\mu(\nabla\psi)}{c_v\rho}\Bigr)
 +\epsilon_K(\gamma-1)\rho\partial_i\psi\cdot\nabla J\partial_iQ_\mu(\nabla\psi)\Bigr]\dd x .
\]

\begin{lemma}\label{l:hheat}
For $\epsilon_K\in(0,1]$,
\begin{equation}\label{e:hheat}
 \int_s^tw_s|\mathcal W_{\rm heat}|\dd r\le C_M\mu_*\bigl(I_s+N_o^2\bigr).
\end{equation}
\end{lemma}

\begin{proof}
By \eqref{e:etst}, $\sum_i\Vert Z^i\Vert_{L^2_w}\le C_MN_o$. By
$|\nabla Q_\mu(\nabla\psi)|\le C|\nabla\psi||\nabla^2\psi|$,
$\nabla\rho=\nabla\phi+(\partial_1\rho_b)e_1$, H\"older's inequality on the cubes
$Q_j$ with $\Vert\nabla\psi\Vert_{H^1(Q_j)}+\Vert\nabla\phi\Vert_{H^1(Q_j)}\le C\mu_*$,
and \eqref{e:gn4},
\[
 \Vert\nabla Q_\mu(\nabla\psi)\Vert_{L^2}
 +\Bigl\Vert\nabla\frac{Q_\mu(\nabla\psi)}{c_v\rho}\Bigr\Vert_{L^2}
 \le C_M\mu_*\Vert\nabla\psi\Vert_{H^2}.
\]
By \eqref{e:Jsym}, the pressure and screened terms of $\mathcal W_{\rm heat}$ are at most
$C_M(\Vert\nabla\Pi\Vert_{L^2}+\Vert\nabla\psi\Vert_{L^2})
\Vert\nabla Q_\mu(\nabla\psi)\Vert_{L^2}$, and the thermal terms are at most
$\sum_i\Vert Z^i\Vert_{L^2}\Vert\nabla(Q_\mu(\nabla\psi)/(c_v\rho))\Vert_{L^2}$.
The Cauchy--Schwarz inequality against $w_s$, the displayed bound and \eqref{e:ha}
give \eqref{e:hheat}.
\end{proof}

\begin{lemma}\label{l:mr}
There are $\epsilon_K,\theta_W,\theta_\eta\in(0,1]$, depending only on $M$ and
the fixed data, such that the energy
\begin{equation}\label{e:men}
\begin{aligned}
 \mathcal E_N={}&\int_{\R^3}e_N(\Pi,V,\zeta)\dd x
 +\epsilon_K\int_{\R^3}\rho V\cdot\nabla J\Pi\dd x
 +\sum_{i=1}^3\int_{\R^3}\bigl[e_N(\partial_i\Pi,\partial_i\psi,\partial_i\zeta)
 +\epsilon_K\rho\partial_i\psi\cdot\nabla J\partial_i\Pi\bigr]\dd x\\
 &+\frac{\theta_W}{2\tau}\int_{\R^3}\rho\zeta^2\dd x
 +\frac{\theta_\eta}2\int_{\R^3}\rho\eta^2\dd x
\end{aligned}
\end{equation}
satisfies
$C_M^{-1}(E_*+\theta_\eta\Vert\eta\Vert_{L^2}^2)\le\mathcal E_N
\le C_M(E_*+\theta_\eta\Vert\eta\Vert_{L^2}^2)$ and $\mathcal E_N(s)\le C_MI_s$, and
the following holds. With $f_p$ of \eqref{e:ps3}, put
\[
\begin{aligned}
 f_{p0}&=f_p,\qquad
 f_{m0}=-\rho h_b\psi_1e_1-\frac{\rho_bA_b\Pi}{\gamma p_b}e_1,\\
 f_{z0}&=-\partial_1\theta_b\,\psi_1-\frac{\rho_bB_b\Pi}{\gamma\rho p_b}
 -\frac{h_b\Pi}{c_v\rho}+\frac{H[\psi]+Q_\mu(\nabla\psi)}{c_v\rho},\qquad
 \dot\psi_0=\frac1\rho\bigl(\mathcal L\psi-\nabla\Pi+f_{m0}\bigr),\\
 f_{pi}&=\partial_if_{p0}-(\partial_iu)\cdot\nabla\Pi-\gamma(\partial_ip)\Div\psi,\qquad
 f_{mi0}=\partial_if_{m0}-(\partial_i\rho)\dot\psi_0-\rho(\partial_iu)\cdot\nabla\psi,\\
 f_{zi0}&=\partial_if_{z0}-(\partial_iu)\cdot\nabla\zeta
 +\partial_i\Bigl(\frac\kappa{c_v\rho}\Bigr)\Delta\zeta-(\gamma-1)(\partial_i\theta)\Div\psi,
\end{aligned}
\]
$\mathbf f_0^0=(f_{p0},f_{m0},f_{z0})$, $\mathbf f_0^i=(f_{pi},f_{mi0},f_{zi0})$,
$(P^0,v^0,z^0)=(\Pi,V,\zeta)$, $(P^i,v^i,z^i)=(\partial_i\Pi,\partial_i\psi,\partial_i\zeta)$
for $1\le i\le3$, and
\begin{equation}\label{e:mf}
\begin{aligned}
 \mathbf F^0&=\mathbf f_0^0+\bigl(-\gamma p\partial_1(b\eta),\
 \rho A_bfe_1+G_\eta,\ B_b(\eta+f)-(\gamma-1)\theta\partial_1(b\eta)\bigr),\\
 \mathbf F^i&=\mathbf f_0^i+\bigl(0,\ \rho\partial_i\{A_b(\eta+f)\}e_1,\
 \partial_i\{B_b(\eta+f)\}\bigr),\qquad 1\le i\le3 .
\end{aligned}
\end{equation}
For a triple $(P,v,z)$ and $\mathbf f=(f_p,f_m,f_z)$ let
$\mathscr R_N(P,v,z;\mathbf f)=\mathscr S_N(P,v,z;\mathbf f)+\mathscr G_N(P,v,z)
+\epsilon_K\mathscr K_N(P,v;f_p,f_m)$, and let
\[
 D_{\mathcal E}=\sum_{i=0}^3\bigl[\mathcal D_N(v^i,z^i)
 +\epsilon_K\Vert\nabla J^{\frac12}P^i\Vert_{L^2}^2\bigr]
 +\frac{\theta_W}{2\tau^2}\int_{\R^3}\rho\zeta^2\dd x
 +\frac{\theta_W\kappa}{c_v\tau}\Vert\nabla\zeta\Vert_{L^2}^2
 +\theta_\eta c\Vert\nabla\eta\Vert_{L^2}^2 .
\]
Then, in the sense of distributions on $(s,T_e)$,
\begin{equation}\label{e:mri}
\begin{aligned}
 \frac{\dd}{\dd t}\mathcal E_N+D_{\mathcal E}
 ={}&\sum_{i=0}^3\bigl[\mathscr R_N(P^i,v^i,z^i;\mathbf F^i)
 +\epsilon_K\mathscr P_N(P^i,v^i,z^i)\bigr]\\
 &+\frac{\theta_W}\tau\int_{\R^3}\rho\zeta\{F_z^0-(\gamma-1)\theta\Div V\}\dd x
 +\theta_\eta\Bigl[c\int_{\R^3}g_T^2\eta^2\dd x+\langle S_\eta,\eta\rangle\Bigr],
\end{aligned}
\end{equation}
and $\frac{\dd}{\dd t}(w_s\mathcal E_N)+w_s(\frac{q_a}\tau\mathcal E_N+D_{\mathcal E})$
equals $w_s$ times the right side. The right side of \eqref{e:mri} equals
\begin{equation}\label{e:mdec}
 \mathcal R_A+\mathcal R_\eta+\mathcal R_{\rm out}+\mathcal T_{\rm in}+\mathcal R_0,
\end{equation}
with $\mathcal R_\eta,\mathcal R_{\rm out}$ of \eqref{e:Reta}, $\mathcal T_{\rm in}$ of
\eqref{e:Tin} and
\[
\begin{aligned}
 \mathcal R_A={}&\sum_{i=0}^3\mathscr R_N(P^i,v^i,z^i;\mathbf f_0^i)
 +\frac{\theta_W}\tau\int_{\R^3}\rho\zeta f_{z0}\dd x
 +\int_{\R^3}\rho A_bF(V_1-d_N\partial_1\zeta)\dd x,\\
 \mathcal R_0={}&\epsilon_K\sum_{i=0}^3\mathscr P_N(P^i,v^i,z^i)
 -\frac{(\gamma-1)\theta_W}\tau\int_{\R^3}\rho\theta\zeta\Div\psi\dd x
 +\theta_\eta\Bigl[c\int_{\R^3}g_T^2\eta^2\dd x+\langle S_\eta,\eta\rangle\Bigr].
\end{aligned}
\]
One has
\begin{equation}\label{e:mrec}
 \int_s^tw_s|\mathcal R_0|\dd r\le\frac12\int_s^tw_sD_{\mathcal E}\dd r
 +\bigl[C_M(\delta+\mu_*)+\omega_b(s)\bigr]\bigl(I_s+N_o^2+\Phi\bigr),
\end{equation}
and, with
\begin{equation}\label{e:works}
\begin{aligned}
 Q_h={}&-\rho V_1^2-\frac{\gamma-1}{2\gamma p}\Pi^2
 +j_N\Bigl[(2-\gamma)V_1\partial_1\zeta-\gamma\sum_{i=2}^3V_i\partial_i\zeta\Bigr]
 +k_N\bigl[(\gamma-1)|\partial_1\zeta|^2+\gamma|\nabla'\zeta|^2\bigr]\\
 ={}&-\rho(V_1-d_N\partial_1\zeta)^2-\frac{\gamma-1}{2\gamma p}\Pi^2
 -\gamma j_NV\cdot\nabla\zeta+\bigl\{(\gamma-1)k_N+\rho d_N^2\bigr\}|\partial_1\zeta|^2
 +\gamma k_N|\nabla'\zeta|^2,\\
 \mathcal R_{\rm princ}={}&\int_{\R^3}h_bQ_h\dd x
 +\int_{\R^3}\frac{\gamma-1}\gamma d_N\frac{\partial_1p_b}{p_b}\Pi\partial_1\zeta\dd x
 +\int_{\R^3}\rho A_bF(V_1-d_N\partial_1\zeta)\dd x,
\end{aligned}
\end{equation}
\begin{equation}\label{e:mr}
 \int_s^tw_s|\mathcal R_A-\mathcal R_{\rm princ}|\dd r
 \le\bigl[C_M(\delta+\mu_*)+\omega_b(s)\bigr]\bigl(I_s+N_o^2+\Phi\bigr).
\end{equation}
\end{lemma}

\begin{proof}
By \eqref{e:espl}, $f_m=f_{m0}+\rho A_b(\eta+f)e_1$ and $f_z=f_{z0}+B_b(\eta+f)$ in
\eqref{e:ps3}, so \eqref{e:peq3} states that $(\Pi,\psi,\zeta)$ solves \eqref{e:fcsy}
with $c=N$ and $\mathbf f=\mathbf f_0^0+(0,\rho A_b(\eta+f)e_1,B_b(\eta+f))$. Since
$D_ub=A_b+h_b\psi_1$ and $\Div V=\Div\psi-\partial_1(b\eta)$, the triple
$(\Pi,V,\zeta)$ solves it with $\mathbf F^0$. By \eqref{e:fcal} with $\alpha=e_i$,
$D_u\psi=\dot\psi_0+A_b(\eta+f)e_1$ and
$\partial_i\{\rho A_b(\eta+f)\}-(\partial_i\rho)A_b(\eta+f)=\rho\partial_i\{A_b(\eta+f)\}$,
the triple $(P^i,v^i,z^i)$ solves it with $\mathbf F^i$. Adding \eqref{e:fcid} and
$\epsilon_K$ times \eqref{e:fcsc} for the four triples, $\theta_W$ times
\eqref{e:fcwv} for $(\Pi,V,\zeta)$ and $\theta_\eta$ times \eqref{e:etaE} gives
\eqref{e:mri}. Since $\mathscr R_N$ is affine in $\mathbf f$,
$F_z^0-(\gamma-1)\theta\Div V=f_{z0}+B_b(\eta+f)-(\gamma-1)\theta\Div\psi$,
$f=F+f_{\rm out}$ and $(\nabla J\Pi)_1=\partial_1J\Pi$, collecting the terms of the
sources containing $\eta$, $f_{\rm out}$ and $F$ gives \eqref{e:mdec}.

By Lemma~\ref{l:FC} and \eqref{e:stb} there is $c_*\in(0,1]$, depending only on
$M$ and the fixed data, with $e_N(P,v,z)\ge c_*(P^2+|v|^2+|\nabla z|^2)$,
$\mathcal D_N(v,z)\ge c_*(\Vert\nabla v\Vert_{L^2}^2+\Vert\nabla^2z\Vert_{L^2}^2)$ and
$\rho\ge c_*$. Since $|\int\rho v\cdot\nabla JP|\le C_M\Vert v\Vert_{L^2}\Vert P\Vert_{L^2}$,
we fix $\epsilon_K$ with $C_M\epsilon_K\le\frac{c_*}2$ and, in addition, so small that
$C_M\epsilon_K\le\frac{c_*}4$ for the constant $C_M$ of the bound of $\mathscr P_N$
below; then $\theta_W$ with $C_M\theta_W\le\frac{c_*}8$ for the constant below; then
$\theta_\eta$ with $C_M\theta_\eta\le\frac{\epsilon_K}4$ for the constant of
\eqref{e:Seta}. By \eqref{e:bsup} we increase $S$ so that
$\Vert b\Vert_{L^\infty}^2\le\frac{\theta_\eta}2$ on $[S,\infty)$. Then
$\frac{c_*}2\Vert V\Vert_{L^2}^2+\frac{c_*\theta_\eta}2\Vert\eta\Vert_{L^2}^2
\ge\frac{c_*}4(\Vert\psi\Vert_{L^2}^2+\theta_\eta\Vert\eta\Vert_{L^2}^2)$, which gives the
equivalence. The bound for $\mathcal E_N(s)$ follows from $\eta(s)=\chi(s)$,
$V(s)=\psi(s)-b\chi(s)e_1$ and $\Vert\chi(s)\Vert_{L^2}+\Vert\Pi(s)\Vert_{H^1}\le C_MI_s^{\frac12}$.

We prove \eqref{e:mrec}. By \eqref{e:Jsym},
$|\int\mathcal Lv\cdot\nabla JP|=|\int\mathbb S(v):\nabla^2JP|
\le C\Vert\nabla v\Vert_{L^2}\Vert\nabla J^{\frac12}P\Vert_{L^2}$, and
$\int\rho v\cdot\nabla Jg=-\int\Div(\rho v)Jg$ with
$\Vert J\{(\gamma-1)\kappa\Delta z-\gamma p\Div v\}\Vert_{L^2}
\le C_M(\Vert\nabla^2z\Vert_{L^2}+\Vert\nabla v\Vert_{L^2})$. Thus
\[
 |\mathscr P_N(P^i,v^i,z^i)|\le\frac18\Vert\nabla J^{\frac12}P^i\Vert_{L^2}^2
 +C_M\bigl(\Vert\nabla v^i\Vert_{L^2}^2+\Vert\nabla^2z^i\Vert_{L^2}^2
 +\Vert(\nabla\rho)\cdot v^i\Vert_{L^2}^2\bigr)
 +\Bigl|\int\rho v^i\cdot[D_u,\nabla J]P^i\Bigr| .
\]
For $i=0$, \eqref{e:etst} bounds $(\nabla\rho)\cdot V$; for $i\ge1$,
$\Vert(\nabla\rho)\cdot\partial_i\psi\Vert_{L^2}\le C(\delta+\omega_b(s))\Vert\nabla\psi\Vert_{L^2}
+C\mu_*\Vert\nabla\psi\Vert_{H^1}$ by \eqref{e:E1} and \eqref{e:prod}. In
$[D_u,\nabla J]=[u_b\partial_1,\nabla J]+[\psi\cdot\nabla,\nabla J]$, \eqref{e:Jcom}
with $a=u_b$, applied with $\Vert\partial_1\Pi\Vert_{L^2}$ for $i=0$ and with
$\Vert P^i\Vert_{L^2}$ for $i\ge1$, gives $\Vert[u_b\partial_1,\nabla J]P^i\Vert_{L^2}
\le C\Vert h_b\Vert_{L^\infty}\Vert\nabla\Pi\Vert_{L^2}$, which with
$\Vert h_b\Vert_{L^2(s,\infty;L^\infty)}^2\le\omega_b(s)$ by \eqref{e:E3} costs
$\omega_b(s)(N_o^2+\Phi)$.
For $i\ge1$, $\psi\cdot\nabla P=\Div(\psi P)-(\Div\psi)P$, \eqref{e:Jsym} and
$\Vert gh\Vert_{H^{-1}}\le C\Vert g\Vert_{H^1_{\rm uloc}}\Vert h\Vert_{L^2}$, the dual form
of \eqref{e:prod}, give $\Vert[\psi\cdot\nabla,\nabla J]P^i\Vert_{L^2}\le C\mu_*\Vert P^i\Vert_{L^2}$.
For $i=0$, integration by parts gives
\[
 \int\rho V\cdot[\psi\cdot\nabla,\nabla J]\Pi
 =-\int\Div(\rho V)\,\psi\cdot\nabla J\Pi-\int\rho V_k(\partial_k\psi_j)\partial_jJ\Pi
 +\int\Div(\rho V)J(\psi\cdot\nabla\Pi),
\]
which is at most $C\mu_*\Vert\Div(\rho V)\Vert_{L^2}\Vert\nabla\Pi\Vert_{L^2}
+C\Vert V\Vert_{L^6}\Vert\nabla\psi\Vert_{L^2}\Vert\nabla J\Pi\Vert_{L^3}$, and
$\Vert\nabla J\Pi\Vert_{L^3}\le C\Vert\Pi\Vert_{L^3}\le C\mu_*$ by \eqref{e:Jloc}. Next,
$|\frac{(\gamma-1)\theta_W}\tau\int\rho\theta\zeta\Div\psi|\le\frac{\theta_W}{8\tau^2}
\int\rho\zeta^2+C_M\theta_W(\Vert\nabla V\Vert_{L^2}^2+\Vert\nabla(b\eta)\Vert_{L^2}^2)$.
Finally \eqref{e:Seta}, \eqref{e:E1Z} and
$\sum_{i=0}^3\Vert\nabla J^{\frac12}P^i\Vert_{L^2}^2=\Vert\nabla\Pi\Vert_{L^2}^2$, which
follows from $\frac{|\xi|^2(1+|\xi|^2)}{1+|\xi|^2}=|\xi|^2$, bound the $\theta_\eta$
term by $\theta_\eta\frac c4\Vert\nabla\eta\Vert_{L^2_w}^2
+\frac{\epsilon_K}4\Vert\nabla\Pi\Vert_{L^2_w}^2
+[C_M(\delta+\mu_*)+\omega_b(s)](N_o^2+\Phi)$. With the choices of
$\epsilon_K,\theta_W,\theta_\eta$, \eqref{e:eta} and \eqref{e:etst}, these bounds give
\eqref{e:mrec}.

We prove \eqref{e:mr}. With
$\mathbf h^i=((\gamma-1)\partial_iQ_\mu(\nabla\psi),0,\partial_i\frac{Q_\mu(\nabla\psi)}{c_v\rho})$,
$\sum_{i=1}^3\{\mathscr S_N(P^i,v^i,z^i;\mathbf h^i)+\epsilon_K\mathscr K_N(P^i,v^i;
(\mathbf h^i)_p,0)\}=\mathcal W_{\rm heat}$, which is bounded by \eqref{e:hheat}. For
$1\le i\le3$, in $\mathscr S_N(P^i,v^i,z^i;\mathbf f_0^i-\mathbf h^i)$ and the screen
$\epsilon_K\mathscr K_N$ the momentum row of $\mathbf f_0^i-\mathbf h^i$ pairs with $W^i$,
the pressure row with $\frac{\partial_i\Pi}{\gamma p}$ and with $\epsilon_K\rho\partial_i\psi$
through $\nabla J$, and the thermal row with $Z^i$. By
$\nabla\rho=\nabla\phi+(\partial_1\rho_b)e_1$, $\nabla p=\nabla\Pi+(\partial_1p_b)e_1$,
$\nabla\theta=\nabla\zeta+(\partial_1\theta_b)e_1$ and
$\partial_iu=\mathbf 1_{\{i=1\}}h_be_1+\partial_i\psi$, each term of the rows is, up to a
factor bounded by $C_M$, either $ag$, with $a$ in the list of \eqref{e:E1Z} or
$a\in\{\partial_1\rho_b,\partial_1\theta_b\}$ and $g$ a derivative of order at most two of
$\Pi,\psi,\zeta$, or a product of two or three factors from these derivatives and
$\nabla\phi$, at most one of them of order two. The terms $ag$ are bounded in $L^2_w$
by \eqref{e:E1}, \eqref{e:E1Z}, \eqref{e:E2} and \eqref{e:E6}; the products by
H\"older's inequality on the cubes $Q_j$, with all factors but one bounded by
$C_M\mu_*$ in $L^\infty$ (order zero) or in $L^3(Q_j)$, $L^6(Q_j)$ (order one), and
by \eqref{e:pbd}, \eqref{e:prod}. For the inertia term $-(\partial_i\phi)\dot\psi_0$ we
use $|\langle\frac{\partial_i\phi}\rho g,W^i\rangle|\le C_M\mu_*\Vert g\Vert_{L^2}
\Vert W^i\Vert_{H^1}$ for $g=\mathcal L\psi-\nabla\Pi$, by \eqref{e:prod}, and for the
screened pressure terms $\Vert\nabla J\partial_ih\Vert_{L^2}\le\Vert h\Vert_{L^2}$. With
\eqref{e:etst}, the weighted integral of these pairings is at most
$[C_M(\delta+\mu_*)+\omega_b(s)](I_s+N_o^2+\Phi)$.

By the equations displayed after \eqref{e:peq3},
\[
 D_u\log p=-\gamma(h_b+\Div\psi)+\frac{\gamma-1}p\bigl\{\kappa(\partial_1^2\theta_b
 +\Delta\zeta)+\nu h_b^2+H[\psi]+Q_\mu(\nabla\psi)\bigr\},
\]
and $\nabla d_N=-d_N\frac{\nabla p}p$, $\nabla j_N=-j_N\frac{\nabla\theta}\theta$,
$\nabla k_N=-k_N(\frac{\nabla p}p+\frac{\nabla\theta}\theta)$. Hence each term of
$\mathscr G_N(P^i,v^i,z^i)$, $1\le i\le3$, is, up to a factor bounded by $C_M$, either
$ag_1g_1'$ with $a\in\{h_b,h_b^2,\partial_1p_b,\partial_1\theta_b,\partial_1^2\theta_b\}$, or
$g_1g_2g_3$, or $g_1g_2g_3g_4$, or $\Delta\zeta|\nabla^2\zeta|^2$, where $g_1,g_1'$ are
components of $\nabla\Pi,\nabla\psi,\nabla^2\psi,\nabla^2\zeta,\nabla^3\zeta$, $g_2$ of
$\nabla\psi,\nabla^2\zeta$, and $g_3,g_4$ of $\nabla\Pi,\nabla\psi,\nabla\zeta$. Here
$\Vert g_1\Vert_{L^2}^2+\Vert g_1'\Vert_{L^2}^2+\Vert g_2\Vert_{H^1}^2\le CD_*$,
$\Vert a\Vert_{L^\infty}\le C(\delta+\omega_b(s))$ by \eqref{e:E1}, \eqref{e:E2},
\eqref{e:E6}, and $\Vert g_3\Vert_{H^1(Q_j)}+\Vert g_4\Vert_{H^1(Q_j)}\le C_M\mu_*$ by
\eqref{e:stb}. H\"older's inequality on the cubes, with $g_1$ in $L^2(Q_j)$, $g_2$ in
$L^6(Q_j)$ and $g_3$ in $L^3(Q_j)$, resp. $g_3,g_4$ in $L^6(Q_j)$, gives
\[
 \int|g_1g_2g_3|+\int|g_1g_2g_3g_4|
 \le C_M\mu_*\sum_j\Vert g_1\Vert_{L^2(Q_j)}\Vert g_2\Vert_{H^1(Q_j)}
 \le C_M\mu_*D_* .
\]
With
$\Vert g\Vert_{L^3(Q_j)}\le C\Vert g\Vert_{L^2(Q_j)}^{\frac12}\Vert g\Vert_{H^1(Q_j)}^{\frac12}$
and $\Vert\nabla^2\zeta\Vert_{L^2(Q_j)}\le\mu_*$,
\[
\begin{aligned}
 \int|\Delta\zeta||\nabla^2\zeta|^2
 &\le C\sum_j\Vert\nabla^2\zeta\Vert_{L^2(Q_j)}^{\frac32}
 \Vert\nabla^2\zeta\Vert_{H^1(Q_j)}^{\frac32}\\
 &\le C\sup_j\Vert\nabla^2\zeta\Vert_{L^2(Q_j)}\sum_j\Vert\nabla^2\zeta\Vert_{H^1(Q_j)}^2
 \le C\mu_*\Vert\nabla^2\zeta\Vert_{H^1}^2 .
\end{aligned}
\]
Hence $\int_s^tw_s|\sum_{i=1}^3\mathscr G_N(P^i,v^i,z^i)|\le[C_M(\delta+\mu_*)+\omega_b(s)]N_o^2$.

For the zero-order work, $\psi_1=V_1+b\eta$ and $\rho d_N=j_N$ give
\[
\begin{aligned}
 (V-d_N\nabla\zeta)\cdot f_{m0}={}&-\rho h_bV_1^2+j_Nh_bV_1\partial_1\zeta
 -\rho h_bb\eta(V_1-d_N\partial_1\zeta)
 -\frac{\rho_bA_b}{\gamma p_b}\Pi(V_1-d_N\partial_1\zeta),\\
 \frac{\Pi f_{p0}}{\gamma p}={}&-\frac{\partial_1p_b}{\gamma p}\Pi(V_1+b\eta)
 -\frac{h_b\Pi^2}p+\frac{\gamma-1}{\gamma p}\Pi\bigl(H[\psi]+Q_\mu(\nabla\psi)\bigr).
\end{aligned}
\]
In \eqref{e:fcg}, the terms containing $D_u\log p$ are
\[
 (D_u\log p)\Bigl[-\frac{\Pi^2}{2\gamma p}+j_NV\cdot\nabla\zeta-k_N|\nabla\zeta|^2\Bigr].
\]
Using $\rho_bA_b=\nu\partial_1h_b-\partial_1p_b$ and
$\nabla d_N=-\frac{d_N}p(\nabla\Pi+(\partial_1p_b)e_1)$, we isolate the
principal work in three groups. First,
\[
 -\nu\int\frac{\partial_1h_b}{\gamma p_b}\Pi(V_1-d_N\partial_1\zeta)
 =\nu\int h_b\partial_1\Bigl[\frac{\Pi(V_1-d_N\partial_1\zeta)}{\gamma p_b}\Bigr];
\]
each term of the derivative is a zero-order quantity of $E_*$ or $b\eta$ times a
first derivative of $\Pi,V$, a term of $\nabla(d_N\partial_1\zeta)$ or
$\partial_1p_b$, and \eqref{e:E3}, \eqref{e:eta}, \eqref{e:pbd} bound its weighted
integral by $\omega_b(s)(N_o^2+\Phi)$. Second, the three terms with the factor
$\partial_1p_b$ combine, by
$\frac1{\gamma p_b}-\frac1{\gamma p}=\frac\Pi{\gamma pp_b}$, to
\[
 \frac{\gamma-1}\gamma d_N\frac{\partial_1p_b}{p_b}\Pi\partial_1\zeta
 -\frac{\partial_1p_b}{\gamma p}\Pi b\eta
 +\frac{\partial_1p_b}{\gamma pp_b}\Pi^2V_1
 -\frac{d_N\partial_1p_b}{p_bp}\Pi^2\partial_1\zeta ;
\]
the first term is in $\mathcal R_{\rm princ}$, and the weighted integrals of the
other three are at most
$C\Vert\partial_1p_b\Vert_{L^2(s,\infty;L^\infty)}\sup_r(\sqrt{w_s}\Vert\Pi\Vert_{L^2})
\Vert b\eta\Vert_{L^2_w}$ and
$C\mu_*\Vert\partial_1p_b\Vert_{L^2(s,\infty;L^\infty)}\Vert\nabla\Pi\Vert_{L^2_w}
\sup_r(\sqrt{w_s}\Vert(V,\nabla\zeta)\Vert_{L^2})$, by
$\int|\Pi^2g|\le\Vert\Pi\Vert_{L^3}\Vert\Pi\Vert_{L^6}\Vert g\Vert_{L^2}$.
Third, the terms with the factor $h_b$ sum to
\[
 h_b\Bigl[-\rho V_1^2-\frac{\gamma-1}{2\gamma p}\Pi^2-\gamma j_NV\cdot\nabla\zeta
 +2j_NV_1\partial_1\zeta+\gamma k_N|\nabla\zeta|^2-k_N|\partial_1\zeta|^2\Bigr]
 =h_bQ_h,
\]
since $-\frac1p+\frac1{2\gamma p}+\frac1{2p}=-\frac{\gamma-1}{2\gamma p}$; the second form of
$Q_h$ in \eqref{e:works} follows from
$-\rho(V_1-d_N\partial_1\zeta)^2=-\rho V_1^2+2j_NV_1\partial_1\zeta-\rho d_N^2|\partial_1\zeta|^2$.

The remaining zero-order terms are $-\rho h_bb\eta(V_1-d_N\partial_1\zeta)$,
$\frac{\gamma-1}{\gamma p}\Pi(H[\psi]+Q_\mu(\nabla\psi))$, $Z^0f_{z0}$,
$\frac{\theta_W}\tau\rho\zeta f_{z0}$, the screen $\epsilon_K\mathscr K_N(\Pi,V;f_{p0},f_{m0})$,
and in $\mathscr G_N(\Pi,V,\zeta)$ the terms with $\Div\psi$, $\partial_i\psi_l$,
$\frac{\gamma-1}p\{\kappa(\partial_1^2\theta_b+\Delta\zeta)+\nu h_b^2+H[\psi]+Q_\mu(\nabla\psi)\}$,
$\frac{d_N}p\Pi\nabla\Pi\cdot\nabla\zeta$ and $\nabla d_N,\nabla j_N,\nabla k_N$.
Their weighted integrals are at most $[C_M(\delta+\mu_*)+\omega_b(s)](I_s+N_o^2+\Phi)$.
Indeed, \eqref{e:E1Z}, \eqref{e:E2} and \eqref{e:gn4} give
$\Vert f_{p0}\Vert_{L^2_w}+\Vert f_{m0}\Vert_{L^2_w}+\Vert f_{z0}\Vert_{L^2_w}
\le[C_M(\delta+\mu_*)+\omega_b(s)]N_o$, which pairs with \eqref{e:etst},
$\Vert\tau^{-1}\zeta\Vert_{L^2_w}\le N_o$ and $\Vert\nabla J\Pi\Vert_{L^2}\le\Vert\nabla\Pi\Vert_{L^2}$;
the screened pressure term is $-\epsilon_K\int\Div(\rho V)Jf_{p0}$, with
$\Vert Jf_{p0}\Vert_{L^2}\le\Vert f_{p0}\Vert_{L^2}$. The coefficient
$\kappa\partial_1^2\theta_b+\nu h_b^2$ is treated by \eqref{e:curv} and \eqref{e:E2} as
in the proof of \eqref{e:Seta}, since
$\Vert\Pi\Vert_{L^2}^2+\Vert V\Vert_{L^2}^2+\Vert\nabla\zeta\Vert_{L^2}^2
\le C_M(E_*+\Vert b\eta\Vert_{L^2}^2)$. The terms with a
factor $h_b$, $\partial_1p_b$ or $\partial_1\theta_b$ and two factors from $q_o$, $b\eta$
and their derivatives are bounded by \eqref{e:E1Z}, \eqref{e:E3} and \eqref{e:eta}, as
in the first group, and the products of three or four such factors by H\"older's
inequality, \eqref{e:stb}, \eqref{e:pbd},
$\Vert\Pi\Vert_{L^4}^2\le\Vert\Pi\Vert_{L^3}\Vert\Pi\Vert_{L^6}$ and
$\Vert b\eta\nabla\zeta\Vert_{L^2_w}\le C_M\mu_*\Phi^{\frac12}$ from the proof of
\eqref{e:etst}. Together with the unscreened term, which is the same in
$\mathcal R_A$ and $\mathcal R_{\rm princ}$, this proves \eqref{e:mr}.
\end{proof}

\subsection{The acoustic correction}

We subtract a correction to the two outgoing sound modes before estimating their work.
With the derivatives of \eqref{e:der}, write $v_\perp=(v_2,v_3)$ for a vector $v$
and $\Delta_m=\partial_y^2+\Delta'$. Norms of functions of $(y,x')$ are taken over
$\R_y\times\R^2$ with $\dd y\dd x'$. Let $s_\sigma,v_\sigma,m_\sigma$ be as in
\eqref{e:mass}, $d_\sigma,\alpha$ as in \eqref{e:ac}, $\vartheta_\sigma$ as in
\eqref{e:CH}, and $Z=(V_b,b_m,\Theta_b)$, $\mathcal K$, $\Omega_\sigma$, $Z_{\sigma,a}$ as in
Lemma~\ref{l:mcoef}. For a triple $f=(P_f,v_f,z_f)$ and a function $F$ on
$\R_y\times\R^2$ put
\begin{equation}\label{e:acg}
\begin{gathered}
 U_\sigma=\frac{\gamma+1}2s_\sigma\vartheta_\sigma,\qquad
 W_\sigma=\frac{U_\sigma}{v_\sigma},\qquad
 \beta_M=\min_{\sigma=\pm1}e^{-\alpha|s_\sigma m_\sigma|},\qquad
 d_{t,\sigma}=\frac{v_\sigma(\mu+(\gamma-1)c)}2,\\
 G_\sigma=\Bigl(\frac{p_*}{P_b}\Bigr)^{\frac14}\Bigl(\frac{v_\sigma}{V_b}\Bigr)^{\frac34}
   \exp\Bigl(-\frac{s_\sigma b_m}{2p_*}\Bigr),\qquad
 \varphi_\sigma[f]=\frac12\Bigl(\frac{P_f\circ X}{\sqrt{\gamma P_b}}
   +\frac{\sigma(v_f)_1\circ X}{\sqrt{V_b}}\Bigr),\\
 J_bF=\bigl(J(F\circ Y)\bigr)\circ X,\qquad
 H_\sigma=\frac{U_\sigma}{V_bG_\sigma}J_b\varphi_\sigma[q_o],\\
 \partial_tu_\sigma^{\rm ac}+s_\sigma\partial_yu_\sigma^{\rm ac}
 -d_\sigma\partial_y^2u_\sigma^{\rm ac}-d_{t,\sigma}\Delta'u_\sigma^{\rm ac}=-H_\sigma,
 \qquad u_\sigma^{\rm ac}(s)=0.
\end{gathered}
\end{equation}
By Lemma~\ref{l:pr}, $|m_\pm|\le C_M$, so $\beta_M\in(0,1]$ depends only on
$M$ and the fixed data. With
\[
\begin{gathered}
 F_\sigma^u=\frac{G_\sigma}2\Bigl[\frac{(\gamma-1)c-\nu}{V_b}\partial_y^2u_\sigma^{\rm ac}
   +V_b\bigl((\gamma-1)c-\mu\bigr)\Delta'u_\sigma^{\rm ac}\Bigr],\qquad
 F_\sigma^e=cG_\sigma\Bigl(\frac{\partial_y^2u_\sigma^{\rm ac}}{V_b}+V_b\Delta'u_\sigma^{\rm ac}\Bigr),\\
 k_\sigma^u=\frac{\sigma(\mu+\lambda)G_\sigma}{\sqrt{\gamma P_bV_b}}\partial_y^2u_\sigma^{\rm ac},\qquad
 k_\sigma^v=\sigma k_\sigma^u+\sqrt{\frac{V_b}{\gamma P_b}}F_\sigma^u,\qquad
 k_\sigma^e=-\sigma\sqrt{\frac{V_b}{\gamma P_b}}F_\sigma^e,
\end{gathered}
\]
the correction is
\begin{equation}\label{e:gc}
\begin{aligned}
 g_c=(P_g,v_g,z_g)=\sum_\sigma\Bigl(&
 \sqrt{\gamma P_b}(G_\sigma\partial_yu_\sigma^{\rm ac}+k_\sigma^u),\
 \sqrt{V_b}(\sigma G_\sigma\partial_yu_\sigma^{\rm ac}+k_\sigma^v)e_1,\\
 &\frac{V_b\sqrt{\gamma P_b}}{c_p}
 (G_\sigma\partial_yu_\sigma^{\rm ac}+k_\sigma^u+k_\sigma^e)\Bigr)\circ Y,
\end{aligned}
\end{equation}
and we put
\begin{equation}\label{e:cor}
 q_s=q_o-g_c,\qquad N_s(s,t)=N[q_s](s,t),
\end{equation}
\begin{equation}\label{e:JH}
\begin{aligned}
 \mathcal J(s,t)=\sum_{\sigma=\pm1}\bigg\{&
 \sup_{s\le r\le t}w_s(r)
 \left(\Vert\nabla_m u_\sigma^{\rm ac}\Vert_{L^2}^2
       +\tau_r^{-1}\Vert u_\sigma^{\rm ac}\Vert_{L^2}^2\right)(r)\\
 &+\int_s^tw_s(r)\left(\Vert\nabla_m^2u_\sigma^{\rm ac}\Vert_{L^2}^2
       +\tau_r^{-1}\Vert\nabla_m u_\sigma^{\rm ac}\Vert_{L^2}^2
       +\tau_r^{-2}\Vert u_\sigma^{\rm ac}\Vert_{L^2}^2\right)(r)\dd r\bigg\},\\
 \mathcal H(s,t)=\sum_{\sigma=\pm1}\bigg\{&
 \sup_{s\le r\le t}w_s(r)\sum_{j=2}^4
                 \Vert\nabla_m^ju_\sigma^{\rm ac}(r)\Vert_{L^2}^2
 +\int_s^tw_s(r)\sum_{j=3}^5
                 \Vert\nabla_m^ju_\sigma^{\rm ac}(r)\Vert_{L^2}^2\dd r\bigg\},
\end{aligned}
\end{equation}
\begin{equation}\label{e:B}
 \mathfrak B(s,t)=I_s+N_o^2(s,t)+N_s^2(s,t)+\Phi(s,t)
                         +\mathcal J(s,t)+\mathcal H(s,t).
\end{equation}

\begin{lemma}\label{l:cor}
One has $u_\sigma^{\rm ac}\in C([s,T_e];H^4)\cap L^2(s,T_e;H^5)$,
$g_c\in C([s,T_e];H^2)$, $(v_g,z_g)\in L^2(s,T_e;H^3)$, $g_c(s)=0$,
\begin{equation}\label{e:gap}
 \tau\partial_yW_\sigma\le1-\beta_M,\qquad
 \Bigl\Vert\frac{U_\sigma}{V_b}-W_\sigma\Bigr\Vert_{L^\infty}\le\frac{C_b}\tau,\qquad
 \partial_tW_\sigma+s_\sigma\partial_yW_\sigma+W_\sigma\partial_yW_\sigma
 =d_\sigma\partial_y^2W_\sigma,
\end{equation}
and
\[
\begin{gathered}
 \mathcal J\le C_MN_o^2,\qquad
 \mathcal H\le\omega_b(s)N_o^2,\qquad
 N[g_c]^2\le C_M\mathcal J+\omega_b(s)N_o^2,\\
 N_o^2+N[g_c]^2+\mathcal J+\mathcal H\le C_MN_s^2,
 \qquad N_s^2\le C_MN_o^2,\\
 E[g_c](r)\le C_M\sum_{\sigma}\sum_{j=1}^4
 \Vert\nabla_m^ju_\sigma^{\rm ac}(r)\Vert_{L^2}^2,\qquad
 \int_s^tw_s\Vert\nabla z_g\Vert_{L^2}^2\dd r
 \le C_M(\mathcal J+\mathcal H)+\omega_b(s)N_o^2.
\end{gathered}
\]
Let $\ell>0$, $P_0\ge0$, and let $\varrho(y,r)\in[0,1]$ be smooth, with
$\varrho=0$ where $|y+h_*-s_\sigma r|\le\ell\tau_r$ and
$\tau_r^j|\partial_y^j\varrho|\le C_\varrho$ for $j\le5$. Then
\[
 \sup_{s\le r\le t}w_s\tau_r^{P_0}
 \Vert\varrho u_\sigma^{\rm ac}(r)\Vert_{H^4}^2
 +\int_s^tw_s\tau_r^{P_0}
 \Vert\varrho u_\sigma^{\rm ac}\Vert_{H^5}^2\dd r
 \le C_{b,\ell,P_0,C_\varrho}e^{-c_\ell s}N_o^2,
\]
with $c_\ell>0$ depending only on $\ell$, $M$ and the fixed data. Finally,
\[
 \sup_{s\le t\le T_e}\bigl(\Vert(P_g,v_g,\nabla z_g)\Vert_{L^3}
 +\Vert(P_g,v_g,\nabla z_g)\Vert_{L^\infty}\bigr)\le C_b\mu_*,
\]
and the same two norms over $\{x:|Y(x_1,t)+h_*|\le2c_0\tau\}$ are at most
$C_be^{-c_bt}\mu_*$.
\end{lemma}

\begin{proof}
Write $u=u_\sigma^{\rm ac}$ when $\sigma$ is fixed. In $x=y-s_\sigma t$, by
\eqref{e:CH}, $p_\sigma=\int_\R\mathsf G(x-z)\Phi(z)\dd z$ with
$\mathsf G(x)=(4\pi d_\sigma\tau)^{-\frac12}e^{-\frac{x^2}{4d_\sigma\tau}}$,
$\Phi=1$ on $(-\infty,0)$ and $\Phi=e^{-\alpha s_\sigma m_\sigma}$ on $(0,\infty)$.
Put $\langle h\rangle=\frac1{p_\sigma}\int_\R h(z)\mathsf G(x-z)\Phi(z)\dd z$.
Since $a_\sigma=-\frac{(p_\sigma)_y}{\alpha p_\sigma}$ and $\alpha d_\sigma=k_\sigma$,
$W_\sigma=\frac{\gamma+1}{2v_\sigma}a_\sigma=-2d_\sigma\frac{(p_\sigma)_y}{p_\sigma}$, and
differentiation under the integral gives
\[
 W_\sigma=\frac{x-\langle z\rangle}\tau,\qquad
 \tau\partial_yW_\sigma=1-\frac{\langle z^2\rangle-\langle z\rangle^2}{2d_\sigma\tau}.
\]
As $e^{-\alpha|s_\sigma m_\sigma|}\max\Phi\le\Phi\le\max\Phi$ and $\int_\R\mathsf G=1$,
\[
 \langle z^2\rangle-\langle z\rangle^2=\min_{\xi\in\R}\langle(z-\xi)^2\rangle
 \ge e^{-\alpha|s_\sigma m_\sigma|}\min_{\xi\in\R}\int_\R(z-\xi)^2\mathsf G(x-z)\dd z
 =2d_\sigma\tau e^{-\alpha|s_\sigma m_\sigma|},
\]
which is the first assertion of \eqref{e:gap}. Multiplying the equation of
$\vartheta_\sigma$ in Lemma~\ref{l:pr} by $2k_\sigma s_\sigma$ gives the third,
since $W_\sigma=2k_\sigma s_\sigma\vartheta_\sigma$. By \eqref{e:CH},
$V_b=V_c(y+h_*,t)-\sum_\iota\vartheta_\iota+(R_U)_1$. One of
$|y+h_*|,|y-s_\sigma t|$ is at least $\frac{|s_\sigma|t}4$ for $t\ge4|h_*|/|s_\sigma|$,
so Lemma~\ref{l:pr} and \eqref{eq:ct-gau} give
$|U_\sigma||V_c(y+h_*,t)-v_\sigma|\le C_Me^{-c_bt}$; with Lemma~\ref{l:asy} this
gives $|U_\sigma||V_b-v_\sigma|\le C_b\tau^{-1}$, the second assertion.

Multiply \eqref{e:acg} by $\frac{w_s}\tau u$, and by $-w_s\Delta_mu$ after
applying $\nabla_m^j$, $j=0,1,2,3$. The transport term integrates to zero,
$\frac{\dd}{\dd t}w_s=-\frac{q_aw_s}\tau$, and with
$d_0=\min_\sigma\min\{d_\sigma,d_{t,\sigma}\}$,
\[
\begin{aligned}
 &\frac{\dd}{\dd t}\Bigl(\frac{w_s}{2\tau}\Vert u\Vert_{L^2}^2\Bigr)
 +\frac{(1+q_a)w_s}{2\tau^2}\Vert u\Vert_{L^2}^2
 +\frac{d_0w_s}\tau\Vert\nabla_mu\Vert_{L^2}^2
 \le\frac{w_s}{\tau}\Bigl|\int H_\sigma u\Bigr|,\\
 &\frac{\dd}{\dd t}\Bigl(\frac{w_s}2\Vert\nabla_m^{j+1}u\Vert_{L^2}^2\Bigr)
 +d_0w_s\Vert\nabla_m^{j+2}u\Vert_{L^2}^2
 \le w_s\Vert\nabla_m^jH_\sigma\Vert_{L^2}\Vert\nabla_m^{j+2}u\Vert_{L^2}.
\end{aligned}
\]
Young's inequality and integration over $[s,t]$ give
\[
 \mathcal J\le C\sum_\sigma\int_s^tw_s\Vert H_\sigma\Vert_{L^2}^2\dd r,\qquad
 \mathcal H\le C\sum_\sigma\sum_{j=1}^3\int_s^tw_s\Vert\nabla_m^jH_\sigma\Vert_{L^2}^2\dd r .
\]
We show
\begin{equation}\label{e:acsrc}
 \Vert H_\sigma\Vert_{L^2}^2\le\frac{C_M}\tau\Vert(\Pi,\psi)\Vert_{L^2}^2,\qquad
 \sum_{j=1}^3\Vert\nabla_m^jH_\sigma\Vert_{L^2}^2
 \le\frac{C_M}\tau\Vert\nabla(\Pi,\psi)\Vert_{L^2}^2
 +\frac{\epsilon_b(t)}\tau\Vert(\Pi,\psi)\Vert_{L^2}^2,
\end{equation}
with $\sup_{r\ge s}\epsilon_b(r)\le\omega_b(s)$. By \eqref{e:chain}, $\nabla_m^j$ of
$(Jh)\circ X$ is a sum of $(\nabla^iJh)\circ X$, $1\le i\le j$, with coefficients
built from $V_b,\partial_yV_b,\partial_y^2V_b$, and
$\Vert\nabla^iJh\Vert_{L^2}\le\Vert\nabla h\Vert_{L^2}$ for $1\le i\le3$ by
\eqref{e:Jsym}. The coefficient $\frac{U_\sigma}{V_bG_\sigma}$ has
$|\partial_y^iU_\sigma|\le C_M\tau^{-\frac{i+1}2}$ by Lemma~\ref{l:pr}; its
factor $\frac1{V_bG_\sigma}$ and the coefficients of $\varphi_\sigma$ are smooth
functions of $Z\in\mathcal K$ whose first two $y$-derivatives tend to zero in
$L^\infty$ by \eqref{e:munw}. The only term with $\partial_y^3Z$ is bounded by
\eqref{eq:planar}:
\[
 \Vert C(Z)\partial_y^3Z\,U_\sigma J_b\varphi_\sigma[q_o]\Vert_{L^2}^2
 \le\frac{C_M}\tau\Vert\partial_y^3Z\Vert_{L^2(\R)}^2
 \Vert J_b\varphi_\sigma[q_o]\Vert_{L^2}\Vert\partial_yJ_b\varphi_\sigma[q_o]\Vert_{L^2}
 \le\frac{C_M}\tau\Vert\partial_y^3Z\Vert_{L^2(\R)}^2\Vert(\Pi,\psi)\Vert_{L^2}^2,
\]
and $\Vert\partial_y^3Z\Vert_{L^2(\R)}\to0$ by \eqref{e:munw}. This proves
\eqref{e:acsrc}. Since $\int_s^tw_s\tau^{-1}\Vert(\Pi,\psi)\Vert_{H^1}^2\le
q_a^{-1}N_o^2$ and $\int_s^tw_s\Vert\nabla(\Pi,\psi)\Vert_{L^2}^2\le N_o^2$, we
obtain $\mathcal J\le C_MN_o^2$ and $\mathcal H\le(C_M\tau_s^{-1}+C\sup_{r\ge s}
\epsilon_b(r))N_o^2\le\omega_b(s)N_o^2$. With
$\mathcal D_\sigma=d_\sigma\partial_y^2+d_{t,\sigma}\Delta'$,
\[
 u(t)=-\int_s^te^{(t-r)\mathcal D_\sigma}H_\sigma(y-s_\sigma(t-r),x',r)\dd r .
\]
The same tests on frequency truncations, with \eqref{e:bgcl}, give
$u\in C([s,T_e];H^4)\cap L^2(s,T_e;H^5)$ and $u(s)=0$ in $H^4$.
The regularity and initial value of $g_c$ follow from \eqref{e:gc}.

Let $g_c^{(2)}$ be \eqref{e:gc} with the terms containing
$\partial_yu_\sigma^{\rm ac}$ omitted, and $g_c^{(1)}=g_c-g_c^{(2)}$. By
\eqref{e:chain}, the derivatives of $P_g,v_g,z_g$ entering $N[g_c]$ are sums of
$C(Z)\nabla_m^{j}u_\sigma^{\rm ac}\circ Y$ and of terms in which some of the
derivatives fall on $C(Z)$ or on $Y$. Here $j\le4$ in $g_c^{(1)}$ and $2\le j\le5$
in $g_c^{(2)}$, the top orders being $\nabla^3z_g$ and $\nabla^2v_g$. By
\eqref{e:E1}, \eqref{e:E2}, \eqref{e:chain} and \eqref{e:munw},
$\sqrt\tau\Vert\partial_yZ\Vert_{L^\infty}\le C_M$ for $t\ge s$ and
$\Vert\partial_yZ\Vert_{L^\infty}+\Vert\partial_y^2Z\Vert_{L^\infty}\to0$. The
third derivatives are bounded by \eqref{eq:planar}: for $j\ge1$,
\[
 \int_s^tw_s\Vert C(Z)\partial_y^3Z\,\nabla_m^ju\Vert_{L^2}^2\dd r
 \le C_M\sup_{s\le r\le t}\bigl(w_s\Vert\nabla_m^ju\Vert_{L^2}
 \Vert\nabla_m^{j+1}u\Vert_{L^2}\bigr)\int_s^\infty\Vert\partial_y^3Z\Vert_{L^2(\R)}^2\dd r,
\]
and the last integral tends to zero by \eqref{e:munw}. Every term of $N[g_c^{(1)}]^2$
with exactly $j$ derivatives on $u$ is then one of the terms of $\mathcal J$
(for $j\le2$, the factor $\tau^{-1}$ being supplied by $|\partial_yZ|^2\le
C_M\tau^{-1}$ when one derivative falls on $C(Z)$) or of $\mathcal H$, up to a
factor $C_M$, plus terms bounded by $\omega_b(s)(\mathcal J+\mathcal H)$; every term
of $N[g_c^{(2)}]^2$ is a term of $\mathcal H$ or carries a vanishing coefficient or
a factor $\tau_s^{-1}$. Hence, by the bounds already proved,
\[
 N[g_c^{(2)}]^2\le\omega_b(s)N_o^2,\qquad
 N[g_c]^2\le C_M\mathcal J+\omega_b(s)N_o^2,\qquad
 \int_s^tw_s\Vert\nabla z_g\Vert_{L^2}^2\dd r\le C_M(\mathcal J+\mathcal H)+\omega_b(s)N_o^2,
\]
and two derivatives of \eqref{e:gc} give the bound of $E[g_c](r)$.

For the reverse bound, \eqref{e:gc} gives
$\frac{P_g\circ X}{\sqrt{\gamma P_b}}=\sum_\iota(G_\iota\partial_yu_\iota^{\rm ac}+k_\iota^u)$
and $\frac{(v_g)_1\circ X}{\sqrt{V_b}}=\sum_\iota(\iota G_\iota\partial_yu_\iota^{\rm ac}+k_\iota^v)$,
so that
\[
 \varphi_\sigma[g_c]=\frac12\sum_{\iota=\pm1}\bigl[(1+\sigma\iota)G_\iota\partial_yu_\iota^{\rm ac}
 +k_\iota^u+\sigma k_\iota^v\bigr]
 =G_\sigma\partial_yu+\frac12\sum_\iota(k_\iota^u+\sigma k_\iota^v).
\]
Inserting $q_o=q_s+g_c$ into $H_\sigma$ and writing $J_b=I+(J_b-I)$ in front of
$G_\sigma\partial_yu$, \eqref{e:acg} becomes
\begin{equation}\label{e:crev}
 \partial_tu+\Bigl(s_\sigma+\frac{U_\sigma}{V_b}\Bigr)\partial_yu
 -d_\sigma\partial_y^2u-d_{t,\sigma}\Delta'u
 =-\frac{U_\sigma}{V_bG_\sigma}\Bigl[J_b\varphi_\sigma[q_s]
 +(J_b-I)(G_\sigma\partial_yu)
 +\frac12J_b\sum_\iota(k_\iota^u+\sigma k_\iota^v)\Bigr].
\end{equation}
Multiply \eqref{e:crev} by $\frac{w_s}\tau u$. By \eqref{e:gap},
\[
 \int\frac{U_\sigma}{V_b}u\,\partial_yu
 =-\frac12\int\partial_yW_\sigma\,u^2+\int\Bigl(\frac{U_\sigma}{V_b}-W_\sigma\Bigr)u\,\partial_yu
 \ge-\frac{1-\beta_M}{2\tau}\Vert u\Vert_{L^2}^2-\frac{C_b}\tau\Vert u\Vert_{L^2}\Vert\partial_yu\Vert_{L^2}.
\]
For the two feedback terms, put $\mathfrak a=\rho_b\bigl(\frac{U_\sigma}{V_bG_\sigma}u\bigr)\circ Y$;
then $|\mathfrak a|\le C_M\tau^{-\frac12}|u\circ Y|$ and, by \eqref{e:munw},
$\Vert\nabla\mathfrak a\Vert_{L^2}\le C_b(\tau^{-\frac12}\Vert\nabla_mu\Vert_{L^2}
+\tau^{-1}\Vert u\Vert_{L^2})$. Since $\dd y=\rho_b\dd x_1$, $J$ is self-adjoint and
$J-I=\Div\nabla J$,
\[
 \int\frac{U_\sigma}{V_bG_\sigma}u\,(J_b-I)(G_\sigma\partial_yu)\dd y\dd x'
 =-\int_{\R^3}\nabla\mathfrak a\cdot\nabla J\bigl((G_\sigma\partial_yu)\circ Y\bigr)\dd x .
\]
Each of $k_\iota^u,k_\iota^v$ is a sum of $c(y,t)\partial_y^2u_\iota^{\rm ac}$ and
$c(y,t)\Delta'u_\iota^{\rm ac}$ with $|c|+\sqrt\tau|\partial_yc|\le C_b$; after
composition with $Y$ it is $\partial_1(\tilde c\,\partial_1u_\iota^E)+
\nabla'\cdot(\tilde c\,\nabla'u_\iota^E)+\tilde c'\partial_1u_\iota^E$, $u_\iota^E=u_\iota^{\rm ac}\circ Y$,
with $|\tilde c|\le C_b$, $|\tilde c'|\le C_b\tau^{-\frac12}$. Moving $J$ onto $\mathfrak a$ and
one derivative onto $J\mathfrak a$, and using \eqref{e:Jsym}, we bound both
feedback pairings by
\[
 C_b\bigl(\tau^{-\frac12}\Vert\nabla_mu\Vert_{L^2}+\tau^{-1}\Vert u\Vert_{L^2}\bigr)
 \sum_\iota\Vert\nabla_mu_\iota^{\rm ac}\Vert_{L^2}.
\]
After the factor $\frac{w_s}\tau$ this is at most
$C_b\tau_s^{-\frac12}\sum_\iota w_s(\tau^{-1}\Vert\nabla_mu_\iota^{\rm ac}\Vert_{L^2}^2
+\tau^{-2}\Vert u_\iota^{\rm ac}\Vert_{L^2}^2)$. The term with $q_s$ is at most
$\frac{q_aw_s}{8\tau^2}\Vert u\Vert_{L^2}^2+\frac{C_Mw_s}\tau\Vert(P_s,v_s)\Vert_{L^2}^2$, and the
drift remainder at most $\frac{d_0w_s}{2\tau}\Vert\partial_yu\Vert_{L^2}^2
+\frac{C_bw_s}{\tau_s\tau^2}\Vert u\Vert_{L^2}^2$. The left side retains
$\frac{(q_a+\beta_M)w_s}{2\tau^2}\Vert u\Vert_{L^2}^2+\frac{d_0w_s}{\tau}\Vert\nabla_mu\Vert_{L^2}^2$.
Summing over $\sigma$, integrating and increasing $S$,
\[
 \sum_\sigma\Bigl[\sup_{s\le r\le t}\frac{w_s}\tau\Vert u_\sigma^{\rm ac}\Vert_{L^2}^2
 +\int_s^tw_s\bigl(\tau^{-2}\Vert u_\sigma^{\rm ac}\Vert_{L^2}^2
 +\tau^{-1}\Vert\nabla_mu_\sigma^{\rm ac}\Vert_{L^2}^2\bigr)\dd r\Bigr]
 \le C_M\int_s^t\frac{w_s}\tau\Vert(P_s,v_s)\Vert_{L^2}^2\dd r\le C_MN_s^2 .
\]
Multiply \eqref{e:crev} by $-w_s\Delta_mu$. The drift term is at most
$C_Mw_s\tau^{-\frac12}\Vert\nabla_mu\Vert_{L^2}\Vert\nabla_m^2u\Vert_{L^2}$, and by
\eqref{e:Jsym} and the same form of $k_\iota^u,k_\iota^v$ the right side of
\eqref{e:crev} has $L^2$ norm at most
$C_M\tau^{-\frac12}(\Vert(P_s,v_s)\Vert_{L^2}+\sum_\iota\Vert\nabla_mu_\iota^{\rm ac}\Vert_{L^2})$.
With the preceding bound this gives $\mathcal J\le C_MN_s^2$. Then
\[
 N_o^2\le2N_s^2+2N[g_c]^2\le C_MN_s^2+\omega_b(s)N_o^2,
\]
and increasing $S$ gives $N_o^2\le C_MN_s^2$, hence the bounds of $N[g_c]^2$ and
$\mathcal H$ by $C_MN_s^2$. Finally $N_s^2\le2N_o^2+2N[g_c]^2\le C_MN_o^2$.

For the exponential bounds put $\Psi_\pm=\pm\lambda(y-s_\sigma t)-3d_\sigma\lambda^2\tau$
with $0<\lambda\le\min\{\frac14,\frac\ell{12d_\sigma}\}$. Then $u_\pm=e^{\Psi_\pm}u$ solves
\[
 \partial_tu_\pm+(s_\sigma\pm2d_\sigma\lambda)\partial_yu_\pm-d_\sigma\partial_y^2u_\pm
 -d_{t,\sigma}\Delta'u_\pm+2d_\sigma\lambda^2u_\pm=-e^{\Psi_\pm}H_\sigma .
\]
Since $\pm\lambda x-\frac{x^2}{8d_\sigma\tau}\le2d_\sigma\lambda^2\tau$, Lemma~\ref{l:pr} gives
$|\partial_y^i(e^{\Psi_\pm}U_\sigma)|\le C_M\tau^{-\frac12}e^{-d_\sigma\lambda^2\tau}$ for
$i\le3$, and the proof of \eqref{e:acsrc} gives
$\Vert e^{\Psi_\pm}H_\sigma\Vert_{H^3}^2\le C_be^{-d_\sigma\lambda^2\tau}\tau^{-1}
\Vert(\Pi,\psi)\Vert_{H^1}^2$. The tests above, the damping $2d_\sigma\lambda^2$ paying the
zero-order terms, give
\[
 \sup_{s\le r\le t}w_s\Vert u_\pm\Vert_{H^4}^2+\int_s^tw_s\Vert u_\pm\Vert_{H^5}^2\dd r
 \le C_{b,\lambda}e^{-d_\sigma\lambda^2\tau_s}N_o^2 .
\]
Where $\varrho\ne0$ and $\tau\ge\frac{2|h_*|}\ell$, one of $\Psi_\pm$ is at least
$\lambda(\ell\tau-|h_*|)-3d_\sigma\lambda^2\tau\ge\frac{\lambda\ell\tau}4$. Since
$e^{\Psi_\pm}\partial^\beta u$ is a combination of $\partial^{\beta'}u_\pm$, $\beta'\le\beta$,
with coefficients bounded by $C\lambda^{|\beta-\beta'|}$, and
$\tau^{P_0}e^{-\frac{\lambda\ell\tau}2}\le C_{P_0,\lambda\ell}e^{-\frac{\lambda\ell\tau}4}$,
the cone bound follows with $c_\ell=\min\{\frac{\lambda\ell}4,d_\sigma\lambda^2\}$.

For the local bounds, \eqref{e:chain}, \eqref{e:Jloc} and \eqref{e:stb}, together
with $\Vert\Pi\Vert_{L^3}+\Vert\psi\Vert_{L^3}\le\mu_*$, give for $p=3,\infty$
\[
 \sum_{j\le2}\Vert\nabla_m^jH_\sigma\Vert_{L^p}\le C_b\mu_*\tau^{-\frac12},\qquad
 \sum_{j\le2}\Vert\nabla_m^j(e^{\Psi_\pm}H_\sigma)\Vert_{L^p}\le C_b\mu_*\tau^{-\frac12},
\]
the derivatives of $J(\varphi_\sigma[q_o]\circ Y)$ being estimated without weight.
Put one derivative on the heat kernel, $\Vert\nabla_me^{\theta\mathcal D_\sigma}\Vert_{L^p\to L^p}
\le C\theta^{-\frac12}$, and the others on the source; since
$\int_s^t(t-r)^{-\frac12}\tau_r^{-\frac12}\dd r\le\pi$, the Duhamel formula gives
$\Vert\nabla_m^ju\Vert_{L^p}+\Vert\nabla_m^ju_\pm\Vert_{L^p}\le C_b\mu_*$ for
$1\le j\le3$, and the damping gives the same for $j=0$ and $u_\pm$. By
\eqref{e:gc}, $(P_g,v_g,\nabla z_g)$ is a sum of bounded coefficients times
$\nabla_m^ju_\sigma^{\rm ac}\circ Y$, $1\le j\le3$, which proves the first local
bound. If $|y+h_*|\le2c_0\tau$, then $|s_\sigma|>8c_0$ gives
$-\sigma(y-s_\sigma t)\ge\frac{|s_\sigma|\tau}2$ for $\tau\ge4(1+\frac{|h_*|}{|s_\sigma|})$, so
$\Psi_{-\sigma}\ge\frac{\lambda|s_\sigma|\tau}4$ for $\lambda\le\frac{|s_\sigma|}{12d_\sigma}$,
which gives the factor $e^{-c_bt}$.
\end{proof}

By \eqref{e:eta} and Lemma~\ref{l:cor},
\begin{equation}\label{e:Bbd}
 \mathfrak B\le C_M\bigl(I_s+N_o^2\bigr).
\end{equation}

For a mass vector $Q=(n,j,j_2,j_3,e)$, with $j_\perp=(j_2,j_3)$, put
\[
\begin{gathered}
 v_Q=j-b_mn,\qquad
 z_Q=\frac{e-b_mj+(b_m^2-E_m)n}{c_v},\qquad
 \pi_Q=\frac{\gamma-1}{V_b}\Bigl(e-b_mj+\frac{b_m^2n}{2}\Bigr),\\
 \mathcal PQ=(\pi_Q,v_Q,j_2,j_3,z_Q)\circ Y,
\end{gathered}
\]
and let $L_bQ=\partial_y\mathcal NQ+\mathcal TQ$, where
\[
\begin{aligned}
 (\mathcal NQ)_n={}&\frac{v_Q}{V_b},\qquad
 (\mathcal NQ)_j=\frac{b_mv_Q}{V_b}+\pi_Q
 -\frac{\nu(v_Q)_y}{V_b}-\lambda\Div'j_\perp,\qquad
 (\mathcal NQ)_\perp=-\mu\Bigl(\frac{(j_\perp)_y}{V_b}+\nabla'v_Q\Bigr),\\
 (\mathcal NQ)_e={}&\Bigl(\frac{E_m}{V_b}+P_b\Bigr)v_Q+b_m\pi_Q
 -b_m\Bigl(\frac{\nu(v_Q)_y}{V_b}+\lambda\Div'j_\perp\Bigr)
 -\frac{\nu(b_m)_yv_Q+\kappa(z_Q)_y}{V_b},\\
 (\mathcal TQ)_n={}&\Div'j_\perp,\qquad
 (\mathcal TQ)_j=b_m\Div'j_\perp-\mu\Div'(j_\perp)_y-\mu V_b\Delta'v_Q,\\
 (\mathcal TQ)_\perp={}&V_b\nabla'\pi_Q-\lambda\nabla'(v_Q)_y
 -\mu V_b\Delta'j_\perp-(\mu+\lambda)V_b\nabla'\Div'j_\perp,\\
 (\mathcal TQ)_e={}&\bigl(E_m+P_bV_b-\lambda(b_m)_y\bigr)\Div'j_\perp
 -b_m\mu\Div'(j_\perp)_y-b_m\mu V_b\Delta'v_Q-\kappa V_b\Delta'z_Q .
\end{aligned}
\]
We write $L_b=L_E-L_V$, where $L_E$ collects the terms of $L_b$ without
$\nu,\mu,\lambda,\kappa$. Put
\[
\begin{gathered}
 R_\sigma=\Bigl(\frac1{\sqrt{\gamma P_b}},
 \frac{b_m}{\sqrt{\gamma P_b}}+\sigma\sqrt{V_b},0,0,
 \frac{E_m+P_bV_b}{\sqrt{\gamma P_b}}+\sigma b_m\sqrt{V_b}\Bigr)^\top,\\
 \mathcal D_{\rm mat}[f]=\sum_\sigma
 \partial_y\Bigl(\frac{R_\sigma U_\sigma}{V_b}J_b\varphi_\sigma[f]\Bigr),\qquad
 (T_{\rm inv})_\perp=-\sum_\sigma V_b\sqrt{\gamma P_b}
 G_\sigma\nabla'\partial_yu_\sigma^{\rm ac},\\
 (T_{\rm inv})_n=(T_{\rm inv})_j=(T_{\rm inv})_e=0,\qquad
 R_g=-(\partial_t+L_b)(\mathcal P^{-1}g_c)-\mathcal D_{\rm mat}[q_s]-T_{\rm inv}.
\end{gathered}
\]
For a mass source $\mathcal F=(\mathcal F_n,\mathcal F_j,\mathcal F_2,\mathcal F_3,\mathcal F_e)$
its physical rows are
\begin{equation}\label{e:conv}
\begin{aligned}
 (\mathfrak P\mathcal F)_p
 &=\Bigl[\frac{\gamma-1}{V_b}
 \Bigl(\mathcal F_e-b_m\mathcal F_j+\frac{b_m^2\mathcal F_n}{2}\Bigr)\Bigr]\circ Y,\qquad
 (\mathfrak P\mathcal F)_m
 =\Bigl[\frac1{V_b}(\mathcal F_j-b_m\mathcal F_n,
 \mathcal F_2,\mathcal F_3)\Bigr]\circ Y,\\
 (\mathfrak P\mathcal F)_z
 &=\Bigl[\frac{\mathcal F_e-b_m\mathcal F_j+(b_m^2-E_m)\mathcal F_n}{c_v}\Bigr]\circ Y,\qquad
 (\mathfrak P\mathcal F)_\chi
 =\frac{(\mathfrak P\mathcal F)_z}{\theta_b}
 -\frac{(\mathfrak P\mathcal F)_p}{c_p\rho_b\theta_b}.
\end{aligned}
\end{equation}
With $H[\psi]$ and $Q_\mu$ of \eqref{e:Lam}, the nonlinear sources are
\begin{equation}\label{e:FN}
\begin{aligned}
 F_p^N={}&-\psi\cdot\nabla\Pi-\gamma\Pi\Div\psi
 +(\gamma-1)Q_\mu(\nabla\psi),\qquad
 F_m^N=\frac\phi\rho(\nabla\Pi-\mathcal L\psi)
 -\rho_b\psi\cdot\nabla\psi
 +\rho_bA_b\frac\phi\rho\frac\Pi{p_b}e_1,\\
 F_z^N={}&-\psi\cdot\nabla\zeta-(\gamma-1)\zeta\Div\psi
 -\frac{\kappa\phi}{c_v\rho_b\rho}\Delta\zeta
 +B_b\frac\phi\rho\frac\Pi{p_b}
 +\frac\phi{c_v\rho_b\rho}(h_b\Pi-H[\psi])
 +\frac{Q_\mu(\nabla\psi)}{c_v\rho},\\
 F_\rho^N={}&\frac{F_p^N-R\rho_bF_z^N}{R\theta_b},
\end{aligned}
\end{equation}
\[
 S_N=\Bigl(\frac{F_\rho^N}{\rho_b},
 \frac{F_{m,1}^N+bF_\rho^N}{\rho_b},
 \frac{F_{m,\perp}^N}{\rho_b},
 c_vF_z^N+\frac{bF_{m,1}^N+(c_v\theta_b+\frac{b^2}2)F_\rho^N}{\rho_b}
 \Bigr)\circ X,
\]
and, in the second convention,
\begin{equation}\label{e:FNt}
\begin{aligned}
 \widetilde F_p^N&=F_p^N,\qquad
 \widetilde F_m^N=\frac\phi\rho(\nabla\Pi-\mathcal L\psi)
 -\rho_b\psi\cdot\nabla\psi
 +\rho_bA_b\frac\phi\rho\frac\Pi{\gamma p_b}e_1,\\
 \widetilde F_z^N&=-\psi\cdot\nabla\zeta-(\gamma-1)\zeta\Div\psi
 -\frac{\kappa\phi}{c_v\rho_b\rho}\Delta\zeta
 +B_b\frac\phi\rho\frac\Pi{\gamma p_b}
 +\frac\phi{c_v\rho_b\rho}(h_b\Pi-H[\psi])
 +\frac{Q_\mu(\nabla\psi)}{c_v\rho};
\end{aligned}
\end{equation}
thus $\widetilde F^N=F^N+(0,\rho_bA_bz_Ne_1,B_bz_N)$ with $z_N$ of \eqref{e:zN}.
Put
\begin{equation}\label{e:Gbrelation}
\begin{gathered}
 G_b=\mathcal L(b\eta e_1)-\rho_bbD_b\eta\,e_1,\\
 G_b=\frac{\rho_b}\rho G_\eta+\frac\phi\rho\mathcal L(b\eta e_1)
 +\rho_bh_b\psi_1\eta\,e_1+\rho_bb(\psi\cdot\nabla\eta)e_1,
\end{gathered}
\end{equation}
where the last identity follows from $D_u=D_b+\psi\cdot\nabla$ and
$\rho=\rho_b+\phi$.

For $j=o,s$ let $q_j=(P_j,v_j,z_j)$ and $V_j=v_j-b\eta e_1$, and put
\begin{equation}\label{e:full}
\begin{aligned}
 f_{p,j}^{(0)}={}&-(\partial_1p_b)(v_j)_1-\gamma h_bP_j
 +(\gamma-1)H[v_j]-\gamma p_b\partial_1(b\eta)+\widetilde F_p^N,\\
 f_{m,j}^{(0)}={}&\Bigl[-\rho_bh_b(v_j)_1
 -\frac{\rho_bA_bP_j}{\gamma p_b}+\rho_bA_bf\Bigr]e_1
 +G_b+\widetilde F_m^N,\\
 f_{z,j}^{(0)}={}&-(\partial_1\theta_b)(v_j)_1
 -\frac{B_bP_j}{\gamma p_b}-\frac{h_bP_j}{c_v\rho_b}
 +\frac{H[v_j]}{c_v\rho_b}+B_bf+B_b\eta
 -(\gamma-1)\theta_b\partial_1(b\eta)+\widetilde F_z^N,\\
 (f_{p,o},f_{m,o},f_{z,o})={}&(f_{p,o}^{(0)},f_{m,o}^{(0)},f_{z,o}^{(0)}),\qquad
 (f_{p,s},f_{m,s},f_{z,s})=(f_{p,s}^{(0)},f_{m,s}^{(0)},f_{z,s}^{(0)})
 +\mathscr D+\mathscr T+F^g-F^\chi,
\end{aligned}
\end{equation}
where
\begin{equation}\label{e:DT}
\begin{gathered}
 \mathscr D=\mathfrak P\mathcal D_{\rm mat}[q_s],\qquad
 \mathscr T=\mathfrak PT_{\rm inv},\qquad F^g=\mathfrak PR_g,\\
 F^\chi=(0,\rho_bA_b\chi_ge_1,B_b\chi_g),\qquad
 \chi_g=\frac1{\theta_b}\Bigl(z_g-\frac{P_g}{c_p\rho_b}\Bigr).
\end{gathered}
\end{equation}

For a triple $f=(P_f,v_f,z_f)$ and a
scalar function $\eta$ put
\begin{equation}\label{e:be}
\begin{gathered}
 V_f=v_f-b\eta e_1,\qquad
 X_f=\frac{P_f}{\sqrt{\gamma p_b}},\qquad
 Y_f=\sqrt{\rho_b}(V_f-d_b\nabla z_f),\qquad
 r_{\sigma,f}=\frac{X_f+\sigma(Y_f)_1}2,\\
 e_b^0(f,\eta)=e_b(P_f,V_f,z_f)
 =\frac12\bigl[X_f^2+|Y_f|^2+(k_b-\rho_bd_b^2)|\nabla z_f|^2\bigr],\qquad
 k_b-\rho_bd_b^2=\frac{c(\nu+c)}{(\gamma-1)\rho_b\theta_b^2}.
\end{gathered}
\end{equation}

With $\epsilon_K,\theta_W$ of Lemma~\ref{l:mr}, for a triple $f=(P_f,v_f,z_f)$,
a source $\mathcal F=(\mathcal F_p,\mathcal F_m,\mathcal F_z)$ and
$\partial_if=(\partial_iP_f,\partial_iv_f,\partial_iz_f)$, put
\begin{equation}\label{e:srcw}
\begin{aligned}
 \mathcal F^{[i]}={}&\Bigl(\partial_i\mathcal F_p,\
 \partial_i\mathcal F_m-\frac{\partial_i\rho_b}{\rho_b}\mathcal F_m,\
 \partial_i\mathcal F_z\Bigr),\\
 \mathscr W_b(f,\eta;\mathcal F)={}&\mathscr S_b(P_f,V_f,z_f;\mathcal F)
 +\epsilon_K\mathscr K_b(P_f,V_f;\mathcal F_p,\mathcal F_m)
 +\frac{\theta_W}\tau\int_{\R^3}\rho_bz_f\mathcal F_z\dd x\\
 &+\sum_{i=1}^3\Bigl[\mathscr S_b(\partial_iP_f,\partial_iv_f,\partial_iz_f;\mathcal F^{[i]})
 +\epsilon_K\mathscr K_b\bigl(\partial_iP_f,\partial_iv_f;(\mathcal F^{[i]})_p,
 (\mathcal F^{[i]})_m\bigr)\Bigr].
\end{aligned}
\end{equation}

\begin{lemma}\label{l:res}
There are $\mathcal V_1,\mathcal V_2,\mathcal V_3,H_g$ such that
\begin{equation}\label{e:rg}
 R_g=\sum_{i=1}^3(\nabla_m)_i\mathcal V_i+H_g,\qquad
 \int_s^tw_s\Bigl[\sum_i\Vert\mathcal V_i\Vert_{L^2}^2
 +\tau\Vert H_g\Vert_{L^2}^2+\Vert R_g\Vert_{H^1}^2\Bigr]\dd r
 \le\omega_b(s)N_o^2,
\end{equation}
and each row of $F^g$ is $\sum_i\partial_i\mathcal V_{\alpha i}+H_\alpha$ with the same
bound in physical variables. One has
\begin{equation}\label{e:qseq}
 (\partial_t+L_b)(\mathcal P^{-1}q_s)=S_N+\mathcal D_{\rm mat}[q_s]+T_{\rm inv}+R_g,
 \qquad q_s(s)=q_o(s).
\end{equation}
For $j=o,s$, $(P_j,V_j,z_j)$ solves \eqref{e:forced} with the sources
\eqref{e:full}, and \eqref{e:wFC} and its fiber form hold for it. Moreover
$\mathscr T_p=\mathscr T_z=(\mathscr T_m)_1=0$,
\[
 (\mathscr T_m)_\perp
 =-\sum_\sigma\bigl(\sqrt{\gamma P_b}G_\sigma
 \nabla'\partial_yu_\sigma^{\rm ac}\bigr)\circ Y,\qquad
 \chi_g=\Bigl[\sum_\sigma\frac{V_b\sqrt{\gamma P_b}}{c_p\Theta_b}
 k_\sigma^e\Bigr]\circ Y,
\]
and $F^\chi$ has the bounds of \eqref{e:rg} in physical variables, also for the
rows $\rho_b\partial_i(A_b\chi_g)e_1$ and $\partial_i(B_b\chi_g)$. For
$\mathcal F\in\{F^g,F^\chi\}$ and every triple $f$ with $N[f]<\infty$,
\[
 \int_s^tw_s\big|\mathscr W_b(f,\eta;\mathcal F)\big|\dd r
 \le\omega_b(s)N_o\bigl(N[f]+\Phi^{\frac12}\bigr).
\]
Let $K\ge1$, let $\mathsf T(t)$ be a real even transverse convolution operator,
strongly measurable in $t$, whose convolution measure has total variation at
most $K$, and let $B_0(x_1,t)$ be measurable, with
$\Vert B_0\Vert_{L^\infty}+\sqrt\tau\Vert\partial_1B_0\Vert_{L^\infty}\le K$. For
$\mathcal F\in\{F^g,F^\chi\}$ and every triple $f$ with $N[f]<\infty$,
\begin{equation}\label{e:rwt}
\begin{aligned}
 \mathscr S_{B_0}^{\mathsf T}(f;\mathcal F)
 ={}&\int_{\R^3}B_0\Bigl[(V_f-d_b\nabla z_f)\cdot\mathsf T\mathcal F_m
              +\frac{P_f\mathsf T\mathcal F_p}{\gamma p_b}\Bigr]\dd x\\
 &+\int_{\R^3}\bigl[\Div(B_0j_bV_f)-\Div(B_0k_b\nabla z_f)\bigr]
                                      \mathsf T\mathcal F_z\dd x,\\
 \int_s^t w_s|\mathscr S_{B_0}^{\mathsf T}(f;\mathcal F)|\dd r
 \le{}&\omega_{b,K}(s)N_o\bigl(N[f]+\Phi^{\frac12}\bigr).
\end{aligned}
\end{equation}
\end{lemma}

\begin{proof}
Fix $(y,t)$ and freeze the coefficients at it. Put
$c_m=\sqrt{\frac{\gamma P_b}{V_b}}$, $c_t=\sqrt{\gamma P_bV_b}$, $\lambda_\sigma=\sigma c_m$,
$\mathcal D_1=V_b^{-1}\partial_y^2+V_b\Delta'$, and represent a mass vector $Q$ by
$\Upsilon=(\Upsilon_p,\Upsilon_1,\Upsilon_\perp,\Upsilon_h)$ through
\[
 \pi_Q=\sqrt{\gamma P_b}\,\Upsilon_p,\qquad
 (v_Q,j_\perp)=\sqrt{V_b}\,(\Upsilon_1,\Upsilon_\perp),\qquad
 z_Q=\frac{V_b\sqrt{\gamma P_b}}{c_p}(\Upsilon_p+\Upsilon_h).
\]
Then $L_E$ and $L_V$ act on $\Upsilon$ as
\[
\begin{aligned}
 \mathcal A\Upsilon={}&\bigl(c_m\partial_y\Upsilon_1+c_t\Div'\Upsilon_\perp,\
 c_m\partial_y\Upsilon_p,\ c_t\nabla'\Upsilon_p,\ 0\bigr),\\
 \mathcal V\Upsilon={}&\Bigl((\gamma-1)c\mathcal D_1(\Upsilon_p+\Upsilon_h),\
 \frac\nu{V_b}\partial_y^2\Upsilon_1+\mu V_b\Delta'\Upsilon_1
 +(\mu+\lambda)\partial_y\Div'\Upsilon_\perp,\\
 &\ \ \mu\mathcal D_1\Upsilon_\perp+(\mu+\lambda)V_b\nabla'\Div'\Upsilon_\perp
 +(\mu+\lambda)\partial_y\nabla'\Upsilon_1,\
 c\mathcal D_1(\Upsilon_p+\Upsilon_h)\Bigr).
\end{aligned}
\]
By \eqref{e:gc}, the $\sigma$ summand of $\mathcal P^{-1}g_c$ is represented by
$G_\sigma\partial_yu\,\mathbf r_\sigma+\mathbf k_\sigma$, with $u=u_\sigma^{\rm ac}$,
$\mathbf r_\sigma=(1,\sigma,0,0)$ and $\mathbf k_\sigma=(k_\sigma^u,k_\sigma^v,0,k_\sigma^e)$;
$R_\sigma$ is represented by $\mathbf r_\sigma$. For a scalar $h$,
$\mathcal A(h\mathbf r_\sigma)=\lambda_\sigma\partial_yh\,\mathbf r_\sigma+(0,0,c_t\nabla'h,0)$, and
with $d_n=\frac{\nu+(\gamma-1)c}{2V_b}$ and
$d_t=\frac{V_b(\mu+(\gamma-1)c)}2$,
\[
\begin{gathered}
 \frac12(1,\sigma,0,0)\cdot\mathcal V(h\mathbf r_\sigma)
       =(d_n\partial_y^2+d_t\Delta')h,\\
 \mathcal V(h\mathbf r_\sigma)
       -\mathbf r_\sigma(d_n\partial_y^2+d_t\Delta')h\\
 \qquad=\bigl(\mathcal Eh,\ -\sigma\mathcal Eh,\
             \sigma(\mu+\lambda)\partial_y\nabla'h,\ c\mathcal D_1h\bigr).
\end{gathered}
\]
where $\mathcal E=\frac12[\frac{(\gamma-1)c-\nu}{V_b}\partial_y^2+V_b((\gamma-1)c-\mu)\Delta']$.
For $h=G_\sigma\partial_yu$ the right side is
$(\partial_yF_\sigma^u,-\sigma\partial_yF_\sigma^u,\sigma(\mu+\lambda)G_\sigma\nabla'\partial_y^2u,
\partial_yF_\sigma^e)$. On the other hand
\[
 (\mathcal A-\lambda_\sigma\partial_y)\mathbf k_\sigma
 =\bigl(c_m\partial_y(k_\sigma^v-\sigma k_\sigma^u),\
 c_m\partial_y(k_\sigma^u-\sigma k_\sigma^v),\ c_t\nabla'k_\sigma^u,\
 -\sigma c_m\partial_yk_\sigma^e\bigr),
\]
and, since $c_m\sqrt{\frac{V_b}{\gamma P_b}}=1$ and $c_t\frac1{\sqrt{\gamma P_bV_b}}=1$, the
definitions of $k_\sigma^u,k_\sigma^v,k_\sigma^e$ give
\[
\begin{gathered}
 c_m\partial_y(k_\sigma^v-\sigma k_\sigma^u)=\partial_yF_\sigma^u,\qquad
 c_m\partial_y(k_\sigma^u-\sigma k_\sigma^v)=-\sigma\partial_yF_\sigma^u,\\
 c_t\nabla'k_\sigma^u=\sigma(\mu+\lambda)G_\sigma\nabla'\partial_y^2u,\qquad
 -\sigma c_m\partial_yk_\sigma^e=\partial_yF_\sigma^e .
\end{gathered}
\]
Hence, with frozen coefficients,
\begin{equation}\label{e:canc}
 (\mathcal A-\lambda_\sigma\partial_y)\mathbf k_\sigma
 =\mathcal V(G_\sigma\partial_yu\,\mathbf r_\sigma)
 -\mathbf r_\sigma(d_n\partial_y^2+d_t\Delta')(G_\sigma\partial_yu).
\end{equation}
Let $\mathcal G_\sigma$ and $K_\sigma$ be the mass vectors represented by
$(F_\sigma^u,-\sigma F_\sigma^u,\sigma(\mu+\lambda)G_\sigma\nabla'\partial_yu,F_\sigma^e)$ and
by $\mathbf k_\sigma$, and $A_\sigma=R_\sigma G_\sigma$; thus
$\mathcal P^{-1}g_c=\sum_\sigma(A_\sigma\partial_yu_\sigma^{\rm ac}+K_\sigma)$.

Now let the coefficients vary. The normal part of $L_E(R_\sigma h)$ is
$\partial_y(\lambda_\sigma R_\sigma h)$, and the transverse part of $L_E(A_\sigma h)$
is the $j_\perp$ row $V_b\sqrt{\gamma P_b}G_\sigma\nabla'h$; hence
$L_E(A_\sigma\partial_yu)=\partial_y(\lambda_\sigma A_\sigma\partial_yu)-T_{{\rm inv},\sigma}$,
$T_{{\rm inv},\sigma}$ being the $\sigma$ summand of $T_{\rm inv}$. With
$\mathcal D_\sigma=d_\sigma\partial_y^2+d_{t,\sigma}\Delta'$, \eqref{e:acg} gives
$\partial_t\partial_yu=-s_\sigma\partial_y^2u+\partial_y\mathcal D_\sigma u-\partial_yH_\sigma$, and
$A_\sigma H_\sigma=\frac{R_\sigma U_\sigma}{V_b}J_b\varphi_\sigma[q_o]$. Therefore
\[
\begin{aligned}
 -(\partial_t+L_b)(A_\sigma\partial_yu)
 ={}&T_{{\rm inv},\sigma}
 +\partial_y\Bigl(\frac{R_\sigma U_\sigma}{V_b}J_b\varphi_\sigma[q_o]\Bigr)
 -(\partial_yA_\sigma)H_\sigma
 -\partial_y\bigl[(\lambda_\sigma-s_\sigma)A_\sigma\partial_yu\bigr]\\
 &-\bigl((\partial_t+s_\sigma\partial_y)A_\sigma\bigr)\partial_yu
 +L_V(A_\sigma\partial_yu)-A_\sigma\partial_y\mathcal D_\sigma u .
\end{aligned}
\]
By the proof of Lemma~\ref{l:cor},
$J_b\varphi_\sigma[q_o]-J_b\varphi_\sigma[q_s]=G_\sigma\partial_yu+(J_b-I)(G_\sigma\partial_yu)
+\frac12J_b\sum_\iota(k_\iota^u+\sigma k_\iota^v)$, and
$\frac{R_\sigma U_\sigma}{V_b}G_\sigma=\frac{U_\sigma}{V_b}A_\sigma$. Thus
\begin{equation}\label{e:Rgx}
\begin{aligned}
 R_g=\sum_\sigma\Bigl\{&\partial_y\Bigl[\Bigl(\frac{U_\sigma}{V_b}-\lambda_\sigma+s_\sigma\Bigr)
 A_\sigma\partial_yu\Bigr]
 +\partial_y\Bigl[\frac{R_\sigma U_\sigma}{V_b}\Bigl((J_b-I)(G_\sigma\partial_yu)
 +\frac12J_b\sum_\iota(k_\iota^u+\sigma k_\iota^v)\Bigr)\Bigr]\\
 &-(\partial_yA_\sigma)H_\sigma
 -\bigl((\partial_t+s_\sigma\partial_y)A_\sigma\bigr)\partial_yu
 +L_V(A_\sigma\partial_yu)-A_\sigma\partial_y\mathcal D_\sigma u
 -(\partial_t+L_b)K_\sigma\Bigr\}.
\end{aligned}
\end{equation}
By \eqref{e:canc}, applied at each point, and $L_b=L_E-L_V$,
\[
\begin{gathered}
 L_V(A_\sigma\partial_yu)-A_\sigma\partial_y\mathcal D_\sigma u
 =\partial_y\mathcal G_\sigma+A_\sigma\partial_y\bigl((d_n-d_\sigma)\partial_y^2u
 +(d_t-d_{t,\sigma})\Delta'u\bigr)+D_\sigma^{(1)},\\
 \partial_y\mathcal G_\sigma-(\partial_t+L_b)K_\sigma
 =-(\partial_t+s_\sigma\partial_y)K_\sigma+(s_\sigma-\lambda_\sigma)\partial_yK_\sigma
 +L_VK_\sigma+D_\sigma^{(2)},
\end{gathered}
\]
where $D_\sigma^{(1)},D_\sigma^{(2)}$ consist of terms in which at least one derivative falls
on a coefficient: smooth functions of $Z$ times $\partial_yZ\,\nabla_m^ku$,
$\partial_y^2Z\,\nabla_m^{k-1}u$ or $(\partial_yZ)^2\nabla_m^{k-1}u$, $k=2,3$. By
\eqref{e:acg}, $(\partial_t+s_\sigma\partial_y)K_\sigma$ is the sum of
$\bigl((\partial_t+s_\sigma\partial_y)C\bigr)\nabla_m^2u$ and $C\,\nabla_m^2(\mathcal D_\sigma u-H_\sigma)$
over the coefficients $C=C(Z)$ of $K_\sigma$.

We estimate the terms of \eqref{e:Rgx} on $\Omega_\sigma$; let $\varrho_\sigma$ be as in
Lemma~\ref{l:cor} with $\ell=\frac{|s_\sigma|}2$ and $\varrho_\sigma=1$ where
$|y+h_*-s_\sigma t|\ge(|s_\sigma|-2c_0)\tau$, so that $1-\varrho_\sigma=0$ outside
$\Omega_\sigma$ for $\tau\ge\frac{|s_\sigma|}{c_0}$. Write $1=\varrho_\sigma+(1-\varrho_\sigma)$ in
every term of \eqref{e:Rgx}. The parts with $\varrho_\sigma$ have the bound \eqref{e:rg} by
the cone bound of Lemma~\ref{l:cor} with $P_0=2$, applied to $\varrho_\sigma$ and to a cutoff
equal to one on the support of $\varrho_\sigma$ and zero where
$|y+h_*-s_\sigma t|\le\frac{|s_\sigma|}4\tau$, and, for the terms with $H_\sigma$, by the
Gaussian factor of $U_\sigma$. The functions $u_\iota^{\rm ac}$, $\iota\ne\sigma$, occur only
with the factor $U_\sigma$, which is at most $C_Me^{-c_bt}$ where $\varrho_\iota\ne1$. On $\Omega_\sigma$, \eqref{e:E10} gives
$\sqrt\tau(|Z-Z_{\sigma,a}|+|\partial_yZ|)\le\omega_b(s)$. Since
$\lambda_\sigma(Z_{\sigma,*})=s_\sigma$ and, with $\lambda_\sigma=\sigma\sqrt{\gamma R}\,\Theta_b^{\frac12}V_b^{-1}$,
\[
 \frac{\dd}{\dd\epsilon}\lambda_\sigma\Bigl(Z_{\sigma,*}+\epsilon\Bigl(-1,s_\sigma,\frac{p_*}{c_v}\Bigr)\Bigr)\Big|_{\epsilon=0}
 =s_\sigma\Bigl(\frac{p_*}{2c_v\theta_\sigma}+\frac1{v_\sigma}\Bigr)
 =\frac{(\gamma+1)s_\sigma}{2v_\sigma},
\]
the definition of $Z_{\sigma,a}$ and \eqref{e:gap} give
$\frac{U_\sigma}{V_b}-\lambda_\sigma+s_\sigma=\frac{U_\sigma}{v_\sigma}-\frac{(\gamma+1)s_\sigma\vartheta_\sigma}{2v_\sigma}
+o_b(\tau^{-\frac12})=o_b(\tau^{-\frac12})$, where $o_b(\tau^{-\frac12})$ denotes a function
bounded by $\omega_b(s)\tau^{-\frac12}$; likewise $d_n-d_\sigma$, $d_t-d_{t,\sigma}$ and
$\lambda_\sigma-s_\sigma$ are $O_b(\tau^{-\frac12})$ and $\partial_yA_\sigma$ is $o_b(\tau^{-\frac12})$.
By \eqref{e:E10} with $a=1$, $\varsigma=s_\sigma$ and $f$ the entries of $A_\sigma$,
$(\partial_t+s_\sigma\partial_y)A_\sigma=\partial_yF_A+H_A$ with $\sqrt\tau|F_A|+\tau|H_A|\le\omega_b(s)$,
and by \eqref{e:E10} with $\varsigma=s_\sigma$ and the coefficients $C$ of $K_\sigma$,
$(\partial_t+s_\sigma\partial_y)C=\partial_yF_C+H_C$ with $\sqrt\tau|F_C|+\tau|H_C|\le\omega_b(s)$.
We place the terms of \eqref{e:Rgx}:
\begin{itemize}
\item the two fluxes in the first line, by
$\Vert(J-I)h\Vert_{L^2}\le\Vert\nabla h\Vert_{L^2}$ and \eqref{e:Jsym}, have $L^2(\Omega_\sigma)$
norms at most $\omega_b(s)\tau^{-\frac12}\Vert\nabla_mu\Vert_{L^2}
+C_b\tau^{-\frac12}\sum_\iota\Vert\nabla_m^2u_\iota^{\rm ac}\Vert_{L^2}$;
\item $(\partial_yA_\sigma)H_\sigma$ is a value with
$\int_s^tw_s\tau\Vert(\partial_yA_\sigma)H_\sigma\Vert_{L^2(\Omega_\sigma)}^2\dd r
\le\omega_b(s)\int_s^tw_s\Vert H_\sigma\Vert_{L^2}^2\dd r\le\omega_b(s)N_o^2$;
\item $(\partial_yF_A+H_A)\partial_yu=\partial_y(F_A\partial_yu)-F_A\partial_y^2u+H_A\partial_yu$, a flux
and two values;
\item $A_\sigma\partial_y((d_n-d_\sigma)\partial_y^2u+(d_t-d_{t,\sigma})\Delta'u)$ is the flux
$A_\sigma(\cdots)$ minus the value $(\partial_yA_\sigma)(\cdots)$;
\item $(\partial_yF_C+H_C)\nabla_m^2u$ is treated like $(\partial_yF_A+H_A)\partial_yu$;
$C\nabla_m^2\mathcal D_\sigma u$, $L_VK_\sigma$ and $(s_\sigma-\lambda_\sigma)\partial_yK_\sigma$ are fluxes
$C\nabla_m^3u$ and values $O_b(\tau^{-\frac12})\nabla_m^3u$, and $C\nabla_m^2H_\sigma$ is the
flux $C\nabla_m H_\sigma$ minus a value $\partial_yC\,\nabla_mH_\sigma$;
\item in $D_\sigma^{(1)},D_\sigma^{(2)}$, $\partial_y^2Z\,\nabla_m^{k-1}u
=\partial_y(\partial_yZ\,\nabla_m^{k-1}u)-\partial_yZ\,\partial_y\nabla_m^{k-1}u$.
\end{itemize}
With $\sqrt\tau|\partial_yZ|\le\omega_b(s)$ on $\Omega_\sigma$, every flux $\mathcal V$ has
$\int_s^tw_s\Vert\mathcal V\Vert_{L^2}^2\le\omega_b(s)(\mathcal J+\mathcal H)+C_b\tau_s^{-1}\mathcal J$
and every value $H$ has $\int_s^tw_s\tau\Vert H\Vert_{L^2}^2\le\omega_b(s)(\mathcal J+\mathcal H)
+C_b\mathcal H+\omega_b(s)N_o^2$, the terms with $\nabla_mH_\sigma$ by \eqref{e:acsrc}. The terms of
\eqref{e:Rgx} and their first derivatives are, with $C$ smooth functions of $Z$ and of
the profiles, of the types
\[
 C\nabla_m^{3,4,5}u,\quad \partial_yC\,\nabla_m^{1,\dots,4}u,\quad
 (\partial_y^2C,\partial_tC)\nabla_m^{1,2,3}u,\quad
 (\partial_y^3C,\partial_t\partial_yC)\nabla_m^{1,2}u,\quad
 C\nabla_m^{2,3}H_\sigma,\quad \partial_yC\,\nabla_mH_\sigma,
\]
where the coefficients of $\nabla_mu$ and $\nabla_m^2u$ are $o_b(\tau^{-\frac12})$, or contain
$\partial_y^2Z$, $\partial_y^3Z$; the latter are bounded by \eqref{eq:planar}
and $\int_{S}^\infty\Vert\partial_y^2Z\Vert_{H^1}^2\dd t<\infty$ from \eqref{e:munw}, as in the proof
of Lemma~\ref{l:cor}. With Lemma~\ref{l:cor} this proves \eqref{e:rg}. For any
matrix $M(Z)$,
\[
 (M\partial_y\mathcal V_1)\circ Y
 =\partial_1\bigl((MV_b\mathcal V_1)\circ Y\bigr)
       -\partial_1\bigl((MV_b)\circ Y\bigr)(\mathcal V_1\circ Y),
\]
whose last coefficient is $O_b(\tau^{-\frac12})$; this gives the physical form of
\eqref{e:rg} for $F^g$.

Since $\Pi=R(\theta\phi+\rho_b\zeta)$,
\[
 \frac{\Pi-R\rho_b\zeta}{R\rho_b\theta_b}-\frac\phi\rho
 =\frac\phi\rho\frac\Pi{p_b}.
\]
Inserting this identity into the three equations of \eqref{e:peq3}, written about
the carrier $(\rho_b,u_be_1,\theta_b)$, gives \eqref{e:FN}, its conservative lift $S_N$ and
$(\partial_t+L_b)\mathcal P^{-1}q_o=S_N$; subtracting the definition of $R_g$ gives
\eqref{e:qseq}. For \eqref{e:full}, write \eqref{e:peq3}, its momentum row multiplied by
$\frac{\rho_b}\rho$, with $D_u=D_b+\psi\cdot\nabla$, $\frac{\rho_b}\rho=1-\frac\phi\rho$ and
$\frac1\rho-\frac1{\rho_b}=-\frac\phi{\rho\rho_b}$. By $p=R\rho\theta$ and \eqref{e:chio},
$\frac p{p_b}\frac\phi\rho=\frac\Pi{\gamma p_b}-\chi_o^b$ with $\chi_o^b=\eta+f+z_N$, so that
\[
 -\rho_bA_b\frac\phi\rho
 =-\frac{\rho_bA_b\Pi}{\gamma p_b}+\rho_bA_b(\eta+f)
 +\rho_bA_b\frac\phi\rho\frac\Pi{\gamma p_b},
\]
and the same with $B_b$ in the temperature equation; the pressure equation has no
such term. Then $v_o=V_o+b\eta e_1$, $D_bb=A_b$ and
$\mathcal L(b\eta e_1)-\rho_bD_b(b\eta)e_1=G_b-\rho_bA_b\eta e_1$ remove the $\eta$ part of the
momentum source, and $\Div v_o=\Div V_o+\partial_1(b\eta)$ gives the terms
$-\gamma p_b\partial_1(b\eta)$ and $-(\gamma-1)\theta_b\partial_1(b\eta)$. This is \eqref{e:full} for
$j=o$. For $j=s$ apply $\mathfrak P$ to \eqref{e:qseq}; $\mathcal D_{\rm mat}[q_s]$, $T_{\rm inv}$ and
$R_g$ give $\mathscr D,\mathscr T,F^g$, and
$\frac1{\theta_b}(z_s-\frac{P_s}{c_p\rho_b})=\chi_o^b-\chi_g$ gives $-F^\chi$.
Since the only nonzero entry of $T_{\rm inv}$ is
$(T_{\rm inv})_\perp$, \eqref{e:conv} gives $\mathscr T$. The identity \eqref{e:wFC} and its
fiber form are those of Lemma~\ref{l:FC} for these sources: $V_j\in C_tL^2\cap L^2_tH^1$ by
\eqref{eq:class}, Lemma~\ref{l:cor} and Lemma~\ref{l:ent}, $G_b\in L^2_tH^{-1}$ and the
temperature source is in $L^2_tL^2$, so that the regularization in the opening paragraph
of this section applies.

The parts $G_\sigma\partial_yu$ and $k^u_\sigma$ of \eqref{e:gc} have $z=\frac{P}{c_p\rho_b}$, which
gives the formula for $\chi_g$. By \eqref{e:bal},
\[
\begin{aligned}
 \rho_bA_b\chi_g
 &=\partial_1(\nu h_b\chi_g)-(\partial_1p_b)\chi_g-\nu h_b\partial_1\chi_g,\\
 B_b\chi_g
 &=\partial_1\Bigl(\frac{\kappa\partial_1\theta_b}{c_v\rho_b}\chi_g\Bigr)
 -\frac{\kappa\partial_1\theta_b}{c_v\rho_b}\partial_1\chi_g
 -\partial_1\Bigl(\frac{\kappa}{c_v\rho_b}\Bigr)(\partial_1\theta_b)\chi_g
 +\frac{\nu h_b^2-p_bh_b}{c_v\rho_b}\chi_g ,
\end{aligned}
\]
and $\chi_g$ is a sum of bounded coefficients times $\nabla_m^2u_\sigma^{\rm ac}\circ Y$; with
\eqref{e:E2}, \eqref{e:E3}, \eqref{e:E6} and Lemma~\ref{l:cor} this gives the bounds of
$F^\chi$, and $\partial_i(\rho_bA_b\chi_g)-(\partial_i\rho_b)A_b\chi_g=\rho_b\partial_i(A_b\chi_g)$.

For a source $\mathcal F=\Div\mathcal V+H$ and a test $T$,
\[
 \int_s^tw_s|\langle T,\mathcal F\rangle|\dd r
 \le\Bigl(\int_s^tw_s\bigl(\Vert\mathcal V\Vert_{L^2}^2
                       +\tau\Vert H\Vert_{L^2}^2\bigr)\dd r\Bigr)^{\frac12}
 \Bigl(\int_s^tw_s\bigl(\Vert\nabla T\Vert_{L^2}^2
                       +\tau^{-1}\Vert T\Vert_{L^2}^2\bigr)\dd r\Bigr)^{\frac12}.
\]
For \eqref{e:rwt}, self-adjointness and commutation of $\mathsf T$ give the tests
\[
\begin{gathered}
 T_m=\mathsf T[B_0(V_f-d_b\nabla z_f)],\qquad
 T_p=\mathsf T\Bigl[\frac{B_0P_f}{\gamma p_b}\Bigr],\\
 T_z=\mathsf T\bigl[\Div(B_0j_bV_f)-\Div(B_0k_b\nabla z_f)\bigr].
\end{gathered}
\]
The measure bound of $\mathsf T$, $\sqrt\tau\Vert\partial_1B_0\Vert_{L^\infty}\le K$,
\eqref{e:E1Z} and \eqref{e:eta} give
$\int_s^tw_s[\sum_{\alpha=p,m}(\Vert\nabla T_\alpha\Vert_{L^2}^2
+\tau^{-1}\Vert T_\alpha\Vert_{L^2}^2)+\Vert T_z\Vert_{L^2}^2]\dd r
\le C_{b,K}(N[f]^2+\Phi)$.
The preceding dual bound for the mechanical rows and the $L^2$ pairing of the heat
row, with \eqref{e:rg}, prove \eqref{e:rwt}; no time derivative of $B_0$ or $\mathsf T$ occurs.

The zero-order unscreened part of $\mathscr W_b(f,\eta;\mathcal F)$ is
$\mathscr S_1^I(f;\mathcal F)$, already covered by \eqref{e:rwt}.
The first-order slots use the $H^1$ bound in \eqref{e:rg}, and the weak temperature
slot uses $z_f/\tau\in L^2_w$. For the zero-order screen,
\[
 \int\rho_bV_f\cdot\nabla J\mathcal F_p
 =-\int\Div(\rho_bV_f)J\Div\mathcal V_p
   +\int\rho_bV_f\cdot\nabla JH_p .
\]
By \eqref{e:Jsym} and \eqref{e:eta}, the remaining test factors are bounded by
$C_b(N[f]+\Phi^{\frac12})$, giving the asserted bound of $\mathscr W_b$.
\end{proof}

A hat denotes the unitary Fourier transform in $x'$, with variable $k\in\R^2$, and
$D_k=(\partial_y,ik_2,ik_3)$; products of complex factors are read as in
Section~\ref{sec:en}. For functions $l(t,k)$, $m(t,k)$ we write $\mathsf L=l(t,D')$,
$\mathsf M=m(t,D')$ for the transverse multipliers with these symbols.

\begin{lemma}\label{l:sp}
Let $l(t,k)\in[0,1]$ be real, even and measurable, and let $m_{\rm sp}(t,k)$ be real,
even and absolutely continuous in $t$, with $1\le m_{\rm sp}\le m_*<\infty$ and
$-\partial_tm_{\rm sp}\ge0$; let $\mathsf L,\mathsf M$ be their multipliers and
$f=q_o-\mathsf Lg_c$. Then, with $\widehat u_\sigma=\widehat{u_\sigma^{\rm ac}}$ and $J_b$ acting on
each Fourier component,
\[
\begin{aligned}
 &\partial_t\widehat u_\sigma
 +\Bigl(s_\sigma+\frac{lU_\sigma}{V_b}\Bigr)\partial_y\widehat u_\sigma
 -d_\sigma\partial_y^2\widehat u_\sigma+d_{t,\sigma}|k|^2\widehat u_\sigma\\
 &\quad=-\frac{U_\sigma}{V_bG_\sigma}
 \Bigl[J_b\widehat{\varphi_\sigma[f]}
 +l(J_b-I)(G_\sigma\partial_y\widehat u_\sigma)
 +\frac l2J_b\sum_{\iota=\pm1}
 (\widehat{k_\iota^u}+\sigma\widehat{k_\iota^v})\Bigr].
\end{aligned}
\]
For $L_0\ge1$ put
\begin{equation}\label{e:spec}
\begin{aligned}
 E_{\rm ac}^{\mathsf M}
 &=\frac12\sum_\sigma\int_{\R\times\R^2}m_{\rm sp}
 \Bigl[L_0\tau^{-1}|\widehat u_\sigma|^2
 +|D_k\widehat u_\sigma|^2
 +\sum_{j=2}^4|D_k^j\widehat u_\sigma|^2\Bigr]\dd y\dd k,\\
 D_{\rm ac}^{\mathsf M}
 &=\sum_\sigma\int_{\R\times\R^2}m_{\rm sp}
 \Bigl[\tau^{-2}|\widehat u_\sigma|^2
 +\tau^{-1}|D_k\widehat u_\sigma|^2
 +|D_k^2\widehat u_\sigma|^2
 +\sum_{j=3}^5|D_k^j\widehat u_\sigma|^2\Bigr]\dd y\dd k,
\end{aligned}
\end{equation}
and let $E_{\rm ac}^{-\mathsf M_t}$ be $E_{\rm ac}^{\mathsf M}$ with $m_{\rm sp}$ replaced by
$-\partial_tm_{\rm sp}$. Let $d_B>0$, let $B_{\rm out}(x_1,t)$ be smooth with
$0\le B_{\rm out}\le1$ and $B_{\rm out}=1$ where $|Y+h_*-s_\sigma t|\le2d_B\tau$ for
$\sigma=\pm1$, and put
\[
 \mathcal P_{m,l}=\int_{\R\times\R^2}
 m_{\rm sp}B_{\rm out}\bigl[l e_b^0(\widehat q_s,\widehat\eta)
 +(1-l)e_b^0(\widehat q_o,\widehat\eta)\bigr]\dd x_1\dd k .
\]
There are $L_0,c_a,C_M>0$, depending only on $M$ and the fixed data, and
$r_s\ge0$ such that
\begin{equation}\label{e:spin}
\begin{gathered}
 \frac{\dd}{\dd t}(w_sE_{\rm ac}^{\mathsf M})
 +w_s\Bigl[c_aD_{\rm ac}^{\mathsf M}
 +E_{\rm ac}^{-\mathsf M_t}+\frac{q_a}{\tau}E_{\rm ac}^{\mathsf M}\Bigr]
 \le\frac{C_Mw_s}{\tau}\mathcal P_{m,l}+w_sr_s,\\
 E_{\rm ac}^{\mathsf M}(s)=0,\qquad
 \int_s^tw_sr_s\dd r\le\omega_{b,m_*,d_B}(s)(N_o^2+\Phi),
\end{gathered}
\end{equation}
after integration over every $[s,t]$. For a triple $f$ and a scalar $\eta$ put
\[
\begin{aligned}
 E_c[f,\eta]&=\Vert(P_f,V_f)\Vert_{L^2}^2
 +\Vert(\nabla P_f,\nabla v_f)\Vert_{L^2}^2
 +\Vert\nabla z_f\Vert_{H^1}^2+\tau^{-1}\Vert z_f\Vert_{L^2}^2,\\
 D_c[f,\eta]&=\Vert\nabla P_f\Vert_{L^2}^2
 +\Vert\nabla V_f\Vert_{L^2}^2+\Vert\nabla^2v_f\Vert_{L^2}^2
 +\Vert\nabla^2z_f\Vert_{H^1}^2
 +\tau^{-1}\Vert\nabla z_f\Vert_{L^2}^2+\tau^{-2}\Vert z_f\Vert_{L^2}^2,
\end{aligned}
\]
let $\mathcal E_{\rm mix}^{\mathsf M}$, $\mathcal D_{\rm mix}^{\mathsf M}$ be the Fourier integrals
of $E_c,D_c$ with the weight $m_{\rm sp}l$ on $(\widehat q_s,\widehat\eta)$ and
$m_{\rm sp}(1-l)$ on $(\widehat q_o,\widehat\eta)$, and for $\theta_a>0$
\[
 \mathcal P_{\theta_a}(s,t)
 =\sup_{s\le r\le t}w_s(\mathcal E_{\rm mix}^{\mathsf M}
 +\theta_aE_{\rm ac}^{\mathsf M})
 +\int_s^tw_s\Bigl[\mathcal D_{\rm mix}^{\mathsf M}
 +\frac{q_a}{\tau}\mathcal E_{\rm mix}^{\mathsf M}
 +\theta_aD_{\rm ac}^{\mathsf M}
 +\frac{\theta_aq_a}{\tau}E_{\rm ac}^{\mathsf M}\Bigr]\dd r.
\]
Then
\begin{equation}\label{e:sprec}
\begin{gathered}
 N_o^2+N_s^2+\theta_a(\mathcal J+\mathcal H)
 \le C_{M,\theta_a}\bigl(\mathcal P_{\theta_a}+m_*N_0(b\eta;s,t)^2\bigr)
 +\omega_{b,m_*,\theta_a}(s)N_o^2,\\
 E[q_o](t)+E[q_s](t)\le C_M\bigl[
 \mathcal E_{\rm mix}^{\mathsf M}(t)+E_{\rm ac}^{\mathsf M}(t)
 +\Vert\mathsf M^{\frac12}(b\eta)(t)\Vert_{L^2}^2\bigr].
\end{gathered}
\end{equation}
\end{lemma}

\begin{proof}
Since $\mathsf L$ commutes with functions of $(y,t)$,
$\varphi_\sigma[q_o]=\varphi_\sigma[f]+\mathsf L\varphi_\sigma[g_c]$, and by the proof of
Lemma~\ref{l:cor}, $\varphi_\sigma[g_c]=G_\sigma\partial_yu_\sigma^{\rm ac}
+\frac12\sum_\iota(k_\iota^u+\sigma k_\iota^v)$. Taking the Fourier transform of
\eqref{e:acg} and writing $J_b=I+(J_b-I)$ in front of $lG_\sigma\partial_y\widehat u_\sigma$,
the term $-\frac{U_\sigma}{V_bG_\sigma}lG_\sigma\partial_y\widehat u_\sigma$ moves to the left and gives
the stated equation; no time derivative of $l$ occurs.

Its first source term satisfies, with $c_B>0$ depending on $d_B$,
\[
 \sum_\sigma\int_{\R^2}m_{\rm sp}
 \Bigl\Vert\frac{U_\sigma}{V_bG_\sigma}
                  J_b\widehat{\varphi_\sigma[f]}\Bigr\Vert_{L^2_y}^2\dd k
 \le\frac{C_M}{\tau}\mathcal P_{m,l}
 +\frac{C_Mm_*}{\tau}\Vert b\Vert_{L^\infty}^2\Vert\eta\Vert_{L^2}^2
 +\frac{C_{b,m_*}}{\tau}e^{-c_Bt}(E[q_o]+E[g_c]),
\]
with $C_M$ independent of $m_*$, $d_B$ and $B_{\rm out}$. Indeed, fix $\sigma$ and put
$I_1=\{y:|y+h_*-s_\sigma t|\le d_B\tau\}$, $I_2=\{y:|y+h_*-s_\sigma t|\le2d_B\tau\}$. On
$\R\setminus I_1$, Lemma~\ref{l:pr} gives $|U_\sigma|\le C_M\tau^{-\frac12}e^{-c_Bt}$ for
$t\ge\frac{2|h_*|}{d_B}$. On $I_1$ write
$\widehat{\varphi_\sigma[f]}=\mathbf 1_{I_2}\widehat{\varphi_\sigma[f]}+\mathbf 1_{\R\setminus I_2}\widehat{\varphi_\sigma[f]}$; in
$x_1$ the sets $X(I_1)$ and $X(\R\setminus I_2)$ are at distance at least $c_Md_B\tau$, so by
Lemma~\ref{l:bessel} the second part contributes at most $e^{-c_Md_B\tau}$ times its norm,
uniformly in $k$, while $J_b$ is a contraction on the first part up to the factor $C_M$
of the change of measure. By \eqref{e:be},
\[
 \varphi_\sigma[f]\circ Y=r_{\sigma,f}+\frac{\sigma\sqrt{\rho_b}}2(d_b\partial_1z_f+b\eta),\qquad
 |\varphi_\sigma[f]\circ Y|^2\le C_M\bigl(e_b^0(f,\eta)+|b\eta|^2\bigr),
\]
and since $\widehat f=l\widehat q_s+(1-l)\widehat q_o$ and $e_b^0(\cdot,\widehat\eta)$ is a convex quadratic
form, $e_b^0(\widehat f,\widehat\eta)\le le_b^0(\widehat q_s,\widehat\eta)+(1-l)e_b^0(\widehat q_o,\widehat\eta)$.
On $I_2$, $B_{\rm out}\circ X=1$; with $|U_\sigma|\le C_M\tau^{-\frac12}$, $m_{\rm sp}\le m_*$ in the
$b\eta$ term and $|\widehat f|\le|\widehat q_o|+|\widehat g_c|$ in the exponentially small terms, this
proves the bound.

Write $\widehat u=\widehat u_\sigma$ and $\int$ for the integral over $\R_y\times\R^2_k$. Multiply the
equation by $m_{\rm sp}\overline{\widehat u}/\tau$ and take real parts:
\[
\begin{aligned}
 &\frac{\dd}{\dd t}\frac1{2\tau}\int m_{\rm sp}|\widehat u|^2
 +\frac1{2\tau}\int(-\partial_tm_{\rm sp})|\widehat u|^2
 +\frac1{2\tau^2}\int m_{\rm sp}|\widehat u|^2
 +\frac1\tau\int m_{\rm sp}\bigl(d_\sigma|\partial_y\widehat u|^2+d_{t,\sigma}|k|^2|\widehat u|^2\bigr)\\
 &\quad=\frac1{2\tau}\int m_{\rm sp}l\,\partial_y\Bigl(\frac{U_\sigma}{V_b}\Bigr)|\widehat u|^2
 +\frac1\tau\operatorname{Re}\int m_{\rm sp}\overline{\widehat u}\,\mathcal S_\sigma,
\end{aligned}
\]
$\mathcal S_\sigma$ being the right side of the equation. Since $0\le l\le1$ and
$1-\beta_M\ge0$, \eqref{e:gap} gives $l\tau\partial_yW_\sigma\le1-\beta_M$, and
$\frac{U_\sigma}{V_b}-W_\sigma$ is handled as in the proof of Lemma~\ref{l:cor}; the left
side therefore retains $\frac{\beta_M}{2\tau^2}\int m_{\rm sp}|\widehat u|^2$. At each $k$ the
two feedback terms of $\mathcal S_\sigma$ have the bounds found in the proof of
Lemma~\ref{l:cor}, with $(J_b-I)$ and $J_b$ replaced by their fiber operators, and are
absorbed for large $s$; the first term of $\mathcal S_\sigma$ is estimated by Young's
inequality and the bound above. Multiply the equation by
$-m_{\rm sp}(\partial_y^2-|k|^2)\overline{\widehat u}$; with $d_0=\min_\sigma\min\{d_\sigma,d_{t,\sigma}\}$ the drift $\frac{lU_\sigma}{V_b}\partial_y\widehat u$ and
the feedback terms cost at most $\frac{d_0}4\int m_{\rm sp}|D_k^2\widehat u|^2
+\frac{C_M}\tau\sum_\iota\int m_{\rm sp}|D_k\widehat u_\iota|^2$, which is paid by
$\frac{L_0d_0}\tau\int m_{\rm sp}|D_k\widehat u|^2$ from $L_0$ times the first test for
$L_0=\frac{4C_M}{d_0}$. For $2\le j\le4$ apply $D_k^{j-1}$ to the Fourier transform of
\eqref{e:acg}, test with $-m_{\rm sp}(\partial_y^2-|k|^2)D_k^{j-1}\overline{\widehat u}$ and move one
derivative of the source onto the test:
\[
 \frac12\frac{\dd}{\dd t}\int m_{\rm sp}|D_k^j\widehat u|^2
 +\frac{d_0}2\int m_{\rm sp}|D_k^{j+1}\widehat u|^2
 +\frac12\int(-\partial_tm_{\rm sp})|D_k^j\widehat u|^2
 \le C\,m_*\Vert\nabla_m^{j-1}H_\sigma\Vert_{L^2}^2 .
\]
Summing, and multiplying by $w_s$, which produces $\frac{q_aw_s}\tau E_{\rm ac}^{\mathsf M}$,
gives \eqref{e:spin} with
\[
 r_s=C_Mm_*\sum_{\sigma}\sum_{j=1}^3\Vert\nabla_m^jH_\sigma\Vert_{L^2}^2
 +\frac{C_Mm_*}{\tau}\Vert b\Vert_{L^\infty}^2\Vert\eta\Vert_{L^2}^2
 +\frac{C_{b,m_*}}{\tau}e^{-c_Bt}(E[q_o]+E[g_c]).
\]
By \eqref{e:acsrc}, \eqref{e:bsup}, $\int_s^tw_s\tau^{-1}\Vert\eta\Vert_{L^2}^2\dd r\le q_a^{-1}\Phi$ and
Lemma~\ref{l:cor}, $\int_s^tw_sr_s\dd r\le\omega_{b,m_*,d_B}(s)(N_o^2+\Phi)$, and $u_\sigma^{\rm ac}(s)=0$ in $H^4$
gives $E_{\rm ac}^{\mathsf M}(s)=0$. A positive truncation of $-\partial_tm_{\rm sp}$ and frequency
truncations justify the integrated inequality.

For Hilbert-space vectors $A,B$ and $0\le l\le1$,
\[
 |A|^2+|B|^2\le3\bigl(l|A|^2+(1-l)|B|^2\bigr)+2|A-B|^2 .
\]
Apply this at each $k$ to the squares in $E[\cdot]$ and $D[\cdot]$, with $A,B$ the
corresponding derivatives of $\widehat q_s,\widehat q_o$. At zero order in the velocity
and at its first derivative in $D$, use
$v_j=V_j+b\eta e_1$ and $\nabla v_j=\nabla V_j+\nabla(b\eta)e_1$. Since $m_{\rm sp}\ge1$, the bound of
$E[g_c]$ in Lemma~\ref{l:cor} by $C_ME_{\rm ac}^{\mathsf M}$ gives the second line of
\eqref{e:sprec}. Integrating, the $b\eta$ terms give $N_0(b\eta;s,t)^2$, the $g_c$ terms give
$N[g_c]^2\le C_M\mathcal J+\omega_b(s)N_o^2$, and $m_{\rm sp}\ge1$, $L_0\ge1$ give
$\theta_a(\mathcal J+\mathcal H)\le4\mathcal P_{\theta_a}$ by Plancherel's identity; this
proves the first line.
\end{proof}

\subsection{Crossing identities}

\begin{lemma}\label{l:cr}
Let $(P,v,z)$ solve \eqref{e:forced}, and put
$h_0=\sqrt{\gamma p_b}$, $c_E=\sqrt{\frac{\gamma p_b}{\rho_b}}$, $X=\frac P{h_0}$,
$Y=\sqrt{\rho_b}(v-d_b\nabla z)$ and $r_\sigma=\frac{X+\sigma Y_1}2$; for
$(P,v,z)=(P_f,V_f,z_f)$ these are $X_f,Y_f,r_{\sigma,f}$ of \eqref{e:be}. Then
\begin{equation}\label{e:acsys}
\begin{gathered}
 D_bX+c_E\Div Y=R_X,\qquad
 D_bY+c_E\nabla X=R_Y,\\
 D_br_\sigma+\sigma c_E\partial_1r_\sigma
 =-\frac{c_E}{2}\Div'Y_\perp+\frac{R_X+\sigma(R_Y)_1}{2},\\
 R_X=\frac{f_p}{h_0}-\frac{XD_bh_0}{h_0}
 +c_E\frac{\partial_1\rho_b}{2\rho_b}Y_1-h_0(\partial_1d_b)\partial_1z,\\
\begin{aligned}
 R_Y={}&\rho_b^{-\frac12}\mathcal Lv
 -\sqrt{\rho_b}d_b\frac{\kappa}{c_v\rho_b}\nabla\Delta z
 +\sqrt{\rho_b}d_b(\gamma-1)\theta_b\nabla\Div v
 -\sqrt{\rho_b}d_b\nabla\Bigl(\frac{\kappa}{c_v\rho_b}\Bigr)\Delta z\\
 &+\sqrt{\rho_b}d_b(\gamma-1)(\nabla\theta_b)\Div v
 +\rho_b^{-\frac12}f_m-\sqrt{\rho_b}d_b\nabla f_z\\
 &-\sqrt{\rho_b}(D_bd_b)\nabla z
 +\sqrt{\rho_b}d_bh_b(\partial_1z)e_1
 -\frac{h_b}{2}Y-\rho_b^{-\frac12}X\nabla h_0 .
\end{aligned}
\end{gathered}
\end{equation}
For every real $B=B(x_1,t)$,
\begin{equation}\label{e:chid}
 (D_bB)e_b(P,v,z)+(\partial_1B)P(v_1-d_b\partial_1z)
 =\sum_{\iota=\pm1}(D_bB+\iota c_E\partial_1B)r_\iota^2
 +\frac{D_bB}2\bigl(|Y_\perp|^2+(k_b-\rho_bd_b^2)|\nabla z|^2\bigr),
\end{equation}
so that the last integral in \eqref{e:wFC} for $c=b$ is the integral of the right side
of \eqref{e:chid} plus $\int_{\R^3}(\partial_1B)(n_b-P(v_1-d_b\partial_1z))\dd x$. Put
\begin{equation}\label{e:asq}
 B_+^{\rm sq}(y,t)=\int_{-\infty}^yW_+^2(\xi,t)\dd\xi,\qquad
 B_-^{\rm sq}(y,t)=\int_y^\infty W_-^2(\xi,t)\dd\xi,
\end{equation}
and let $I_\sigma h$ denote $\int_{-\infty}^yh\dd\xi$ for $\sigma=1$ and $\int_y^\infty h\dd\xi$ for
$\sigma=-1$. Then
\[
\begin{gathered}
 D_b(B_\sigma^{\rm sq}\circ Y)
 =(-|s_\sigma|W_\sigma^2+R_{B,\sigma})\circ Y,
 \qquad\partial_1(B_\sigma^{\rm sq}\circ Y)
 =\sigma\rho_b(W_\sigma^2\circ Y),\\
 R_{B,\sigma}=\sigma d_\sigma\partial_y(W_\sigma^2)
 -\frac{2\sigma}{3}W_\sigma^3-2d_\sigma I_\sigma\bigl((\partial_yW_\sigma)^2\bigr),\qquad
 0\le B_\sigma^{\rm sq}\le C_M\tau^{-\frac12},\qquad
 \Vert R_{B,\sigma}\Vert_{L^\infty}\le C_M\tau^{-\frac32}.
\end{gathered}
\]
With $w_\sigma=W_\sigma\circ Y$,
\begin{equation}\label{e:cross}
\begin{aligned}
 \frac{\dd}{\dd t}\int_{\R^3}w_\sigma r_{-\sigma}^2\dd x
 ={}&-2s_\sigma\int_{\R^3}((\partial_yW_\sigma)\circ Y)r_{-\sigma}^2\dd x
 -\int_{\R^3}w_\sigma c_Er_{-\sigma}\Div'Y_\perp\dd x+\operatorname{Err}_1,\\
 \frac{\dd}{\dd t}\int_{\R^3}w_\sigma r_\sigma r_{-\sigma}\dd x
 ={}&-2s_\sigma\int_{\R^3}((\partial_yW_\sigma)\circ Y)r_\sigma r_{-\sigma}\dd x
 -2\sigma\int_{\R^3}w_\sigma c_E(\partial_1r_\sigma)r_{-\sigma}\dd x\\
 &-\frac12\int_{\R^3}w_\sigma c_E(r_\sigma+r_{-\sigma})\Div'Y_\perp\dd x+\operatorname{Err}_2,
\end{aligned}
\end{equation}
where
\[
\begin{aligned}
 a_\sigma^{\rm err}
 &=(s_\sigma-\sigma\rho_bc_E)((\partial_yW_\sigma)\circ Y)
 +\bigl(-W_\sigma\partial_yW_\sigma+d_\sigma\partial_y^2W_\sigma\bigr)\circ Y
 +(h_b-\sigma\partial_1c_E)w_\sigma,\\
 \operatorname{Err}_1
 &=\int_{\R^3}\bigl[a_\sigma^{\rm err}r_{-\sigma}^2
 +w_\sigma r_{-\sigma}(R_X-\sigma(R_Y)_1)\bigr]\dd x,\\
 \operatorname{Err}_2
 &=\int_{\R^3}\Bigl[a_\sigma^{\rm err}r_\sigma r_{-\sigma}
 +\frac{w_\sigma}{2}\bigl\{(r_\sigma+r_{-\sigma})R_X
 +\sigma(r_{-\sigma}-r_\sigma)(R_Y)_1\bigr\}\Bigr]\dd x.
\end{aligned}
\]
The identities \eqref{e:chid} and \eqref{e:cross} are quadratic in $(P,v,z)$; they hold
for the associated symmetric bilinear forms of two solutions, hence for Fourier
components and with a self-adjoint transverse multiplier independent of $t$
inserted in one slot. The weights $B_\sigma^{\rm sq}\circ Y$, $w_\sigma$ and $w_\sigma^2$ satisfy
the hypothesis of \eqref{e:wFC} on $[s,T_e]$.
\end{lemma}

\begin{proof}
By $D_b\rho_b=-\rho_bh_b$ in \eqref{e:bal},
\[
 D_b\sqrt{\rho_b}=-\frac{h_b}2\sqrt{\rho_b},\qquad
 D_b\nabla z=\nabla D_bz-h_b(\partial_1z)e_1,\qquad
 h_0d_b=\frac{(\gamma-1)\kappa}{h_0},
\]
the last by \eqref{e:FC}. Hence the first row of
\eqref{e:forced} gives $D_bX=\frac{f_p}{h_0}-h_0\Div v+h_0d_b\Delta z-\frac{XD_bh_0}{h_0}$, while
$c_E\sqrt{\rho_b}=h_0$ gives $c_E\Div Y=h_0(\Div v-d_b\Delta z-\partial_1d_b\,\partial_1z)
+c_E\frac{\partial_1\rho_b}{2\rho_b}Y_1$; this is the equation of $X$. The second row, divided by
$\sqrt{\rho_b}$, and the gradient of the third row give the equation of $Y$, with
$\rho_b^{-\frac12}\nabla P=c_E\nabla X+\rho_b^{-\frac12}X\nabla h_0$. Half the sum of the
equation of $X$ and $\sigma$ times the first component of that of $Y$ is the equation of
$r_\sigma$, since $\frac{c_E}2(\partial_1Y_1+\sigma\partial_1X)=\sigma c_E\partial_1r_\sigma$.

Since $r_+^2+r_-^2=\frac{X^2+Y_1^2}2$, $r_+^2-r_-^2=XY_1$, $P(v_1-d_b\partial_1z)=c_EXY_1$ and
$e_b(P,v,z)=\frac12[X^2+|Y|^2+(k_b-\rho_bd_b^2)|\nabla z|^2]$ by \eqref{e:be}, both sides of
\eqref{e:chid} equal $\frac{D_bB}2(X^2+|Y|^2+(k_b-\rho_bd_b^2)|\nabla z|^2)+c_E(\partial_1B)XY_1$.

Multiplying the equation of $W_\sigma$ in \eqref{e:gap} by $2W_\sigma$,
\[
 \partial_t(W_\sigma^2)+s_\sigma\partial_y(W_\sigma^2)+\frac23\partial_y(W_\sigma^3)
 =d_\sigma\partial_y^2(W_\sigma^2)-2d_\sigma(\partial_yW_\sigma)^2 ;
\]
integration over $(-\infty,y)$ for $\sigma=1$ and over $(y,\infty)$ for $\sigma=-1$, with
$\sigma s_\sigma=|s_\sigma|$, gives $\partial_tB_\sigma^{\rm sq}=-|s_\sigma|W_\sigma^2+R_{B,\sigma}$, and
$\partial_yB_\sigma^{\rm sq}=\sigma W_\sigma^2$. By \eqref{e:chain}, $D_b(F\circ Y)=(\partial_tF)\circ Y$ and
$\partial_1(F\circ Y)=\rho_b(\partial_yF)\circ Y$. The bounds of Lemma~\ref{l:pr} give
$\Vert W_\sigma\Vert_{L^2}^2\le C_M\tau^{-\frac12}$,
$\Vert W_\sigma\partial_yW_\sigma\Vert_{L^\infty}+\Vert W_\sigma^3\Vert_{L^\infty}\le C_M\tau^{-\frac32}$ and
$\Vert\partial_yW_\sigma\Vert_{L^2}^2\le C_M\tau^{-\frac32}$, hence the bounds of
$B_\sigma^{\rm sq}$ and $R_{B,\sigma}$; with \eqref{e:gap} they also give the boundedness of the
weights and of their $D_b$ and $\partial_1$ derivatives required by \eqref{e:wFC}.

The two crossing identities follow together from \eqref{e:acsys}:
\[
\begin{aligned}
 (D_b-\sigma c_E\partial_1)r_{-\sigma}^2
 &=-c_Er_{-\sigma}\Div'Y_\perp+r_{-\sigma}(R_X-\sigma(R_Y)_1),\\
 (D_b-\sigma c_E\partial_1)(r_\sigma r_{-\sigma})
 &=-2\sigma c_E(\partial_1r_\sigma)r_{-\sigma}
 -\frac{c_E}2(r_\sigma+r_{-\sigma})\Div'Y_\perp\\
 &\quad+\frac12\bigl[(r_\sigma+r_{-\sigma})R_X
 +\sigma(r_{-\sigma}-r_\sigma)(R_Y)_1\bigr].
\end{aligned}
\]
Multiply by $w_\sigma$ and integrate. Reynolds' identity and one integration
by parts give the same coefficient for both products; by \eqref{e:gap},
\[
 D_bw_\sigma+h_bw_\sigma-\sigma\partial_1(c_Ew_\sigma)
 =-2s_\sigma((\partial_yW_\sigma)\circ Y)+a_\sigma^{\rm err}.
\]
This gives \eqref{e:cross} with the stated errors. All operations are bilinear,
and their coefficients depend only on $(x_1,t)$, so the identities commute with
transverse multipliers.
\end{proof}

For a triple $f=(P_f,v_f,z_f)$, put
\begin{equation}\label{e:btest}
\begin{aligned}
 T_0[f]&=v_f-d_b\nabla z_f+\epsilon_K\nabla JP_f,&
 Z_0[f]&=\Div(j_bv_f)-\Div(k_b\nabla z_f),\\
 T_i[f]&=\partial_iv_f-d_b\nabla\partial_iz_f+\epsilon_K\nabla J\partial_iP_f,&
 Z_i[f]&=\Div(j_b\partial_iv_f)-\Div(k_b\nabla\partial_iz_f)\quad(1\le i\le3).
\end{aligned}
\end{equation}

With $T_0[f],Z_0[f],T_i[f],Z_i[f]$ of \eqref{e:btest} and $G_b$ of \eqref{e:Gbrelation}, put, for a
triple $f$ and a scalar $z_0$,
\begin{equation}\label{e:Bent}
\begin{aligned}
 \mathscr B_\eta[f]={}&\int_{\R^3}
 \Bigl[T_0[f]\cdot G_b+(\partial_1P_f)b\eta
 +Z_0[f]\bigl\{B_b\eta-(\gamma-1)\theta_b\partial_1(b\eta)\bigr\}\Bigr]\dd x\\
 &-\epsilon_K\int_{\R^3}\rho_bv_f\cdot
 \nabla J\bigl(\gamma p_b\partial_1(b\eta)\bigr)\dd x
 +\sum_{i=1}^3\int_{\R^3}
 \bigl[\rho_b(T_i[f])_1\partial_i(A_b\eta)
 +Z_i[f]\partial_i(B_b\eta)\bigr]\dd x\\
 &+\frac{\theta_W}{\tau}\int_{\R^3}\rho_bz_fB_b\eta\dd x,\\
 \mathscr T_\chi[f;z_0]={}&\int_{\R^3}
 \bigl[\rho_bA_bz_0(T_0[f])_1+B_bz_0Z_0[f]\bigr]\dd x
 +\sum_{i=1}^3\int_{\R^3}
 \bigl[\rho_b(T_i[f])_1\partial_i(A_bz_0)
 +Z_i[f]\partial_i(B_bz_0)\bigr]\dd x\\
 &+\frac{\theta_W}{\tau}\int_{\R^3}\rho_bz_fB_bz_0\dd x.
\end{aligned}
\end{equation}

We first record the pairing needed for the nonlinear inertia term. Let $K\ge1$
and let a vector $T$ satisfy
\[
 N_0(T;s,t)\le K\mathfrak B^{\frac12},\qquad
 \Vert T\Vert_{L^3}+\Vert T\Vert_{L^\infty}\le K\mu_*,\qquad
 \tau|T|\le K\mu_*\quad\text{on }\{|Y+h_*|\le\tfrac32c_0\tau\}.
\]
Then
\begin{equation}\label{e:inertial-work}
\begin{aligned}
 \Bigl\Vert T\nabla\frac\phi\rho\Bigr\Vert_{L^2_w}
 &\le C_M\mu_*\bigl(KN_o(s,t)+N_0(T;s,t)\bigr),\\
 \int_s^tw_s\Bigl|\int_{\R^3}\frac\phi\rho
 (\nabla\Pi-\mathcal L\psi)\cdot T\dd x\Bigr|\dd r
 &\le C_{M,K}\mu_*\mathfrak B.
\end{aligned}
\end{equation}
Indeed, \eqref{e:rhoD}, \eqref{e:state} and
$\Vert\nabla\zeta\Vert_{L^6}\le C\Vert\nabla^2\zeta\Vert_{L^2}$ give
\[
 \Bigl\Vert T\nabla\frac\phi\rho\Bigr\Vert_{L^2}
 \le C_MK\mu_*\bigl(\Vert\nabla\Pi\Vert_{L^2}+\Vert\nabla^2\zeta\Vert_{L^2}\bigr)
 +C_M\mu_*\Vert(|\partial_1\rho_b|+|\partial_1\theta_b|)T\Vert_{L^2}.
\]
Here the factors $\phi,\zeta$ multiplying background derivatives are bounded
by \eqref{e:stb}. Integration and \eqref{e:E1Z} prove the first line;
the bound $N_0(T;s,t)\le K\mathfrak B^{\frac12}$ is needed only for the second.
For the pressure work, use \eqref{e:espl} in the form
$\phi/\rho=\rho_b\Pi/(\gamma\rho p_b)-\eta-f$ and $f=F$ on
$\{|Y+h_*|\le\tfrac32c_0\tau\}$ to obtain
\[
 \Bigl|\int_{\R^3}\frac\phi\rho\nabla\Pi\cdot T\dd x\Bigr|
 \le C_MK\mu_*\Vert\nabla\Pi\Vert_{L^2}
 \Bigl(\Vert\nabla T\Vert_{L^2}+\Vert\nabla\eta\Vert_{L^2}
 +\tau^{-1}\Vert F\Vert_{L^2}
 +\Vert\mathbf1_{\{|Y+h_*|>\frac32c_0\tau\}}f\Vert_{L^2}\Bigr).
\]
The last two factors are controlled in $L^2_w$ by \eqref{e:Ch} and
\eqref{e:ext}, respectively. For the viscous work,
\[
 \int_{\R^3}\frac\phi\rho\mathcal L\psi\cdot T\dd x
 =-\int_{\R^3}\mathbb S(\psi):
 \Bigl(\frac\phi\rho\nabla T+T\otimes\nabla\frac\phi\rho\Bigr)\dd x.
\]
The first line of \eqref{e:inertial-work} and
$\Vert\phi/\rho\Vert_{L^\infty}\le C_M\mu_*$ now give the second line.

\begin{lemma}\label{l:sd}
Let $K\ge1$, let $\mathsf L(t)$ be a real even transverse convolution operator whose
kernel has $L^1$ norm at most $K$ for every $t$, or $\mathsf L=I$, and put
$g_L=\mathsf Lg_c$. Then
\[
\begin{gathered}
 \int_s^tw_s|\mathscr W_b(g_L,0;F^N)|\dd r
 +\int_s^tw_s|\mathscr T_\chi[g_L;z_N]|\dd r
 \le C_{b,K}\mu_*\mathfrak B,\\
 \int_s^tw_s|\mathscr B_\eta[g_L]|\dd r
 \le\bigl(C_{M,K}\delta+C_{b,K}\mu_*+\omega_{b,K}(s)\bigr)\mathfrak B,
\end{gathered}
\]
and
\[
\begin{aligned}
 &\sup_{s\le r\le t}w_s\Vert T_0[g_L]\Vert_{L^2}^2
 +\int_s^tw_s\bigl[\Vert\nabla T_0[g_L]\Vert_{L^2}^2
 +\tau^{-1}\Vert T_0[g_L]\Vert_{L^2}^2
 +\Vert Z_0[g_L]\Vert_{L^2}^2\bigr]\dd r\le C_{M,K}\mathfrak B,\\
 &\sum_{i=1}^3\int_s^tw_s\bigl[\Vert T_i[g_L]\Vert_{H^1}^2
 +\Vert Z_i[g_L]\Vert_{L^2}^2\bigr]\dd r
 +\int_s^tw_s\Vert\nabla^2(\mathsf LP_g)\Vert_{L^2}^2\dd r
 \le C_{b,K}\mathfrak B.
\end{aligned}
\]
\end{lemma}

\begin{proof}
Young's inequality for convolutions in $x'$ gives
$\Vert\mathsf LF\Vert_{L^p}\le K\Vert F\Vert_{L^p}$, $2\le p\le\infty$, and $\mathsf L$ commutes with
$\partial_i$, with $J$ and with functions of $(x_1,t)$. Hence the test bounds follow from
Lemma~\ref{l:cor} and from
$\int_s^tw_s\Vert\nabla^2P_g\Vert_{L^2}^2\dd r\le C_b(\mathcal J+\mathcal H)+\omega_b(s)N_o^2$, which
the computation of $N[g_c]$ in the proof of Lemma~\ref{l:cor} gives. Put $T=T_0[g_L]$. By
Lemma~\ref{l:cor} and $\nabla K_J\in L^1$,
$\Vert T\Vert_{L^3}+\Vert T\Vert_{L^\infty}\le C_{b,K}\mu_*$.

The sources \eqref{e:FN} satisfy $\Vert F^N\Vert_{L^2_w}\le C_b\mu_*N_o$ and
$\Vert\nabla F_z^N\Vert_{L^2_w}\le C_b\mu_*(I_s+N_o^2)^{\frac12}$: in
\[
\begin{aligned}
 \nabla F_z^N={}&-\nabla(\psi\cdot\nabla\zeta)
 -(\gamma-1)\nabla(\zeta\Div\psi)
 -\nabla\Bigl(\frac{\kappa\phi}{c_v\rho_b\rho}\Delta\zeta\Bigr)
 +\nabla\Bigl(B_b\frac{\phi}{\rho}\frac{\Pi}{p_b}\Bigr)\\
 &+\nabla\Bigl(\frac{\phi}{c_v\rho_b\rho}(h_b\Pi-H[\psi])\Bigr)
 +\frac{\nabla Q_\mu(\nabla\psi)}{c_v\rho}
 -\frac{Q_\mu(\nabla\psi)\nabla\rho}{c_v\rho^2}
\end{aligned}
\]
the first two groups use \eqref{e:gn4} and \eqref{e:prod}, the third uses
$\nabla^3\zeta\in L^2_w$, the $B_b$ and $h_b$ groups use \eqref{e:E3}, and the last two use
$|\nabla Q_\mu(\nabla\psi)|\le C|\nabla\psi||\nabla^2\psi|$, \eqref{e:rhoD} and \eqref{e:ha}.

In the zero slots of $\mathscr W_b(g_L,0;F^N)$,
\[
\begin{gathered}
 \Bigl|\int P_{g_L}\psi\cdot\nabla\Pi\dd x\Bigr|
 \le C\Vert\psi\Vert_{L^3}\Vert\nabla\Pi\Vert_{L^2}\Vert\nabla P_{g_L}\Vert_{L^2},\\
 \Bigl|\int\rho_b(\psi\cdot\nabla\psi)\cdot T\dd x\Bigr|
 \le C\Vert\psi\Vert_{L^3}\Vert\nabla\psi\Vert_{L^2}\Vert\nabla T\Vert_{L^2},\\
 \Bigl|\int\rho_bA_b\frac\phi\rho\frac\Pi{p_b}T_1\dd x\Bigr|
 \le C_M\Vert\Pi\Vert_{L^3}\Bigl\Vert A_b\frac\phi\rho\Bigr\Vert_{L^2}\Vert\nabla T\Vert_{L^2},
\end{gathered}
\]
the other terms of $F_p^N$ are bounded in the same way, and
$\Vert A_b\frac\phi\rho\Vert_{L^2_w}\le C_b\mathfrak B^{\frac12}$ by \eqref{e:E3}, \eqref{e:ABF} and
\eqref{e:ext}, since \eqref{e:espl} gives
\[
 \frac\phi\rho=\frac{\frac\Pi{\gamma p_b}-\eta-f}{1+\frac\Pi{\gamma p_b}} .
\]
To apply \eqref{e:inertial-work}, it remains to check the bound of $T$ near the contact.
The sets $\{|Y+h_*|\le\tfrac32c_0\tau\}$ and $\{|Y+h_*|>2c_0\tau\}$ are at
distance at least $c_M\tau$, since $\partial_1Y=\rho_b\le C_M$.
Young's inequality, the exponential tail of $\nabla K_J$ in Lemma~\ref{l:bessel}
and Lemma~\ref{l:cor} give, for $p=3,\infty$,
\[
 \Vert\mathbf1_{\{|Y+h_*|\le\frac32c_0\tau\}}\nabla JP_g\Vert_{L^p}
 \le C\Vert\mathbf1_{\{|Y+h_*|\le2c_0\tau\}}P_g\Vert_{L^p}
 +Ce^{-c_M\tau}\Vert P_g\Vert_{L^p}
 \le C_be^{-c_br}\mu_* .
\]
These indicators commute with $\mathsf L$, which acts only in $x'$.
The same local bound for $v_g,\nabla z_g$ in Lemma~\ref{l:cor} therefore gives
$\tau|T|\le C_{b,K}\mu_*$ on $\{|Y+h_*|\le\tfrac32c_0\tau\}$.
Thus \eqref{e:inertial-work} controls the remaining momentum term
$\frac\phi\rho(\nabla\Pi-\mathcal L\psi)$ and the product
$T\nabla(\phi/\rho)$ used below.
In the first slots the derivatives of $F_p^N,F_m^N$ are moved onto the tests, keeping
$-\frac{\partial_i\rho_b}{\rho_b}F_m^N$, and the thermal slots use $Z_i[g_L]$ and
$\nabla F_z^N$. In the screened slot
$\int\rho_bv_{g_L}\cdot\nabla JF_p^N\dd x=-\int\Div(\rho_bv_{g_L})JF_p^N\dd x$ and
$\Vert JF_p^N\Vert_{L^2_w}\le C_b\mu_*N_o$, since $J:L^{\frac65}\to L^2$ and
$\Vert Q_\mu(\nabla\psi)\Vert_{L^{6/5}}\le C\Vert\psi\Vert_{L^3}\Vert\nabla^2\psi\Vert_{L^2}$; the weak value
uses $\frac{z_{g_L}}\tau$. This proves the bound of $\mathscr W_b(g_L,0;F^N)$.

For $z_N=-\frac{(\gamma-1)\phi\Pi}{\gamma\rho p_b}$,
\[
 \Bigl|\int\rho_bA_bz_NT_1\dd x\Bigr|
 \le C_M\Vert\Pi\Vert_{L^3}\Bigl\Vert A_b\frac\phi\rho\Bigr\Vert_{L^2}\Vert\nabla T\Vert_{L^2},
\]
and \eqref{e:E3}, \eqref{e:E6} and \eqref{e:ehp} give
$\Vert B_bz_N\Vert_{L^2_w}+\Vert\nabla(A_bz_N,B_bz_N)\Vert_{L^2_w}\le C_b\mu_*N_o$. Pairing with
the test bounds proves the bound of $\mathscr T_\chi[g_L;z_N]$.

By \eqref{e:Gbrelation},
\[
 \int T\cdot G_b\dd x=\Bigl\langle G_\eta,\frac{\rho_b}\rho T\Bigr\rangle
 +\int T\cdot\Bigl[\frac\phi\rho\mathcal L(b\eta e_1)
 +\rho_bh_b\psi_1\eta e_1+\rho_bb(\psi\cdot\nabla\eta)e_1\Bigr]\dd x .
\]
Since $\nabla\frac{\rho_b}\rho=-\nabla\frac\phi\rho$, the bounds above give
$N_0(\frac{\rho_b}\rho T;s,t)\le C_{b,K}\mathfrak B^{\frac12}$, and \eqref{e:Getb} gives
$\int_s^tw_s|\langle G_\eta,\frac{\rho_b}\rho T\rangle|\dd r\le\omega_{b,K}(s)\mathfrak B$. In the
other three terms,
$\int\frac\phi\rho\mathcal L(b\eta e_1)\cdot T\dd x=-\int\mathbb S(b\eta e_1):[\frac\phi\rho\nabla T
+T\otimes\nabla\frac\phi\rho]\dd x$ is bounded with \eqref{e:eta}, and
\[
 \Bigl|\int T_1\rho_b\bigl(h_b\psi_1\eta+b\psi\cdot\nabla\eta\bigr)\dd x\Bigr|
 \le C\Vert\psi\Vert_{L^3}\Vert\nabla\eta\Vert_{L^2}
 \bigl(\Vert h_b\Vert_{L^\infty}\Vert T\Vert_{L^2}+\Vert b\Vert_{L^\infty}\Vert\nabla T\Vert_{L^2}\bigr)
\]
with \eqref{e:E2} and \eqref{e:bsup}. Hence
$\int_s^tw_s|\int T\cdot G_b\dd x|\dd r\le(C_{b,K}\mu_*+\omega_{b,K}(s))\mathfrak B$. By
\eqref{e:E3}, \eqref{e:E6} and \eqref{e:eta},
\[
 \Vert\nabla(b\eta)\Vert_{L^2_w}
 +\Vert\nabla(A_b\eta,B_b\eta)\Vert_{L^2_w}
 +\Vert B_b\eta\Vert_{L^2_w}\le\omega_b(s)\Phi^{\frac12},
\]
which bounds the zero thermal slot, the first slots and the weak value of $\mathscr B_\eta[g_L]$.
The zero pressure and screened terms use
$\Vert b\eta\Vert_{L^2_w}^2\le(C_M\delta^2+\omega_b(s))(I_s+N_o^2)$ from \eqref{e:eta}, the
second after $p_b\partial_1(b\eta)=\partial_1(p_bb\eta)-(\partial_1p_b)b\eta$ and the transfer of $\partial_1$
onto $\rho_bv_{g_L}$:
\[
 \int_s^tw_s\Bigl|\int(\partial_1P_{g_L})b\eta\dd x\Bigr|\dd r
 +\epsilon_K\int_s^tw_s
 \Bigl|\int\rho_bv_{g_L}\cdot\nabla J\bigl(\gamma p_b\partial_1(b\eta)\bigr)\dd x\Bigr|\dd r
 \le\bigl(C_{M,K}\delta+\omega_{b,K}(s)\bigr)\mathfrak B .
\]
This proves the bound of $\mathscr B_\eta[g_L]$.
\end{proof}

With $\epsilon_K,\theta_W$ of Lemma~\ref{l:mr}, for a triple $f=(P_f,v_f,z_f)$ and a
scalar $\eta$, put
\begin{equation}\label{e:Eb}
\begin{aligned}
 \mathcal E_b(f,\eta)={}&\int_{\R^3}e_b^0(f,\eta)\dd x
 +\epsilon_K\int_{\R^3}\rho_bV_f\cdot\nabla JP_f\dd x\\
 &+\sum_{i=1}^3\int_{\R^3}\bigl[e_b^0(\partial_if,0)
 +\epsilon_K\rho_b\partial_iv_f\cdot\nabla J\partial_iP_f\bigr]\dd x
 +\frac{\theta_W}{2\tau}\int_{\R^3}\rho_bz_f^2\dd x.
\end{aligned}
\end{equation}

For $f=q_o,q_s$, the background identity for $(f,\eta)$ is the sum of \eqref{e:fcid} and
$\epsilon_K$ times \eqref{e:fcsc}, with $c=b$, for $(P_f,V_f,z_f)$ and for
$\partial_if$, $1\le i\le3$, and of $\theta_W$ times
\eqref{e:fcwv} for $(P_f,V_f,z_f)$, with the sources of Lemma~\ref{l:res} and of
\eqref{e:fcal}. With $(P_f^0,v_f^0,z_f^0)=(P_f,V_f,z_f)$ and $(P_f^i,v_f^i,z_f^i)=\partial_if$, its
left side is $\frac{\dd}{\dd t}\mathcal E_b(f,\eta)+D_b[f]$, where $D_b[f]$ is $D_{\mathcal E}$ of
Lemma~\ref{l:mr} with $\mathcal D_b$, $\rho_b$, these triples and without the term with $\eta$, and
we put
\[
 \mathcal R_{0,b}[f]=\epsilon_K\sum_{i=0}^3\mathscr P_b(P_f^i,v_f^i,z_f^i)
 -\frac{(\gamma-1)\theta_W}\tau\int_{\R^3}\rho_b\theta_bz_f\Div v_f\dd x .
\]
For a triple $f$ put
\[
\begin{gathered}
 H_b[f]=-\gamma j_bV_f\cdot\nabla z_f
 +\bigl((\gamma-1)k_b+\rho_bd_b^2\bigr)|\partial_1z_f|^2+\gamma k_b|\nabla'z_f|^2,\\
 Q_h^b[f]=-\rho_b\bigl((V_f)_1-d_b\partial_1z_f\bigr)^2-\frac{\gamma-1}{2\gamma p_b}P_f^2+H_b[f]
 =-(Y_f)_1^2-\frac{\gamma-1}2X_f^2+H_b[f],
\end{gathered}
\]
and write $X=X_{q_s}$, $Y=Y_{q_s}$, $r_\sigma=r_{\sigma,q_s}$ and $w_\sigma=W_\sigma\circ Y$.

\begin{lemma}\label{l:cd}
One has
\[
\begin{gathered}
 \frac1{\theta_b}\Bigl(z_s-\frac{P_s}{c_p\rho_b}\Bigr)=\chi_o^b-\chi_g,\\
 \Bigl|\mathcal E_N-\mathcal E_b(q_o,\eta)-\frac{\theta_\eta}2\int_{\R^3}\rho\eta^2\dd x\Bigr|
 \le C_b\mu_*\bigl[E[q_o]+\Vert b\Vert_{L^\infty}^2\Vert\eta\Vert_{L^2}^2\bigr],\\
 \bigl|D_{\mathcal E}-D_b[q_o]-\theta_\eta c\Vert\nabla\eta\Vert_{L^2}^2\bigr|
 \le C_b\mu_*\bigl[D_*+\Vert\nabla(b\eta)\Vert_{L^2}^2\bigr],
\end{gathered}
\]
and the background identity for $(q_s,\eta)$ minus that for $(q_o,\eta)$ reads
\begin{equation}\label{e:cdif}
 \frac{\dd}{\dd t}\bigl[\mathcal E_b(q_s,\eta)-\mathcal E_b(q_o,\eta)\bigr]+D_b[q_s]-D_b[q_o]
 =\mathcal R_{0,b}[q_s]-\mathcal R_{0,b}[q_o]+\mathcal R_b
 +\mathscr W_b(q_s,\eta;\mathscr D+\mathscr T)+\mathcal O_b,
\end{equation}
with $\int_s^tw_s|\mathcal O_b|\dd r\le(C_M\delta+C_b\mu_*+\omega_b(s))\mathfrak B$ and
\[
\begin{aligned}
 \mathcal R_b={}&\int_{\R^3}h_b\bigl\{Q_h^b[q_s]-Q_h^b[q_o]\bigr\}\dd x
 +\frac{\gamma-1}\gamma\int_{\R^3}d_b\frac{\partial_1p_b}{p_b}(P_s\partial_1z_s-\Pi\partial_1\zeta)\dd x\\
 &+\int_{\R^3}\rho_bA_bF\bigl\{(V_s)_1-d_b\partial_1z_s-V_1+d_b\partial_1\zeta\bigr\}\dd x .
\end{aligned}
\]
With $Q_h$ and $\mathcal R_{\rm princ}$ of \eqref{e:works} and
$\mathcal R_h=-\int_{\R^3}h_b\bigl[Y_1^2+\frac{\gamma-1}2X^2\bigr]\dd x$,
\[
\begin{aligned}
 \mathcal R_{\rm princ}+\mathcal R_b={}&\mathcal R_h+\int_{\R^3}h_bH_b[q_s]\dd x
 +\frac{\gamma-1}\gamma\int_{\R^3}d_b\frac{\partial_1p_b}{p_b}P_s\partial_1z_s\dd x
 +\int_{\R^3}\rho A_bF(V_1-d_N\partial_1\zeta)\dd x\\
 &+\int_{\R^3}\rho_bA_bF\bigl\{(V_s)_1-d_b\partial_1z_s-V_1+d_b\partial_1\zeta\bigr\}\dd x+\mathcal C_o,
\end{aligned}
\]
$\int_s^tw_s|\mathcal C_o|\dd r\le(C_b\mu_*+\omega_b(s))\mathfrak B$, where $\mathcal C_o=\int h_b(Q_h-Q_h^b[q_o])\dd x+\frac{\gamma-1}\gamma\int(d_N-d_b)\frac{\partial_1p_b}{p_b}\Pi\partial_1\zeta\dd x$.
The tests of $\mathscr D$ and $\mathscr T$ satisfy
\begin{equation}\label{e:Dbd}
 \int_s^tw_s\Vert\mathscr D\Vert_{H^1}^2\dd r\le\varepsilon_{D,b}(s)N_s^2,\qquad
 \varepsilon_{D,b}(s)=C_b\Bigl[\tau_s^{-1}+\int_s^\infty\tau^{-1}
 \Vert\partial_1^2Q_b\Vert_{L^2(\R)}^2\dd r\Bigr]\longrightarrow0,
\end{equation}
\begin{equation}\label{e:Dslots}
\begin{aligned}
 \int_s^tw_s\bigg[&\Bigl|\int_{\R^3}Z_{0,s}\mathscr D_z\dd x\Bigr|
 +\epsilon_K|\mathscr K_b(P_s,V_s;\mathscr D_p,\mathscr D_m)|
 +\frac{\theta_W}\tau\Bigl|\int_{\R^3}\rho_bz_s\mathscr D_z\dd x\Bigr|\\
 &+\sum_{i=1}^3\bigl(|\mathscr S_b(\partial_iP_s,\partial_iv_s,\partial_iz_s;\mathscr D^{[i]})|
 +\epsilon_K|\mathscr K_b(\partial_iP_s,\partial_iv_s;\mathscr D^{[i]}_p,\mathscr D^{[i]}_m)|\bigr)
 \bigg]\dd r\le\omega_b(s)\mathfrak B,
\end{aligned}
\end{equation}
with $Z_{0,s}=\Div(j_bV_s)-\Div(k_b\nabla z_s)$ and $\mathscr D^{[i]}$ as in \eqref{e:srcw}, and
\[
 \int_s^tw_s\Bigl|\mathscr W_b(q_s,\eta;\mathscr T)
 -\int_{\R^3}(\mathscr T_m)_\perp\cdot(T_0[q_s])_\perp\dd x\Bigr|\dd r\le\omega_b(s)\mathfrak B .
\]
With
$S_{\mathscr D}=\int_{\R^3}[(V_s-d_b\nabla z_s)\cdot\mathscr D_m+\frac{P_s\mathscr D_p}{\gamma p_b}]\dd x$,
$\beta_{\rm cr}=\frac{2(\gamma-3)}{\gamma+1}$,
\[
\begin{gathered}
 E_{\rm th}=-2\sum_\sigma\int_{\R^3}\partial_1\Bigl(\frac{r_\sigma}{\rho_b}\Bigr)\rho_b(U_\sigma\circ Y)
 J\Bigl(\frac{\sigma\sqrt{\rho_b}d_b}2\partial_1z_s\Bigr)\dd x,\\
 \mathcal R_n=-\sum_\sigma\int_{\R^3}((\partial_yW_\sigma)\circ Y)
 \bigl(r_{-\sigma}^2+\beta_{\rm cr}r_+r_-\bigr)\dd x,
\end{gathered}
\]
one has
\[
 S_{\mathscr D}=\sum_\sigma\int_{\R^3}(\partial_1(U_\sigma\circ Y))r_\sigma^2\dd x+E_{\rm th}+E_{\rm sel},
 \qquad
 \mathcal R_h+S_{\mathscr D}=\mathcal R_n+E_{\rm th}+E_{\rm n},
\]
\[
 \int_s^tw_s|E_{\rm sel}|\dd r\le\omega_b(s)\mathfrak B,\qquad
 \int_s^tw_s|E_{\rm n}|\dd r\le\bigl(C_M\delta+\omega_b(s)\bigr)\mathfrak B,
\]
and, for every $\varepsilon>0$,
\[
 \int_s^tw_s|E_{\rm th}|\dd r
 \le\varepsilon\sum_\sigma\int_s^tw_s\Vert\nabla r_\sigma\Vert_{L^2}^2\dd r
 +C_\varepsilon\sum_\sigma\int_s^tw_s\int_{\R^3}w_\sigma^2|\nabla z_s|^2\dd x\dd r
 +\omega_b(s)\mathfrak B .
\]
For every symmetric $\R$-bilinear form $\langle\cdot,\cdot\rangle$, all $X,Y_1$ in its domain,
$r_\pm=\frac{X\pm Y_1}2$ and $\sigma=\pm1$,
\begin{equation}\label{e:bcr}
 -\frac2{\gamma+1}\Bigl[\langle Y_1,Y_1\rangle+\frac{\gamma-1}2\langle X,X\rangle\Bigr]
 +\langle r_\sigma,r_\sigma\rangle
 =-\langle r_{-\sigma},r_{-\sigma}\rangle-\beta_{\rm cr}\langle r_+,r_-\rangle .
\end{equation}
With
\[
 C_{\rm cr}=-\sum_\sigma\frac1{2s_\sigma}\int_{\R^3}w_\sigma
 \bigl(r_{-\sigma}^2+\beta_{\rm cr}r_+r_-\bigr)\dd x,
\]
\begin{equation}\label{e:crj}
\begin{aligned}
 J_{\rm cr}={}&\sum_\sigma\frac1{2s_\sigma}\int_{\R^3}w_\sigma c_E
 \Bigl[\Bigl(r_{-\sigma}+\frac{\beta_{\rm cr}}2(r_\sigma+r_{-\sigma})\Bigr)\Div'Y_\perp
 +2\sigma\beta_{\rm cr}(\partial_1r_\sigma)r_{-\sigma}\Bigr]\dd x,\\
 E_{\rm cr}={}&-\sum_\sigma\frac{\operatorname{Err}_1+\beta_{\rm cr}\operatorname{Err}_2}{2s_\sigma}
\end{aligned}
\end{equation}
and $\operatorname{Err}_1,\operatorname{Err}_2$ of \eqref{e:cross} for $(P_s,V_s,z_s)$, one has
\[
\begin{gathered}
 \frac{\dd}{\dd t}C_{\rm cr}=-\mathcal R_n+J_{\rm cr}+E_{\rm cr},\qquad
 \int_s^tw_s|E_{\rm cr}|\dd r\le\bigl(C_M\delta+C_b\mu_*+\omega_b(s)\bigr)\mathfrak B,\\
 |J_{\rm cr}|\le\varepsilon\Vert\nabla(X,Y)\Vert_{L^2}^2
 +C_\varepsilon\sum_\sigma\int_{\R^3}w_\sigma^2\bigl(r_{-\sigma}^2+|Y_\perp|^2\bigr)\dd x
 \qquad(\varepsilon>0),
\end{gathered}
\]
and
\[
\begin{gathered}
 |C_{\rm cr}(t)|\le C_M\tau^{-\frac12}E[q_s](t)+\omega_b(s)\Vert\eta(t)\Vert_{L^2}^2,\qquad
 |C_{\rm cr}(s)|\le C_M\tau_s^{-\frac12}I_s,\\
 0\le\sum_\sigma\int_{\R^3}(B_\sigma^{\rm sq}\circ Y)e_b^0(q_s,\eta)\dd x
 \le C_M\tau^{-\frac12}\bigl[E[q_s]+\Vert b\Vert_{L^\infty}^2\Vert\eta\Vert_{L^2}^2\bigr],\\
 \int_s^tw_s\frac{\dd C_{\rm cr}}{\dd r}\dd r
 =w_s(t)C_{\rm cr}(t)-C_{\rm cr}(s)+q_a\int_s^t\frac{w_s}\tau C_{\rm cr}\dd r .
\end{gathered}
\]
\end{lemma}

\begin{proof}
The first identity follows from $q_s=q_o-g_c$ and \eqref{e:DT}. In \eqref{e:men} and in
$\mathcal E_b(q_o,\eta)$ the triples coincide, and the coefficients differ by $\rho-\rho_b=\phi$,
$p-p_b=\Pi$ and the differences of $d,j,k$ of \eqref{e:FC}, all bounded by $C_b\mu_*$ by
\eqref{e:stb}; with $\Vert V\Vert_{L^2}\le\Vert\psi\Vert_{L^2}+\Vert b\Vert_{L^\infty}\Vert\eta\Vert_{L^2}$ this
gives the comparison of the energies, and the same argument with \eqref{e:fcd} gives that of
$D_{\mathcal E}$ and $D_b[q_o]$.

By Lemma~\ref{l:res}, $(P_j,V_j,z_j)$ solves \eqref{e:forced} with the sources \eqref{e:full},
and $(P_j,v_j,z_j)$ solves it with the sources \eqref{e:full} plus
$(\gamma p_b\partial_1(b\eta),\rho_bA_b\eta e_1-G_b,(\gamma-1)\theta_b\partial_1(b\eta))$, since
$D_bb=A_b$; the triples $\partial_iq_j$ solve it with the sources \eqref{e:fcal} of these. In
\eqref{e:fcwv} for $(P_j,V_j,z_j)$, $\Div V_j=\Div v_j-\partial_1(b\eta)$, so its right side pairs
the temperature source in the second system with $\theta_W\rho_bz_j/\tau$,
plus the last term of $\mathcal R_{0,b}[q_j]$. All
source terms of the background identities are linear, $q_s-q_o=-g_c$, and the terms of
\eqref{e:full} that contain neither $P_j$ nor $v_j$ are the same for $j=o,s$, except
$\mathscr D+\mathscr T+F^g-F^\chi$. By $\widetilde F^N=F^N+(0,\rho_bA_bz_Ne_1,B_bz_N)$, the source
terms of the difference are therefore
\[
 -\mathscr W_b(g_c,0;F^N)-\mathscr T_\chi[g_c;z_N]-\mathscr B_\eta[g_c]
 -\mathscr W_b\bigl(g_c,0;(0,\rho_bA_bfe_1,B_bf)\bigr)
 +\mathscr W_b(q_s,\eta;\mathscr D+\mathscr T+F^g-F^\chi),
\]
and the other terms of the difference, besides $\mathcal R_{0,b}[q_s]-\mathcal R_{0,b}[q_o]$, are
the terms of the two identities with derivatives of $\rho_b,u_b,\theta_b$ in their
coefficients, including those of the sources linear in $(P_j,v_j)$. Lemma~\ref{l:sd} with
$\mathsf L=I$ and $K=1$ bounds the first three terms, and Lemma~\ref{l:res} gives
$\int_s^tw_s|\mathscr W_b(q_s,\eta;F^g-F^\chi)|\le\omega_b(s)\mathfrak B$. In
$\mathscr W_b(g_c,0;(0,\rho_bA_bfe_1,B_bf))$, $-\int\rho_bA_bF((v_g)_1-d_b\partial_1z_g)\dd x$ is the
last term of $\mathcal R_b$; its other terms contain $f_{\rm out}$, bounded by \eqref{e:ext}, or
$A_bF$, $B_bF$ or their first derivatives paired with $\epsilon_K\partial_1JP_g$, $Z_0[g_c]$,
$T_i[g_c]$, $Z_i[g_c]$ or $\frac{z_g}\tau$, bounded by \eqref{e:ABF}, \eqref{e:Ch} and
Lemma~\ref{l:cor}. The terms with coefficient derivatives are grouped as in the proof of
\eqref{e:mr}: the terms with the factor $\partial_1p_b$ give
$\frac{\gamma-1}\gamma\int d_b\frac{\partial_1p_b}{p_b}P_j\partial_1z_j\dd x-\int\frac{\partial_1p_b}{\gamma p_b}P_jb\eta\dd x$,
by \eqref{e:bal}, the zero-order terms with the factor $h_b$ give
$\int h_bQ_h^b[q_j]\dd x$ plus terms with the factor $b\eta$, and all the others are bounded as the
corresponding terms of $\mathcal R_A-\mathcal R_{\rm princ}$ there, with $N_s\le C_MN_o$. This
proves \eqref{e:cdif}. In $\mathcal C_o$ the coefficients of $Q_h-Q_h^b[q_o]$ are differences of
functions of $(\rho,\theta)$ and of $(\rho_b,\theta_b)$; writing $h_b=\partial_1b$ and
integrating by parts once, \eqref{e:bsup}, \eqref{e:rhoD}, \eqref{e:stb}, \eqref{e:prod} and
\eqref{e:pbd} bound the first term, and $d_N-d_b=-\frac{(\gamma-1)\kappa\Pi}{\gamma pp_b}$,
$\Vert\Pi\Vert_{L^4}^2\le\Vert\Pi\Vert_{L^3}\Vert\Pi\Vert_{L^6}$ and \eqref{e:E2} bound the second. The
identity for $\mathcal R_{\rm princ}+\mathcal R_b$ is the sum of the definitions.

By Lemma~\ref{l:pr} and $c_0<\frac18\min\{c_-,c_+\}$,
$|\partial_y^jW_\sigma|\le C_Me^{-c_M\tau}$ on $\R\setminus\Omega_\sigma$, and
$|\partial_y^jW_\sigma|\le C_M\tau^{-\frac{j+1}2}$; with \eqref{e:munw} and \eqref{e:E10}, for every
first $y$-derivative $a$ of a smooth function of $Z\in\mathcal K$ and $j=0,1$,
\begin{equation}\label{e:Wloc}
 |a\,\partial_y^jW_\sigma|\le\omega_b(s)\tau^{-\frac{j+2}2},\qquad
 \Bigl(\Bigl|1-\frac{v_\sigma}{V_b}\Bigr|
 +\Bigl|s_\sigma-\sigma\sqrt{\frac{\gamma P_b}{V_b}}\Bigr|\Bigr)|\partial_y^jW_\sigma|
 \le C_M\tau^{-\frac{j+2}2},
\end{equation}
since the two factors in brackets vanish at $Z_{\sigma,*}$ and $|Z_{\sigma,a}-Z_{\sigma,*}|\le
C_M\tau^{-\frac12}$. By \eqref{e:chain}, $\partial_1(G\circ Y)=\rho_b(\partial_yG)\circ Y$, so
$h_b$, $\partial_1\rho_b$, $\partial_1c_E$ and $\partial_1h_0$ are of the form $\rho_b(a\circ Y)$; also
$\rho_bc_E=\sqrt{\gamma P_b/V_b}\circ Y$ and $\partial_1(U_\sigma\circ Y)=\rho_bv_\sigma((\partial_yW_\sigma)\circ Y)$.

By \eqref{e:conv}, $\mathfrak PR_\sigma=C_\sigma^D$ with
$C_\sigma^D=(h_0,\sigma\sqrt{\rho_b}e_1,\frac{h_0}{c_p\rho_b})$, using $P_bV_b=R\Theta_b$ and
$(\gamma-1)c_p=\gamma R$. With $f_\sigma=\varphi_\sigma[q_s]\circ Y$, so that $(J_b\varphi_\sigma[q_s])\circ Y=Jf_\sigma$,
and $(\partial_yG)\circ Y=\rho_b^{-1}\partial_1(G\circ Y)$,
\[
\begin{aligned}
 \mathscr D&=\sum_\sigma
       \bigl(A_\sigma^D\partial_1Jf_\sigma+B_\sigma^DJf_\sigma\bigr),\\
 A_\sigma^D&=(U_\sigma\circ Y)C_\sigma^D,\\
 B_\sigma^D&=\frac{C_\sigma^D}{\rho_b}
                   \partial_1(\rho_b(U_\sigma\circ Y))
       +\rho_b(U_\sigma\circ Y)\,\mathfrak P(\partial_yR_\sigma).
\end{aligned}
\]
where $\mathfrak P(\partial_yR_\sigma)$ are the physical rows of $\partial_yR_\sigma$. By
Lemma~\ref{l:pr}, \eqref{e:Wloc} and \eqref{e:unwE},
\[
 |A_\sigma^D|\le C_b\tau^{-\frac12},\qquad
 |\partial_1A_\sigma^D|+|B_\sigma^D|\le C_b\tau^{-1},\qquad
 |\partial_1B_\sigma^D|\le C_b\bigl(\tau^{-\frac32}+\tau^{-\frac12}|\partial_1^2Q_b|\bigr),
\]
and $\Vert\nabla f_\sigma\Vert_{L^2}\le C_b(\Vert\nabla(P_s,v_s)\Vert_{L^2}+\tau^{-\frac12}\Vert(P_s,v_s)\Vert_{L^2})$.
The symbols of $\partial_1J$ and $\nabla\partial_1J$ are bounded by $|\xi|$, and that of $\partial_1J$ by
$\frac12$; hence $\Vert\partial_1Jf_\sigma\Vert_{L^2}+\Vert\nabla\partial_1Jf_\sigma\Vert_{L^2}\le2\Vert\nabla f_\sigma\Vert_{L^2}$
and, by \eqref{eq:planar}, $\sup_{x_1}\Vert Jf_\sigma(x_1,\cdot)\Vert_{L^2(\R^2)}^2\le2\Vert Jf_\sigma\Vert_{L^2}\Vert\partial_1Jf_\sigma\Vert_{L^2}
\le\Vert f_\sigma\Vert_{L^2}^2$; the product with $\partial_1^2Q_b$ uses its $L^2(\R)$ norm and
this trace bound. Hence \eqref{e:Dbd}, since $\int_s^\infty\Vert\partial_1^2Q_b\Vert_{H^1}^2\dd r\to0$ by
\eqref{e:unwE}. In \eqref{e:Dslots}, $\Div(\rho_bV_s)=\rho_b\Div v_s-\rho_b\partial_1(b\eta)
+(\partial_1\rho_b)(V_s)_1$ and $Z_{0,s}$ contain only $\nabla v_s$, $\partial_1(b\eta)$, $\nabla^2z_s$ and
first derivatives of $\rho_b,j_b,k_b$, so that, by Lemma~\ref{l:cor} and \eqref{e:eta},
\[
 \Vert\nabla(P_s,v_s)\Vert_{L^2_w}+\Vert\nabla^2z_s\Vert_{L^2_w}
 +\Vert Z_{0,s}\Vert_{L^2_w}+\sum_{i=1}^3\Vert Z_i[q_s]\Vert_{L^2_w}
 +\Vert\Div(\rho_bV_s)\Vert_{L^2_w}+\Bigl\Vert\frac{z_s}\tau\Bigr\Vert_{L^2_w}\le C_b\mathfrak B^{\frac12};
\]
with $\int\rho_bV_s\cdot\nabla J\mathscr D_p\dd x=-\int\Div(\rho_bV_s)J\mathscr D_p\dd x$, the
Cauchy--Schwarz inequality in each slot and \eqref{e:Dbd} give \eqref{e:Dslots}. By
Lemma~\ref{l:res}, $\mathscr T_p=\mathscr T_z=(\mathscr T_m)_1=0$, so the zero slots of
$\mathscr W_b(q_s,\eta;\mathscr T)$ reduce to $\int(\mathscr T_m)_\perp\cdot(T_0[q_s])_\perp\dd x$, and its
first slots are $\int T_i[q_s]\cdot\mathscr T_m^{[i]}\dd x$ with
$\mathscr T_m^{[i]}=\partial_i\mathscr T_m-\frac{\partial_i\rho_b}{\rho_b}\mathscr T_m$. By \eqref{e:chain},
\[
 \sum_{i=1}^3\Vert\mathscr T_m^{[i]}\Vert_{L^2}^2
 \le C_M\sum_\sigma\Vert\nabla_m^3u_\sigma^{\rm ac}\Vert_{L^2}^2
 +C_b\Vert\partial_1Q_b\Vert_{L^\infty}^2\sum_\sigma\Vert\nabla_m^2u_\sigma^{\rm ac}\Vert_{L^2}^2,
\]
and Lemma~\ref{l:cor} with \eqref{e:unwE} gives $\sum_i\Vert\mathscr T_m^{[i]}\Vert_{L^2_w}^2\le
\omega_b(s)N_o^2$; this proves the bound of $\mathscr W_b(q_s,\eta;\mathscr T)$.

The zero tests of $\mathscr D$ are $\frac{P_s}{\gamma p_b}=\frac X{h_0}$ and
$V_s-d_b\nabla z_s=\frac Y{\sqrt{\rho_b}}$, and $\frac X{h_0}h_0+\frac{Y_1}{\sqrt{\rho_b}}\sigma\sqrt{\rho_b}=2r_\sigma$.
Hence
\[
 S_{\mathscr D}=\sum_\sigma\Bigl[\,2\int_{\R^3}\frac{r_\sigma}{\rho_b}\partial_1(\rho_b(U_\sigma\circ Y)Jf_\sigma)\dd x
 +\int_{\R^3}\rho_b(U_\sigma\circ Y)\Bigl(\frac X{h_0}\mathfrak P(\partial_yR_\sigma)_p
 +\frac{Y}{\sqrt{\rho_b}}\cdot\mathfrak P(\partial_yR_\sigma)_m\Bigr)Jf_\sigma\dd x\Bigr],
\]
and one integration by parts turns the first integral into
$-2\int\partial_1(\frac{r_\sigma}{\rho_b})\rho_b(U_\sigma\circ Y)Jf_\sigma\dd x$. By \eqref{e:acg} and \eqref{e:be},
\[
 f_\sigma=r_\sigma+\frac{\sigma\sqrt{\rho_b}}2\bigl(d_b\partial_1z_s+b\eta\bigr),\qquad
 Jf_\sigma=r_\sigma+(J-I)r_\sigma+J\Bigl(\frac{\sigma\sqrt{\rho_b}}2b\eta\Bigr)
 +J\Bigl(\frac{\sigma\sqrt{\rho_b}d_b}2\partial_1z_s\Bigr),
\]
and
$-2\int\partial_1(\frac{r_\sigma}{\rho_b})\rho_b(U_\sigma\circ Y)r_\sigma\dd x
=\int(\partial_1(U_\sigma\circ Y))r_\sigma^2\dd x+2\int\frac{\partial_1\rho_b}{\rho_b}(U_\sigma\circ Y)r_\sigma^2\dd x$.
The last term of $Jf_\sigma$ gives $E_{\rm th}$, and
\[
\begin{aligned}
 E_{\rm sel}=\sum_\sigma\Bigl[&\,2\int_{\R^3}\frac{\partial_1\rho_b}{\rho_b}(U_\sigma\circ Y)r_\sigma^2\dd x
 -2\int_{\R^3}\partial_1\Bigl(\frac{r_\sigma}{\rho_b}\Bigr)\rho_b(U_\sigma\circ Y)
 \Bigl\{(J-I)r_\sigma+J\Bigl(\frac{\sigma\sqrt{\rho_b}}2b\eta\Bigr)\Bigr\}\dd x\\
 &+\int_{\R^3}\rho_b(U_\sigma\circ Y)\Bigl(\frac X{h_0}\mathfrak P(\partial_yR_\sigma)_p
 +\frac{Y}{\sqrt{\rho_b}}\cdot\mathfrak P(\partial_yR_\sigma)_m\Bigr)Jf_\sigma\dd x\Bigr].
\end{aligned}
\]
By \eqref{e:Wloc}, $|\partial_1\rho_b\,(U_\sigma\circ Y)|+|(U_\sigma\circ Y)((\partial_yR_\sigma)\circ Y)|\le\omega_b(s)\tau^{-1}$; the symbol
of $J-I$ is $-\frac{|\xi|^2}{1+|\xi|^2}$ and $\frac{|\xi|^2}{1+|\xi|^2}\le\frac{|\xi|}2$, so
$\Vert(J-I)r_\sigma\Vert_{L^2}\le\frac12\Vert\nabla r_\sigma\Vert_{L^2}$;
$|(U_\sigma\circ Y)|\le C_M\tau^{-\frac12}$, $\Vert Jg\Vert_{L^2}\le\Vert g\Vert_{L^2}$ and \eqref{e:eta} bound the
$b\eta$ term. With $\int_s^tw_s\tau^{-1}\Vert(X,Y,f_\sigma)\Vert_{L^2}^2\dd r+\int_s^tw_s\Vert\nabla(X,Y)\Vert_{L^2}^2\dd r
\le C_M\mathfrak B$, which follows from \eqref{e:B}, \eqref{e:E1Z} and \eqref{e:eta}, this gives the
bound of $E_{\rm sel}$.

By \eqref{e:Um} and \eqref{e:CH}, $b\circ X=b_m=H_c(y+h_*,t)+\sum_\sigma s_\sigma\vartheta_\sigma+(R_U)_2$, and
$s_\sigma\partial_y\vartheta_\sigma=\frac2{\gamma+1}\partial_yU_\sigma=\frac{2v_\sigma}{\gamma+1}\partial_yW_\sigma$ by
\eqref{e:acg}. Hence
\[
\begin{aligned}
 \mathcal R_h={}&-\frac2{\gamma+1}\sum_\sigma\int_{\R^3}((\partial_yW_\sigma)\circ Y)\Bigl(Y_1^2+\frac{\gamma-1}2X^2\Bigr)\dd x\\
 &-\frac2{\gamma+1}\sum_\sigma\int_{\R^3}(\rho_bv_\sigma-1)((\partial_yW_\sigma)\circ Y)\Bigl(Y_1^2+\frac{\gamma-1}2X^2\Bigr)\dd x\\
 &-\int_{\R^3}\rho_b(\partial_yH_c)(Y+h_*,t)\Bigl(Y_1^2+\frac{\gamma-1}2X^2\Bigr)\dd x\\
 &+\int_{\R^3}((R_U)_2\circ Y)\partial_1\Bigl(Y_1^2+\frac{\gamma-1}2X^2\Bigr)\dd x .
\end{aligned}
\]
By \eqref{e:Wloc} the second integral is at most $C_M\tau^{-\frac32}\Vert(X,Y)\Vert_{L^2}^2$; by
\eqref{eq:pr-gau}, $|\partial_yH_c|\le C\delta\tau^{-1}$, so the third is at most
$C_M\delta\tau^{-1}\Vert(X,Y)\Vert_{L^2}^2$; by Lemma~\ref{l:asy},
$\sqrt\tau\Vert(R_U)_2\Vert_{L^\infty}\le\omega_b(s)$, so the fourth is at most
$\omega_b(s)\tau^{-\frac12}\Vert(X,Y)\Vert_{L^2}\Vert\nabla(X,Y)\Vert_{L^2}$. In the same way
$\sum_\sigma\int(\partial_1(U_\sigma\circ Y))r_\sigma^2\dd x=\sum_\sigma\int((\partial_yW_\sigma)\circ Y)r_\sigma^2\dd x$ up to a
term at most $C_M\tau^{-\frac32}\Vert(X,Y)\Vert_{L^2}^2$. Now \eqref{e:bcr} with
$\langle g,h\rangle=\int((\partial_yW_\sigma)\circ Y)gh\dd x$ gives
\[
 \sum_\sigma\int_{\R^3}((\partial_yW_\sigma)\circ Y)\Bigl[r_\sigma^2-\frac2{\gamma+1}\Bigl(Y_1^2+\frac{\gamma-1}2X^2\Bigr)\Bigr]\dd x
 =\mathcal R_n,
\]
which proves the identity for $\mathcal R_h+S_{\mathscr D}$ with the stated bound of $E_{\rm n}$.
The algebraic identity \eqref{e:bcr} follows from
\[
 \langle Y_1,Y_1\rangle+\frac{\gamma-1}2\langle X,X\rangle
 =\frac{\gamma+1}2\bigl(\langle r_+,r_+\rangle+\langle r_-,r_-\rangle\bigr)+(\gamma-3)\langle r_+,r_-\rangle,
\]
since $X=r_++r_-$ and $Y_1=r_+-r_-$.

For $E_{\rm th}$, $[a,J]F(x)=\int K_J(x-z)(a(x)-a(z))F(z)\dd z$ and $\int|x|K_J(x)\dd x=2$ give
$\Vert[a,J]F\Vert_{L^2}\le2\Vert\nabla a\Vert_{L^\infty}\Vert F\Vert_{L^2}$. Since $(U_\sigma\circ Y)=v_\sigma w_\sigma$ and
$|\partial_1(U_\sigma\circ Y)|\le C_M\tau^{-1}$,
\[
 \Bigl\Vert (U_\sigma\circ Y)J\Bigl(\frac{\sigma\sqrt{\rho_b}d_b}2\partial_1z_s\Bigr)\Bigr\Vert_{L^2}
 \le C_M\Vert w_\sigma\nabla z_s\Vert_{L^2}+C_M\tau^{-1}\Vert\nabla z_s\Vert_{L^2},
\]
and $\partial_1(\frac{r_\sigma}{\rho_b})\rho_b=\partial_1r_\sigma-\frac{\partial_1\rho_b}{\rho_b}r_\sigma$ with
$|\partial_1\rho_b\,(U_\sigma\circ Y)|\le\omega_b(s)\tau^{-1}$; Young's inequality gives the bound of $E_{\rm th}$.

By Lemma~\ref{l:res}, $(P_s,V_s,z_s)$ solves \eqref{e:forced} with the sources
$\mathbf f_s$ of \eqref{e:full}, so Lemma~\ref{l:cr} applies. Since $r_\sigma r_{-\sigma}=r_+r_-$,
\eqref{e:cross} gives $\frac{\dd}{\dd t}C_{\rm cr}=-\mathcal R_n+J_{\rm cr}+E_{\rm cr}$. In $J_{\rm cr}$,
$\int w_\sigma c_Eg\Div'Y_\perp\dd x=-\int w_\sigma c_E\nabla'g\cdot Y_\perp\dd x$ for $g=r_{\pm\sigma}$, and
$|\nabla r_{\pm\sigma}|\le|\nabla X|+|\nabla Y|$; Young's inequality gives the bound of $J_{\rm cr}$.
In $\operatorname{Err}_1,\operatorname{Err}_2$, \eqref{e:Wloc} gives $|a_\sigma^{\rm err}|\le\omega_b(s)\tau^{-1}$;
the second-order terms of $R_X,R_Y$ are integrated by parts once onto $w_\sigma r_{\pm\sigma}$, with
$\sqrt\tau|w_\sigma|+\tau|\partial_1w_\sigma|\le C_M$; the coefficient terms are bounded with
\eqref{e:Wloc}, \eqref{e:bal}, \eqref{e:E2} and \eqref{e:E3}; and the sources $\mathbf f_s$ are bounded as
in the proof of Lemma~\ref{l:sd} and above, with \eqref{e:Gbrelation} and \eqref{e:Getb} for $G_b$,
\eqref{e:ABF} for $A_bF,B_bF$, \eqref{e:Dbd} for $\mathscr D$ and Lemma~\ref{l:res} for $F^g,F^\chi$.
Here $\mathscr T$ does not enter, since only $R_X$ and $(R_Y)_1$ occur and
$\mathscr T_p=\mathscr T_z=(\mathscr T_m)_1=0$. This gives the bound of $E_{\rm cr}$.

By Lemma~\ref{l:pr} and Lemma~\ref{l:cr}, $\Vert w_\sigma\Vert_{L^\infty}+\Vert B_\sigma^{\rm sq}\Vert_{L^\infty}\le
C_M\tau^{-\frac12}$; with $e_b^0\ge0$, $r_\pm^2+e_b^0(q_s,\eta)\le C_M(P_s^2+|v_s|^2+|\nabla z_s|^2+|b\eta|^2)$ and
$C_M\tau^{-\frac12}\Vert b\Vert_{L^\infty}^2\le C_M\tau^{-\frac32}$ by \eqref{e:bsup}, this gives the bounds at
$t$. At $t=s$, $g_c(s)=0$ gives $q_s(s)=q_o(s)$, and $\eta(s)=\chi(s)$ with
$\Vert\chi(s)\Vert_{L^2}^2+E[q_o](s)\le C_MI_s$. The last identity follows from
$\partial_rw_s=-\frac{q_aw_s}\tau$.
\end{proof}

\subsection{The contact work}\label{s:work}

Put $a_c=\frac{\alpha_G}2$ with $\alpha_G$ of Lemma~\ref{l:hes}. For a
function $g$ on $\R\times[s,t]$, $g\circ Y$ is the function
$x\mapsto g(Y(x_1,r),r)$ at the time $r$, and $A^{a_c}_g$,
$\mathscr N^{a_c}_{s,t}(g)$ are the observations of \eqref{e:nobs}.

\begin{lemma}\label{l:cbg}
With $b_a,p_a$ of \eqref{e:cmp} and $\widehat p$ of \eqref{e:Uhat},
\begin{equation}\label{e:chsplit}
 h_b\circ X=a_h+r_h,\qquad
 a_h=\frac{(\partial_yH_c)(y+h_*,t)}{V_b},\qquad
 r_h=\frac{\partial_y(b_m-b_a)+\sum_\sigma s_\sigma\partial_y\vartheta_\sigma}{V_b},
\end{equation}
\begin{equation}\label{e:thsplit}
 (\partial_1p_b)\circ X=a_p+r_p,\qquad
 a_p=\frac{(\partial_y\widehat p)(y+h_*,t)}{V_b},\qquad
 r_p=\frac{\partial_y(P_b-p_a)
   +\partial_y\bigl(p_a-\widehat p(y+h_*,t)\bigr)}{V_b},
\end{equation}
and
\[
 |a_h|\le\frac{C\delta}\tau,\qquad
 |a_p|\le\frac{C\delta^2}{\tau^{3/2}},\qquad
 \tau\bigl(A^{a_c}_{a_h}+A^{a_c}_{a_p}\bigr)\le C_M\delta^2,\qquad
 \int_s^\infty\bigl(A^{a_c}_{r_h}+A^{a_c}_{r_p}\bigr)\dd t\le\omega_b(s).
\]
For every function $g$ on $\R\times[s,t]$ and almost every $r\in[s,t]$,
at the time $r$,
\begin{equation}\label{e:cspace}
 \Vert(g\circ Y)F\Vert_{L^2}^2\le C_ME_CA^{a_c}_g,\qquad
 \Vert(g\circ Y)\nabla F\Vert_{L^2}^2\le C_MD_CA^{a_c}_g .
\end{equation}
\end{lemma}

\begin{proof}
By \eqref{e:chain}, $(\partial_1a)\circ X=\frac{\partial_y(a\circ X)}{V_b}$
for $a=a(x_1,t)$; since $b\circ X=b_m$ and $p_b\circ X=P_b$,
$h_b\circ X=\frac{\partial_yb_m}{V_b}$ and
$(\partial_1p_b)\circ X=\frac{\partial_yP_b}{V_b}$. By \eqref{e:cmp} and
$r_\sigma=(-1,s_\sigma,p_*)$ of \eqref{e:mass}, $b_a=H_c(y+h_*,t)+\sum_\sigma s_\sigma\vartheta_\sigma$,
and the value of $p_a$ at $\vartheta_+=\vartheta_-=0$ is $\widehat p(y+h_*,t)$.
Thus \eqref{e:chsplit} and \eqref{e:thsplit} hold. By \eqref{eq:pr-gau},
$|\partial_yH_c|\le C\delta\tau^{-1}$ and
$|\partial_y\widehat p|\le C\delta^2\tau^{-\frac32}$; with
$V_b\ge\frac1{2\rho_-}$ (Lemma~\ref{l:mcoef}) and
$\sqrt\tau\int_\R e^{-\frac{a_c(y+h_*)^2}{2\tau}}\dd y=\sqrt{\frac{2\pi}{a_c}}\,\tau$,
this gives the bounds for $a_h$ and $a_p$. Every term of
$p_a-\widehat p(y+h_*,t)$ contains a factor $\vartheta_\sigma$, so
Lemma~\ref{l:pr} and \eqref{eq:pr-gau} give
\[
 \bigl|\partial_y\bigl(p_a-\widehat p(y+h_*,t)\bigr)\bigr|
 +\sum_\sigma|\partial_y\vartheta_\sigma|
 \le C_M\tau^{-1}\sum_\sigma e^{-\frac{(y-s_\sigma t)^2}{8d_\sigma\tau}} .
\]
Since $s_\sigma\ne0$, completing the square in $y$ gives $A^{a_c}\le C_be^{-c_bt}$
for these terms. Lemma~\ref{l:bulk} for $\partial_y(b_m-b_a)$ and
$\partial_y(P_b-p_a)$ gives the last bound. By $\dd x_1=V_b\dd y$,
\eqref{e:ctr} and $a_c\le\alpha_G$,
\[
 \Vert(g\circ Y)F\Vert_{L^2}^2
 =\int_\R V_b|g|^2\Vert(F\circ X)(y,\cdot)\Vert_{L^2(\R^2)}^2\dd y
 \le C_ME_C\sqrt\tau\int_\R e^{-\frac{a_c(y+h_*)^2}{2\tau}}|g|^2\dd y,
\]
which is the first bound of \eqref{e:cspace}; the second line of
\eqref{e:ctr} gives the second.
\end{proof}

For $\mathsf L=l(t,D')$ with $l$ in the class of Lemma~\ref{l:sp},
$g_c=(P_g,v_g,z_g)$ of \eqref{e:gc} and $d_N,d_b$ of \eqref{e:FC}, put
\begin{equation}\label{e:WA}
 Z_{\mathsf L}=\frac\rho{\rho_b}(V_1-d_N\partial_1\zeta)
 -\mathsf L\bigl((v_g)_1-d_b\partial_1z_g\bigr),\qquad
 \mathcal W_A(s,t)=\int_s^tw_s\Bigl|\int_{\R^3}\rho_bA_bF\,Z_{\mathsf L}\dd x\Bigr|\dd r .
\end{equation}

\begin{lemma}\label{l:cw}
Let $\epsilon_A>0$. There are $c_1,c_2>0$, depending only on $M$,
$\epsilon_A$ and the fixed data, such that the following holds after
$\delta_1(M)$ is decreased and $S$ is increased, both depending also on
$\epsilon_A$. For every $\mathsf L$ as in \eqref{e:WA} and every
decomposition $P_b-p_*=g_1+g_2$ in distributions on $\R\times[s,T_e]$ with
$A^{a_c}_{g_2}\in L^1_{\rm loc}([s,T_e])$ and
$\mathscr N^{a_c}_{s,T_e}(g_1)<\infty$,
\begin{equation}\label{e:cwo}
 \mathcal W_A(s,t)\le\bigl[\epsilon_A+c_1\delta\,\mathscr N^{a_c}_{s,t}(g_1)\bigr]
 \mathfrak B(s,t)
 +c_2\delta^2\int_s^tA^{a_c}_{g_2}(r)\,\mathfrak B(s,r)\dd r .
\end{equation}
\end{lemma}

\begin{proof}
Let $s\le r\le t$. Since $\rho=\rho_b+\phi$,
$\theta=\theta_b+\zeta$ and $\gamma R=(\gamma-1)c_p$,
$\frac\rho{\rho_b}d_N=\frac c{\rho_b\theta}=d_b\frac{\theta_b}\theta$, hence
\[
 \frac\rho{\rho_b}(V_1-d_N\partial_1\zeta)
 =V_1-d_b\partial_1\zeta+\frac\phi{\rho_b}V_1+d_b\frac\zeta\theta\partial_1\zeta .
\]
By $V=\psi-b\eta e_1$, $\partial_1d_b=-d_b\frac{\partial_1p_b}{p_b}$,
\eqref{e:eta} and \eqref{e:E1Z},
$N_0(V_1-d_b\partial_1\zeta;s,r)\le C_M(N_o+\Phi^{\frac12})(s,r)$. Next,
\[
\begin{aligned}
 \nabla\Bigl(\frac\phi{\rho_b}V_1\Bigr)
 &=\frac{V_1\nabla\phi}{\rho_b}+\frac\phi{\rho_b}\nabla V_1
 -\frac{\phi\partial_1\rho_b}{\rho_b^2}V_1e_1,\\
 \nabla\Bigl(d_b\frac\zeta\theta\partial_1\zeta\Bigr)
 &=d_b\frac\zeta\theta\nabla\partial_1\zeta
 +d_b\frac{\theta_b}{\theta^2}(\partial_1\zeta)\nabla\zeta
 -d_b\Bigl(\frac{\zeta\partial_1\theta_b}{\theta^2}
 +\frac{\partial_1p_b}{p_b}\frac\zeta\theta\Bigr)(\partial_1\zeta)e_1 .
\end{aligned}
\]
By \eqref{e:state}, $\Vert\psi\Vert_{L^\infty}\le C_M\mu_*$ and \eqref{e:pbd},
\[
 \Vert\psi\nabla\phi\Vert_{L^2}\le C_M\mu_*\bigl(\Vert\nabla\Pi\Vert_{L^2}
 +\Vert\nabla^2\zeta\Vert_{L^2}
 +\Vert(|\partial_1\rho_b|+|\partial_1\theta_b|)\psi\Vert_{L^2}\bigr);
\]
on the cubes $Q_j$, $\Vert\nabla\phi\Vert_{L^3(Q_j)}\le C_M\mu_*$ and
$\Vert b\eta\Vert_{L^6(Q_j)}\le C\Vert b\eta\Vert_{H^1(Q_j)}$ give
$\Vert(b\eta)\nabla\phi\Vert_{L^2}\le C_M\mu_*\Vert b\eta\Vert_{H^1}$, and
$\Vert b\eta\Vert_{L^2_w}^2\le C_M\Phi$ by \eqref{e:cbg} and \eqref{e:N0};
\eqref{e:pbd} also gives
$\Vert(\partial_1\zeta)\nabla\zeta\Vert_{L^2}\le C\mu_*\Vert\nabla^2\zeta\Vert_{L^2}$.
With $|\phi|+|\zeta|\le C_M\mu_*$ from \eqref{e:stb}, \eqref{e:E1Z} and \eqref{e:eta}, the two
displays give
$N_0(\frac\phi{\rho_b}V_1)+N_0(d_b\frac\zeta\theta\partial_1\zeta)
\le C_M\mu_*(N_o+\Phi^{\frac12})$. Since $0\le l\le1$, $\mathsf L$ commutes
with $\nabla$ and $\Vert\mathsf LH\Vert_{L^2}\le\Vert H\Vert_{L^2}$, so
$N_0(\mathsf LH)\le N_0(H)$; moreover
$N_0((v_g)_1-d_b\partial_1z_g)\le C_MN[g_c]$, and
$N[g_c]^2\le C_M\mathcal J+\omega_b(s)N_o^2$ by Lemma~\ref{l:cor}. By
\eqref{e:B},
\begin{equation}\label{e:cZL}
 N_0(Z_{\mathsf L};s,r)^2\le C_M\mathfrak B(s,r)\qquad(s\le r\le t).
\end{equation}

The hypotheses give $g_1,g_2\in L^1_{\rm loc}(\R\times[s,T_e])$ and, for
almost every $r$, $g_1(\cdot,r),\partial_yg_1(\cdot,r),g_2(\cdot,r)\in L^2_{\rm loc}(\R)$,
hence $g_2(\cdot,r)=P_b(\cdot,r)-p_*-g_1(\cdot,r)\in H^1_{\rm loc}(\R)$. By
\eqref{e:chain}, $\partial_1(g\circ Y)=\rho_b(\partial_yg)\circ Y$. With
$\rho_bA_b=\nu\partial_1h_b-\partial_1p_b$ from \eqref{e:bal} and
$p_b=P_b\circ Y$, for almost every $r$,
\[
 \int_{\R^3}\rho_bA_bFZ_{\mathsf L}\dd x
 =-\nu\int_{\R^3}h_b\partial_1(FZ_{\mathsf L})\dd x
 +\int_{\R^3}(g_2\circ Y)\partial_1(FZ_{\mathsf L})\dd x
 -\int_{\R^3}\rho_b\bigl((\partial_yg_1)\circ Y\bigr)FZ_{\mathsf L}\dd x ;
\]
the integrations by parts hold for $FZ_{\mathsf L}\chi(\frac{x_1}R)$ with a cut-off
$\chi$, and \eqref{e:cspace} gives the limit $R\to\infty$. By
\eqref{e:cspace} and the Cauchy--Schwarz inequality, for a function $g$
on $\R\times[s,t]$,
\begin{equation}\label{e:cpb}
\begin{aligned}
 \Bigl|\int_{\R^3}(g\circ Y)\partial_1(FZ_{\mathsf L})\dd x\Bigr|
 &\le C_M\sqrt{A^{a_c}_g}\bigl(\sqrt{D_C}\Vert Z_{\mathsf L}\Vert_{L^2}
 +\sqrt{E_C}\Vert\partial_1Z_{\mathsf L}\Vert_{L^2}\bigr),\\
 \Bigl|\int_{\R^3}\rho_b\bigl((\partial_yg)\circ Y\bigr)FZ_{\mathsf L}\dd x\Bigr|
 &\le C_M\sqrt{A^{a_c}_{\partial_yg}}\sqrt{E_C}\Vert Z_{\mathsf L}\Vert_{L^2} .
\end{aligned}
\end{equation}

By \eqref{e:chsplit}, $h_b=(a_h+r_h)\circ Y$. For $a_h$, the first line of
\eqref{e:cpb}, $\tau A^{a_c}_{a_h}\le C_M\delta^2$, $E_C\le\tau D_C$ and the
Cauchy--Schwarz inequality in time give the bound
$C_M\delta(\int_s^tw_sD_C)^{\frac12}
[(\int_s^t\frac{w_s}\tau\Vert Z_{\mathsf L}\Vert_{L^2}^2)^{\frac12}
+(\int_s^tw_s\Vert\partial_1Z_{\mathsf L}\Vert_{L^2}^2)^{\frac12}]
\le C_M\delta\,\mathcal C_C^{\frac12}N_0(Z_{\mathsf L})$; for $r_h$, the bound is
$C_M(\int_s^\infty A^{a_c}_{r_h})^{\frac12}\mathcal C_C^{\frac12}N_0(Z_{\mathsf L})$.
By \eqref{e:Ch}, $N_0(\psi)\le2N_o$ and \eqref{e:cZL},
\[
 \int_s^tw_s\Bigl|\nu\int_{\R^3}h_b\partial_1(FZ_{\mathsf L})\dd x\Bigr|\dd r
 \le\bigl(C_M\delta^2+\omega_b(s)\bigr)\mathfrak B(s,t).
\]

By \eqref{e:cZL} and \eqref{e:Ch} on $[s,r]$, with
$N_0(\psi;s,r)\le2N_o(s,r)$, the running essential suprema give, for almost every
$r\in[s,t]$,
\begin{equation}\label{e:cpref}
 w_s(r)\Vert Z_{\mathsf L}(r)\Vert_{L^2}^2\le C_M\mathfrak B(s,r),\qquad
 w_s(r)E_C(r)\le C_M\delta^2\mathfrak B(s,r).
\end{equation}

For $\lambda\in(0,1]$, Young's inequality in \eqref{e:cpb} gives
\begin{equation}\label{e:cyoung}
 w_s\Bigl|\int_{\R^3}(g_2\circ Y)\partial_1(FZ_{\mathsf L})\dd x\Bigr|
 \le\lambda\Bigl(\frac{w_sD_C}{\delta^2}+w_s\Vert\partial_1Z_{\mathsf L}\Vert_{L^2}^2\Bigr)
 +\frac{C_M}\lambda A^{a_c}_{g_2}\bigl(\delta^2w_s\Vert Z_{\mathsf L}\Vert_{L^2}^2+w_sE_C\bigr).
\end{equation}
Integrating over $[s,t]$ and using \eqref{e:Ch}, \eqref{e:cZL} and
\eqref{e:cpref},
\[
 \int_s^tw_s\Bigl|\int_{\R^3}(g_2\circ Y)\partial_1(FZ_{\mathsf L})\dd x\Bigr|\dd r
 \le\lambda C_M\mathfrak B(s,t)
 +\frac{C_M}\lambda\delta^2\int_s^tA^{a_c}_{g_2}(r)\mathfrak B(s,r)\dd r .
\]
By the second line of \eqref{e:cpb} and \eqref{e:cpref},
$w_s|\int_{\R^3}\rho_b((\partial_yg_1)\circ Y)FZ_{\mathsf L}\dd x|
\le C_M\delta\sqrt{A^{a_c}_{\partial_yg_1}(r)}\,\mathfrak B(s,r)$ for almost every
$r$. Since $\mathfrak B(s,\cdot)$ is nondecreasing and
$\sqrt{A^{a_c}_{\partial_yg_1}}
=\tau^{\frac14}\Vert e^{-\frac{a_c(y+h_*)^2}{4\tau}}\partial_yg_1\Vert_{L^2(\R)}$ is
the first term of the integrand of $\mathscr N^{a_c}$ in \eqref{e:nobs},
\begin{equation}\label{e:cweak}
 \int_s^tw_s\Bigl|\int_{\R^3}\rho_b\bigl((\partial_yg_1)\circ Y\bigr)FZ_{\mathsf L}\dd x\Bigr|\dd r
 \le C_M\delta\int_s^t\sqrt{A^{a_c}_{\partial_yg_1}(r)}\,\mathfrak B(s,r)\dd r
 \le C_M\delta\,\mathscr N^{a_c}_{s,t}(g_1)\,\mathfrak B(s,t).
\end{equation}

Adding the three bounds,
\[
 \mathcal W_A(s,t)\le\bigl[C_M\delta^2+\omega_b(s)+\lambda C_M
 +C_M\delta\,\mathscr N^{a_c}_{s,t}(g_1)\bigr]\mathfrak B(s,t)
 +\frac{C_M}\lambda\delta^2\int_s^tA^{a_c}_{g_2}(r)\mathfrak B(s,r)\dd r .
\]
Choose $\lambda$, $\delta_1(M)$ and $S$ so that
$\lambda C_M+C_M\delta_1(M)^2+\sup_{s\ge S}\omega_b(s)\le\epsilon_A$.
Then \eqref{e:cwo} holds with $c_1=C_M$ and $c_2=C_M/\lambda$.
\end{proof}

\subsection{Thermal and nonlinear works}\label{s:heat}

Let $K\ge1$, let $l_A,l_S$ be real, even and measurable functions of
$(t,k)$ with $|l_A|+|l_S|\le K$, put $\mathsf A=l_A(t,D')$,
$\mathsf S=l_S(t,D')$ and $z_{\mathsf A}=\zeta-\mathsf Az_g$, and let
$B_0=B_0(x_1,t)$ be measurable with $|B_0|\le K$. For a triple
$f=(P_f,v_f,z_f)$ with $V_f$ of \eqref{e:be} put, with $V$ of \eqref{e:cut},
the coefficients of \eqref{e:FC},
\begin{equation}\label{e:Tth}
\begin{aligned}
 \mathcal T_N={}&\int_{\R^3}h_b\bigl[-\gamma j_NV\cdot\nabla\zeta
 +\{(\gamma-1)k_N+\rho d_N^2\}|\partial_1\zeta|^2
 +\gamma k_N|\nabla'\zeta|^2\bigr]\dd x\\
 &+\int_{\R^3}\frac{\gamma-1}\gamma d_N\frac{\partial_1p_b}{p_b}\Pi\partial_1\zeta\dd x,\\
 \mathcal T_b^{B_0,\mathsf S}[f]={}&\int_{\R^3}B_0h_b\Bigl[-\gamma j_bV_f\cdot\mathsf S\nabla z_f
 +\{(\gamma-1)k_b+\rho_bd_b^2\}\partial_1z_f\,\mathsf S\partial_1z_f\Bigr]\dd x\\
 &+\int_{\R^3}B_0h_b\gamma k_b\sum_{i=2}^3\partial_iz_f\,\mathsf S\partial_iz_f\dd x
 +\int_{\R^3}B_0\frac{\gamma-1}\gamma d_b\frac{\partial_1p_b}{p_b}P_f\,\mathsf S\partial_1z_f\dd x .
\end{aligned}
\end{equation}
By \eqref{e:works}, and \eqref{e:WA} with $l=0$,
\[
 \mathcal R_{\rm princ}=-\int_{\R^3}h_b\Bigl[\rho(V_1-d_N\partial_1\zeta)^2
 +\frac{\gamma-1}{2\gamma p}\Pi^2\Bigr]\dd x+\mathcal T_N
 +\int_{\R^3}\rho_bA_bF\,Z_{\mathsf L}\dd x .
\]

\begin{lemma}\label{l:tw}
Let $K,\mathsf A,\mathsf S,z_{\mathsf A},B_0$ be as in \eqref{e:Tth}.
\begin{enumerate}
\item[(i)] For $a_0\in\{h_b,\partial_1p_b\}$, $1\le i\le3$, a measurable scalar
$C$ with $|C|\le K$ and a scalar $Z$ with $N_0(Z;s,t)<\infty$,
\begin{equation}\label{e:twm}
 \int_s^tw_s\Bigl|\int_{\R^3}Ca_0Z\,\partial_iz_{\mathsf A}\dd x\Bigr|\dd r
 \le\bigl(C_{M,K}\delta^2+\omega_{b,K}(s)\bigr)N_0(Z;s,t)\,\mathfrak B(s,t)^{\frac12}.
\end{equation}
\item[(ii)] One has
\begin{equation}\label{e:twb}
 \int_s^tw_s\bigl(|\mathcal T_N|+|\mathcal T_b^{B_0,\mathsf S}[q_o-\mathsf Ag_c]|\bigr)\dd r
 \le\bigl(C_{M,K}\delta^2+\omega_{b,K}(s)\bigr)\mathfrak B(s,t).
\end{equation}
\item[(iii)] Let $T_p,T_z$ be scalar and $T_v,T_\theta$ vector functions on
$\R^3\times[s,t]$ with
\begin{equation}\label{e:srct}
\begin{gathered}
 \Vert T_v\Vert_{L^3}+\Vert T_v\Vert_{L^\infty}\le K\mu_*,\qquad
 \tau|T_v|\le K\mu_*\ \text{ on }\{|Y+h_*|<2c_0\tau\}
 \qquad\text{for a.e. }r\in[s,t],\\
 N_0(T_v;s,t)+N_0(T_p;s,t)+\Vert T_\theta\Vert_{L^2_w}
 +\Vert\nabla T_\theta\Vert_{L^2_w}+\Vert T_z\Vert_{L^2_w}
 \le K\mathfrak B(s,t)^{\frac12}.
\end{gathered}
\end{equation}
Then, with $\widetilde F^N$ and $f^{(0)}_{z,j}$ of
\eqref{e:full}, $j=o,s$, and $\mathfrak B=\mathfrak B(s,t)$,
\begin{equation}\label{e:srcrec}
\begin{gathered}
 \Vert\widetilde F^N\Vert_{L^2_w}\le C_M\mu_*N_o,\qquad
 \Vert f^{(0)}_{z,j}\Vert_{L^2_w}
 \le\bigl[C_M(\delta+\mu_*)+\omega_b(s)\bigr]\mathfrak B^{\frac12},\\
 \int_s^tw_s\Bigl|\int_{\R^3}\bigl(T_p\widetilde F^N_p+(T_v+T_\theta)\cdot\widetilde F^N_m
 +T_z\widetilde F^N_z\bigr)\dd x\Bigr|\dd r\le C_{M,K}\mu_*\mathfrak B,\\
 \int_s^tw_s\Bigl(|\langle G_b,T_v+T_\theta\rangle|
 +\Bigl|\int_{\R^3}\gamma p_b\partial_1(b\eta)T_p\dd x\Bigr|\Bigr)\dd r
 \le\bigl[C_{M,K}(\delta+\mu_*)+\omega_{b,K}(s)\bigr]\mathfrak B .
\end{gathered}
\end{equation}
If the first and second lines of \eqref{e:srct} hold with separate
constants $K_\infty$ and $K_E$, respectively, the coefficient in
the last bound of \eqref{e:srcrec} can be taken as
\[
 C_{M,K_E}\delta+C_{M,K_\infty,K_E}\mu_*
                         +\omega_{b,K_\infty,K_E}(s).
\]
\end{enumerate}
\end{lemma}

\begin{proof}
(i) Let $(a,r_0)$ be $(a_h,r_h)$ for $a_0=h_b$ and $(a_p,r_p)$ for
$a_0=\partial_1p_b$, so that $a_0=(a+r_0)\circ Y$ by \eqref{e:chsplit} and
\eqref{e:thsplit}. By \eqref{e:cspace} and Lemma~\ref{l:cbg},
$\Vert(a\circ Y)\nabla F\Vert_{L^2}\le C_M\delta\tau^{-\frac12}\sqrt{D_C}$ and
$\Vert(a\circ Y)F\Vert_{L^2}^2\le C_M\delta^2\frac{E_C}\tau\le C_M\delta^2D_C$, while
the terms with $r_0$ are bounded by $\sup_rw_s\Vert Z\Vert_{L^2}^2$, $\sup_rw_sE_C$
and $\int_s^\infty A^{a_c}_{r_0}\le\omega_b(s)$. The Cauchy--Schwarz inequality in
time and \eqref{e:N0} give
\begin{equation}\label{e:thF}
 \int_s^tw_s\int_{\R^3}|a_0Z\nabla F|\dd x\dd r
 \le\bigl(C_M\delta+\omega_b(s)\bigr)\mathcal C_C^{\frac12}N_0(Z),\qquad
 \Vert a_0F\Vert_{L^2_w}\le\bigl(C_M\delta+\omega_b(s)\bigr)\mathcal C_C^{\frac12}.
\end{equation}
The last five terms of \eqref{e:tz7} have total $L^2_w$ norm at most
$C_M(N_o+\Phi^{\frac12})$, by \eqref{e:N0}, \eqref{e:E1Z},
\eqref{e:ehp} and \eqref{e:ext}. By \eqref{e:E1Z},
$\Vert a_0Z\Vert_{L^2_w}\le\omega_b(s)N_0(Z)$ and
$\Vert(\partial_1\theta_b)Z\Vert_{L^2_w}\le C(\delta+\omega_b(s))N_0(Z)$. Pairing
$Ca_0Z$ with the first two terms and the remainder of \eqref{e:tz7}, using \eqref{e:thF} and
\eqref{e:Ch},
\begin{equation}\label{e:tw0}
 \int_s^tw_s\Bigl|\int_{\R^3}Ca_0Z\,\partial_i\zeta\dd x\Bigr|\dd r
 \le\bigl(C_{M,K}\delta^2+\omega_{b,K}(s)\bigr)N_0(Z)\bigl(N_o+\Phi^{\frac12}\bigr).
\end{equation}
Since $\mathsf A$ commutes with $\partial_i$ and
$\Vert\mathsf AH\Vert_{L^2}\le K\Vert H\Vert_{L^2}$, Lemma~\ref{l:cor} gives
$\Vert\partial_i\mathsf Az_g\Vert_{L^2_w}^2\le K^2(C_M(\mathcal J+\mathcal H)+\omega_b(s)N_o^2)$,
hence
$\int_s^tw_s|\int Ca_0Z\partial_i\mathsf Az_g|\le K\Vert a_0Z\Vert_{L^2_w}
\Vert\partial_i\mathsf Az_g\Vert_{L^2_w}\le\omega_{b,K}(s)N_0(Z)\mathfrak B^{\frac12}$.
With \eqref{e:tw0} and \eqref{e:B} this proves \eqref{e:twm}.

(ii) Let $f=q_o-\mathsf Ag_c$; then $P_f=\Pi-\mathsf AP_g$, $V_f=V-\mathsf Av_g$
and $z_f=z_{\mathsf A}$. Since $N[g_c]\le N_o+N_s$, $N_0(\mathsf AH)\le KN_0(H)$
and \eqref{e:eta},
\[
 N_0(V)+N_0(P_f)+N_0(V_f)+N_0(\nabla z_f)\le C_K\mathfrak B^{\frac12}.
\]
All coefficients of $\mathcal T_N$ in \eqref{e:Tth} are bounded by $C_M$.
Apply (i) with $\mathsf A=0$ and $Z\in\{V_i,\partial_i\zeta,\Pi\}$.
For $\mathcal T_b^{B_0,\mathsf S}[f]$, its coefficients depend only on $(x_1,t)$,
so self-adjointness moves $\mathsf S$ onto the other factor.
Apply (i) with $Z\in\{\mathsf S(V_f)_i,\mathsf S\partial_iz_f,\mathsf SP_f\}$;
here $N_0(\mathsf SZ)\le KN_0(Z)$ since $\mathsf S$ commutes with derivatives.
The displayed bound then gives \eqref{e:twb}.

(iii) By \eqref{e:stb}, \eqref{e:gn4},
$\Vert\psi\cdot\nabla\zeta\Vert_{L^2}\le\Vert\psi\Vert_{L^3}\Vert\nabla\zeta\Vert_{L^6}$
and $\Vert h_b\Vert_{L^\infty}\le\omega_b(s)\le1$ from \eqref{e:E6},
\[
 \Vert\widetilde F^N\Vert_{L^2}\le C_M\mu_*\bigl(\Vert\nabla\Pi\Vert_{L^2}
 +\Vert\nabla\psi\Vert_{H^1}+\Vert\nabla^2\zeta\Vert_{L^2}
 +\Vert(|A_b|+|B_b|+|h_b|)\Pi\Vert_{L^2}\bigr),
\]
and \eqref{e:E1Z} gives the first bound of \eqref{e:srcrec}. The other terms of
$f^{(0)}_{z,j}$ are $-(\partial_1\theta_b)(v_j)_1$,
$-\frac{B_bP_j}{\gamma p_b}-\frac{h_bP_j}{c_v\rho_b}$, $\frac{H[v_j]}{c_v\rho_b}$,
$B_b\eta-(\gamma-1)\theta_b\partial_1(b\eta)$ and $B_bf=B_bF+B_bf_{\rm out}$; by
\eqref{e:E1Z}, \eqref{e:E6}, \eqref{e:eta}, \eqref{e:ABF}, \eqref{e:Ch} and
\eqref{e:ext}, their $L^2_w$ norms are at most
\[
 C\bigl(\delta+\omega_b(s)\bigr)N[q_j],\qquad \omega_b(s)N[q_j],\qquad
 \omega_b(s)N[q_j],\qquad \omega_b(s)\Phi^{\frac12},\qquad
 \bigl(C_M\delta^2+\omega_b(s)\bigr)N_o,
\]
with $N[q_o]=N_o$ and $N[q_s]=N_s$. This gives the second bound.

In the second line of \eqref{e:srcrec}, the terms with $T_\theta$ and $T_z$
are bounded by the first line and \eqref{e:srct}. With $\widetilde F^N_p=F^N_p$
of \eqref{e:FN}, $\Vert(\Pi,\psi)\Vert_{L^3}\le\mu_*$,
$\Vert T_p\Vert_{L^6}\le C\Vert\nabla T_p\Vert_{L^2}$ and \eqref{e:pbd},
\[
 \Bigl|\int_{\R^3}T_p\widetilde F^N_p\dd x\Bigr|
 \le C\mu_*\bigl(\Vert\nabla\Pi\Vert_{L^2}+\Vert\nabla\psi\Vert_{H^1}\bigr)
 \Vert\nabla T_p\Vert_{L^2}.
\]
The hypotheses \eqref{e:srct} allow \eqref{e:inertial-work} with $T=T_v$,
which controls $\frac\phi\rho(\nabla\Pi-\mathcal L\psi)$ in the momentum row.
For its other two terms,
\[
 \Bigl|\int_{\R^3}\rho_b(\psi\cdot\nabla\psi)\cdot T_v\dd x\Bigr|
 +\Bigl|\int_{\R^3}\rho_bA_b\frac\phi\rho\frac\Pi{\gamma p_b}(T_v)_1\dd x\Bigr|
 \le C\bigl(\Vert\psi\Vert_{L^3}\Vert\nabla\psi\Vert_{L^2}
 +\Vert A_b\tfrac\phi\rho\Vert_{L^2}\Vert\Pi\Vert_{L^3}\bigr)\Vert T_v\Vert_{L^6},
\]
where
$\Vert A_b\frac\phi\rho\Vert_{L^2_w}\le C_M(N_o+\Phi^{\frac12})$ by \eqref{e:espl},
\eqref{e:E1Z}, \eqref{e:ABF}, \eqref{e:E6} and \eqref{e:ext}.
Integration, with \eqref{e:srct}, gives the second line of \eqref{e:srcrec}.

By \eqref{e:Gbrelation},
\[
 \langle G_b,T_v\rangle=\Bigl\langle G_\eta,\frac{\rho_b}\rho T_v\Bigr\rangle
 +\int_{\R^3}\frac\phi\rho\mathcal L(b\eta e_1)\cdot T_v\dd x
 +\int_{\R^3}\rho_b\bigl(h_b\psi_1\eta+b\psi\cdot\nabla\eta\bigr)(T_v)_1\dd x .
\]
Since $\nabla\frac{\rho_b}\rho=-\nabla\frac\phi\rho$, the first line of
\eqref{e:inertial-work} and \eqref{e:srct} give
$N_0(\frac{\rho_b}\rho T_v)\le C_M(N_0(T_v)+K\mu_*N_o)$, and \eqref{e:Getb} bounds
the first term. For the second, integrate by parts and use the first line of
\eqref{e:inertial-work}, with
$N_0(b\eta)^2\le\omega_b(s)\Phi$ from \eqref{e:eta}. The third
is at most
$C\Vert h_b\Vert_{L^\infty}\Vert\psi\Vert_{L^2}\Vert\eta\Vert_{L^6}\Vert T_v\Vert_{L^3}
+C\Vert b\Vert_{L^\infty}\Vert\psi\Vert_{L^3}\Vert\nabla\eta\Vert_{L^2}
\Vert T_v\Vert_{L^6}$, and \eqref{e:E3}, \eqref{e:cbg} apply. Next,
$\langle G_b,T_\theta\rangle=-\int\mathbb S(b\eta e_1):\nabla T_\theta
-\int\rho_bbD_b\eta(T_\theta)_1$; the first term is at most
$C\Vert\nabla(b\eta)\Vert_{L^2}\Vert\nabla T_\theta\Vert_{L^2}$, and
Lemma~\ref{l:mat} with $T_{\rm th}=(T_\theta)_1$ bounds the second, since
$\mathfrak B\le C_M(I_s+N_o^2+\Phi)$ by Lemma~\ref{l:cor}. Finally
$\int\gamma p_b\partial_1(b\eta)T_p=-\int\gamma b\eta\,\partial_1(p_bT_p)$, and
$\Vert b\eta\Vert_{L^2_w}^2\le(C_M\delta^2+\omega_b(s))(I_s+N_o^2)$ by \eqref{e:eta},
while \eqref{e:E1Z} bounds $(\partial_1p_b)T_p$. This proves the third line of \eqref{e:srcrec}.
Keeping the two test constants separate gives the stated refinement:
the surviving $\delta$ term uses only $N_0(T_p)$ from the second
line of \eqref{e:srct}; the first-line bounds occur with $\mu_*$
or a vanishing background tail.
\end{proof}

\subsection{Transverse frequencies}\label{s:freq}

For a bounded real even symbol $a(t,k)$, put $\mathsf A=a(t,D')$.
If $\mathcal Q$ is an integral of quadratic or bilinear terms with coefficients
depending only on $(x_1,t)$, then $\mathcal Q^{\mathsf A}$ denotes insertion of
$\mathsf A$ in one factor. The choice of factor is immaterial, since $\mathsf A$
is self-adjoint and commutes with the coefficients, $\nabla$ and $J$.
By Plancherel's identity, this is the corresponding fiber integral with weight
$a(t,k)$, $\partial_j=ik_j$ for $j=2,3$, and products read as
$\operatorname{Re}\widehat g\,\overline{\widehat h}$.
It is linear in $a$; if $a\ge0$, it equals $\mathcal Q$ with
$\mathsf A^{1/2}$ applied to both factors. For $a$ absolutely continuous in time,
the fiber identities of Lemma~\ref{l:FC}, \eqref{e:chid} and \eqref{e:cross}
acquire the term $\mathcal Q^{\partial_t\mathsf A}$ on the right.

Fix a smooth nonincreasing $\chi_{\rm sp}\colon[0,\infty)\to[0,1]$, constant near $0$,
with $\chi_{\rm sp}=1$ on $[0,1]$, $\chi_{\rm sp}=0$ on $[2,\infty)$ and
$\sqrt{\chi_{\rm sp}}$, $\sqrt{1-\chi_{\rm sp}}$ smooth, and smooth
$\zeta_1,\zeta_2\colon[0,\infty)\to[0,1]$ with $\zeta_1=0$ on $[0,\frac14]$,
$\zeta_1=1$ on $[\frac12,\infty)$, $\zeta_2=1$ on $[0,3]$ and $\zeta_2=0$ on $[4,\infty)$.
For $0<\varepsilon_L<1$, $A_{\rm sp}\ge4$ and $L_{\rm sp}>0$ put
\begin{equation}\label{e:split}
\begin{gathered}
 r_k=|k|\tau^{\frac34},\qquad
 l=\chi_{\rm sp}\Bigl(\frac{r_k}{\varepsilon_L}\Bigr),\qquad
 h_H=1-\chi_{\rm sp}\Bigl(\frac{r_k}{A_{\rm sp}}\Bigr),\qquad
 \zeta_{\rm sp}(v)=\zeta_1\Bigl(\frac v{\varepsilon_L}\Bigr)\zeta_2\Bigl(\frac v{A_{\rm sp}}\Bigr),\\
 m_{\rm sp}=\exp\Bigl\{L_{\rm sp}\int_{r_k}^\infty\frac{\zeta_{\rm sp}(v)}v\dd v\Bigr\},\qquad
 Q_H=m_{\rm sp}h_H,\qquad v_{\rm mid}=m_{\rm sp}(1-l-h_H),\\
 m_*=\Bigl(\frac{16A_{\rm sp}}{\varepsilon_L}\Bigr)^{L_{\rm sp}},
\end{gathered}
\end{equation}
and let $\mathsf L,\mathsf M,\mathsf Q_H$ be the multipliers with the symbols
$l,m_{\rm sp},Q_H$. Since $\partial_tr_k=\frac{3r_k}{4\tau}$,
\[
 -\partial_tm_{\rm sp}=\frac{3L_{\rm sp}}{4\tau}\zeta_{\rm sp}(r_k)m_{\rm sp},\qquad
 0\le-\partial_tl\le\frac{C_\chi}\tau,\qquad 0\le\partial_th_H\le\frac{C_\chi}\tau,\qquad
 1\le m_{\rm sp}\le m_*,
\]
with $C_\chi=\frac32\sup|\chi_{\rm sp}'|$; $l+h_H\le1$, $m_{\rm sp}\le4^{L_{\rm sp}}$ where
$r_k\ge A_{\rm sp}$, and $\zeta_{\rm sp}(r_k)=1$ on the supports of $\partial_tl$,
$\partial_th_H$ and $1-l-h_H$, which lie in $\{\varepsilon_L\le r_k\le2A_{\rm sp}\}$. Thus
$l$ and $m_{\rm sp}$ belong to the class of Lemma~\ref{l:sp}. Each of $l$, $m_{\rm sp}$,
$Q_H$, $\varepsilon_L^{-1}\tau^{\frac34}kl$ and their products is $g(\tau^{\frac34}k)$ with $g$ smooth
on $\R^2$ and $g-c_0$ compactly supported for a constant $c_0$, and $\tau\partial_t$ of it is
$\frac34(\xi\cdot\nabla g)(\xi)$ at $\xi=\tau^{\frac34}k$. Hence each of these
multipliers $\mathsf T$, and $\tau\partial_t\mathsf T$, is a constant multiple of $I$ plus the
convolution in $x'$ with the dilation by $\tau^{-\frac34}$ of a fixed integrable kernel, and
\begin{equation}\label{e:kern}
 \Vert\mathsf Tg\Vert_{L^p(\R^3)}\le C_*\Vert g\Vert_{L^p(\R^3)}\qquad(1\le p\le\infty)
\end{equation}
for each of them. Let $\epsilon_{\rm rec}>0$. Here and below, the parameters
$\theta_a,\varepsilon_L,L_{\rm sp},
A_{\rm sp},\lambda_F,\theta_H,\Lambda_G,L_c,\epsilon_{\rm rec}$ are fixed in the proof of
Lemma~\ref{l:cl}, depending only on $M$ and the fixed data; $C_{M,*}$, $C_{b,*}$,
$\omega_{b,*}$ may depend also on them, and $C_M$ does not.

For $L_c\ge1$ and $\sigma=\pm1$ let $\beta_\sigma(\xi)=0$ for $\xi\le4c_0$,
$\beta_\sigma(\xi)=1$ for $\xi\ge|s_\sigma|$, and
\begin{equation}\label{e:cone}
\begin{gathered}
 \beta_\sigma(\xi)=1-\chi_{\rm sp}\Bigl(1+\frac1{L_c}\log\frac{|s_\sigma|-4c_0}{|s_\sigma|-\xi}\Bigr)
 \quad(4c_0<\xi<|s_\sigma|),\\
 B_{\rm out}=\sum_{\sigma=\pm1}\beta_\sigma(\varkappa_\sigma),\qquad
 \varkappa_\sigma=\frac{\sigma(Y+h_*)}\tau .
\end{gathered}
\end{equation}
Since $\chi_{\rm sp}$ is constant on $[0,1]$ and on $[2,\infty)$ and $4c_0<\frac12|s_\sigma|$,
$\beta_\sigma$ is smooth and nondecreasing; with
$d_B=\frac14e^{-L_c}\min_\sigma(|s_\sigma|-4c_0)$, $\beta_\sigma=1$ for
$\xi\ge|s_\sigma|-4d_B$, and
\[
 0\le(|s_\sigma|-\xi)\beta_\sigma'(\xi)\le\frac{\sup|\chi_{\rm sp}'|}{L_c},\qquad
 \frac{\beta_\sigma'(\xi)^2}{\beta_\sigma(\xi)}
 =\frac{4\bigl((\sqrt{1-\chi_{\rm sp}})'\bigr)^2}{L_c^2(|s_\sigma|-\xi)^2}
 \le\frac{\sup|(\sqrt{1-\chi_{\rm sp}})'|^2}{4L_c^2d_B^2},
\]
\[
 D_b\beta_\sigma(\varkappa_\sigma)=-\frac{\varkappa_\sigma}\tau\beta_\sigma'(\varkappa_\sigma),\qquad
 \partial_1\beta_\sigma(\varkappa_\sigma)=\frac{\sigma\rho_b}\tau\beta_\sigma'(\varkappa_\sigma),
\]
by $D_bY=0$ and $\partial_1Y=\rho_b$. The two summands of $B_{\rm out}$ have disjoint
supports, so $0\le B_{\rm out}\le1$; $B_{\rm out}=0$ on $\{|Y+h_*|\le4c_0\tau\}$, which
contains $\operatorname{supp}F$ by \eqref{e:cut}, and $B_{\rm out}=1$ on
$\{|Y+h_*-s_\sigma t|\le2d_B\tau\}$ once $2d_B\tau\ge|s_\sigma|$. After increasing $S$,
$B_{\rm out}$ is a weight of the class of Lemma~\ref{l:sp}.

With $L_0$ of Lemma~\ref{l:sp}, $\theta_a\in(0,1]$ and the fiber energy
$e^0_b(\widehat f,\widehat\eta)$ of \eqref{e:be}, put
\begin{equation}\label{e:ann}
\begin{aligned}
 e_{\rm ac}&=\frac12\sum_\sigma\int_\R\Bigl[L_0\tau^{-1}|\widehat u_\sigma|^2+|D_k\widehat u_\sigma|^2
 +\sum_{j=2}^4|D_k^j\widehat u_\sigma|^2\Bigr]\dd y,\qquad
 e_{\rm mix}=l\,e_b^0(\widehat q_s,\widehat\eta)+(1-l)e_b^0(\widehat q_o,\widehat\eta),\\
 \mathcal D_{\rm ann}&=\int_{\R^2}(-\partial_tm_{\rm sp})\Bigl[\int_\R B_{\rm out}e_{\rm mix}\dd x_1
 +\theta_ae_{\rm ac}\Bigr]\dd k,\qquad
 E^{\rm mix}_{\rm out}=\int_{\R\times\R^2}B_{\rm out}(m_{\rm sp}-1)e_{\rm mix}\dd x_1\dd k,\\
 \Xi_{\rm out}&=\int_{\R\times\R^2}B_{\rm out}(m_{\rm sp}-1)\partial_tl
 \bigl\{e_b^0(\widehat q_s,\widehat\eta)-e_b^0(\widehat q_o,\widehat\eta)\bigr\}\dd x_1\dd k,\\
 \mathcal W_{\rm sw}&=\mathcal E_b^{\partial_t\mathsf L}(q_s,\eta)-\mathcal E_b^{\partial_t\mathsf L}(q_o,\eta)
 +\Xi_{\rm out},
\end{aligned}
\end{equation}
so that $E_{\rm ac}^{\mathsf M}=\int m_{\rm sp}e_{\rm ac}\dd k$ and
$E_{\rm ac}^{-\mathsf M_t}=\int(-\partial_tm_{\rm sp})e_{\rm ac}\dd k$ in \eqref{e:spec}. For a
triple $f$ put $\mathcal Q_h(\widehat f)=|(\widehat Y_f)_1|^2+\frac{\gamma-1}2|\widehat X_f|^2$
and $R_h^{\mathsf A}[f]=-\int_{\R\times\R^2}h_ba\,\mathcal Q_h(\widehat f)\dd x_1\dd k$; thus
$R_h^I[q_s]$ is $\mathcal R_h$ of Lemma~\ref{l:cd}.

\begin{lemma}\label{l:low}
Let $\Lambda_G\ge1$, $\mathsf Q_L=\mathsf M\mathsf L$,
$\rho_J=1+\Lambda_G\sum_\sigma B_\sigma^{\rm sq}\circ Y$ with $B_\sigma^{\rm sq}$ of \eqref{e:asq},
$Y_L=(v_s)_\perp-d_b\nabla'z_s$, $F_T=(\mathscr T_m)_\perp$, and
\begin{equation}\label{e:LC}
\begin{gathered}
 C_\sigma^m=\frac{V_b\sqrt{\gamma P_b}G_\sigma}{s_\sigma},\qquad
 W_L=\Bigl(\sum_\sigma C_\sigma^m\nabla'\mathsf Lu_\sigma^{\rm ac}\Bigr)\circ Y,\qquad
 C_L=\int_{\R^3}\rho_b\rho_JY_L\cdot\mathsf MW_L\dd x,\\
 \mathscr C_L=C_{\rm cr}^{\mathsf Q_L}+\Lambda_G\sum_\sigma\int_{\R\times\R^2}
 (B_\sigma^{\rm sq}\circ Y)\,m_{\rm sp}l\,e_b^0(\widehat q_s,\widehat\eta)\dd x_1\dd k-C_L,\\
 \mathcal R_L=R_h^{\mathsf Q_L}[q_s]+S_{\mathscr D}^{\mathsf Q_L}
 +\int_{\R^3}F_T\cdot\mathsf Q_LY_L\dd x
 +\epsilon_K\int_{\R^3}F_T\cdot\mathsf L\nabla'JP_s\dd x,
\end{gathered}
\end{equation}
with $C_{\rm cr}$ and $S_{\mathscr D}$ of Lemma~\ref{l:cd}. There are $C_M$, independent of the
parameters, and $\Lambda_*$, depending only on $M$, $\epsilon_{\rm rec}$ and $m_*$, such that for
$\Lambda_G\ge\Lambda_*$
\begin{equation}\label{e:lowc}
\begin{aligned}
 \int_s^tw_s\mathcal R_L\dd r+w_s(t)\mathscr C_L(t)-\mathscr C_L(s)
 \le{}&C_M\varepsilon_L^2\int_s^t\frac{w_s}\tau\Bigl[\int_{\R\times\R^2}
 \bigl(1+B_{\rm out}(m_{\rm sp}-1)\bigr)e_{\rm mix}\dd x_1\dd k+E_{\rm ac}^{\mathsf M}\Bigr]\dd r\\
 &+\bigl(\epsilon_{\rm rec}+C_{M,*}\delta+C_{b,*}\mu_*+\omega_{b,*}(s)\bigr)\mathfrak B,
\end{aligned}
\end{equation}
\[
 |\mathscr C_L(t)|\le C_{M,*}\tau^{-\frac14}\bigl(E[q_s]+E_{\rm ac}^{\mathsf M}\bigr)(t)
 +\omega_{b,*}(s)\Vert\eta(t)\Vert_{L^2}^2,\qquad
 |\mathscr C_L(s)|\le C_{M,*}\tau_s^{-\frac14}I_s .
\]
\end{lemma}

\begin{proof}
Write $u_\sigma=u_\sigma^{\rm ac}$ and $r_\sigma=r_{\sigma,q_s}$. Since $\mathsf Q_L$ is
self-adjoint and nonnegative, $\langle g,h\rangle=\int_{\R^2}g\,\mathsf Q_Lh\dd x'$ is a symmetric
bilinear form, \eqref{e:bcr} holds for it, and
$|\langle g,h\rangle|\le\langle g,g\rangle^{\frac12}\langle h,h\rangle^{\frac12}$. With
$\mathsf Q_L$ inserted, the identities of Lemma~\ref{l:cd} and \eqref{e:cross} hold with the
terms $\mathcal Q^{\partial_t\mathsf Q_L}$ added, and by \eqref{e:kern} their bounds hold with
$C_{M,*},C_{b,*},\omega_{b,*}$; the bounds of $J_{\rm cr}$ and $E_{\rm th}$ hold with their
squares replaced by the $\mathsf Q_L$ forms. Hence
\[
 R_h^{\mathsf Q_L}[q_s]+S_{\mathscr D}^{\mathsf Q_L}+\frac{\dd}{\dd t}C_{\rm cr}^{\mathsf Q_L}
 =J_{\rm cr}^{\mathsf Q_L}+E_{\rm cr}^{\mathsf Q_L}+E_{\rm th}^{\mathsf Q_L}
 +E_{\rm sel}^{\mathsf Q_L}+E_{\rm n}^{\mathsf Q_L}+C_{\rm cr}^{\partial_t\mathsf Q_L},
\]
where $\int_s^tw_s(|E_{\rm cr}^{\mathsf Q_L}|+|E_{\rm sel}^{\mathsf Q_L}|+|E_{\rm n}^{\mathsf Q_L}|)
\le(C_{M,*}\delta+C_{b,*}\mu_*+\omega_{b,*}(s))\mathfrak B$ and, since
$|\partial_t(m_{\rm sp}l)|\le\frac{C_*}\tau$, $\int_s^tw_s|C_{\rm cr}^{\partial_t\mathsf Q_L}|
\le C_*(\tau_s^{-\frac12}+\omega_b(s))\mathfrak B$. For $\varepsilon>0$,
\[
\begin{aligned}
 \int_s^tw_s\bigl(|J_{\rm cr}^{\mathsf Q_L}|+|E_{\rm th}^{\mathsf Q_L}|\bigr)\dd r
 \le{}&C_M\varepsilon m_*N_s^2+\omega_{b,*}(s)\mathfrak B\\
 &+C_\varepsilon\sum_\sigma\int_s^tw_s\Bigl[\int_{\R^3}w_\sigma^2
 \bigl(r_{-\sigma}^2+|(Y_{q_s})_\perp|^2+|\nabla z_s|^2\bigr)\dd x\Bigr]^{\mathsf Q_L}\dd r .
\end{aligned}
\]

Put $B=\rho_J-1$. By \eqref{e:asq}, $0\le B\le C_M\Lambda_G\tau^{-\frac12}$,
$|\partial_1B|\le C_M\Lambda_G\tau^{-1}$, $D_bB=\Lambda_G\sum_\sigma(-|s_\sigma|W_\sigma^2+R_{B,\sigma})\circ Y$
and $\partial_1B=\Lambda_G\sum_\sigma\sigma\rho_bw_\sigma^2$; by Lemma~\ref{l:pr},
$B\le C_*e^{-c_M\tau}$ on $\{|Y+h_*|\le3c_0\tau\}$. Since $\rho_bc_E=\sqrt{\gamma P_b/V_b}\circ Y$,
\eqref{e:Wloc} gives $|\rho_bc_E-|s_\sigma||\,w_\sigma^2\le C_M\tau^{-\frac32}$, so \eqref{e:chid}
reads
\[
\begin{aligned}
 &(D_bB)e_b(P_s,V_s,z_s)+(\partial_1B)P_s\bigl((V_s)_1-d_b\partial_1z_s\bigr)\\
 &\quad=-\Lambda_G\sum_\sigma|s_\sigma|w_\sigma^2\Bigl[2r_{-\sigma}^2
 +\frac{|(Y_{q_s})_\perp|^2+(k_b-\rho_bd_b^2)|\nabla z_s|^2}2\Bigr]
 +O\bigl(\Lambda_G\tau^{-\frac32}\bigr)e_b(P_s,V_s,z_s).
\end{aligned}
\]
The fiber form of \eqref{e:wFC} for $(P_s,V_s,z_s)$ with the sources \eqref{e:full}, the weight
$B$ and $\mathsf Q_L$ inserted gives, with
$E_{\rm sq}=\int_{\R\times\R^2}B\,m_{\rm sp}l\,e_b^0(\widehat q_s,\widehat\eta)\dd x_1\dd k$,
\[
\begin{gathered}
 \frac{\dd}{\dd t}E_{\rm sq}+\mathcal O_{\rm sq}
 =\int_{\R^3}BF_T\cdot\mathsf Q_LY_L\dd x+\mathcal Q^{\partial_t\mathsf Q_L}_{\rm sq}+E'_{\rm sq},\\
 \mathcal O_{\rm sq}=\Lambda_G\sum_\sigma|s_\sigma|\Bigl[\int_{\R^3}w_\sigma^2\Bigl(2r_{-\sigma}^2
 +\frac{|(Y_{q_s})_\perp|^2+(k_b-\rho_bd_b^2)|\nabla z_s|^2}2\Bigr)\dd x\Bigr]^{\mathsf Q_L},
\end{gathered}
\]
where $\mathcal Q^{\partial_t\mathsf Q_L}_{\rm sq}=\int B\,\partial_t(m_{\rm sp}l)e_b^0\le0$ and $E'_{\rm sq}$
collects $-\mathcal D_b^{B,\mathsf Q_L}$, $\mathscr G_b^{B,\mathsf Q_L}$, the flux
$\int(\partial_1B)[n_b-P_s((V_s)_1-d_b\partial_1z_s)]^{\mathsf Q_L}$, the $O(\Lambda_G\tau^{-3/2})$ term
and the works of all sources except $F_T$; here $\mathscr T_p=\mathscr T_z=(\mathscr T_m)_1=0$ by
Lemma~\ref{l:res} and $(V_s-d_b\nabla z_s)_\perp=Y_L$. Each term of $E'_{\rm sq}$ has the factor $B$
or $\partial_1B$. Hence $|\mathcal D_b^{B,\mathsf Q_L}|\le C_*\tau^{-\frac12}(D[q_s]+\Vert\nabla(b\eta)\Vert_{L^2}^2)$,
the other quadratic terms are bounded by
$C_*\tau^{-1}(E[q_s]+\Vert\eta\Vert_{L^2}^2)^{\frac12}(E[q_s]+D[q_s]+\Vert\nabla(b\eta)\Vert_{L^2}^2)^{\frac12}$,
using $|\partial_1(\rho_b,\theta_b,u_b,p_b)|\le C_M\tau^{-\frac12}$ by \eqref{e:E1}, \eqref{e:E2}, and the flux
with $f_{z,s}$ by $C_*\tau^{-1}E[q_s]^{\frac12}\Vert f_{z,s}\Vert_{L^2}$, with
$\Vert f_{z,s}\Vert_{L^2_w}\le C_b\mathfrak B^{\frac12}$ by \eqref{e:srcrec}, \eqref{e:Dbd}, \eqref{e:rg}
and \eqref{e:eta}. The works of $\mathscr D$ and $F^g-F^\chi$ are bounded by \eqref{e:Dbd} and
\eqref{e:rwt}. For the other sources, Lemma~\ref{l:tw}(iii) applies to the tests
\[
\begin{gathered}
 T_v=B\mathsf Q_L(\psi-v_g),\qquad T_\theta=-B\mathsf Q_L(b\eta e_1+d_b\nabla z_s),\\
 T_p=\frac{B\mathsf Q_LP_s}{\gamma p_b},\qquad
 T_z=B\mathsf Q_L\{\Div(j_bV_s)-\Div(k_b\nabla z_s)\},
\end{gathered}
\]
by Lemma~\ref{l:cor} and $B\le C_*e^{-c_M\tau}$ on $\{|Y+h_*|<2c_0\tau\}$, with $K=C_{b,*}$ in the
first line of \eqref{e:srct} and $K=C_{M,*}$ in the second; in \eqref{e:srcrec} the first constant
multiplies only $\mu_*$. The terms with $A_bF,B_bF$ are bounded by $B\le C_*e^{-c_M\tau}$ on
$\operatorname{supp}F$, those with $f_{\rm out}$ by \eqref{e:ext}. Therefore
\[
 \int_s^tw_s|E'_{\rm sq}|\dd r\le\bigl(C_{M,*}\delta+C_{b,*}\mu_*
 +C_*\tau_s^{-\frac12}+\omega_{b,*}(s)\bigr)\mathfrak B .
\]

Choose $\varepsilon$ with $C_M\varepsilon m_*=\frac{\epsilon_{\rm rec}}2$ and then $\Lambda_*$ such that
$\Lambda_*|s_\sigma|\min\{1,k_b-\rho_bd_b^2\}\ge4C_\varepsilon$; the last term of the bound of
$J_{\rm cr}^{\mathsf Q_L},E_{\rm th}^{\mathsf Q_L}$ is then at most $\frac12\int_s^tw_s\mathcal O_{\rm sq}$.
By $\int_s^tw_s\frac{\dd C}{\dd r}=w_s(t)C(t)-C(s)+q_a\int_s^t\frac{w_s}\tau C$ and
$|C_{\rm cr}^{\mathsf Q_L}|+E_{\rm sq}\le C_*\tau^{-\frac12}(E[q_s]+\Vert\eta\Vert_{L^2}^2)$
(Lemma~\ref{l:cd}), the two identities give
\begin{equation}\label{e:lowcr}
\begin{aligned}
 &\int_s^tw_s\bigl(R_h^{\mathsf Q_L}[q_s]+S_{\mathscr D}^{\mathsf Q_L}\bigr)\dd r
 +\bigl[w_s(C_{\rm cr}^{\mathsf Q_L}+E_{\rm sq})\bigr]_s^t\\
 &\qquad\le\int_s^tw_s\int_{\R^3}BF_T\cdot\mathsf Q_LY_L\dd x\dd r
 +\Bigl(\frac{\epsilon_{\rm rec}}2+C_{M,*}\delta+C_{b,*}\mu_*+C_*\tau_s^{-\frac12}
 +\omega_{b,*}(s)\Bigr)\mathfrak B .
\end{aligned}
\end{equation}

By \eqref{e:full}, Lemma~\ref{l:res}, $(G_b)_\perp=(\mathcal L(b\eta e_1))_\perp$ and
$(\mathscr D_m)_\perp=0$, with $f_z=f_{z,s}+(\gamma-1)\theta_b\partial_1(b\eta)$,
\[
 \rho_bD_b(v_s)_\perp+\nabla'P_s-(\mathcal Lv_s)_\perp=F_T+(\widetilde F^N_m)_\perp+(F^g_m)_\perp,\qquad
 D_bz_s-\frac\kappa{c_v\rho_b}\Delta z_s+(\gamma-1)\theta_b\Div v_s=f_z,
\]
and $\Vert f_z\Vert_{L^2_w}\le C_b\mathfrak B^{\frac12}$. Hence
\[
\begin{aligned}
 D_bY_L={}&\frac{-\nabla'P_s+(\mathcal Lv_s)_\perp+F_T+(\widetilde F^N_m)_\perp+(F^g_m)_\perp}{\rho_b}\\
 &-d_b\nabla'\Bigl(\frac\kappa{c_v\rho_b}\Delta z_s-(\gamma-1)\theta_b\Div v_s+f_z\Bigr)
 -(D_bd_b)\nabla'z_s .
\end{aligned}
\]
By \eqref{e:acg}, \eqref{e:LC}, $\frac{s_\sigma C_\sigma^m}{V_b}=\sqrt{\gamma P_b}G_\sigma$ and
Lemma~\ref{l:res}, $\rho_bD_bW_L=\mathsf LF_T+R_L$ with
\[
 R_L=\Bigl[\sum_\sigma\frac1{V_b}\bigl\{(\partial_tC_\sigma^m)\nabla'\mathsf Lu_\sigma
 +C_\sigma^m\bigl(d_\sigma\nabla'\mathsf L\partial_y^2u_\sigma+d_{t,\sigma}\nabla'\Delta'\mathsf Lu_\sigma
 -\nabla'\mathsf LH_\sigma+\nabla'(\partial_t\mathsf L)u_\sigma\bigr)\bigr\}\Bigr]\circ Y,
\]
the derivative of $C_\sigma^m$ taken at fixed $y$. Since
$\frac{\dd}{\dd t}\int_{\R^3}\rho_bg\dd x=\int_{\R^3}\rho_bD_bg\dd x$, $\mathsf M$ commutes with $D_b$ up
to $\partial_t\mathsf M$, and the transverse derivatives may be integrated by parts, the
substitution of $D_bY_L$ and $\rho_bD_bW_L$ in $\frac{\dd}{\dd t}C_L$ gives
\begin{equation}\label{e:lowid}
 w_sT_L=\frac{\dd}{\dd t}(w_sC_L)+w_s(P_L+A_L+E_L),\qquad C_L(s)=0,
\end{equation}
where $T_L=\int\rho_JF_T\cdot\mathsf Q_LY_L\dd x+\epsilon_K\int F_T\cdot\mathsf L\nabla'JP_s\dd x$,
\[
\begin{aligned}
 P_L={}&-\sum_\sigma\int_{\R^3}\rho_JP_s\bigl(C_\sigma^m\mathsf M\Delta'\mathsf Lu_\sigma\bigr)\circ Y\dd x,\\
 A_L={}&\sum_{\sigma,\iota=\pm1}\int_{\R^3}\rho_J\bigl(\sqrt{\gamma P_b}G_\sigma C_\iota^m
 \nabla'\partial_yu_\sigma\cdot\mathsf M\nabla'\mathsf Lu_\iota\bigr)\circ Y\dd x,\\
 E_L={}&-\int_{\R^3}\bigl[\rho_b(D_b\rho_J)Y_L\cdot\mathsf MW_L+\rho_b\rho_JY_L\cdot(\partial_t\mathsf M)W_L
 +\rho_JY_L\cdot\mathsf MR_L\bigr]\dd x\\
 &-\int_{\R^3}\rho_b\rho_Jd_b\Bigl(\frac\kappa{c_v\rho_b}\Delta z_s-(\gamma-1)\theta_b\Div v_s+f_z\Bigr)
 \Div'(\mathsf MW_L)\dd x\\
 &-\int_{\R^3}\rho_J\mathsf MW_L\cdot\bigl((\mathcal Lv_s)_\perp+(\widetilde F^N_m)_\perp+(F^g_m)_\perp\bigr)\dd x
 +\int_{\R^3}\rho_b\rho_J(D_bd_b)\nabla'z_s\cdot\mathsf MW_L\dd x\\
 &+\epsilon_K\int_{\R^3}F_T\cdot\mathsf L\nabla'JP_s\dd x+\frac{q_a}\tau C_L ;
\end{aligned}
\]
$P_L$ comes from $\nabla'P_s$, $A_L$ from $F_T$ in $D_bY_L$, and $C_L(s)=0$ since $u_\sigma(s)=0$.

Since $|k|\le2\varepsilon_L\tau^{-\frac34}$ on the support of $l$, the symbols of $\nabla'\mathsf L$,
$\Delta'\mathsf L$ and $\nabla'\partial_t\mathsf L$ are bounded by $2\varepsilon_L\tau^{-\frac34}$,
$4\varepsilon_L^2\tau^{-\frac32}$ and $C\varepsilon_L\tau^{-\frac74}$. With $|\partial_1C_\sigma^m|\le C_b\tau^{-\frac12}$,
\eqref{e:JH}, Lemma~\ref{l:cor}, $1\le\rho_J\le2$, $|\partial_1\rho_J|+|D_b\rho_J|\le C_M\Lambda_G\tau^{-1}$
(after increasing $S$) and \eqref{e:kern},
\begin{equation}\label{e:ltst}
 \sup_{[s,t]}w_s\Vert Z\Vert_{L^2}^2+\int_s^tw_s\bigl(\Vert\nabla Z\Vert_{L^2}^2
 +\tau^{-1}\Vert Z\Vert_{L^2}^2+\Vert\Div'Z\Vert_{L^2}^2\bigr)\dd r
 \le C_{b,*}\tau_s^{-\frac12}\mathcal J,\qquad Z\in\{W_L,\rho_J\mathsf MW_L\}.
\end{equation}
By \eqref{e:acg}, $\Vert U_\sigma\Vert_{L^\infty}\le C_M\tau^{-\frac12}$, the integrable kernel of $J$
and $\Vert(\Pi,\psi)\Vert_{L^3\cap L^\infty}\le C\mu_*$, the heat kernel gives
$\Vert u_\sigma\Vert_{L^3\cap L^\infty}\le\int_s^t\Vert H_\sigma\Vert_{L^3\cap L^\infty}\le C_b\mu_*\sqrt\tau$,
with an extra factor $C_be^{-c_bt}$ on $\{|y+h_*|\le2c_0\tau\}$ by the Gaussian factor of
$U_\sigma$ (Lemma~\ref{l:pr}). Hence $T_v=(0,\rho_J\mathsf MW_L)$ satisfies \eqref{e:srct} with
$K=C_{b,*}$, $\Vert T_v\Vert_{L^3\cap L^\infty}\le C_{b,*}\mu_*\varepsilon_L\tau^{-\frac14}$. On the
outgoing cone $\Omega_\sigma$, \eqref{e:E10} with $\varsigma=0$, $a=\frac1{V_b}$, $f=C_\sigma^m$ gives
$\frac1{V_b}\partial_tC_\sigma^m=\partial_yF_\sigma+H_\sigma'$ with $\sqrt\tau|F_\sigma|\le C_b$,
$\tau|H_\sigma'|\le\omega_b(s)$; off a fixed smaller cone, Lemma~\ref{l:cor} bounds $u_\sigma$ with
the factor $e^{-c_*s}$. Writing $(\partial_yF_\sigma)g=\partial_y(F_\sigma g)-F_\sigma\partial_yg$, the
normal diffusion as $\partial_y$ of a flux plus a value and the transverse diffusion as a transverse
divergence, and using \eqref{e:acsrc}, we obtain
\begin{equation}\label{e:lflux}
 R_L=\sum_{i=1}^3\partial_i\mathcal V_i+\mathcal H_L,\qquad
 \int_s^tw_s\Bigl(\sum_i\Vert\mathcal V_i\Vert_{L^2}^2+\tau\Vert\mathcal H_L\Vert_{L^2}^2\Bigr)\dd r
 \le\bigl(C_{b,*}\tau_s^{-\frac12}+\omega_{b,*}(s)\bigr)\mathfrak B .
\end{equation}
With $\Vert\nabla Y_L\Vert_{L^2_w}+\Vert\tau^{-\frac12}Y_L\Vert_{L^2_w}\le C_bN_s$, the terms of $E_L$ are
bounded as follows: those with $D_b\rho_J$, $\partial_t\mathsf M$ and $\frac{q_a}\tau C_L$ by
$C_*\tau^{-\frac74}\Vert Y_L\Vert_{L^2}\Vert u_\sigma\Vert_{L^2}$; the $R_L$ term by \eqref{e:lflux}
against $\rho_J\mathsf MY_L$; the Lam\'e term after one integration by \eqref{e:ltst}; the thermal and
$f_z$ terms by $\Vert(\Delta z_s,\Div v_s,f_z)\Vert_{L^2_w}$ and \eqref{e:ltst}; the
$(\widetilde F^N_m)_\perp$ term by \eqref{e:srcrec} with $T_v=(0,\rho_J\mathsf MW_L)$; the $(F^g_m)_\perp$
term by \eqref{e:rg}; the $D_bd_b$ term by \eqref{e:E10} with $\varsigma=0$, $a=1$, $f=d_b$, one
normal integration and \eqref{e:ltst}; and the screen, written in mass measure and integrated once
in $y$, by $C_b\tau^{-\frac32}\Vert u_\sigma\Vert_{L^2}(\Vert\nabla P_s\Vert_{L^2}+\tau^{-\frac12}\Vert P_s\Vert_{L^2})$.
Thus $\int_s^tw_s|E_L|\le(C_{b,*}\mu_*+C_{b,*}\tau_s^{-\frac14}+\omega_{b,*}(s))\mathfrak B$.

On the support of $m_{\rm sp}l$, $|k|^2\le4\varepsilon_L^2\tau^{-\frac32}$. In $P_L$ write
$P_s=B_{\rm out}P_s+(1-B_{\rm out})P_s$. By Young's inequality, in mass coordinates and Fourier
variables,
\[
\begin{aligned}
 m_{\rm sp}l|k|^2|B_{\rm out}\widehat P_s|\,|\widehat u_\sigma|
 &\le\frac{2\varepsilon_L^2m_{\rm sp}l}\tau\bigl(B_{\rm out}|\widehat P_s|^2+\tau^{-1}|\widehat u_\sigma|^2\bigr),\\
 m_{\rm sp}l|k|^2|\partial_y\widehat u_\sigma|\,|\widehat u_\iota|
 &\le\frac{2\varepsilon_L^2m_{\rm sp}l}\tau\bigl(|\partial_y\widehat u_\sigma|^2+\tau^{-1}|\widehat u_\iota|^2\bigr),
\end{aligned}
\]
and $|\widehat P_s|^2\le2\gamma p_be_b^0(\widehat q_s,\widehat\eta)$,
$m_{\rm sp}lB_{\rm out}\le l(1+B_{\rm out}(m_{\rm sp}-1))$, $le_b^0(\widehat q_s,\widehat\eta)\le e_{\rm mix}$.
Lemma~\ref{l:cor} with $\varrho=1-B_{\rm out}$ and $\ell=2d_B$ bounds the part with
$(1-B_{\rm out})P_s$ by $C_{b,*}e^{-c_*s}\mathfrak B$. Hence
\[
 \int_s^tw_s(|P_L|+|A_L|)\dd r\le C_M\varepsilon_L^2\int_s^t\frac{w_s}\tau\Bigl[\int_{\R\times\R^2}
 \bigl(1+B_{\rm out}(m_{\rm sp}-1)\bigr)e_{\rm mix}\dd x_1\dd k+E_{\rm ac}^{\mathsf M}\Bigr]\dd r
 +\omega_{b,*}(s)\mathfrak B,
\]
with $C_M$ depending only on the bounds of $C_\sigma^m$, $G_\sigma$, $p_b$ and $\rho_J\le2$. Adding
\eqref{e:lowcr}, \eqref{e:lowid} integrated over $[s,t]$ and these bounds, and increasing $S$ so that
$C_*\tau_s^{-\frac14}\le\frac{\epsilon_{\rm rec}}2$, proves \eqref{e:lowc}.

By Plancherel's identity in mass coordinates and $|k|\le2\varepsilon_L\tau^{-\frac34}$ on the support of $l$,
$|C_L|\le C_M\varepsilon_L\tau^{-\frac14}\int m_{\rm sp}l(|\widehat Y_L|^2+\tau^{-1}\sum_\sigma|\widehat u_\sigma|^2)$;
with the bounds of $C_{\rm cr}$ and of the square
potential in Lemma~\ref{l:cd}, \eqref{e:bsup}, $\eta(s)=\chi(s)$ and $q_s(s)=q_o(s)$, this gives the
bounds of $\mathscr C_L$.
\end{proof}

\begin{lemma}\label{l:mix}
Put, with $\mathscr S_{B_0}^{\mathsf T}$ of \eqref{e:rwt},
\[
 R_h^{\rm out}=-\int_{\R\times\R^2}h_bB_{\rm out}(m_{\rm sp}-1)
 \bigl[l\mathcal Q_h(\widehat q_s)+(1-l)\mathcal Q_h(\widehat q_o)\bigr]\dd x_1\dd k,\qquad
 \mathscr S_{\rm out}=\mathscr S_{B_{\rm out}}^{(\mathsf M-I)\mathsf L}(q_s;\mathscr D+\mathscr T).
\]
There is $L_*$, depending only on $M$, $\epsilon_{\rm rec}$ and $m_*$, such that for $L_c\ge L_*$
\begin{equation}\label{e:mix}
\begin{aligned}
 &w_s(t)E_{\rm out}^{\rm mix}(t)+\int_s^tw_s\Bigl[\int_{\R\times\R^2}B_{\rm out}(-\partial_tm_{\rm sp})
 e_{\rm mix}\dd x_1\dd k+\frac{q_a}\tau E_{\rm out}^{\rm mix}\Bigr]\dd r\\
 &\quad\le E_{\rm out}^{\rm mix}(s)+\int_s^tw_s\bigl(R_h^{\rm out}+\mathscr S_{\rm out}+\Xi_{\rm out}\bigr)\dd r
 +\bigl(\epsilon_{\rm rec}+C_{M,*}\delta+C_{b,*}\mu_*+\omega_{b,*}(s)\bigr)\mathfrak B,
\end{aligned}
\end{equation}
and $0\le E_{\rm out}^{\rm mix}(s)\le C_{M,*}I_s$. There is $C_M$, independent of the parameters, such
that for $L_{\rm sp}\ge C_MC_\chi\max\{1,\theta_a^{-1}\}$
\begin{equation}\label{e:msw}
\begin{gathered}
 \int_s^tw_s|\mathcal W_{\rm sw}|\dd r\le\frac18\int_s^tw_s\mathcal D_{\rm ann}\dd r+\omega_{b,*}(s)\mathfrak B,\\
 \Bigl|\int_s^tw_sR_h^{v_{\rm mid}}[q_o]\dd r\Bigr|\le\frac18\int_s^tw_s\mathcal D_{\rm ann}\dd r
 +\bigl(C_{M,*}\delta+\omega_{b,*}(s)\bigr)\mathfrak B .
\end{gathered}
\end{equation}
\end{lemma}

\begin{proof}
Apply the fiber form of \eqref{e:wFC}, $c=b$, with the weight $B_{\rm out}$,
to $(P_s,V_s,z_s)$ with $(\mathsf M-I)\mathsf L$ inserted and to $(\Pi,V,\zeta)$
with $(\mathsf M-I)(I-\mathsf L)$ inserted, with the sources
$\mathbf f_s,\mathbf f_o$ of \eqref{e:full}, and add; both multipliers are nonnegative and at most $m_*$. By
$\partial_t((m_{\rm sp}-1)l)=(\partial_tm_{\rm sp})l+(m_{\rm sp}-1)\partial_tl$ and
$\partial_t((m_{\rm sp}-1)(1-l))=(\partial_tm_{\rm sp})(1-l)-(m_{\rm sp}-1)\partial_tl$, the terms with the time
derivatives of the symbols add up to $-\int B_{\rm out}(-\partial_tm_{\rm sp})e_{\rm mix}+\Xi_{\rm out}$.
In \eqref{e:wFC} the term $B_{\rm out}\{\Div(j_bV_j)-\Div(k_b\nabla z_j)\}f_{z,j}$ of $\mathscr S_b^{B_{\rm out}}$
and the part $(\partial_1B_{\rm out})(j_b(V_j)_1-k_b\partial_1z_j)f_{z,j}$ of the flux add up to
$\{\Div(B_{\rm out}j_bV_j)-\Div(B_{\rm out}k_b\nabla z_j)\}f_{z,j}$; thus the source terms are
$\mathscr S_{B_{\rm out}}^{(\mathsf M-I)\mathsf L}(q_s;\mathbf f_s)+\mathscr S_{B_{\rm out}}^{(\mathsf M-I)(I-\mathsf L)}(q_o;\mathbf f_o)$.

On the support of $\beta_\sigma'(\varkappa_\sigma)$, $\varkappa_\sigma>4c_0$, which lies in $\Omega_\sigma$; there
$\rho_bc_E=\sqrt{\gamma P_b/V_b}\circ Y$ and \eqref{e:E10} give $|\rho_bc_E-|s_\sigma||\le C_b\tau^{-\frac12}$.
For $B=\beta_\sigma(\varkappa_\sigma)$, \eqref{e:chid} and the derivatives in \eqref{e:cone} give
\[
\begin{aligned}
 &(D_bB)e_b(P,v,z)+(\partial_1B)P(v_1-d_b\partial_1z)\\
 &\quad=\frac{\beta_\sigma'}\tau\Bigl[(\rho_bc_E-\varkappa_\sigma)r_\sigma^2-(\rho_bc_E+\varkappa_\sigma)r_{-\sigma}^2
 -\frac{\varkappa_\sigma}2\bigl(|Y_\perp|^2+(k_b-\rho_bd_b^2)|\nabla z|^2\bigr)\Bigr],
\end{aligned}
\]
and $\beta_\sigma'(\rho_bc_E-\varkappa_\sigma)\le(|s_\sigma|-\varkappa_\sigma)\beta_\sigma'+C_{b,*}\tau^{-\frac12}\le
\frac{\sup|\chi_{\rm sp}'|}{L_c}+C_{b,*}\tau^{-\frac12}$. With the two multipliers inserted, the time
integral of these terms is at most
$\frac{Cm_*}{q_aL_c}\mathfrak B+C_{b,*}\tau_s^{-\frac12}\mathfrak B$, since
$\int_s^tw_s\frac{q_a}\tau(E[q_o]+E[q_s]+\Vert b\eta\Vert_{L^2}^2)\le C\mathfrak B$ by \eqref{e:bsup}.
Choose $L_*$ with $\frac{Cm_*}{q_aL_*}\le\frac{\epsilon_{\rm rec}}2$.

The curl--div identity and the square formula of Lemma~\ref{l:FC} give, on each fiber,
\[
\begin{aligned}
 \mathcal D_b^{B_{\rm out}}(v,z)={}&\int_{\R^3}B_{\rm out}\Bigl[\mu|\operatorname{curl}v|^2
 +\frac\nu\gamma(\Div v)^2\\
 &\qquad+\frac{\kappa k_b}{c_v\rho_b}
 \Bigl(\Delta z-\frac{(\gamma-1)c_v\rho_b\theta_b}{\kappa}\Div v\Bigr)^2\Bigr]\dd x\\
 &-\mu\int_{\R^3}(\partial_1B_{\rm out})\sum_{l=1}^3(v_l\partial_lv_1-v_1\partial_lv_l)\dd x .
\end{aligned}
\]
Hence, with the nonnegative multipliers inserted, the first integral is nonnegative; the second, the flux in
\eqref{e:wFC} except its source part, and all terms with $\partial_1B_{\rm out}$ below are bounded
by $|\partial_1B_{\rm out}|\le C_*\tau^{-1}$ and give $C_{b,*}\tau_s^{-\frac12}\mathfrak B$.

The terms with coefficient derivatives, those of $\mathscr G_b^{B_{\rm out}}$ and those of
\eqref{e:full} linear in $(P_j,v_j)$, are grouped as in the proof of \eqref{e:cdif}, with the factor
$B_{\rm out}$: by \eqref{e:bal}, $\rho_bA_b+\partial_1p_b=\nu\partial_1h_b$ and
$\frac{B_b}{\theta_b}=-(\gamma-1)h_b+\frac{\kappa\partial_1^2\theta_b+\nu h_b^2}{c_v\rho_b\theta_b}$; the terms
with $\partial_1h_b$ and $\partial_1^2\theta_b$ are integrated once by parts,
\[
 \int_{\R^3}B_{\rm out}(\partial_1h_b)Q\dd x=-\int_{\R^3}h_b\partial_1(B_{\rm out}Q)\dd x,\qquad
 \int_{\R^3}B_{\rm out}\frac{\partial_1^2\theta_b}{\rho_b\theta_b}Q\dd x
 =-\int_{\R^3}\partial_1\theta_b\,\partial_1\Bigl(\frac{B_{\rm out}Q}{\rho_b\theta_b}\Bigr)\dd x
\]
for the quadratic forms $Q$ that occur. They give $R_h^{\rm out}$,
$\mathcal T_b^{B_{\rm out},(\mathsf M-I)\mathsf L}[q_s]+\mathcal T_b^{B_{\rm out},(\mathsf M-I)(I-\mathsf L)}[q_o]$ of
\eqref{e:Tth}, the terms $-\int B_{\rm out}\frac{\partial_1p_b}{\gamma p_b}P_jb\eta$ with the two multipliers
inserted, and terms with one first
coefficient derivative times products of $f_j$ and $\nabla f_j$. By \eqref{e:twb} with $\mathsf A=I,0$,
\eqref{e:E1Z}, \eqref{e:E2}, \eqref{e:bsup} and \eqref{e:twm}, all but $R_h^{\rm out}$ cost
$(C_{M,*}\delta+\omega_{b,*}(s))\mathfrak B$.

The other sources of \eqref{e:full} are common to $j=s,o$, except $\mathscr D+\mathscr T+F^g-F^\chi$.
Since $\mathsf Lq_s+(I-\mathsf L)q_o=q_o-\mathsf Lg_c$, the common sources are tested by
\[
\begin{gathered}
 T_v=B_{\rm out}(\mathsf M-I)(\psi-\mathsf Lv_g),\qquad
 T_\theta=-B_{\rm out}(\mathsf M-I)\bigl(b\eta e_1+d_b\nabla(\zeta-\mathsf Lz_g)\bigr),\\
 T_p=\frac{B_{\rm out}}{\gamma p_b}(\mathsf M-I)(\Pi-\mathsf LP_g),\\
 T_z=\Div\bigl[B_{\rm out}j_b(\mathsf M-I)(V-\mathsf Lv_g)\bigr]
 -\Div\bigl[B_{\rm out}k_b(\mathsf M-I)\nabla(\zeta-\mathsf Lz_g)\bigr].
\end{gathered}
\]
By \eqref{e:kern}, Lemma~\ref{l:cor}, \eqref{e:eta} and $B_{\rm out}=0$ on $\{|Y+h_*|\le4c_0\tau\}$, they
satisfy \eqref{e:srct} with $K=C_{b,*}$ in its first line and $K=C_{M,*}$ in its second, and
Lemma~\ref{l:tw}(iii) bounds their works by $(C_{M,*}\delta+C_{b,*}\mu_*+\omega_{b,*}(s))\mathfrak B$; the first
constant in \eqref{e:srct} multiplies only $\mu_*$ in \eqref{e:srcrec}. The terms with $A_bF$, $B_bF$
vanish, since $B_{\rm out}F=0$; those with $f_{\rm out}$ are bounded by \eqref{e:ext} and \eqref{e:E3}.
The part of $f_{z,j}^{(0)}$ linear in $(P_j,v_j)$ is paired with the divergence test of $q_j$ by the
first line of \eqref{e:srcrec}. The work of $F^g-F^\chi$ is $\omega_{b,*}(s)\mathfrak B$ by \eqref{e:rwt}
with $B_0=B_{\rm out}$ and $\mathsf T=(\mathsf M-I)\mathsf L$, and that of $\mathscr D+\mathscr T$ is $\mathscr S_{\rm out}$.
Adding these bounds, multiplying by $w_s$ and integrating proves \eqref{e:mix}; $q_s(s)=q_o(s)$ and
$m_{\rm sp}\le m_*$ give the bound of $E_{\rm out}^{\rm mix}(s)$.

Write $e_j=e_b^0(\widehat q_j,\widehat\eta)$ and $e_g=e_b^0(\widehat g_c,0)$.
The Hilbert-space inequality at the end of Lemma~\ref{l:sp}, applied to the squares in
\eqref{e:be}, gives
\[
 |e_s-e_o|\le e_s+e_o\le3e_{\rm mix}+2e_g,
\]
since $q_o-q_s=g_c$. Moreover, $\int_\R e_g\dd x_1\le C_Me_{\rm ac}$ by \eqref{e:gc}. The zero-order part of $\mathcal W_{\rm sw}$ is
$\int\partial_tl\,(m_{\rm sp}B_{\rm out}+1-B_{\rm out})(e_s-e_o)\dd x_1\dd k$. On the support of $\partial_tl$,
$|\partial_tl|\le\frac{C_\chi}\tau=\frac{4C_\chi}{3L_{\rm sp}}\frac{-\partial_tm_{\rm sp}}{m_{\rm sp}}$, so its part with
$m_{\rm sp}B_{\rm out}$ is at most
$\frac{4CC_\chi}{3L_{\rm sp}}\int(-\partial_tm_{\rm sp})(B_{\rm out}e_{\rm mix}+C_Me_{\rm ac})\le\frac18\mathcal D_{\rm ann}$
for $L_{\rm sp}\ge C_MC_\chi\max\{1,\theta_a^{-1}\}$; its part with $1-B_{\rm out}$ contains $Z_{g_c}$ and is
$\omega_{b,*}(s)\mathfrak B$ by Lemma~\ref{l:cor} with $\varrho=1-B_{\rm out}$, $\ell=2d_B$. The
first-order blocks, the screens and the weak temperature value of $\mathcal W_{\rm sw}$ have the factor
$|\partial_tl|\le\frac{C_\chi}\tau$ and at least one factor of $g_c$; by
$V_s\cdot\nabla JP_s-V\cdot\nabla J\Pi=-v_g\cdot\nabla J\Pi-V\cdot\nabla JP_g+v_g\cdot\nabla JP_g$,
$z_s^2-\zeta^2=-2\zeta z_g+z_g^2$, the Cauchy--Schwarz inequality in time and Lemma~\ref{l:cor}, they
cost $C_*\tau_s^{-\frac12}\mathfrak B$. This proves the first bound of \eqref{e:msw}.

By Lemma~\ref{l:asy}, $b\circ X=H_c(y+h_*,t)+\sum_\sigma s_\sigma\vartheta_\sigma+(R_U)_2$, and
$s_\sigma\vartheta_\sigma=\frac{2v_\sigma W_\sigma}{\gamma+1}$ by \eqref{e:acg}; since
$\partial_1(G\circ Y)=\rho_b(\partial_yG)\circ Y$,
\begin{equation}\label{e:hdec}
 h_b=\rho_b\bigl(\partial_yH_c\bigr)(Y+h_*,t)+\frac{2\rho_b}{\gamma+1}\sum_\sigma v_\sigma(\partial_yW_\sigma)\circ Y
 +\partial_1\bigl((R_U)_2\circ Y\bigr),
\end{equation}
with $|\partial_yH_c|\le C\delta\tau^{-1}$ by \eqref{eq:pr-gau}, $|\partial_yW_\sigma|\le C_M\tau^{-1}$ and
$\sup_{r\ge s}\sqrt\tau\Vert(R_U)_2\Vert_{L^\infty}\le\omega_b(s)$. In $R_h^{v_{\rm mid}}[q_o]$, the first term
of \eqref{e:hdec} costs $C_{M,*}\delta\mathfrak B$; the last, after one integration by parts in $x_1$ and the
Cauchy--Schwarz inequality, costs $\omega_{b,*}(s)\mathfrak B$; the second is exponentially small where
$B_{\rm out}<1$ by Lemma~\ref{l:pr}, and where $B_{\rm out}=1$ it is bounded by
$\frac{C_M}\tau v_{\rm mid}\mathcal Q_h(\widehat q_o)\le\frac{C_M}{L_{\rm sp}}(-\partial_tm_{\rm sp})B_{\rm out}e_{\rm mix}$,
since $\zeta_{\rm sp}=1$ on the support of $1-l-h_H$ and $(1-l)e_b^0(\widehat q_o,\widehat\eta)\le e_{\rm mix}$.
This proves the second bound of \eqref{e:msw}.
\end{proof}

\begin{lemma}\label{l:tf}
Let $a_H=(16\max_\sigma d_\sigma)^{-1}$, $\xi_\sigma=\frac{Y-s_\sigma t}{\sqrt\tau}$, $\lambda_F,\theta_H>0$,
$\mathcal A_\perp=\Div'(-\Delta')^{-1}$ on functions $\R^3\to\R^2$, and
\[
\begin{gathered}
 \varphi_\sigma^G=e^{-a_H\xi_\sigma^2},\qquad \psi_\sigma^G=(1+\xi_\sigma^2)\varphi_\sigma^G,\qquad
 B_+^G=\int_{-\infty}^{\xi_+}e^{-a_Hv^2}\dd v,\qquad B_-^G=\int_{\xi_-}^\infty e^{-a_Hv^2}\dd v,\\
 C_F=\lambda_F\sum_\sigma\int_{\R\times\R^2}B_\sigma^GQ_He_b^0(\widehat q_o,\widehat\eta)\dd x_1\dd k,\qquad
 C_I=-\frac{\theta_H}\tau\sum_\sigma\int_{\R^3}\psi_\sigma^Gr_{\sigma,q_o}\mathsf Q_H\mathcal A_\perp(Y_{q_o})_\perp\dd x .
\end{gathered}
\]
There are $\lambda_*>0$, depending only on $M$, $\epsilon_{\rm rec}$ and $L_{\rm sp}$, and $\theta_*$, depending
only on $M$, such that for $\lambda_F\le\lambda_*$ and $\theta_H\ge\theta_*$ there is $A_*$, depending only
on $M,\epsilon_{\rm rec},L_{\rm sp},\lambda_F,\theta_H$, such that for $A_{\rm sp}\ge A_*$, with $C_H=C_F+C_I$,
\begin{equation}\label{e:high}
 \int_s^tw_sR_h^{Q_H}[q_o]\dd r+w_s(t)C_H(t)-C_H(s)
 \le\bigl(\epsilon_{\rm rec}+C_{M,*}\delta+C_{b,*}\mu_*+\omega_{b,*}(s)\bigr)\mathfrak B,
\end{equation}
$C_F\ge0$, $|C_I|\le C_*\theta_HA_{\rm sp}^{-1}\tau^{-\frac14}(E[q_o]+\Vert b\Vert_{L^\infty}^2\Vert\eta\Vert_{L^2}^2)$ and
$|C_H(s)|\le C_*I_s$.
\end{lemma}

\begin{proof}
Write $X,Y,r_\sigma$ for $X_{q_o},Y_{q_o},r_{\sigma,q_o}$. The symbol of $\mathcal A_\perp\mathsf Q_H$ is
$\frac{ik}{|k|^2}Q_H=\tau^{\frac34}A_{\rm sp}^{-1}g(\frac{\tau^{3/4}k}{A_{\rm sp}})$ with
$g(\xi)=\frac{i\xi}{|\xi|^2}Q_H(A_{\rm sp}\xi\tau^{-3/4})$, a smooth function independent of $A_{\rm sp}$ and
$\tau$, zero near $0$ and equal to $\frac{i\xi}{|\xi|^2}$ for $|\xi|\ge4$; its inverse transform is
$O(|x'|^{-1})$ near $0$ and, since $|\partial^\alpha g|\le C|\xi|^{-1-|\alpha|}$, $O(|x'|^{-3})$ at infinity.
Hence $\Vert\mathcal A_\perp\mathsf Q_Hf\Vert_{L^p}\le C_{L_{\rm sp}}A_{\rm sp}^{-1}\tau^{\frac34}\Vert f\Vert_{L^p}$,
$1\le p\le\infty$, and the same holds for $\mathcal A_\perp\partial_t\mathsf Q_H$ with $C_{L_{\rm sp}}\tau^{-1}$
in addition. By Lemma~\ref{l:pr}, $|\partial_y^jW_\sigma|\le C_M\tau^{-\frac{j+1}2}e^{-2a_H\xi_\sigma^2}$ and the
analogous bounds for $\vartheta_\sigma$.

Apply the fiber form of \eqref{e:wFC}, $c=b$, with the weight $\lambda_FB_\sigma^G$ and $\mathsf Q_H$
inserted, to $(\Pi,V,\zeta)$ with the sources $\mathbf f_o$. Since
$\partial_1B_\sigma^G=\sigma\rho_b\tau^{-\frac12}\varphi_\sigma^G$ and
$D_bB_\sigma^G=-|s_\sigma|\tau^{-\frac12}\varphi_\sigma^G-\frac{\sigma\xi_\sigma}{2\tau}\varphi_\sigma^G$,
\eqref{e:chid} and $r_\sigma^2-r_{-\sigma}^2=\sigma XY_1$ give
\[
\begin{aligned}
 &(D_bB_\sigma^G)e_b+(\partial_1B_\sigma^G)\Pi(V_1-d_b\partial_1\zeta)
 =-|s_\sigma|\tau^{-\frac12}\varphi_\sigma^G\Bigl[2r_{-\sigma}^2
 +\frac{|Y_\perp|^2+(k_b-\rho_bd_b^2)|\nabla\zeta|^2}2\Bigr]\\
 &\qquad+\sigma(\rho_bc_E-|s_\sigma|)\tau^{-\frac12}\varphi_\sigma^GXY_1
 -\frac{\sigma\xi_\sigma}{2\tau}\varphi_\sigma^Ge_b,
\end{aligned}
\]
with $e_b=e_b(\Pi,V,\zeta)$, and by \eqref{e:E10} the last two terms are at most
$C_M\tau^{-1}\psi_\sigma^G(r_\sigma^2+r_{-\sigma}^2+e_b)$ plus terms exponentially small in $\tau$. The
other terms of \eqref{e:wFC} are treated as in the proof of Lemma~\ref{l:mix}, with $B_{\rm out}$ replaced
by $\lambda_FB_\sigma^G$: $0\le B_\sigma^G\le\sqrt{\pi/a_H}$, $|\partial_1B_\sigma^G|\le C\tau^{-\frac12}$, and
$B_\sigma^G$ and $\partial_1B_\sigma^G$ are exponentially small on $\{|Y+h_*|\le3c_0\tau\}$; the diffusion is
nonnegative up to the curl flux, and the fluxes cost $C_*\lambda_F\mathfrak B$; the coefficient terms give
$-\lambda_F\int B_\sigma^Gh_bQ_H\mathcal Q_h(\widehat q_o)$, bounded by \eqref{e:hdec} as $R_h^{Q_H}[q_o]$ below,
and $(C_{M,*}\delta+\omega_{b,*}(s))\mathfrak B$; the sources are received by Lemma~\ref{l:tw}(iii), \eqref{e:Dbd},
\eqref{e:rwt}, \eqref{e:ext} with the tests of the proof of Lemma~\ref{l:mix}, $B_{\rm out}(\mathsf M-I)$
replaced by $\lambda_FB_\sigma^G\mathsf Q_H$ and $\mathsf Lg_c$ by $0$. Of the time derivative of the symbol,
$\partial_tQ_H=(\partial_tm_{\rm sp})h_H+m_{\rm sp}\partial_th_H$, the first part is nonpositive; the second is
nonnegative, supported in $\{A_{\rm sp}\le r_k\le2A_{\rm sp}\}$ where $m_{\rm sp}\le4^{L_{\rm sp}}$, and at most
$\frac{C_\chi}\tau4^{L_{\rm sp}}$. Hence, with
$\mathcal O_F=\lambda_F\sum_\sigma|s_\sigma|\int\tau^{-\frac12}\varphi_\sigma^GQ_H[2|\widehat r_{-\sigma}|^2
+\frac12(|\widehat Y_\perp|^2+(k_b-\rho_bd_b^2)|D_k\widehat\zeta|^2)]\dd x_1\dd k$,
\begin{equation}\label{e:front}
\begin{aligned}
 &w_s(t)C_F(t)-C_F(s)+\int_s^tw_s\mathcal O_F\dd r
 \le C_M\lambda_F\sum_\sigma\int_s^tw_s\tau^{-1}\int\psi_\sigma^GQ_H
 e_b^0(\widehat q_o,\widehat\eta)\dd x_1\dd k\dd r\\
 &\qquad+\Bigl(\frac{C_{L_{\rm sp}}\lambda_F}{q_a}+C_{M,*}\delta+C_{b,*}\mu_*+\omega_{b,*}(s)\Bigr)\mathfrak B .
\end{aligned}
\end{equation}

By \eqref{e:acsys} for $(\Pi,V,\zeta)$ with the sources $\mathbf f_o$, $\mathcal A_\perp\nabla'=-I$ and
$\partial_1c_E$ independent of $x'$, we have
$D_b(\mathcal A_\perp Y_\perp)=c_EX+\mathcal A_\perp(R_Y)_\perp$, and
$D_br_\sigma=-\sigma c_E\partial_1r_\sigma-\frac{c_E}2\Div'Y_\perp
+\frac12(R_X+\sigma(R_Y)_1)$.
Since $\frac{\dd}{\dd t}\int g=\int(D_bg+h_bg)$ and $X=r_\sigma+r_{-\sigma}$,
\[
\begin{aligned}
 &\frac{\dd}{\dd t}C_I+\frac{\theta_H}\tau\sum_\sigma\bigl[\psi_\sigma^Gc_Er_\sigma^2\bigr]^{\mathsf Q_H}
 =-\frac{\theta_H}\tau\sum_\sigma\Bigl[\Bigl(D_b\psi_\sigma^G+h_b\psi_\sigma^G-\frac{\psi_\sigma^G}\tau\Bigr)r_\sigma\mathcal A_\perp Y_\perp\\
 &\quad+\psi_\sigma^G\Bigl(c_Er_\sigma r_{-\sigma}-\sigma c_E(\partial_1r_\sigma)\mathcal A_\perp Y_\perp
 -\frac{c_E}2(\Div'Y_\perp)\mathcal A_\perp Y_\perp\\
 &\qquad\qquad+\frac12(R_X+\sigma(R_Y)_1)\mathcal A_\perp Y_\perp
 +r_\sigma\mathcal A_\perp(R_Y)_\perp\Bigr)\Bigr]^{\mathsf Q_H}
 +C_I^{\partial_t\mathsf Q_H},
\end{aligned}
\]
where $[\cdot]^{\mathsf Q_H}$ denotes the integral with $\mathsf Q_H$ inserted. On the support of $Q_H$,
$|\widehat{\mathcal A_\perp Y_\perp}|\le|k|^{-1}|\widehat Y_\perp|\le A_{\rm sp}^{-1}\tau^{\frac34}|\widehat Y_\perp|$, and
$(\psi_\sigma^G)^2\le C\varphi_\sigma^G$. For $\epsilon>0$:
\[
 \frac{\theta_H}\tau\bigl|\bigl[\psi_\sigma^Gc_E(\partial_1r_\sigma)\mathcal A_\perp Y_\perp\bigr]^{\mathsf Q_H}\bigr|
 \le\epsilon\bigl[|\partial_1r_\sigma|^2\bigr]^{\mathsf Q_H}
 +\frac{C_\epsilon\theta_H^2}{A_{\rm sp}^2}\tau^{-\frac12}\bigl[\varphi_\sigma^G|Y_\perp|^2\bigr]^{\mathsf Q_H};
\]
$|D_b\psi_\sigma^G|\le C\tau^{-\frac12}(1+\xi_\sigma^4)\varphi_\sigma^G$, so the first term on the right has
the same bound with $\epsilon\tau^{-1}[r_\sigma^2]^{\mathsf Q_H}$ in place of the first square;
$\frac{|k\cdot\widehat Y_\perp|^2}{|k|^2}\le|\widehat Y_\perp|^2$ bounds the $\Div'Y_\perp$ term by
$\frac{C\theta_H}\tau[\psi_\sigma^G|Y_\perp|^2]^{\mathsf Q_H}$; and
$\frac{\theta_H}\tau|\psi_\sigma^Gc_Er_\sigma r_{-\sigma}|\le\frac{\theta_H}{4\tau}\psi_\sigma^Gc_Er_\sigma^2
+\frac{C\theta_H}\tau\psi_\sigma^Gr_{-\sigma}^2$. In $R_X,R_Y$ each term with a first or material
coefficient derivative $c_1$ has, by \eqref{e:floc} and the Cauchy--Schwarz inequality in time,
\[
 \int_s^t\frac{w_s\theta_H}{\tau^{1/4}A_{\rm sp}}\int|\psi_\sigma^Gc_1fg|\dd x\dd r
 \le\epsilon\int_s^t\frac{w_s}\tau\Vert f\Vert_{L^2}^2\dd r
 +\frac{C_\epsilon\theta_H^2}{A_{\rm sp}^2}\sup_{[s,t]}w_s\Vert g\Vert_{L^2}^2
 \int_s^\infty\sqrt\tau\Vert\psi_\sigma^Gc_1\Vert_{L^\infty}^2\dd r,
\]
the last integral at most $\omega_b(s)$; this includes $h_b\psi_\sigma^G$. The viscous and heat terms of
$R_Y$ are integrated once onto the test and cost $C_*\theta_HA_{\rm sp}^{-1}\tau_s^{-\frac14}\mathfrak B$; in
$-\sqrt{\rho_b}d_b\nabla f_{z,o}$ the normal derivative is integrated once and
$\mathcal A_\perp\nabla'f_{z,o}=-f_{z,o}$, so only $\Vert f_{z,o}\Vert_{L^2_w}\le C_b\mathfrak B^{\frac12}$ enters. The
other sources are paired with the tests $\frac{\theta_H}\tau\psi_\sigma^G(\mathcal A_\perp\mathsf Q_H)^{(*)}$
applied to $(\Pi,\psi)$ and $Y_\perp$. To verify their temperature parts, note that
$F$ is supported in $\{|Y+h_*|\le3c_0\tau\}$, whereas $|s_\sigma|>8c_0$.
After increasing $S$, the acoustic Gaussians are exponentially small on $\operatorname{supp}F$.
Thus \eqref{e:tz7}, \eqref{e:eta}, \eqref{e:E1Z}, \eqref{e:ehp},
\eqref{e:ext} and \eqref{e:Ch} give
\[
 \int_s^t w_s\bigl\|e^{-a_H\xi_\sigma^2/4}\nabla\zeta\bigr\|_{L^2}^2\dd r
 \le C_M\mathfrak B.
\]
The transverse operators commute with this Gaussian. Their multiplier bound and
the $\nabla^2\zeta$ term in $D_*$ therefore give, for the representative test
$T=\tau^{-1}\psi_\sigma^G(\mathcal A_\perp\mathsf Q_H)^*\partial_1\zeta$,
\[
 \Vert T\Vert_{L^2_w}^2+\Vert\nabla T\Vert_{L^2_w}^2
 \le C_*A_{\rm sp}^{-2}\tau_s^{-1/2}\mathfrak B.
\]
In the mechanical tests, place the terms containing $b\eta$ or $\zeta$ in $T_\theta$;
the remaining part $T_v$ contains only $(\Pi,\psi)$. The preceding bound,
$\mathcal A_\perp\nabla'=-I$ and \eqref{e:eta} control $T_\theta$.
The multiplier kernel bound gives the $L^3$ and $L^\infty$ bounds for $T_v$;
the acoustic Gaussians give its central pointwise bound.
Together with the same bounds for $T_p,T_z$, this verifies \eqref{e:srct}
with $K=C_*$ in the first line and $K=C_{M,*}$ in the second.
Lemma~\ref{l:tw}(iii) applies. Finally
$|C_I^{\partial_t\mathsf Q_H}|\le C_*\theta_HA_{\rm sp}^{-1}\tau^{-\frac54}(E[q_o]+\Vert b\Vert_{L^\infty}^2\Vert\eta\Vert_{L^2}^2)$.

By \eqref{e:hdec}, $|\partial_yW_\sigma|\le C_M\tau^{-1}\varphi_\sigma^G$ and $Q_H\le4^{L_{\rm sp}}$,
\[
 \int_s^tw_sR_h^{Q_H}[q_o]\dd r\le C_M\sum_\sigma\int_s^tw_s\tau^{-1}\bigl[\varphi_\sigma^G(r_\sigma^2+r_{-\sigma}^2)\bigr]^{\mathsf Q_H}\dd r
 +\bigl(C_{M,*}\delta+\omega_{b,*}(s)\bigr)\mathfrak B .
\]
On $\{|\xi_\sigma|\le\tau^{\frac18}\}$, $\tau^{-1}\psi_\sigma^G\le2\tau^{-\frac14}\tau^{-\frac12}\varphi_\sigma^G$; on its
complement $\psi_\sigma^G\le Ce^{-a_H\tau^{1/4}/2}$. Hence every term with $\tau^{-1}\psi_\sigma^G$ times
$r_{-\sigma}^2$, $|Y_\perp|^2$ or $|\nabla\zeta|^2$ is at most $\frac14\mathcal O_F$ after increasing $S$. Choose
$\epsilon$ and $\lambda_*$ with $C_M4^{L_{\rm sp}}\epsilon+\frac{C_{L_{\rm sp}}\lambda_*}{q_a}\le\frac{\epsilon_{\rm rec}}2$,
then $\theta_*$ with $\theta_*\min c_E\ge8C_M$, so that $\frac{\theta_H}{4\tau}[\psi_\sigma^Gc_Er_\sigma^2]^{\mathsf Q_H}$
dominates the terms with $\tau^{-1}\psi_\sigma^Gr_\sigma^2$, and then $A_*$ with
$\frac{C_\epsilon\theta_H^2}{A_*^2}\le\frac{\lambda_F}8\min_\sigma|s_\sigma|$, so that the terms with
$\frac{\theta_H^2}{A_{\rm sp}^2}\tau^{-\frac12}\varphi_\sigma^G|Y_\perp|^2$ are at most $\frac14\mathcal O_F$. Adding
\eqref{e:front} and the identity of $C_I$, multiplied by $w_s$ and integrated, with
$q_a\int_s^t\frac{w_s}\tau C_F\ge0$ and $q_a\int_s^t\frac{w_s}\tau|C_I|\le C_*\theta_HA_{\rm sp}^{-1}\tau_s^{-\frac14}\mathfrak B$,
proves \eqref{e:high}. The endpoint bounds follow from $|\widehat{\mathcal A_\perp Y_\perp}|\le A_{\rm sp}^{-1}\tau^{\frac34}|\widehat Y_\perp|$,
$Q_H\le4^{L_{\rm sp}}$, $q_o(s)$ and $\eta(s)=\chi(s)$.
\end{proof}

\subsection{Closing the weighted estimate}\label{s:close}

\begin{lemma}\label{l:cl}
Let $M\ge1$. There are $C_1,C_2>0$ and $\delta_W(M)>0$, chosen before the datum, depending
only on $M$, the fixed physical data and the fixed test parameters, and independent of $u_*$,
such that each background with $0<\delta\le\delta_W$ admits $S_W\ge1$ and $\mu_b>0$ with the
following property. Let $S\ge S_W$, let $(\rho,u,\theta)$ be as in the first paragraph of
Section~\ref{sec:en} on $\R^3\times[S,T_e]$ with
$\mu_*\le\mu_b$, and let
\begin{equation}\label{e:pc}
 P_b-p_*=g_1+g_2\ \text{ in distributions on }\R\times[S,T_e],\qquad
 A^{a_c}_{g_2}\in L^1_{\rm loc}([S,T_e]),\qquad
 C_1\delta\,\mathscr N^{a_c}_{S,T_e}(g_1)\le\frac12,
\end{equation}
with $a_c$ of Lemma~\ref{l:cw}. Put
\begin{equation}\label{e:bp}
 b_P(\ell)=C_2\delta^2A^{a_c}_{g_2}(\ell).
\end{equation}
Then, for $S\le t\le T_e$,
\begin{equation}\label{e:cl}
 \mathfrak B(S,t)\le C_{0,b}I_S+\int_S^tb_P(\ell)\,\mathfrak B(S,\ell)\dd\ell,
\end{equation}
where $C_{0,b}$ does not depend on $T_e$, the solution, the decomposition or its observation
sizes. If moreover, with $C_{P,b},\eta_P\ge0$,
\begin{equation}\label{e:bpbd}
 \int_S^tb_P(\ell)\dd\ell\le C_{P,b}+\eta_P\log\frac{\tau_t}{\tau_S},\qquad S\le t\le T_e,
\end{equation}
then, for $S\le t\le T_e$,
\begin{equation}\label{e:grow}
 N_o(S,t)^2\le\mathfrak B(S,t)-I_S,\qquad
 \mathfrak B(S,t)\le C_{0,b}e^{C_{P,b}}I_S\Bigl(\frac{\tau_t}{\tau_S}\Bigr)^{\eta_P}.
\end{equation}
\end{lemma}

\begin{proof}
Let $\mathsf L$ be the multiplier of \eqref{e:split}, and $\mathcal W_A$ that of \eqref{e:WA} for
this $\mathsf L$. There is $C_J$, depending only on $M$ and the fixed data, such that each
background with $0<\delta\le\delta_1(M)$ admits $S_J$, $\mu_J$ and $C_{I,b}$ with
\begin{equation}\label{e:raw}
 \mathfrak B(s,t)\le C_{I,b}I_s+C_J\mathcal W_A(s,t)
\end{equation}
for every $s\ge S_J$, every solution as in the statement on $[s,T_e]$ with $\mu_*\le\mu_J$, and
$s\le t<T_e$.

Assume \eqref{e:raw}. Let $\epsilon_A=\frac1{4C_J}$, let $c_1,c_2$ be the constants of
Lemma~\ref{l:cw} for this $\epsilon_A$, put $C_1=2C_Jc_1$, $C_2=2C_Jc_2$ and
$C_{0,b}=2C_{I,b}$, and let $\delta_W\le\delta_1(M)$, $S_W\ge S_J$ and $\mu_b\le\mu_J$ be such that
Lemma~\ref{l:cw} applies for $\delta\le\delta_W$ and $S\ge S_W$. Let
$S\le t<T_e$. By \eqref{e:nobs}, $\mathscr N^{a_c}_{S,t}(g_1)\le\mathscr N^{a_c}_{S,T_e}(g_1)$, so
\eqref{e:pc} gives $C_Jc_1\delta\mathscr N^{a_c}_{S,t}(g_1)\le\frac14$, and \eqref{e:raw} with $s=S$
and \eqref{e:cwo} give
\[
 \mathfrak B(S,t)\le C_{I,b}I_S+\frac12\mathfrak B(S,t)
 +C_Jc_2\delta^2\int_S^tA^{a_c}_{g_2}(\ell)\,\mathfrak B(S,\ell)\dd\ell .
\]
Absorbing $\frac12\mathfrak B(S,t)$ gives \eqref{e:cl}.
Continuity of the energies and monotone convergence of the time integrals
extend it to $t=T_e$. Gronwall's inequality then gives
\[
 \mathfrak B(S,t)\le C_{0,b}I_S\exp\Bigl\{\int_S^tb_P(\ell)\dd\ell\Bigr\}.
\]
Together with \eqref{e:B} and \eqref{e:bpbd}, this proves \eqref{e:grow}.

We prove \eqref{e:raw}. Let $s\le t'\le t<T_e$ and $\mathfrak B=\mathfrak B(s,t)$. By \eqref{e:mri},
\eqref{e:mdec} and \eqref{e:etaE},
\begin{equation}\label{e:mfl}
\begin{aligned}
 &\frac{\dd}{\dd t}\Bigl[\mathcal E_N-\frac{\theta_\eta}2\int_{\R^3}\rho\eta^2\dd x\Bigr]
 +D_{\mathcal E}-\theta_\eta c\Vert\nabla\eta\Vert_{L^2}^2\\
 &\quad=\mathcal R_A+\mathcal R_\eta+\mathcal R_{\rm out}+\mathcal T_{\rm in}
 +\epsilon_K\sum_{i=0}^3\mathscr P_N(P^i,v^i,z^i)
 -\frac{(\gamma-1)\theta_W}\tau\int_{\R^3}\rho\theta\zeta\Div\psi\dd x .
\end{aligned}
\end{equation}
By the first paragraph of Section~\ref{s:freq}, the proof of \eqref{e:cdif} with $\mathsf L$
inserted gives
\[
\begin{aligned}
 &\frac{\dd}{\dd t}\bigl[\mathcal E_b^{\mathsf L}(q_s,\eta)-\mathcal E_b^{\mathsf L}(q_o,\eta)\bigr]
 +D_b^{\mathsf L}[q_s]-D_b^{\mathsf L}[q_o]
 =\mathcal R_{0,b}^{\mathsf L}[q_s]-\mathcal R_{0,b}^{\mathsf L}[q_o]+\mathcal R_b^{\mathsf L}\\
 &\qquad+\mathscr W_b(\mathsf Lq_s,\mathsf L\eta;\mathscr D+\mathscr T)
 +\mathcal E_b^{\partial_t\mathsf L}(q_s,\eta)-\mathcal E_b^{\partial_t\mathsf L}(q_o,\eta)
 +\mathcal O_b^{\mathsf L},
\end{aligned}
\]
with $\int_s^tw_s|\mathcal O_b^{\mathsf L}|\dd r\le(C_M\delta+C_b\mu_*+\omega_b(s))\mathfrak B$. Indeed, in that
proof Lemma~\ref{l:sd} applies with $g_L=\mathsf Lg_c$ and $K$ the $L^1$ norm of the inverse
Fourier transform of $\chi_{\rm sp}(|\xi|)$ on $\R^2$, which is the $L^1$ norm of the kernel of
$\mathsf L$, and the other bounds there are bounds of pairings of sources with tests, which hold
for the tests $\mathsf LT$, since $\mathsf L$ commutes with $\nabla$, $J$ and functions of
$(x_1,t)$ and $\Vert\mathsf LT\Vert_{L^p}\le K\Vert T\Vert_{L^p}$.

Put
\[
\begin{aligned}
 \mathscr E_J={}&\mathcal E_N-\frac{\theta_\eta}2\int_{\R^3}\rho\eta^2\dd x
 +\mathcal E_b^{\mathsf L}(q_s,\eta)-\mathcal E_b^{\mathsf L}(q_o,\eta)
 +E_{\rm out}^{\rm mix}+\theta_aE_{\rm ac}^{\mathsf M}+\mathscr C_L+C_H,\\
 \mathscr E_+={}&\mathcal E_b^{\mathsf L}(q_s,\eta)+\mathcal E_b^{I-\mathsf L}(q_o,\eta)
 +E_{\rm out}^{\rm mix}+\theta_aE_{\rm ac}^{\mathsf M},\\
 \mathscr P_J={}&\sup_{s\le r\le t}w_s\mathscr E_+
 +\int_s^tw_s\Bigl[D_b^{\mathsf L}[q_s]+D_b^{I-\mathsf L}[q_o]+\theta_ac_aD_{\rm ac}^{\mathsf M}
 +\mathcal D_{\rm ann}+\frac{q_a}\tau\mathscr E_+\Bigr]\dd r,
\end{aligned}
\]
with $\mathscr C_L$ of Lemma~\ref{l:low}, $C_H$ of Lemma~\ref{l:tf} and $c_a$ of \eqref{e:spin}. By
Lemma~\ref{l:FC} for $c=b$ on each fiber and the choices of $\epsilon_K,\theta_W$ in the proof of
Lemma~\ref{l:mr}, whose constant $c_*$ serves also for $(\rho_b,\theta_b)$ by \eqref{e:stb} and \eqref{e:E0}, for $0\le a(t,k)\le1$ and $\mathsf A=a(t,D')$ the integral of
$a\,e_b^0(\widehat f,\widehat\eta)$ over $\R\times\R^2$ is at most $2\mathcal E_b^{\mathsf A}(f,\eta)$, and
the Fourier integrals of $E_c[f,\eta]$, $D_c[f,\eta]$ of Lemma~\ref{l:sp} with the weight $a$ are at
most $C_M\mathcal E_b^{\mathsf A}(f,\eta)$, $C_MD_b^{\mathsf A}[f]$. Hence, by $m_{\rm sp}B_{\rm out}\le
1+B_{\rm out}(m_{\rm sp}-1)$,
\begin{equation}\label{e:Pml}
 \mathcal P_{m,l}\le\int_{\R\times\R^2}\bigl(1+B_{\rm out}(m_{\rm sp}-1)\bigr)e_{\rm mix}\dd x_1\dd k
 \le2\mathscr E_+ ,
\end{equation}
and \eqref{e:sprec} with $m_{\rm sp}=1$, $m_*=1$, whose budget $\mathcal P_{\theta_a}$ is then at most
$C_M\mathscr P_J$, since the acoustic terms of \eqref{e:spec} increase with $m_{\rm sp}$, together with $\Phi\le C_M(I_s+N_o^2)$ and
$N_0(b\eta)^2\le\omega_b(s)\Phi$ of \eqref{e:eta}, gives
\begin{equation}\label{e:coer}
 \mathfrak B\le C_R\bigl(I_s+\mathscr P_J\bigr)+\omega_b(s)\mathfrak B,
\end{equation}
with $C_R$ depending only on $M$ and $\theta_a$.

Multiply \eqref{e:mfl} and the identity above by $w_s$, integrate over $[s,t']$, and add
\eqref{e:mix}, $\theta_a$ times \eqref{e:spin}, \eqref{e:lowc} and \eqref{e:high}, all on $[s,t']$.
The terms with time derivatives of the multipliers add up to $\mathcal W_{\rm sw}$ of \eqref{e:ann} on
the right and $\mathcal D_{\rm ann}$ on the left. By the comparisons of Lemma~\ref{l:cd},
\[
\begin{gathered}
 \mathcal E_N-\frac{\theta_\eta}2\int_{\R^3}\rho\eta^2\dd x
 +\mathcal E_b^{\mathsf L}(q_s,\eta)-\mathcal E_b^{\mathsf L}(q_o,\eta)
 \ge\mathcal E_b^{\mathsf L}(q_s,\eta)+\mathcal E_b^{I-\mathsf L}(q_o,\eta)
 -C_b\mu_*\bigl[E[q_o]+\Vert b\Vert_{L^\infty}^2\Vert\eta\Vert_{L^2}^2\bigr],\\
 D_{\mathcal E}-\theta_\eta c\Vert\nabla\eta\Vert_{L^2}^2+D_b^{\mathsf L}[q_s]-D_b^{\mathsf L}[q_o]
 \ge D_b^{\mathsf L}[q_s]+D_b^{I-\mathsf L}[q_o]
 -C_b\mu_*\bigl[D_*+\Vert\nabla(b\eta)\Vert_{L^2}^2\bigr],
\end{gathered}
\]
by linearity in the symbol. The work terms on the right are grouped as follows.

With $\mathcal C_o$ of Lemma~\ref{l:cd}, $Q_h^b[f]=-(Y_f)_1^2-\frac{\gamma-1}2X_f^2+H_b[f]$,
$(V_s,z_s)=(V-v_g,\zeta-z_g)$, \eqref{e:Tth}, \eqref{e:WA},
$1+B_{\rm out}(m_{\rm sp}-1)=m_{\rm sp}-(1-B_{\rm out})(m_{\rm sp}-1)$ and
$m_{\rm sp}(1-l)=v_{\rm mid}+Q_H$, the normal works of \eqref{e:works}, of $\mathcal R_b^{\mathsf L}$ and
of \eqref{e:mix} add up to
\[
\begin{aligned}
 \mathcal R_{\rm princ}+\mathcal R_b^{\mathsf L}+R_h^{\rm out}
 ={}&R_h^{\mathsf Q_L}[q_s]+R_h^{v_{\rm mid}}[q_o]+R_h^{Q_H}[q_o]
 +\mathcal T_b^{1,\mathsf L}[q_s]+\mathcal T_b^{1,I-\mathsf L}[q_o]\\
 &+\int_{\R^3}\rho_bA_bF\,Z_{\mathsf L}\dd x+\mathcal C_o\\
 &+\int_{\R\times\R^2}h_b(1-B_{\rm out})(m_{\rm sp}-1)
 \bigl[l\mathcal Q_h(\widehat q_s)+(1-l)\mathcal Q_h(\widehat q_o)\bigr]\dd x_1\dd k .
\end{aligned}
\]
By \eqref{e:twb} with $(B_0,\mathsf S,\mathsf A)=(1,\mathsf L,I)$ and $(1,I-\mathsf L,0)$, the
$\mathcal T_b$ terms cost $(C_M\delta^2+\omega_b(s))\mathfrak B$, and $\mathcal C_o$ costs
$(C_b\mu_*+\omega_b(s))\mathfrak B$ by Lemma~\ref{l:cd}; the contact term gives $\mathcal W_A(s,t)$. In
the last integral, the part of \eqref{e:hdec} with $\partial_yH_c$ costs $C_{M,*}\delta\mathfrak B$; the part
with $\partial_yW_\sigma$ is at most $C_{b,*}e^{-c_*s}\mathfrak B$, since $|Y+h_*-s_\sigma t|>2d_B\tau$ on
the support of $1-B_{\rm out}$ and Lemma~\ref{l:pr} applies; the part with $(R_U)_2$, after one
integration by parts in $x_1$ with $|\partial_1B_{\rm out}|\le C_*\tau^{-1}$, costs
$\omega_{b,*}(s)\mathfrak B$.

Since $\mathscr T_p=\mathscr T_z=(\mathscr T_m)_1=0$, $(T_0[q_s])_\perp=Y_L+\epsilon_K\nabla'JP_s$ with
$Y_L$ of Lemma~\ref{l:low}, and
$\mathsf L+B_{\rm out}(\mathsf M-I)\mathsf L=\mathsf Q_L-(1-B_{\rm out})(\mathsf M-I)\mathsf L$, the
zero-order mechanical part of $\mathscr W_b(\mathsf Lq_s,\mathsf L\eta;\mathscr D+\mathscr T)$ and the mechanical
rows of $\mathscr S_{\rm out}$ add up to
\[
\begin{aligned}
 &S_{\mathscr D}^{\mathsf Q_L}+\int_{\R^3}F_T\cdot\mathsf Q_LY_L\dd x
 +\epsilon_K\int_{\R^3}F_T\cdot\mathsf L\nabla'JP_s\dd x\\
 &\quad-\int_{\R^3}(1-B_{\rm out})\Bigl[(V_s-d_b\nabla z_s)\cdot(\mathsf M-I)\mathsf L
 (\mathscr D_m+\mathscr T_m)+\frac{P_s(\mathsf M-I)\mathsf L\mathscr D_p}{\gamma p_b}\Bigr]\dd x .
\end{aligned}
\]
The first line is $\mathcal R_L-R_h^{\mathsf Q_L}[q_s]$ of \eqref{e:LC}. On the support of
$1-B_{\rm out}$ the coefficients of $\mathscr D$ in the proof of Lemma~\ref{l:cd}, which contain
$U_\sigma\circ Y$, are at most $C_be^{-c_*\tau}$ by Lemma~\ref{l:pr}, and Lemma~\ref{l:cor} with
$\varrho=1-B_{\rm out}\circ X$, $\ell=2d_B$ bounds $(1-B_{\rm out})\mathscr T_m$; with
$m_{\rm sp}\le m_*$ the second line costs $C_{b,*}e^{-c_*s}\mathfrak B$. The other terms of
$\mathscr W_b(\mathsf Lq_s,\mathsf L\eta;\mathscr D+\mathscr T)$ are those of \eqref{e:Dslots} and of
$\mathscr W_b(q_s,\eta;\mathscr T)-\int(\mathscr T_m)_\perp\cdot(T_0[q_s])_\perp\dd x$ with $\mathsf L$
inserted, which cost $\omega_b(s)\mathfrak B$ as in Lemma~\ref{l:cd}; the heat row of
$\mathscr S_{\rm out}$ costs $C_{b,*}\varepsilon_{D,b}(s)^{\frac12}\mathfrak B$ by \eqref{e:Dbd}, since
$\Vert\Div(B_{\rm out}j_bV_s)-\Div(B_{\rm out}k_b\nabla z_s)\Vert_{L^2_w}\le C_{b,*}\mathfrak B^{\frac12}$.

Since the triples $(P^i,v^i,z^i)$ of Lemma~\ref{l:mr} are those of $q_o$ in $\mathcal R_{0,b}[q_o]$,
$v_{q_o}=\psi$ and $\rho\theta-\rho_b\theta_b=\frac\Pi R$, linearity in the symbol gives
\[
\begin{aligned}
 &\epsilon_K\sum_{i=0}^3\mathscr P_N(P^i,v^i,z^i)
 -\frac{(\gamma-1)\theta_W}\tau\int_{\R^3}\rho\theta\zeta\Div\psi\dd x
 +\mathcal R_{0,b}^{\mathsf L}[q_s]-\mathcal R_{0,b}^{\mathsf L}[q_o]\\
 &\quad=\epsilon_K\sum_{i=0}^3(\mathscr P_N-\mathscr P_b)(P^i,v^i,z^i)
 -\frac{(\gamma-1)\theta_W}{R\tau}\int_{\R^3}\Pi\zeta\Div\psi\dd x
 +\mathcal R_{0,b}^{\mathsf L}[q_s]+\mathcal R_{0,b}^{I-\mathsf L}[q_o],
\end{aligned}
\]
and by \eqref{e:fck} and $D_u-D_b=\psi\cdot\nabla$, for a triple $(P,v,z)$,
\[
\begin{aligned}
 (\mathscr P_N-\mathscr P_b)(P,v,z)={}&\int_{\R^3}\phi v\cdot\nabla J\bigl[(\gamma-1)\kappa\Delta z
 -\gamma p\Div v\bigr]\dd x-\gamma\int_{\R^3}\rho_bv\cdot\nabla J(\Pi\Div v)\dd x\\
 &+\int_{\R^3}\phi v\cdot[D_u,\nabla J]P\dd x+\int_{\R^3}\rho_bv\cdot[\psi\cdot\nabla,\nabla J]P\dd x .
\end{aligned}
\]
For $i\ge1$ each term is at most $C_b\mu_*D_*$, by $\int gv\cdot\nabla Jh=-\int\Div(gv)Jh$,
\eqref{e:stb}, \eqref{e:prod} and the commutator bounds in the proof of \eqref{e:mrec}. For $i=0$ the
first two terms equal $-\int\Div(\phi V)J[(\gamma-1)\kappa\Delta\zeta-\gamma p\Div V]$ and
$\gamma\int\Div(\rho_bV)J(\Pi\Div V)$, where $\Div(\phi V)=\phi\Div V+(\nabla\rho)\cdot V-(\partial_1\rho_b)V_1$
is bounded in $L^2_w$ by \eqref{e:stb}, \eqref{e:etst} and \eqref{e:E1Z}, and the commutator terms are
bounded as in the proof of \eqref{e:mrec}, with $\rho$ replaced by $\phi$ and by $\rho_b$. With
$|\Pi|\le C_b\mu_*$ for the weak value, the first two terms on the right cost
$(C_M\delta+C_b\mu_*+\omega_b(s))\mathfrak B$. Let $(\mathsf A,f)$ be $(\mathsf L,q_s)$ or
$(I-\mathsf L,q_o)$. Then $\mathcal R_{0,b}^{\mathsf A}[f]$ and $D_b^{\mathsf A}[f]$ are $\mathcal R_{0,b}$ and
$D_b$ for the triple $\mathsf A^{\frac12}f$ with $\eta$ replaced by $\mathsf A^{\frac12}\eta$, and
$N[\mathsf A^{\frac12}f]\le N[f]$, $N_0(\mathsf A^{\frac12}(b\eta))\le N_0(b\eta)$. The proof of
\eqref{e:mrec} for the carrier $b$, where $(\nabla\rho_b)\cdot v=(\partial_1\rho_b)v_1$ is bounded by
\eqref{e:E1Z}, the commutator is $[u_b\partial_1,\nabla J]$ and there is no $\theta_\eta$ term, and with
$\Div v_f=\Div V_f+\partial_1(b\eta)$ in the weak value, gives
\[
 \int_s^tw_s\bigl|\mathcal R_{0,b}^{\mathsf A}[f]\bigr|\dd r
 \le\frac12\int_s^tw_sD_b^{\mathsf A}[f]\dd r+\bigl(C_M\delta^2+\omega_b(s)\bigr)\mathfrak B .
\]

The remaining terms are bounded by \eqref{e:mr}, \eqref{e:eret}, \eqref{e:Tib}, the bound of
$\mathcal O_b^{\mathsf L}$, and \eqref{e:spin}, \eqref{e:Pml} with $\theta_aC_M\le\frac{q_a}{32}$:
\[
 \theta_a\int_s^{t'}w_s\Bigl(\frac{C_M}\tau\mathcal P_{m,l}+r_s\Bigr)\dd r
 \le\frac{q_a}{16}\int_s^{t'}\frac{w_s}\tau\mathscr E_+\dd r+\omega_{b,*}(s)\mathfrak B ;
\]
$\mathcal W_{\rm sw}$ and $R_h^{v_{\rm mid}}[q_o]$ by \eqref{e:msw}; $\mathcal R_L$ by \eqref{e:lowc}, whose
first term is at most $\frac{q_a}4\int_s^{t'}\frac{w_s}\tau\mathscr E_+\dd r$ by \eqref{e:Pml},
$E_{\rm ac}^{\mathsf M}\le\theta_a^{-1}\mathscr E_+$ and $12C_M\varepsilon_L^2\le q_a\theta_a$; and
$R_h^{Q_H}[q_o]$ by \eqref{e:high}. Collecting, with $C_F\ge0$,
\begin{equation}\label{e:Jab}
\begin{aligned}
 &w_s(t')\mathscr E_+(t')+\frac12\int_s^{t'}w_s\Bigl[D_b^{\mathsf L}[q_s]+D_b^{I-\mathsf L}[q_o]
 +\theta_ac_aD_{\rm ac}^{\mathsf M}+\mathcal D_{\rm ann}+\frac{q_a}\tau\mathscr E_+\Bigr]\dd r\\
 &\quad\le\mathscr E_J(s)+w_s(t')\bigl(|\mathscr C_L|+|C_I|\bigr)(t')+\mathcal W_A(s,t)
 +\bigl(3\epsilon_{\rm rec}+C_{M,*}(\delta+\delta^2)+C_{b,*}\mu_*+\omega_{b,*}(s)\bigr)\mathfrak B .
\end{aligned}
\end{equation}

Since $q_s(s)=q_o(s)$, $E_{\rm ac}^{\mathsf M}(s)=0$ and $\mathcal E_N(s)\le C_MI_s$, the endpoint bounds of
Lemmas~\ref{l:low}, \ref{l:mix} and~\ref{l:tf} give $\mathscr E_J(s)\le C_{M,*}I_s$, and with
$w_sE_{\rm ac}^{\mathsf M}\le C_*(\mathcal J+\mathcal H)$ they give
$\sup_{[s,t]}w_s(|\mathscr C_L|+|C_I|)\le(C_{M,*}\tau_s^{-\frac14}+\omega_{b,*}(s))\mathfrak B$. Taking the
supremum over $t'$ in \eqref{e:Jab},
\[
 \mathscr P_J\le3C_{M,*}I_s+3\mathcal W_A(s,t)
 +3\bigl(3\epsilon_{\rm rec}+C_{M,*}(\delta+\delta^2)+C_{b,*}\mu_*+\omega_{b,*}(s)\bigr)\mathfrak B ,
\]
and \eqref{e:coer} gives \eqref{e:raw} with $C_J=6C_R$ and $C_{I,b}=2C_R(1+3C_{M,*})$, provided
\[
 9C_R\epsilon_{\rm rec}\le\frac18,\qquad 3C_RC_{M,*}(\delta+\delta^2)\le\frac18,\qquad
 3C_R\bigl(C_{b,*}\mu_*+\omega_{b,*}(s)\bigr)+\omega_b(s)\le\frac14 .
\]
The parameters are chosen in the order
\[
\begin{gathered}
 \epsilon_K,\theta_W,\theta_\eta\to\theta_a\to C_R,C_J\to\epsilon_{\rm rec}\to\varepsilon_L\to L_{\rm sp}
 \to\lambda_F\to\theta_H\to A_{\rm sp}\\
 \to\Lambda_G\to L_c\to\epsilon_A,c_1,c_2\to\delta_W\to S_W,\mu_b ,
\end{gathered}
\]
each depending only on $M$, the fixed data and the preceding ones, except $S_W,\mu_b$, which depend
also on the background: $\epsilon_K,\theta_W,\theta_\eta$ of Lemma~\ref{l:mr}; $\theta_a\le1$ with
$\theta_aC_M\le\frac{q_a}{32}$ for $C_M$ of \eqref{e:spin}; $C_R$ of \eqref{e:coer} and $C_J=6C_R$;
$\epsilon_{\rm rec}$ with the first condition above; $\varepsilon_L$ with $12C_M\varepsilon_L^2\le q_a\theta_a$
for $C_M$ of \eqref{e:lowc}; $L_{\rm sp}$ as in \eqref{e:msw}; $\lambda_F,\theta_H,A_{\rm sp}$ as in
Lemma~\ref{l:tf}, which fixes $m_*$ of \eqref{e:split}; $\Lambda_G\ge\Lambda_*$ of Lemma~\ref{l:low} and
$L_c\ge L_*$ of Lemma~\ref{l:mix}, which depend on $m_*$; $\epsilon_A,c_1,c_2$ of Lemma~\ref{l:cw};
$\delta_W$ below the thresholds of the lemmas used, with the second condition above; and
$\mu_b\le\mu_0$ and $S_W$ above the entrances of the lemmas used, with the third condition for all
$s\ge S_W$ and $\mu_*\le\mu_b$, which is possible since each $\omega_{b,*}(s)$ tends to $0$. Each
condition on $s$ then holds for every $s\ge S_W$, so that \eqref{e:raw} and Lemma~\ref{l:cw} hold
for every $s\ge S_W$.
\end{proof}

\section{Decay and continuation}\label{s:dec}

We use the assumptions and notation of Section~\ref{sec:en}.

\subsection{The time-weighted first-order estimate}\label{s:age}

Fix a real radial $\varphi\in C_c^\infty(\R^3)$ with $0\le\varphi\le1$,
$\varphi=1$ on $\{|\xi|\le1\}$ and $\varphi=0$ on $\{|\xi|\ge2\}$. Set
\begin{equation}\label{e:low}
 \widehat{P_1f}=\varphi\widehat f.
\end{equation}
We also write
$\Vert f\Vert_{L^p_{\rm uloc}}=\sup_{z\in\R^3}\Vert f\Vert_{L^p(z+(-1,1)^3)}$.

\begin{lemma}\label{l:prd}
Let $f,a,G$ be real functions on $\R^3$ with $f\in L^3$, $\nabla f\in H^1$,
$a\in L^\infty$ and $\nabla a\in H^2$. Whenever the right-hand sides are
finite,
\begin{align}
 \Vert(\nabla f)G\Vert_{L^2}
 &\le C\bigl(\Vert f\Vert_{L^3}+\Vert f\Vert_{H^2_{\rm uloc}}\bigr)
       \Vert\nabla G\Vert_{L^2}
   +C\Vert G\Vert_{L^2_{\rm uloc}}\Vert\nabla^2f\Vert_{L^2},\label{e:prd1}\\
 \Vert(\nabla f)G\Vert_{L^2}
 &\le C\Vert f\Vert_{L^3}\Vert G\Vert_{L^2}
   +C\Vert G\Vert_{L^2_{\rm uloc}}\Vert\nabla^2f\Vert_{H^1},\label{e:prd2}\\
 \Vert\nabla f\Vert_{L^3}^3
 &\le C\Vert f\Vert_{L^3}\Vert\nabla^2f\Vert_{L^2}^2,\label{e:prd3}
\end{align}
\begin{equation}\label{e:prd4}
 \Vert\nabla a\Vert_{L^6}^3\le C\Vert a\Vert_{L^\infty}^2\Vert\nabla^3a\Vert_{L^2},
 \qquad
 \Vert\nabla^2a\Vert_{L^3}^3\le C\Vert a\Vert_{L^\infty}\Vert\nabla^3a\Vert_{L^2}^2,
 \qquad
 \Vert\nabla a\Vert_{L^{12}}^2\le C\Vert a\Vert_{L^\infty}\Vert\nabla^3a\Vert_{L^2}.
\end{equation}
\end{lemma}

\begin{proof}
The kernels of $P_1$ and $\nabla P_1$ are Schwartz functions; hence
$\Vert\nabla P_1f\Vert_{L^3}+\Vert\nabla P_1f\Vert_{L^\infty}
\le C\Vert f\Vert_{L^3}$, $\Vert P_1G\Vert_{L^\infty}\le C\Vert G\Vert_{L^2_{\rm uloc}}$
and $\Vert(1-P_1)f\Vert_{H^2_{\rm uloc}}\le C\Vert f\Vert_{H^2_{\rm uloc}}$.
Since the symbol of $1-P_1$ vanishes on $\{|\xi|\le1\}$,
$\Vert\nabla(1-P_1)f\Vert_{H^1}\le C\Vert\nabla^2f\Vert_{L^2}$,
$\Vert\nabla(1-P_1)f\Vert_{H^2}\le C\Vert\nabla^2f\Vert_{H^1}$ and
$\Vert(1-P_1)G\Vert_{H^1}\le C\Vert\nabla G\Vert_{L^2}$. In
\[
 (\nabla f)G=(\nabla P_1f)G+\nabla(1-P_1)f\,P_1G
   +\nabla(1-P_1)f\,(1-P_1)G
\]
the three terms are bounded in $L^2$ by
$\Vert\nabla P_1f\Vert_{L^3}\Vert G\Vert_{L^6}\le C\Vert f\Vert_{L^3}\Vert\nabla G\Vert_{L^2}$,
by $\Vert P_1G\Vert_{L^\infty}\Vert\nabla(1-P_1)f\Vert_{L^2}$ and, by the
first inequality of \eqref{e:prod}, by
$C\Vert\nabla(1-P_1)f\Vert_{H^1_{\rm uloc}}\Vert(1-P_1)G\Vert_{H^1}$; this
gives \eqref{e:prd1}. For \eqref{e:prd2},
$\Vert(\nabla P_1f)G\Vert_{L^2}\le\Vert\nabla P_1f\Vert_{L^\infty}\Vert G\Vert_{L^2}$
and, on each $Q_j$,
$\Vert\nabla(1-P_1)f\,G\Vert_{L^2(Q_j)}
\le C\Vert\nabla(1-P_1)f\Vert_{H^2(Q_j)}\Vert G\Vert_{L^2_{\rm uloc}}$; we sum the
squares over $j$.

The Gagliardo--Nirenberg inequality between $L^3$ and $\dot H^2$
gives \eqref{e:prd3}; between $L^\infty$ and $\dot H^3$ it gives
\eqref{e:prd4}. Approximation yields the stated regularity classes.
\end{proof}

By \eqref{e:prd4}, $\Vert\nabla a\Vert_{L^6}\Vert\nabla^2a\Vert_{L^3}
\le C\Vert a\Vert_{L^\infty}\Vert\nabla^3a\Vert_{L^2}$.

\begin{lemma}\label{l:age}
There is $S_C\ge S$ such that every solution on $[S_C,T_e]$ satisfies, with
\[
 \omega(t)=w_{S_C}(t)\tau^{\frac14}=\tau_{S_C}^{q_a}\tau^{\frac14-q_a},
\]
\begin{equation}\label{e:ptm}
 \sup_{S_C\le r\le t}\omega(r)\Vert\nabla(\Pi,\psi)(r)\Vert^2\le C_{b,S_C}\bigl(I_{S_C}+N_o(S_C,t)^2\bigr)
 \qquad(S_C\le t<T_e).
\end{equation}
\end{lemma}

\begin{proof}
Fix $t\in(S_C,T_e)$ and write $w=w_{S_C}$.
All time suprema and integrals in this proof are over $[S_C,t]$. Then
\[
 \omega=w\tau^{\frac14},\qquad
 \omega'=\left(\frac14-q_a\right)\frac\omega\tau\le\frac14w,\qquad
 \omega(S_C)=\tau_{S_C}^{\frac14},
\]
since $\tau^{\frac14}\le\tau$.

\emph{Lower-order bounds.} Put
\[
 g_\rho=\pa_1\log\rho_b,\qquad g_p=\pa_1\log p_b,
\]
with $g_T=\pa_1\log\theta_b$ as in \eqref{e:bc}. Then
\[
\begin{aligned}
 \nabla\log\frac{\rho}{\rho_b}&=\frac{\nabla\phi-g_\rho\phi e_1}\rho,
 &\nabla\log\frac{\theta}{\theta_b}&=\frac{\nabla\zeta-g_T\zeta e_1}\theta,\\
 \nabla\log\frac{p}{p_b}&=\frac{\nabla\Pi-g_p\Pi e_1}p,
 &1-\frac{\rho_b}{\rho}&=\frac\phi\rho,\qquad 1-\frac{p_b}{p}=\frac\Pi p.
\end{aligned}
\]
These formulas, their derivatives, $\Pi=R(\theta\phi+\rho_b\zeta)$,
\eqref{e:ha} with $s=S_C$, \eqref{e:trace}, \eqref{e:E1}, \eqref{e:gn4}, and
\eqref{e:curv} with \eqref{e:hpl} for the terms with $\pa_1^2(\rho_b,\theta_b)$,
give
\begin{equation}\label{e:agl}
\begin{aligned}
 &\sup w\bigl(E_*+\Vert\nabla(\phi,\log\frac{\rho}{\rho_b},\log\frac{\theta}{\theta_b},\log\frac{p}{p_b})\Vert^2+\Vert\nabla^2\psi\Vert^2\bigr)\\
 &\quad+\int w\bigl(D_*+\Vert\nabla\log\frac{p}{p_b}\Vert^2
   +\Vert\nabla^2(\phi,\Pi,\log\frac{\rho}{\rho_b},\log\frac{\theta}{\theta_b},\log\frac{p}{p_b})\Vert^2+\Vert\nabla^3\psi\Vert^2\bigr)
 \le C_M\bigl(I_{S_C}+N_o(S_C,t)^2\bigr) .
\end{aligned}
\end{equation}
Let $\mathcal Y$ be the set of functions $Y$ on $\R^3\times[S_C,t]$ with
$\sup w\Vert Y\Vert^2\le C_M\bigl(I_{S_C}+N_o(S_C,t)^2\bigr)$, and $\mathcal Z$ the set of $Z$ with
$\int w\Vert Z\Vert^2\le C_M\bigl(I_{S_C}+N_o(S_C,t)^2\bigr)$. By \eqref{e:agl}, $\mathcal Y$ contains $\Pi$,
$\psi$, $\nabla(\phi,\psi,\zeta,\Pi,\log\frac{\rho}{\rho_b},\log\frac{\theta}{\theta_b},\log\frac{p}{p_b})$, $\nabla^2(\psi,\zeta)$ and
$\tau^{-\frac12}(\phi,\zeta)$, and $\mathcal Z$ contains $\nabla(\Pi,\psi,\log\frac{p}{p_b})$,
$\nabla^2(\phi,\Pi,\psi,\zeta,\log\frac{\rho}{\rho_b},\log\frac{\theta}{\theta_b},\log\frac{p}{p_b})$, $\nabla^3(\psi,\zeta)$ and
$\tau^{-\frac12}(\Pi,\psi,\nabla\zeta)$.

\emph{Coefficients.} Put $L_T=\frac{\kappa\pa_1^2\theta_b}{c_v\rho_b\theta_b}$,
\[
\begin{gathered}
 \mathcal C_1=\{u_b-u_*,h_b,\pa_1\rho_b,\pa_1\theta_b,\pa_1p_b,g_\rho,g_T,g_p\},\\
 \mathcal C_2=\{h_b,\pa_1p_b,A_b,B_b,\pa_1h_b,\pa_1^2p_b,
   \pa_1^2\rho_b,\pa_1^2\theta_b,\pa_1g_\rho,\pa_1g_T,\pa_1g_p\}.
\end{gathered}
\]
By \eqref{e:E8}, and by \eqref{e:curv}, \eqref{e:E1} and
$\pa_1g_\rho=\frac{\pa_1^2\rho_b}{\rho_b}-g_\rho^2$,
$\pa_1g_T=\frac{\pa_1^2\theta_b}{\theta_b}-g_T^2$, $g_p=g_\rho+g_T$ for the last
five members of $\mathcal C_2$,
\[
 \sup_{r\ge S_C}\tau\Vert a\Vert_{L^\infty}^2\le C_b\quad(a\in\mathcal C_1),\qquad
 \int_{S_C}^\infty\tau^{\frac12}\Vert a\Vert_{L^\infty}^2\dd r\le C_b\quad(a\in\mathcal C_2),
\]
\[
 \int_{S_C}^\infty\Bigl(\tau\Vert(\pa_1A_b,\pa_1B_b,\pa_1L_T)\Vert_{L^2(\R)}^2
   +\tau^{\frac14}\Vert\pa_1^3\log\rho_b\Vert_{L^2(\R)}^2\Bigr)\dd r\le C_b .
\]
For $a=a(x_1,r)$ and $F\in\{\phi,\zeta\}$, resp. $F\in\{\psi,\Pi\}$,
$\Vert aF\Vert^2\le\Vert a\Vert_{L^2(\R)}^2\sup_{x_1}\Vert F(x_1,\cdot)\Vert_{L^2(\R^2)}^2$,
and the last factor is at most $C_M\sqrt\tau E_*$, resp. $2E_*$, by
\eqref{e:trace} and \eqref{eq:planar}. With $\omega=w\tau^{\frac14}$ and the
Cauchy--Schwarz inequality in $r$, these bounds give, for $Y\in\mathcal Y$ and
$Z,Z'\in\mathcal Z$,
\begin{equation}\label{e:agp}
\begin{gathered}
 \int\omega\Vert a\Vert_{L^\infty}\Vert Y\Vert\Vert Z\Vert
 +\int\omega\Vert a\Vert_{L^\infty}^2\Vert Y\Vert^2\le C_b\bigl(I_{S_C}+N_o(S_C,t)^2\bigr)
 \quad(a\in\mathcal C_2),\\
 \int\omega\Vert a\Vert_{L^\infty}\Vert Z\Vert\Vert Z'\Vert
 +\int\omega\Vert a\Vert_{L^\infty}^2\Vert Z\Vert^2\le C_b\bigl(I_{S_C}+N_o(S_C,t)^2\bigr)
 \quad(a\in\mathcal C_1),\\
 \int\omega\bigl(\Vert aF\Vert^2+\Vert aF\Vert\Vert Z\Vert\bigr)\le C_b\bigl(I_{S_C}+N_o(S_C,t)^2\bigr)
 \quad(a\in\{\pa_1A_b,\pa_1B_b,\pa_1L_T\},\ F\in\{\phi,\zeta\}),\\
 \int\omega\Vert(\pa_1^3\log\rho_b)\psi\Vert^2\le C_b\bigl(I_{S_C}+N_o(S_C,t)^2\bigr) .
\end{gathered}
\end{equation}
Moreover, by \eqref{e:E0} and \eqref{e:stb},
$\Vert(\rho-\rho_-,\theta-\theta_-,p-p_*)\Vert_{L^\infty}\le C_0\delta+\omega_b(S_C)
+C_M\mu_*$.

\emph{Products.} Put
\[
\begin{gathered}
 \mathcal T^2=\Vert\nabla^2(\Pi,\psi,\log\frac{p}{p_b})\Vert^2+\Vert\nabla^3(\psi,\zeta,\log\frac{\theta}{\theta_b})\Vert^2.
\end{gathered}
\]
As for $\Pi$ in the proof of Lemma~\ref{l:st}, the $L^\infty$ and
$H^2_{\rm uloc}$ norms of $\phi,\psi,\zeta,\Pi,\log\frac{\rho}{\rho_b},\log\frac{\theta}{\theta_b},\log\frac{p}{p_b}$ are at most
$C_M\mu_*$, and so are the $L^3$ norms of $\psi,\Pi,\log\frac{p}{p_b}$, since
$|\log\frac{p}{p_b}|\le C_M|\Pi|$. Hence \eqref{e:gn4} and Lemma~\ref{l:prd} give, for
$f\in\{\Pi,\log(p/p_b),\psi_1,\psi_2,\psi_3\}$ and
$a,d\in\{\zeta,\log(\theta/\theta_b),\psi_1,\psi_2,\psi_3\}$,
\begin{equation}\label{e:agn}
\begin{gathered}
 \Vert\nabla f\Vert_{L^3}+\Vert\nabla^2a\Vert_{L^3}
       \le C_M\mu_*^{\frac13}\mathcal T^{\frac23},\\
 \Vert\nabla f\Vert_{L^4}+\Vert\nabla a\Vert_{L^{12}}
       \le C_M\mu_*^{\frac12}\mathcal T^{\frac12},\qquad
 \Vert\nabla a\Vert_{L^6}\le C_M\mu_*^{\frac23}\mathcal T^{\frac13},\\
 \Vert(\nabla f)G\Vert\le C_M\mu_*\mathcal T\quad
 \bigl(G\in\{\nabla^2\psi,\nabla^2\zeta,\nabla^2\log\frac{\theta}{\theta_b}\}\bigr),\\
 \Vert(\nabla\psi)\nabla^2\log\frac{p}{p_b}\Vert+\Vert(\nabla f)\,\nabla a\otimes\nabla d\Vert
 \le C_M\mu_*\mathcal T,
\end{gathered}
\end{equation}
The product bounds use \eqref{e:prd1}. For $G=\pa_ia\,\pa_jd$,
$\Vert\nabla G\Vert\le\Vert\nabla^2a\Vert_{L^3}\Vert\nabla d\Vert_{L^6}
+\Vert\nabla a\Vert_{L^6}\Vert\nabla^2d\Vert_{L^3}$ and
$\Vert G\Vert_{L^2_{\rm uloc}}\le C_M\mu_*^2$.
The product $(\nabla\psi)\nabla^2\log\frac{p}{p_b}$ uses \eqref{e:prd2}.
Each factor $F$ in the interpolation bounds of
\eqref{e:agn}, each of $\nabla f,\nabla^2a$ in $L^6$ and of the factors of
$\mathcal T$ in $L^2$, each of $\psi,\Pi,\log\frac{p}{p_b}$ in $L^3$ and each of
$\phi,\psi,\zeta,\Pi,\log\frac{\rho}{\rho_b},\log\frac{\theta}{\theta_b},\log\frac{p}{p_b}$ in $L^\infty$ satisfies
$\Vert F\Vert_{L^{p_F}}\le C_M\mu_*^{1-t_F}\mathcal T^{t_F}$ with
$0\le t_F\le1$. We expand the gradients of the state variables by
\begin{equation}\label{e:aglg}
\begin{gathered}
 \nabla\log\rho=\frac{\nabla\Pi+(\pa_1p_b)e_1}p
                       -\frac{\nabla\zeta+(\pa_1\theta_b)e_1}\theta,\\
 \nabla\log\theta=\nabla \log\frac{\theta}{\theta_b}+g_Te_1,\qquad
 \nabla\log p=\nabla\log\frac{p}{p_b}+g_pe_1.
\end{gathered}
\end{equation}
$\nabla \log\frac{\rho}{\rho_b}=\nabla\log\frac{p}{p_b}-\nabla \log\frac{\theta}{\theta_b}$ and \eqref{e:state}. Then every nonlinear
term below is a bounded function times either a product of $k\ge3$ such
factors with $\sum_F\frac1{p_F}=1$ and $\sum_Ft_F=2$, whose integral is at
most $C_M\mu_*^{k-2}\mathcal T^2$ by H\"older's inequality, or a factor of
$\mathcal T$ times one of the products in \eqref{e:agn}; in both
cases its integral is at most $C_M\mu_*\mathcal T^2$. For instance,
\[
\begin{gathered}
 \int|\nabla\psi||\nabla\Pi|^2\dd x
 \le\Vert\nabla\psi\Vert_{L^3}\Vert\nabla\Pi\Vert_{L^3}^2,\qquad
 \int|\psi||\nabla\psi||\nabla^2\psi|\dd x
 \le\Vert\psi\Vert_{L^3}\Vert\nabla\psi\Vert_{L^6}\Vert\nabla^2\psi\Vert,\\
 \int|\nabla\Pi||\nabla^2\zeta||\nabla^3\zeta|\dd x
 \le\Vert(\nabla\Pi)\nabla^2\zeta\Vert\Vert\nabla^3\zeta\Vert,\qquad
 \int|\nabla\psi|^2|\nabla^2\zeta|^2\dd x
 \le\Vert\nabla\psi\Vert_{L^6}^2\Vert\nabla^2\zeta\Vert_{L^3}^2 .
\end{gathered}
\]
A term with one derivative less and an additional factor $a\in\mathcal C_1$
is bounded in the same way, with \eqref{e:gn4}, by
$C_M\mu_*\Vert a\Vert_{L^\infty}\Vert Z\Vert(\Vert Z'\Vert+\mathcal T)$ with
$Z,Z'\in\mathcal Z$, and then by \eqref{e:agp} and Young's inequality.

\emph{The second-order estimate in logarithmic variables.} Put
$b=u_b-u_*$ and $D_0=\pa_t+u_*\pa_1$.
By \eqref{eq:NSF},
\[
 D_u\log\rho+\Div u=0,\qquad D_uu+R\theta\nabla\log p=\frac{\mathcal Lu}\rho,
 \qquad c_vD_u\log\theta+R\Div u=\frac{\kappa\Delta\theta+Q_\mu(\nabla u)}{\rho\theta},
\]
and \eqref{e:bal} gives the same equations for $(\rho_b,u_be_1,\theta_b)$ with
$D_b$ in place of $D_u$. We subtract them and use
$D_uf=D_bf+\psi_1\pa_1f$ for $f=f(x_1,t)$,
$\frac{\Delta\theta}\theta=\Delta\log\theta+|\nabla\log\theta|^2$,
$Q_\mu(\nabla u)=\nu h_b^2+H[\psi]+Q_\mu(\nabla\psi)$ and $R=(\gamma-1)c_v$:
\begin{equation}\label{e:lgeq}
\begin{aligned}
 D_0\log\frac{\rho}{\rho_b}+\Div\psi&=\mathcal N_\rho,\\
 D_0\psi-\rho_-^{-1}\mathcal L\psi+R\theta_-\nabla\log\frac{p}{p_b}&=\mathcal N_\psi,\\
 D_0\log\frac{\theta}{\theta_b}-\frac\kappa{c_v\rho_-}\Delta \log\frac{\theta}{\theta_b}+(\gamma-1)\Div\psi&=\mathcal N_T,
\end{aligned}
\end{equation}
where
\[
\begin{aligned}
 \mathcal N_\rho={}&-\psi\cdot\nabla \log\frac{\rho}{\rho_b}-b\pa_1\log\frac{\rho}{\rho_b}-g_\rho\psi_1,\\
 \mathcal N_\psi={}&-\psi\cdot\nabla\psi-b\pa_1\psi
   +(\rho^{-1}-\rho_-^{-1})\mathcal L\psi-R(\theta-\theta_-)\nabla\log\frac{p}{p_b}\\
 &-h_b\psi_1e_1-\frac{A_b\phi}{\rho}e_1-\frac{\Pi}{\rho}g_pe_1,\\
 \mathcal N_T={}&-\psi\cdot\nabla \log\frac{\theta}{\theta_b}-b\pa_1\log\frac{\theta}{\theta_b}
   +\Bigl(\frac\kappa{c_v\rho}-\frac\kappa{c_v\rho_-}\Bigr)\Delta \log\frac{\theta}{\theta_b}\\
 &+\frac\kappa{c_v\rho}\bigl(|\nabla \log\frac{\theta}{\theta_b}|^2+2g_T\pa_1\log\frac{\theta}{\theta_b}\bigr)
   +\frac{H[\psi]+Q_\mu(\nabla\psi)}{c_v\rho\theta}\\
 &-g_T\psi_1-\frac{L_T\phi}{\rho}
      -\frac{\nu h_b^2\Pi}{c_v\rho_b\theta_bp}.
\end{aligned}
\]
Let
\begin{equation}\label{e:log}
\begin{gathered}
 E_{2\log}=\frac12\int\bigl(R\theta_-|\nabla^2\log\frac{\rho}{\rho_b}|^2+|\nabla^2\psi|^2
   +c_v\theta_-|\nabla^2\log\frac{\theta}{\theta_b}|^2\bigr)\dd x,\\
 \Gamma=\sum_{i=1}^3\int\pa_i\psi\cdot\nabla\pa_i\log\frac{p}{p_b}\dd x,\\
 D_H=\Vert\nabla^2\log\frac{p}{p_b}\Vert^2+\Vert\nabla^3\psi\Vert^2+\Vert\nabla^3\log\frac{\theta}{\theta_b}\Vert^2 .
\end{gathered}
\end{equation}
We apply $\pa_i\pa_j$ to \eqref{e:lgeq}, multiply the rows by
$R\theta_-\pa_{ij}\log\frac{\rho}{\rho_b}$, $\pa_{ij}\psi$, $c_v\theta_-\pa_{ij}\log\frac{\theta}{\theta_b}$, integrate and
sum over $i,j$. Since $p/p_b=(\rho/\rho_b)(\theta/\theta_b)$ and
$R=(\gamma-1)c_v$, the terms with $\Div\psi$ and $\nabla\log\frac{p}{p_b}$ cancel after
one integration by parts, and
\[
\begin{aligned}
 &E_{2\log}'+\frac1{\rho_-}\sum_{i,j}\bigl(\mu\Vert\nabla\pa_{ij}\psi\Vert^2
   +(\mu+\lambda)\Vert\Div\pa_{ij}\psi\Vert^2\bigr)
 +\frac{\kappa\theta_-}{\rho_-}\Vert\nabla^3\log\frac{\theta}{\theta_b}\Vert^2\\
 &\qquad=\sum_{i,j}\int\bigl(R\theta_-\pa_{ij}\log\frac{\rho}{\rho_b}\,\pa_{ij}\mathcal N_\rho
   +\pa_{ij}\psi\cdot\pa_{ij}\mathcal N_\psi
   +c_v\theta_-\pa_{ij}\log\frac{\theta}{\theta_b}\,\pa_{ij}\mathcal N_T\bigr)\dd x,
\end{aligned}
\]
where $\mu+\lambda>0$ by \eqref{eq:cf-hyp}. By \eqref{e:lgeq},
$D_0\log\frac{p}{p_b}=-\gamma\Div\psi+\frac\kappa{c_v\rho_-}\Delta \log\frac{\theta}{\theta_b}+\mathcal N_\rho+\mathcal N_T$; with
$\int\pa_1(\pa_i\psi\cdot\nabla\pa_i\log\frac{p}{p_b})\dd x=0$ this gives
\[
\begin{aligned}
 \Gamma'+R\theta_-\Vert\nabla^2\log\frac{p}{p_b}\Vert^2
 ={}&\sum_i\int\Bigl(\frac{\mathcal L\pa_i\psi\cdot\nabla\pa_i\log\frac{p}{p_b}}{\rho_-}
   +\gamma(\Div\pa_i\psi)^2-\frac\kappa{c_v\rho_-}\Div\pa_i\psi\,\Delta\pa_i\log\frac{\theta}{\theta_b}\Bigr)\dd x\\
 &+\sum_i\int\bigl(\pa_i\mathcal N_\psi\cdot\nabla\pa_i\log\frac{p}{p_b}
   -\Div\pa_i\psi\,\pa_i(\mathcal N_\rho+\mathcal N_T)\bigr)\dd x .
\end{aligned}
\]
We fix $\varepsilon_\Gamma\in(0,1]$, depending only on the fixed data, such
that Young's inequality bounds $\varepsilon_\Gamma$ times the first integral
by $\frac{R\theta_-\varepsilon_\Gamma}2\Vert\nabla^2\log\frac{p}{p_b}\Vert^2
+\frac\mu{2\rho_-}\Vert\nabla^3\psi\Vert^2
+\frac{\kappa\theta_-}{2\rho_-}\Vert\nabla^3\log\frac{\theta}{\theta_b}\Vert^2+C\Vert\nabla^2\psi\Vert^2$.
Then $E_H=E_{2\log}+\varepsilon_\Gamma\Gamma$ satisfies, with $c_H>0$ and
$C_H$ depending only on the fixed data,
\begin{equation}\label{e:lgid}
 E_H'+2c_HD_H\le C_H\Vert\nabla^2\psi\Vert^2+\mathcal S_\rho+\mathcal S_\psi+\mathcal S_T,
 \qquad
 |\varepsilon_\Gamma\Gamma|\le\frac14E_{2\log}+C\Vert\nabla\psi\Vert^2,
\end{equation}
where
\[
\begin{aligned}
 \mathcal S_\rho&=R\theta_-\sum_{i,j}\int\pa_{ij}\log\frac{\rho}{\rho_b}\,\pa_{ij}\mathcal N_\rho\dd x
   -\varepsilon_\Gamma\sum_i\int\Div\pa_i\psi\,\pa_i\mathcal N_\rho\dd x,\\
 \mathcal S_\psi&=\sum_{i,j}\int\pa_{ij}\psi\cdot\pa_{ij}\mathcal N_\psi\dd x
   +\varepsilon_\Gamma\sum_i\int\pa_i\mathcal N_\psi\cdot\nabla\pa_i\log\frac{p}{p_b}\dd x,\\
 \mathcal S_T&=c_v\theta_-\sum_{i,j}\int\pa_{ij}\log\frac{\theta}{\theta_b}\,\pa_{ij}\mathcal N_T\dd x
   -\varepsilon_\Gamma\sum_i\int\Div\pa_i\psi\,\pa_i\mathcal N_T\dd x .
\end{aligned}
\]

In $\mathcal S_\psi$ and $\mathcal S_T$ we move one derivative from
$\pa_{ij}\mathcal N_\psi$, $\pa_{ij}\mathcal N_T$; thus
$|\mathcal S_\psi|+|\mathcal S_T|\le C\mathcal T(\Vert\nabla\mathcal N_\psi\Vert
+\Vert\nabla\mathcal N_T\Vert)$. We write
$\bigl(\nabla\mathcal N_\psi,\nabla\mathcal N_T\bigr)=\mathcal P+\mathcal B$, where
$\mathcal B$ collects the terms with a factor $b$ or a derivative of the
background, for instance $g_\rho\rho^{-1}\mathcal L\psi$ from
$\nabla\rho^{-1}=-\rho^{-1}(\nabla \log\frac{\rho}{\rho_b}+g_\rho e_1)$. The coefficient
differences in $\mathcal P$ are at most $C_0\delta+\omega_b(S_C)+C_M\mu_*$,
and the products in $\mathcal P$ have at most third derivatives of $\psi$,
$\log\frac{\theta}{\theta_b}$ and second derivatives of $\log\frac{p}{p_b}$; by \eqref{e:agn},
\[
 \Vert\mathcal P\Vert\le C_M\bigl(\delta+\omega_b(S_C)+\mu_*\bigr)\mathcal T;
\]
for instance $\Vert\nabla\psi\nabla \log\frac{\theta}{\theta_b}\Vert\le\Vert\nabla\psi\Vert_{L^3}
\Vert\nabla \log\frac{\theta}{\theta_b}\Vert_{L^6}$, $\Vert\psi\nabla^2\log\frac{\theta}{\theta_b}\Vert\le\Vert\psi\Vert_{L^3}
\Vert\nabla^2\log\frac{\theta}{\theta_b}\Vert_{L^6}$ and $\Vert\nabla \log\frac{\theta}{\theta_b}\Vert_{L^6}^3\le C_M\mu_*^2\mathcal T$.
Each term of $\mathcal B$ is one of the products of \eqref{e:agp}, possibly
with an additional factor bounded by $C_M\mu_*$ as above; in particular
$\pa_1A_b(1-\frac{\rho_b}{\rho})$ and $\pa_1L_T(\frac{\rho_b}{\rho}-1)$ are trace terms with
$|1-\frac{\rho_b}{\rho}|=\frac{|\phi|}\rho$, and $h_b$, $\pa_1h_b$, $A_b$, $L_T$, $\pa_1g_p$,
$\pa_1g_T$ multiply members of $\mathcal Y$. Hence
$\int\omega\Vert\mathcal B\Vert^2\le C_b\bigl(I_{S_C}+N_o(S_C,t)^2\bigr)$.

In the first part of $\mathcal S_\rho$ put $w_*=\psi+be_1$, so that
$\mathcal N_\rho=-w_*\cdot\nabla \log\frac{\rho}{\rho_b}-g_\rho\psi_1$, and
\[
\begin{aligned}
 &-\sum_{i,j}\int\pa_{ij}\log\frac{\rho}{\rho_b}\,
              \pa_{ij}(w_*\cdot\nabla\log\frac{\rho}{\rho_b})\dd x\\
 &\quad=\int\frac{\Div w_*}2
                     |\nabla^2\log\frac{\rho}{\rho_b}|^2\dd x\\
 &\qquad-2\sum_{i,j,k}\int\pa_iw_{*k}\,
             \pa_{jk}\log\frac{\rho}{\rho_b}\,
             \pa_{ij}\log\frac{\rho}{\rho_b}\dd x\\
 &\qquad-\sum_{i,j,k}\int\pa_{ij}w_{*k}\,
             \pa_k\log\frac{\rho}{\rho_b}\,
             \pa_{ij}\log\frac{\rho}{\rho_b}\dd x.
\end{aligned}
\]
With $\nabla \log\frac{\rho}{\rho_b}=\nabla\log\frac{p}{p_b}-\nabla \log\frac{\theta}{\theta_b}$, the part with $\psi$ consists of
cubic terms, at most $C_M\mu_*\mathcal T^2$; for instance
$\int|\nabla\psi||\nabla^2\log\frac{p}{p_b}|^2\le\Vert(\nabla\psi)\nabla^2\log\frac{p}{p_b}\Vert
\Vert\nabla^2\log\frac{p}{p_b}\Vert$ and $\int|\nabla^2\psi||\nabla \log\frac{\theta}{\theta_b}||\nabla^2\log\frac{\theta}{\theta_b}|
\le\Vert\nabla^2\psi\Vert\Vert\nabla \log\frac{\theta}{\theta_b}\Vert_{L^6}\Vert\nabla^2\log\frac{\theta}{\theta_b}\Vert_{L^3}$.
The part with $be_1$ equals
$\int(\frac{h_b}2|\nabla^2\log\frac{\rho}{\rho_b}|^2-2h_b\sum_j|\pa_{1j}\log\frac{\rho}{\rho_b}|^2-\pa_1h_b\,\pa_1\log\frac{\rho}{\rho_b}\,\pa_1^2\log\frac{\rho}{\rho_b})$
and is bounded by \eqref{e:agp}, since $\nabla^2\log\frac{\rho}{\rho_b}\in\mathcal Z$ and
$\nabla \log\frac{\rho}{\rho_b}\in\mathcal Y$. For $F=g_\rho\psi_1$, the substitution
$\log\frac{\rho}{\rho_b}=\log\frac{p}{p_b}-\log\frac{\theta}{\theta_b}$ and one integration by parts in the part with $\log\frac{\theta}{\theta_b}$ give
\[
 -\sum_{i,j}\int\pa_{ij}\log\frac{\rho}{\rho_b}\,\pa_{ij}F\dd x
 =-\sum_{i,j}\int\bigl(\pa_{ij}\log\frac{p}{p_b}\,\pa_{ij}F+\pa_j\pa_{ij}\log\frac{\theta}{\theta_b}\,\pa_iF\bigr)\dd x
 \le D_H^{\frac12}\bigl(\Vert\nabla^2F\Vert+\Vert\nabla F\Vert\bigr),
\]
and $\nabla^2F$ consists of $g_\rho\nabla^2\psi_1$, $\pa_1g_\rho\nabla\psi_1$ and
$(\pa_1^3\log\rho_b)\psi_1$, so that $\int\omega(\Vert\nabla^2F\Vert^2
+\Vert\nabla F\Vert^2)\le C_b\bigl(I_{S_C}+N_o(S_C,t)^2\bigr)$ by \eqref{e:agp}. In the second part of
$\mathcal S_\rho$, $\nabla\mathcal N_\rho$ is bounded in the same way as $\nabla\mathcal N_\psi$.
Collecting these bounds, we obtain for every $\sigma\in(0,1]$
\[
 \int\omega\bigl(|\mathcal S_\rho|+|\mathcal S_\psi|+|\mathcal S_T|\bigr)
 \le\bigl(\sigma+C_M(\delta+\omega_b(S_C)+\mu_*)\bigr)\int\omega\mathcal T^2
   +C_{\sigma,b}\bigl(I_{S_C}+N_o(S_C,t)^2\bigr) .
\]
We multiply \eqref{e:lgid} by $\omega$ and integrate over $[S_C,t]$. By
$\omega'\le\frac14w$, $|E_H|\le C(E_{2\log}+\Vert\nabla\psi\Vert^2)$,
\eqref{e:agl} and $\omega(S_C)E_H(S_C)\le C_M\tau_{S_C}^{\frac14}I_{S_C}$,
\begin{equation}\label{e:logb}
\begin{aligned}
 &\omega(t)E_H(t)+2c_H\int\omega D_H\\
 &\quad\le C_{\sigma,b,S_C}\bigl(I_{S_C}+N_o(S_C,t)^2\bigr)
       +C_H\int\omega\Vert\nabla^2\psi\Vert^2\\
 &\qquad+\bigl(\sigma+C_M(\delta+\omega_b(S_C)+\mu_*)\bigr)\int\omega\mathcal T^2 .
\end{aligned}
\end{equation}

Differentiating $\Pi=p_b(\frac{p}{p_b}-1)$ twice,
\[
\begin{aligned}
 \nabla^2\Pi={}&\frac{\Pi\,\pa_1^2p_b}{p_b}e_1\otimes e_1\\
 &+\frac{p\,\pa_1p_b}{p_b}
     \left(e_1\otimes\nabla\log\frac{p}{p_b}
                  +\nabla\log\frac{p}{p_b}\otimes e_1\right)\\
 &+p\left(\nabla^2\log\frac{p}{p_b}
                  +\nabla\log\frac{p}{p_b}\otimes\nabla\log\frac{p}{p_b}\right),
\end{aligned}
\]
with $|\frac{p}{p_b}-1|=\frac{|\Pi|}{p_b}$ and
$\Vert\nabla\log\frac{p}{p_b}\Vert_{L^4}^2\le C_M\mu_*\Vert\nabla^2\log\frac{p}{p_b}\Vert$; hence, by
\eqref{e:agp},
\begin{equation}\label{e:plog}
\begin{aligned}
 \int\omega\Vert\nabla^2\Pi\Vert^2
 &\le C_M\int\omega\Vert\nabla^2\log\frac{p}{p_b}\Vert^2
       +C_b\bigl(I_{S_C}+N_o(S_C,t)^2\bigr),\\
 \int\omega\mathcal T^2
 &\le C_M\int\omega\bigl(D_H+\Vert\nabla^2\psi\Vert^2
   +\Vert\nabla^3\zeta\Vert^2\bigr)\\
 &\quad+C_b\bigl(I_{S_C}+N_o(S_C,t)^2\bigr).
\end{aligned}
\end{equation}

\emph{The first-order estimate.} Let
\[
 D=\Div\psi,\qquad Z=\Delta\zeta,\qquad
 \chi_d=\frac\kappa{c_v\rho},\qquad
 c_\theta=(\gamma-1)\theta,\qquad l_T=\frac{k_N}\rho.
\]
We use $d_N$ from \eqref{e:FC}.
The gradient of the first row of \eqref{e:peq3}, the divergence of the second
row divided by $\rho$, and the Laplacian of the third row give, with
$\Div\mathcal L\psi=\nu\Delta D$,
\begin{equation}\label{e:dc}
\begin{aligned}
 D_u\nabla\Pi+\gamma p\nabla D-(\gamma-1)\kappa\nabla Z&=F_G,\\
 D_uD-\frac\nu\rho\Delta D+\Div(\rho^{-1}\nabla\Pi)&=F_D,\\
 D_uZ-\chi_d\Delta Z+c_\theta\Delta D&=F_Z,
\end{aligned}
\end{equation}
where
\[
\begin{aligned}
 F_G={}&\nabla f_p-\gamma D\nabla p-(\nabla u)^\top \nabla\Pi,\qquad
 F_D=\nabla\rho^{-1}\cdot\mathcal L\psi+\Div(\rho^{-1}f_m)
   -\sum_{j,k}\pa_ju_k\,\pa_k\psi_j,\\
 F_Z={}&\Delta f_z+\Div(Z\nabla\chi_d-D\nabla c_\theta)
   +\nabla\chi_d\cdot\nabla Z-\nabla c_\theta\cdot\nabla D
   -(\Delta u)\cdot\nabla\zeta-2\sum_j\pa_ju\cdot\nabla\pa_j\zeta .
\end{aligned}
\]
We multiply the rows of \eqref{e:dc} by $\frac{\nabla\Pi}{\gamma p\rho}$, $D-d_NZ$, $l_TZ-d_ND$ and
integrate. Since $(\gamma-1)\kappa/(\gamma p\rho)=d_N/\rho$,
the pressure terms cancel up to $\int\rho^{-1}Z\nabla\Pi\cdot\nabla d_N$; since
$c_\theta l_T=d_N(\frac\nu\rho+\chi_d)$, the coupling terms of second order
combine. With $K_T=\frac{D_u\theta}\theta$, $D_u\bigl(\frac1{\gamma p\rho}\bigr)=(2\Div u-K_T)\frac1{\gamma p\rho}$,
$D_ud_N=(\Div u-K_T)d_N$, $D_ul_T=2(\Div u-K_T)l_T$ and
$\int afD_uf\dd x=\frac{\dd}{\dd t}\int\frac{af^2}2\dd x
-\int\frac{D_ua+a\Div u}2f^2\dd x$, we obtain
\begin{equation}\label{e:dcid}
 E_{\rm dc}'+D_{\rm dc}=\int\bigl(\frac{\nabla\Pi}{\gamma p\rho}\cdot F_G+(D-d_NZ)F_D
   +(l_TZ-d_ND)F_Z\bigr)\dd x+\mathcal R_{\rm sp}+\mathcal R_{\rm tm},
\end{equation}
where
\[
\begin{aligned}
 E_{\rm dc}={}&\frac12\int\bigl(\frac{|\nabla\Pi|^2}{\gamma p\rho}+D^2-2d_NDZ+l_TZ^2\bigr)\dd x,\\
 D_{\rm dc}={}&\int\bigl((\frac\nu\rho+d_Nc_\theta)|\nabla D|^2+l_T\chi_d|\nabla Z|^2
   -2l_Tc_\theta\nabla D\cdot\nabla Z\bigr)\dd x,\\
 \mathcal R_{\rm sp}={}&-\int\rho^{-1}Z\nabla\Pi\cdot\nabla d_N\dd x
  +\int\bigl(\frac\nu\rho Z\nabla d_N-(D-d_NZ)\nabla\Bigl(\frac\nu\rho\Bigr)\bigr)\cdot\nabla D\dd x\\
 &-\int\bigl(\chi_dZ\nabla l_T-\chi_dD\nabla d_N+(l_TZ-d_ND)\nabla\chi_d\bigr)
   \cdot\nabla Z\dd x\\
 &+\int\bigl(c_\theta Z\nabla l_T-c_\theta D\nabla d_N+(l_TZ-d_ND)\nabla c_\theta\bigr)
   \cdot\nabla D\dd x,\\
 \mathcal R_{\rm tm}={}&\int\Bigl(\frac1{2\gamma p\rho}(3\Div u-K_T)|\nabla\Pi|^2+\frac{\Div u}2D^2\Bigr)\dd x\\
 &-\int\Bigl(d_N(2\Div u-K_T)DZ-\frac{l_T}2(3\Div u-2K_T)Z^2\Bigr)\dd x .
\end{aligned}
\]
The expressions of $e_c$ and
$\mathcal D_c$ in Lemma~\ref{l:FC} for the carrier $N$ give, pointwise,
\[
\begin{gathered}
 \frac{|\nabla\Pi|^2}{\gamma p\rho}+D^2-2d_NDZ+l_TZ^2
 =\frac{|\nabla\Pi|^2}{\gamma p\rho}+(D-d_NZ)^2+\frac{k_N-\rho d_N^2}\rho Z^2,\\
 \rho\bigl((\frac\nu\rho+d_Nc_\theta)|\nabla D|^2+l_T\chi_d|\nabla Z|^2
   -2l_Tc_\theta\nabla D\cdot\nabla Z\bigr)
 =\frac\nu\gamma|\nabla D|^2+\frac{\kappa k_N}{c_v\rho}
   \Bigl|\nabla Z-\frac{(\gamma-1)c_v\rho\theta}\kappa\nabla D\Bigr|^2,
\end{gathered}
\]
the second by the computation of $\mathcal D_c$ there, since
$\rho(\frac\nu\rho+d_Nc_\theta)=\nu+(\gamma-1)j_N\theta$,
$\rho l_T\chi_d=\frac{\kappa k_N}{c_v\rho}$ and $\rho l_Tc_\theta=(\gamma-1)k_N\theta$.
With $\Vert D\Vert^2+\Vert\operatorname{curl}\psi\Vert^2=\Vert\nabla\psi\Vert^2$,
$\Vert Z\Vert=\Vert\nabla^2\zeta\Vert$, $\Vert\nabla Z\Vert=\Vert\nabla^3\zeta\Vert$
and \eqref{e:stb},
\[
 E_{\rm dc}\ge\frac1{C_M}\bigl(\Vert \nabla\Pi\Vert^2+\Vert D\Vert^2+\Vert Z\Vert^2\bigr),
 \qquad
 D_{\rm dc}\ge\frac1{C_M}\bigl(\Vert\nabla D\Vert^2+\Vert\nabla^3\zeta\Vert^2\bigr).
\]

Let $c_\Omega=\nabla\log\rho\times\psi$ and
$\Omega=\rho^{-1}\operatorname{curl}(\rho\psi)=\operatorname{curl}\psi+c_\Omega$. We
write the second row of \eqref{e:peq3} as
$D_u(\rho\psi)+\rho\psi\Div u+\nabla\Pi-\mathcal L\psi=f_m$ and apply
$\operatorname{curl}$; with
$\operatorname{curl}((u\cdot\nabla)m)=(u\cdot\nabla)\operatorname{curl}m
+\sum_k\nabla u_k\times\pa_km$,
$\operatorname{curl}f_m=-\nabla(\rho h_b\psi_1)\times e_1-A_b\nabla\phi\times e_1$,
$\sum_k\nabla\psi_k\times\pa_k\psi=(\Div\psi)\operatorname{curl}\psi
-(\operatorname{curl}\psi\cdot\nabla)\psi$ and, for a matrix $B$ and vectors $a,d$,
$(B^\top a)\times d+a\times(B^\top d)+B(a\times d)=({\rm tr}B)\,a\times d$, used
with $B=(\pa_j\psi_i)_{i,j}$, $a=\nabla\log\rho$, $d=\psi$, this gives
\[
 D_u\Omega+(\Div u)\Omega-(\Omega\cdot\nabla)u-\frac\mu\rho\Delta\Omega
 =-\frac\mu\rho\Delta c_\Omega-\frac{A_b}\rho\nabla\phi\times e_1
  +\nabla\log\rho\times\nabla\frac{|\psi|^2}2-\nabla(\Div u)\times\psi,
\]
and, after one integration by parts in the term with $\Delta c_\Omega$,
\[
\begin{aligned}
 \frac{\dd}{\dd t}\frac{\Vert\Omega\Vert^2}2+\mu\int\frac{|\nabla\Omega|^2}\rho\dd x
 ={}&\int\Bigl(\Omega\cdot(\Omega\cdot\nabla)u-\frac{\Div u}2|\Omega|^2
   +\mu\sum_k\frac{\pa_k\log\rho}\rho\,\Omega\cdot\pa_k\Omega\Bigr)\dd x\\
 &+\mu\sum_k\int\pa_k\Bigl(\frac\Omega\rho\Bigr)\cdot\pa_kc_\Omega\dd x\\
 &+\int\Omega\cdot\Bigl(-\frac{A_b}\rho\nabla\phi\times e_1
   +\nabla\log\rho\times\nabla\frac{|\psi|^2}2-\nabla(\Div u)\times\psi\Bigr)\dd x .
\end{aligned}
\]
By \eqref{e:aglg}, \eqref{e:agn} and $\Vert\psi\nabla\zeta\Vert\le\Vert\psi\Vert_{L^3}
\Vert\nabla\zeta\Vert_{L^6}$,
\[
 \Vert c_\Omega\Vert\le C_M\mu_*(\Vert\nabla\Pi\Vert+\Vert\nabla^2\zeta\Vert)
   +C_b\tau^{-\frac12}\Vert\psi\Vert,\qquad
 \Vert\nabla c_\Omega\Vert\le C_M\mu_*\mathcal T+\Vert\mathcal B_c\Vert,
\]
where $\mathcal B_c$ consists of $(\pa_1g_\rho)\psi$ and of members of
$\mathcal C_1$ times $\nabla\psi$, so that $\int\omega\Vert\mathcal B_c\Vert^2
\le C_b\bigl(I_{S_C}+N_o(S_C,t)^2\bigr)$; in particular $\Omega$ belongs to $\mathcal Y$ and $\mathcal Z$,
with $C_b$ in place of $C_M$.

All terms on the right of \eqref{e:dcid} and of the identity for $\Omega$
are now bounded by the rules above. The nonlinear terms give
$C_M\mu_*\int\omega\mathcal T^2$. In
$\int(l_TZ-d_ND)F_Z$ we move one derivative from $\Delta f_z$ and from
$\Div(Z\nabla\chi_d-D\nabla c_\theta)$; then $\nabla f_z$ and $\nabla c_\Omega$ are
paired with $\nabla(l_TZ-d_ND)$ and $\nabla\Omega$, and the background parts
of $\nabla f_z$, namely $\pa_1^2\theta_b\psi_1$, $\pa_1\theta_b\nabla\psi_1$,
$\pa_1B_b\frac\phi\rho$, $B_b\nabla\frac\phi\rho$ and the derivatives of
$\frac{h_b\Pi}{c_v\rho}$ and $\frac{H[\psi]}{c_v\rho}$, have
$\int\omega\Vert\cdot\Vert^2\le C_b\bigl(I_{S_C}+N_o(S_C,t)^2\bigr)$ by \eqref{e:agp}. The other
background terms are paired with $\nabla\Pi$, $D$, $Z$ or $\Omega$, which belong to
$\mathcal Z$, and are bounded by \eqref{e:agp}; for instance
$\Div(\rho^{-1}f_m)=-\pa_1h_b\psi_1-h_b\pa_1\psi_1-\pa_1A_b\frac\phi\rho
-A_b\pa_1\frac\phi\rho$ in $F_D$, the term
$-\frac{A_b}\rho\nabla\phi\times e_1$ and the term
$\frac{\rho_bB_b}{\rho\theta}$ of
\[
 K_T=\frac{\rho_bB_b}{\rho\theta}
   +\frac{\kappa Z+H[\psi]+Q_\mu(\nabla\psi)-h_b\Pi}{c_v\rho\theta}
   -(\gamma-1)D
\]
in $\mathcal R_{\rm tm}$. Moreover
$\Vert\nabla^2\psi\Vert^2=\Vert\nabla D\Vert^2+\Vert\nabla\operatorname{curl}\psi\Vert^2$
and $\Vert\nabla\operatorname{curl}\psi\Vert^2\le2\Vert\nabla\Omega\Vert^2
+2\Vert\nabla c_\Omega\Vert^2$. Put $E_{DC}=E_{\rm dc}+\frac12\Vert\Omega\Vert^2$. We
multiply the sum of \eqref{e:dcid} and of the identity for $\Omega$ by
$\omega$ and integrate over $[S_C,t]$; with $\omega'\le\frac14w$,
$\nabla\Pi,D,Z,\Omega\in\mathcal Z$ and $\omega(S_C)E_{DC}(S_C)\le
C_b\tau_{S_C}^{\frac14}I_{S_C}$, there is $c_D>0$, depending only on $M$ and
the fixed data, such that for every $\sigma\in(0,1]$
\begin{equation}\label{e:dcest}
 \omega(t)E_{DC}(t)+2c_D\int\omega\bigl(\Vert\nabla^2\psi\Vert^2
   +\Vert\nabla^3\zeta\Vert^2\bigr)
 \le C_{\sigma,b,S_C}\bigl(I_{S_C}+N_o(S_C,t)^2\bigr)+(\sigma+C_M\mu_*)\int\omega\mathcal T^2 .
\end{equation}

\emph{Conclusion.} By the bounds of $E_{\rm dc}$ and $c_\Omega$,
$\frac12\Vert\Omega\Vert^2\ge\frac14\Vert\operatorname{curl}\psi\Vert^2-\frac12\Vert c_\Omega\Vert^2$
and $\mu_*\le\mu_0$ small, with $E_1=\Vert\nabla\Pi\Vert^2+\Vert\nabla\psi\Vert^2
+\Vert\nabla^2\zeta\Vert^2$,
\[
 E_{DC}\ge2c_EE_1-C_b\tau^{-1}\Vert\psi\Vert^2,\qquad
 \omega\tau^{-1}\Vert\psi\Vert^2\le w\Vert\psi\Vert^2\le \bigl(I_{S_C}+N_o(S_C,t)^2\bigr),
\]
with $c_E>0$ depending only on $M$ and the fixed data. By \eqref{e:lgid},
$E_H\ge\frac34E_{2\log}-C\Vert\nabla\psi\Vert^2$. We fix $\varepsilon_H\in(0,1]$
with $\varepsilon_HC\le c_E$ and $\varepsilon_HC_H\le c_D$ and add
\eqref{e:dcest} and $\varepsilon_H$ times \eqref{e:logb}:
\[
\begin{aligned}
 &\omega\bigl(c_EE_1+\tfrac34\varepsilon_HE_{2\log}\bigr)(t)
 +c_D\int\omega\bigl(\Vert\nabla^2\psi\Vert^2+\Vert\nabla^3\zeta\Vert^2\bigr)
 +2\varepsilon_Hc_H\int\omega D_H\\
 &\qquad\le C_{\sigma,b,S_C}\bigl(I_{S_C}+N_o(S_C,t)^2\bigr)
   +\bigl(\sigma+C_M(\delta+\omega_b(S_C)+\mu_*)\bigr)\int\omega\mathcal T^2 .
\end{aligned}
\]
By \eqref{e:plog} the last term is at most $C_M(\sigma+C_M(\delta+\omega_b(S_C)
+\mu_*))\int\omega(D_H+\Vert\nabla^2\psi\Vert^2+\Vert\nabla^3\zeta\Vert^2)
+C_b\bigl(I_{S_C}+N_o(S_C,t)^2\bigr)$. We fix $\sigma$ with $C_M\sigma\le\frac14\min\{c_D,\varepsilon_Hc_H\}$,
then decrease $\delta_1(M)$ and $\mu_0$ and choose $S_C\ge S$ such that
$C_M^2(\delta+\omega_b(S_C)+\mu_*)\le\frac14\min\{c_D,\varepsilon_Hc_H\}$. Then
\[
 \omega(t)\bigl(E_1+E_{2\log}\bigr)(t)
 +\int_{S_C}^t\omega\bigl(\Vert\nabla^2\psi\Vert^2+\Vert\nabla^3\zeta\Vert^2+D_H\bigr)\dd r
 \le C_{b,S_C}\bigl(I_{S_C}+N_o(S_C,t)^2\bigr) .
\]
Taking the supremum over $r\in[S_C,t]$ and using
$\Vert\nabla(\Pi,\psi)\Vert^2\le E_1$ and the monotonicity of
$N_o(S_C,\cdot)$ proves \eqref{e:ptm}.
\end{proof}

\subsection{Local damping}\label{s:damp}

Put $\ell_{\rm low}(t)=\Vert P_1q(t)\Vert_{L^\infty}$.
The kernel $K_1$ of $P_1$ is a Schwartz function, and
$\ell_{\rm low}\le\Vert K_1\Vert_{L^1}\Vert q\Vert_{L^\infty}\le C\mu_*$.

\begin{lemma}\label{l:damp}
\begin{enumerate}
\item[(a)] Let $d>0$, $K\ge0$ and $s\le t<T_e$. If
$\ell_{\rm low}(r)\le K\tau_r^{-d}$ for $s\le r\le t$, then
\[
 \Vert q(t)\Vert_{H^2_{\rm uloc}}+\Vert q(t)\Vert_{L^\infty}
 \le C_{b,s,d}\bigl(I_s^{\frac12}+K\bigr)\tau_t^{-d}.
\]
\item[(b)] For $s\le t<T_e$,
$\sup_{s\le r\le t}\Vert q(r)\Vert_{H^2}^2
+\int_s^t\Vert\nabla(\psi,\zeta)\Vert_{H^2}^2\dd r\le C_be^{C_b(t-s)}I_s$.
\item[(c)] Let $(\rho,u,\theta)$ solve \eqref{eq:NSF} on $\R^3\times[0,T)$,
$S<T<\infty$, with $q$ in \eqref{eq:class}, $\rho,\theta$ bounded below by
positive constants on every $[0,T']$, $T'<T$, and $\mu_*\le\mu_0$ on every
$[S,T']$, $T'<T$. Then for some $\ell>0$ the solution extends to
$\R^3\times[0,T+\ell]$, with $q$ in \eqref{eq:class} and $\rho,\theta$
bounded below by a positive constant on $[0,T+\ell]$.
\end{enumerate}
\end{lemma}

\begin{proof}
Let $P_0$ be the Fourier multiplier with a real even symbol in
$C_c^\infty(\{|\xi|<1\})$ equal to $1$ on $\{|\xi|\le\frac12\}$; then
$P_0P_1=P_0$. Put $q^h=(1-P_0)q$ and, for $z\in\R^3$ and $L\ge1$,
\[
 \omega_z(x,t)=\Bigl(1+\frac{|x-z-u_*(t-s)e_1|^2}{L^2}\Bigr)^{-2},
 \qquad
 \Vert f\Vert_{j,z}^2=\sum_{k=0}^j\int\omega_z|\nabla^kf|^2\dd x,
\]
\[
\begin{gathered}
 E_z^0=\frac12\sum_{k=0}^2\binom2k\int\omega_z
   \Bigl(\frac{R\theta_-}{\rho_-}|\nabla^k\phi^h|^2+\rho_-|\nabla^k\psi^h|^2
   +\frac{c_v\rho_-}{\theta_-}|\nabla^k\zeta^h|^2\Bigr)\dd x,\\
 J_z=\int\omega_z\Bigl(\psi^h\cdot\nabla\phi^h
   +\sum_i\pa_i\psi^h\cdot\nabla\pa_i\phi^h\Bigr)\dd x,\qquad
 E_z=E_z^0+\epsilon_cJ_z,\\
 D_z=\Vert\nabla\phi^h\Vert_{1,z}^2+\Vert\nabla(\psi^h,\zeta^h)\Vert_{2,z}^2 .
\end{gathered}
\]
We show that there are $\epsilon_c,L,c_d,c_2>0$, depending only on $M$ and
the fixed data, such that on $(s,T_e)$
\begin{equation}\label{e:ldif}
 \frac1C\Vert q^h\Vert_{2,z}^2\le E_z\le C\Vert q^h\Vert_{2,z}^2,\qquad
 E_z'+c_dE_z+c_2D_z\le C_L\ell_{\rm low}^2 .
\end{equation}

By \eqref{e:E0} and \eqref{e:E6}, for a given $\beta_*>0$ we may decrease
$\delta_1(M)$ and increase $S$ so that, for $t\ge S$,
\[
 \Vert(\rho_b-\rho_-,u_b-u_*,\theta_b-\theta_-)\Vert_{L^\infty}
 +\Vert\pa_1(\rho_b,u_b,\theta_b)\Vert_{H^2(\R)}
 +\Vert(A_b,B_b)\Vert_{H^1(\R)}\le\beta_* .
\]
With $\varphi_z=\omega_z^{1/2}$,
\[
 (\pa_t+u_*\pa_1)\omega_z=0,\qquad
 |\nabla^j\varphi_z|\le C_jL^{-j}\varphi_z,\qquad
 \frac{\varphi_z(x)}{\varphi_z(x-v)}\le C\Bigl(1+\frac{|v|}L\Bigr)^2 .
\]
Hence convolution with a Schwartz function is bounded in
$\Vert\cdot\Vert_{j,z}$, uniformly in $z$ and $L$; for a Fourier multiplier
$B_0$ with symbol in $C_c^\infty$,
$\Vert[B_0,\varphi_z]f\Vert\le CL^{-1}\Vert\varphi_zf\Vert$; by
\eqref{eq:planar} applied to $\varphi_zf$,
$\Vert a(x_1)f\Vert_{0,z}\le C\Vert a\Vert_{L^2(\R)}\Vert f\Vert_{1,z}$; and
$\sup_Q\omega_z\le C\inf_Q\omega_z$ on every cube $Q$ of side $2$. Let the
symbol of $B_0$ equal $1$ near $0$ and vanish for $|\xi|\ge\frac14$. Then
$B_0q^h=0$, so $B_0(\varphi_zq^h)=[B_0,\varphi_z]q^h$, and
$\Vert(1-B_0)g\Vert\le C\Vert\nabla g\Vert$; for $L$ large,
\[
 \Vert q^h\Vert_{0,z}\le C\Vert\nabla q^h\Vert_{0,z},\qquad
 \Vert q^h\Vert_{2,z}^2\le CD_z,\qquad
 \Vert P_0q\Vert_{j,z}\le C_{j,L}\ell_{\rm low},
\]
the last since $\nabla^kP_0q=(\nabla^kK_0)*P_1q$ for the kernel $K_0$ of
$P_0$ and $\int\omega_z\dd x\le CL^3$. Moreover
$\Vert q^h\Vert_{2,v}\le C_L\Vert q\Vert_{H^2_{\rm uloc}}$, since
$\sum_j\sup_{Q_j}\omega_v\le CL^3$, and
$\Vert q\Vert_{H^2(Q)}\le C\Vert q\Vert_{2,v}$ for the cube $Q$ of side $2$
centred at $v+u_*(t-s)e_1$; thus
\[
 m=\beta_*+\Vert q\Vert_{H^2_{\rm uloc}}+\sup_v\Vert q^h\Vert_{2,v}
 \le\beta_*+C_L\mu_* .
\]

Since $1-P_0$ commutes with the constant coefficients in \eqref{e:cpeq},
$q^h$ satisfies \eqref{e:cpeq} with the sources $(1-P_0)(N_\phi,N_\psi,N_\zeta)$
of \eqref{e:nphys}. Write $N_\phi=-w_*\cdot\nabla\phi+N_\phi'$ with
$w_*=(u_b-u_*)e_1+\psi$. By the products on the cubes $Q_j$ with
$H^2(Q_j)\subset L^\infty(Q_j)$, $H^1(Q_j)\subset L^6(Q_j)$, by the weighted
trace bound for the terms with $\pa_1^3\rho_b$, $\pa_1^2u_b$, $\pa_1A_b$,
$\pa_1B_b$, and by $\Vert q\Vert_{2,z}\le\Vert q^h\Vert_{2,z}+C_L\ell_{\rm low}$,
\[
\begin{gathered}
 \Vert(N_\phi,N_\psi,N_\zeta)\Vert_{0,z}+\Vert N_\phi\Vert_{1,z}
   \le Cm\bigl(\Vert q^h\Vert_{2,z}+C_L\ell_{\rm low}\bigr),\\
 \Vert N_\phi'\Vert_{2,z}+\Vert(N_\psi,N_\zeta)\Vert_{1,z}
   \le Cm\bigl(\Vert q^h\Vert_{2,z}+\Vert\nabla(\psi^h,\zeta^h)\Vert_{2,z}
   +C_L\ell_{\rm low}\bigr);
\end{gathered}
\]
for instance
$\Vert\frac{Q_\mu(\nabla\psi)}\rho\Vert_{1,z}\le C\Vert q\Vert_{H^2_{\rm uloc}}
\Vert\nabla\psi\Vert_{2,z}$ and
$\Vert(\pa_1^3\rho_b)\psi_1\Vert_{0,z}\le C\beta_*\Vert\psi\Vert_{1,z}$.
For the density transport, $(1-P_0)(w_*\cdot\nabla\phi)
=w_*\cdot\nabla\phi^h-P_0(w_*\cdot\nabla\phi)+w_*\cdot\nabla P_0\phi$, and for
$|\alpha|=2$, $\alpha=e_i+e_j$,
\[
 \pa^\alpha(w_*\cdot\nabla\phi^h)-w_*\cdot\nabla\pa^\alpha\phi^h
 =\pa_iw_*\cdot\nabla\pa_j\phi^h+\pa_jw_*\cdot\nabla\pa_i\phi^h
   +\pa_{ij}w_*\cdot\nabla\phi^h .
\]
On each $Q_j$,
$\int_{Q_j}|\nabla\psi||\nabla^2\phi^h|^2\le C\Vert\nabla\psi\Vert_{H^2(Q_j)}
\Vert\phi^h\Vert_{H^2(Q_j)}\Vert\nabla^2\phi^h\Vert_{L^2(Q_j)}$ and
$\int_{Q_j}|\nabla^2\psi||\nabla\phi^h||\nabla^2\phi^h|\le
\Vert\nabla^2\psi\Vert_{L^3(Q_j)}\Vert\nabla\phi^h\Vert_{L^6(Q_j)}
\Vert\nabla^2\phi^h\Vert_{L^2(Q_j)}$; we multiply by $\sup_{Q_j}\omega_z$,
sum over $j$ and use the Cauchy--Schwarz inequality. With
$\int\omega_zf\,w_*\cdot\nabla f=-\frac12\int\Div(\omega_zw_*)f^2$, the
terms with $u_b-u_*$, and
$\Vert\nabla\psi\Vert_{2,z}\le\Vert\nabla\psi^h\Vert_{2,z}+C_L\ell_{\rm low}$,
\[
 \Bigl|\sum_{|\alpha|\le2}\int\omega_z\pa^\alpha\phi^h\,
   \pa^\alpha(w_*\cdot\nabla\phi^h)\dd x\Bigr|
 \le C\Bigl(m+\frac1L\Bigr)\bigl(\Vert q^h\Vert_{2,z}^2
   +\Vert\nabla\psi^h\Vert_{2,z}^2\bigr)+C_L\ell_{\rm low}^2 .
\]

We apply $\pa^\alpha$, $|\alpha|\le2$, to the equations of $q^h$, multiply
by $\omega_z$ times the coefficients of $E_z^0$ and $\pa^\alpha q^h$, and
integrate. By $(\pa_t+u_*\pa_1)\omega_z=0$ the transport by $u_*$ gives
$(E_z^0)'$; the first-order couplings cancel up to terms with $\nabla\omega_z$;
in the parabolic rows we move one derivative from the sources to the test
functions. Differentiating $J_z$ by the same equations gives
$-\frac{R\theta_-}{\rho_-}\Vert\nabla\phi^h\Vert_{1,z}^2$ and terms with at
most one derivative of $N_\phi$. Hence, with $c_1,c_r,C_r$ depending only on
the fixed data,
\[
\begin{aligned}
 (E_z^0)'+c_1\Vert\nabla(\psi^h,\zeta^h)\Vert_{2,z}^2
 &\le C\Bigl(m+\frac1L\Bigr)\bigl(\Vert q^h\Vert_{2,z}^2+D_z\bigr)
   +C_L\ell_{\rm low}^2,\\
 J_z'+c_r\Vert\nabla\phi^h\Vert_{1,z}^2
 &\le C_r\Vert\nabla(\psi^h,\zeta^h)\Vert_{2,z}^2
   +C\Bigl(m+\frac1L\Bigr)\bigl(\Vert q^h\Vert_{2,z}^2+D_z\bigr)+C_L\ell_{\rm low}^2 .
\end{aligned}
\]
We choose $\epsilon_c$ with $\epsilon_cC_r\le\frac{c_1}2$ and
$\epsilon_c|J_z|\le\frac12E_z^0$, then $L$ large, and then $\beta_*$ and
$\mu_0$ small, so that $C(m+\frac1L)(1+C)\le\frac14\min\{c_1,\epsilon_cc_r\}$
with the constant of $\Vert q^h\Vert_{2,z}^2\le CD_z$. This and
$\Vert q^h\Vert_{2,z}^2\le CD_z$ give \eqref{e:ldif}.

By \eqref{e:ldif} and $E_z(s)\le CI_s$,
\[
 E_z(t)\le Ce^{-c_d(t-s)}I_s+C_L\int_s^te^{-c_d(t-r)}\ell_{\rm low}(r)^2\dd r,
\]
and $\Vert q\Vert_{H^2_{\rm uloc}}+\Vert q\Vert_{L^\infty}
\le C(\sup_zE_z^{1/2}+C_L\ell_{\rm low})$. If
$\ell_{\rm low}(r)\le K\tau_r^{-d}$, splitting the integral at $\frac t2$ gives
\[
 \int_s^te^{-c_d(t-r)}\tau_r^{-2d}\dd r+e^{-c_d(t-s)}\le C_{s,d}\tau_t^{-2d},
\]
which proves (a). The same tests on $q$ with $\omega_z$ replaced by $1$ and
$P_0$ by $0$ give, for
\[
 E_{\rm gl}=\frac12\sum_{k=0}^2\binom2k\int\Bigl(\frac{R\theta_-}{\rho_-}
   |\nabla^k\phi|^2+\rho_-|\nabla^k\psi|^2+\frac{c_v\rho_-}{\theta_-}
   |\nabla^k\zeta|^2\Bigr)\dd x
 +\epsilon_c\int\Bigl(\psi\cdot\nabla\phi+\sum_i\pa_i\psi\cdot\nabla\pa_i\phi\Bigr)\dd x,
\]
$\frac1C\Vert q\Vert_{H^2}^2\le E_{\rm gl}\le C\Vert q\Vert_{H^2}^2$ and
$E_{\rm gl}'+c_{\rm gl}\Vert\nabla(\psi,\zeta)\Vert_{H^2}^2\le C_bE_{\rm gl}$, where
the terms with $m$ are absorbed as before and the other background terms,
which need not be small here, are at most $C_bE_{\rm gl}$; Gronwall's
inequality gives (b). In (c), (b) on $[S,T']$ and
\eqref{e:stb} give \eqref{e:cont} with $T_{\rm fin}=T$, and
Lemma~\ref{l:cont} with $T_f=T+1$ gives the extension.
\end{proof}

\subsection{The temperature at low frequency}\label{s:adj}

Fix the adjoint smoothing exponents by
\begin{equation}\label{e:decpar}
\begin{gathered}
 r_a=\frac{128}{127},\qquad
 \varsigma=\frac32\left(\frac1{r_a}-\frac12\right)=\frac{189}{256},\\
 d_*=\frac14\left(\varsigma-\frac12\right)-\frac{q_a}2
     =\frac{57}{1024}.
\end{gathered}
\end{equation}
\begin{lemma}\label{l:adj}
Let $s\le t<T_e$ and $\mathcal X\ge0$ with $N_o(s,t)\le\mathcal X$ and
$\tau_r^{\frac14-q_a}\Vert\nabla(\Pi,\psi)(r)\Vert^2\le\mathcal X^2$ for $s\le r\le t$. Then
\begin{equation}\label{e:thout}
 \sup_{\Vert J\Vert_{L^\infty}\le1}
 \Bigl\Vert P_1\Bigl(J\frac{\zeta(t)}{\theta_b(t)}\Bigr)\Bigr\Vert_{L^\infty}
 \le C_{b,s}\bigl(I_s^{\frac12}+\mathcal X+\mathcal X^2\bigr)\tau_t^{-d_*},
\end{equation}
where the supremum is over real $J\in L^\infty(\R^3)$.
\end{lemma}

\begin{proof}
By \eqref{eq:phzeta},
$D_u\theta_b=B_b+\psi_1\pa_1\theta_b$, $\phi+\frac{\rho\zeta}{\theta_b}=\frac\Pi{R\theta_b}$ and
$\rho(1+\frac\zeta{\theta_b})=\frac{\rho\theta}{\theta_b}=\rho_b+\frac\Pi{R\theta_b}$,
\[
 \rho D_u\Bigl(\frac\zeta{\theta_b}\Bigr)-\frac\kappa{c_v\theta_b}\Delta\zeta
 =F_\nabla+F_c+F_R+F_Q,
\]
where, with $k_{1,2}$ the second component of $k_1$ in \eqref{e:E1},
$g_{T,c}=\frac{(\Theta_h)_y}{V_b\Theta_b}\circ Y$ and
$g_{T,R}=\frac{k_{1,2}}{\theta_b}$, so that $g_T=g_{T,c}+g_{T,R}$,
\[
\begin{gathered}
 F_\nabla=\frac{H[\psi]-p\Div\psi}{c_v\theta_b},\qquad
 F_c=-\Bigl(\rho_b+\frac\Pi{R\theta_b}\Bigr)g_{T,c}\psi_1,\qquad
 F_Q=\frac{Q_\mu(\nabla\psi)}{c_v\theta_b},\\
 F_R=-\Bigl(\rho_b+\frac\Pi{R\theta_b}\Bigr)g_{T,R}\psi_1
   -\frac{B_b\Pi}{R\theta_b^2}-\frac{h_b\Pi}{c_v\theta_b}.
\end{gathered}
\]
By \eqref{e:E0}, \eqref{e:E1}, \eqref{e:E5} and \eqref{e:stb}, for $r\ge s$,
$\theta_b\ge\frac{\theta_-}2$, $\rho\ge\frac{\rho_-}4$,
\[
\begin{gathered}
 \Vert g_{T,c}\Vert_{L^\infty}\le\frac{2C_0\delta}{\theta_-}\tau^{-\frac12},\qquad
 \int_s^\infty\Vert g_{T,R}\Vert_{L^\infty}^2\dd r\le C_M,\\
 \int_s^\infty c_R^2\dd r\le C_M,\qquad
 c_R=\Vert(g_{T,R},B_b,h_b)\Vert_{L^2(\R)} .
\end{gathered}
\]

Fix $x\in\R^3$ and $J$ with $\Vert J\Vert_{L^\infty}\le1$, and let $\Psi$
solve, for $s\le r\le t$,
\[
 -D_u\Psi=\frac{\kappa\theta_b}{c_v\rho}\Delta\frac\Psi{\theta_b},\qquad
 \Psi(t)=\frac{JK_1(x-\cdot)}{\rho(t)},
\]
with $\Psi\in C([s,t];L^2)\cap L^2(s,t;H^1)$. To construct it, let
$j_\varepsilon$ be a nonnegative spatial mollifier and set
\[
 \rho^\varepsilon=j_\varepsilon*\rho,\qquad
 u^\varepsilon=\frac{j_\varepsilon*(\rho u)}{\rho^\varepsilon},\qquad
 J^\varepsilon=j_\varepsilon*J.
\]
Then $\pa_r\rho^\varepsilon+\Div(\rho^\varepsilon u^\varepsilon)=0$,
the density bounds are preserved, and $\Vert J^\varepsilon\Vert_\infty\le1$.
Galerkin approximation for the regularized backward equations, with
terminal value $J^\varepsilon K_1(x-\cdot)/\rho^\varepsilon(t)$,
followed by the weighted $L^2$ energy estimate and the maximum principle,
gives solutions $\Psi^\varepsilon$ satisfying
\[
 \sup_{s\le r\le t}
   \bigl(\Vert\Psi^\varepsilon(r)\Vert+
         \Vert\Psi^\varepsilon(r)\Vert_\infty\bigr)
 +\Vert\nabla\Psi^\varepsilon\Vert_{L^2(s,t;L^2)}\le C_{b,s,t},
\]
uniformly in $\varepsilon$. Moreover,
\[
 \pa_r(\rho^\varepsilon\Psi^\varepsilon)
 =-\Div(\rho^\varepsilon u^\varepsilon\Psi^\varepsilon)
  -\frac\kappa{c_v}\theta_b
       \Delta\frac{\Psi^\varepsilon}{\theta_b}.
\]
The multiplier bounds for $\rho^\varepsilon$ give uniform bounds for
$\rho^\varepsilon\Psi^\varepsilon$ in $L^2(s,t;H^1)$ and its time
derivative in $L^2(s,t;H^{-1})$. The compactness criterion
\cite[Corollary~4]{Simon1987}, applied on bounded spatial domains,
and the strong convergence of the coefficients give a weak solution
with the stated terminal value. Its weighted unknown belongs to $C([s,t];L^2)$;
since $\rho\in C([s,t];L^\infty)$ and $\inf\rho>0$, so does $\Psi$.
The energy and smoothing calculations below hold for these approximants
with uniform bounds and pass to the limit. In particular, the $L^1$
heating source is paired using the uniform $L^\infty$ bound and weak-star convergence of
$\Psi^\varepsilon$. Since
$\pa_t(\frac{\rho\zeta}{\theta_b})=\rho D_u\Bigl(\frac\zeta{\theta_b}\Bigr)-\Div(\frac{\rho u\zeta}{\theta_b})$, the forward--backward
weak pairing, justified by space-time smoothing and spatial cutoffs, gives
$\frac{\dd}{\dd r}\int\frac{\rho\zeta\Psi}{\theta_b}\dd x=\int(F_\nabla+F_c+F_R+F_Q)\Psi\dd x$,
and $\int\frac{\rho\zeta\Psi}{\theta_b}\dd x=P_1(J\frac{\zeta(t)}{\theta_b(t)})(x)$ at $r=t$. Hence
\begin{equation}\label{e:dual}
 |P_1(J\frac{\zeta(t)}{\theta_b(t)})(x)|\le C\Vert\zeta(s)\Vert\Vert\Psi(s)\Vert
   +\int_s^t\bigl(\Vert F_\nabla+F_c+F_R\Vert\Vert\Psi\Vert
   +\Vert F_Q\Vert_{L^1}\Vert\Psi\Vert_{L^\infty}\bigr)\dd r .
\end{equation}

We show, with $K(t,r)=(\frac{\tau_t}{\tau_r})^{\frac1{64}}(1+t-r)^{-\varsigma}$,
\begin{equation}\label{e:adj}
 \Vert\Psi(r)\Vert+\Vert\Psi(r)\Vert_{L^\infty}\le C_bK(t,r),
 \qquad s\le r\le t .
\end{equation}
Since $\theta_b\Delta\frac\Psi{\theta_b}=\Delta\Psi-2g_T\pa_1\Psi
+(g_T^2-\pa_1g_T)\Psi$ and $\pa_r\rho+\Div(\rho u)=0$, the function
$Y_a(j)=\int\rho|\Psi|^a\dd x$ at $r=t-j$ satisfies, for $1<a<\infty$,
\[
 \frac1a\frac{\dd Y_a}{\dd j}
 =\frac\kappa{c_v}\int\Bigl(-(a-1)|\Psi|^{a-2}|\nabla\Psi|^2
   +(a-2)g_T|\Psi|^{a-2}\Psi\pa_1\Psi+g_T^2|\Psi|^a\Bigr)\dd x .
\]
By Young's inequality the middle term is at most
$\frac{a-1}2\int|\Psi|^{a-2}|\nabla\Psi|^2+\frac{(2-a)^2}{2(a-1)}\int g_T^2|\Psi|^a$,
and $\int|\Psi|^{a-2}|\nabla\Psi|^2=\frac4{a^2}\Vert\nabla|\Psi|^{\frac a2}\Vert^2$.
Put $C_*=\frac{64\kappa C_0^2}{c_v\theta_-^2\rho_-}$ and decrease $\delta_1(M)$
so that
\[
 C_*\Bigl(1+\frac{(2-r_a)^2}{2(r_a-1)}\Bigr)\delta^2\le\frac{r_a}{64} .
\]
With $g_T^2\le2g_{T,c}^2+2g_{T,R}^2$ and $\int|\Psi|^a\le\frac4{\rho_-}Y_a$ this
gives, for $a\in\{r_a,2\}$ and $d_a=\frac{2\kappa(a-1)}{c_va}$,
\[
 \frac{\dd Y_a}{\dd j}+d_a\Vert\nabla|\Psi|^{\frac a2}\Vert^2
 \le\Bigl(\frac{r_a}{64\tau_{t-j}}+b_R(t-j)\Bigr)Y_a,\qquad
 b_R=C_M\Vert g_{T,R}\Vert_{L^\infty}^2 .
\]
Let $W(j)=\exp\{-\frac1{r_a}\int_{t-j}^t(\frac{r_a}{64\tau_v}+b_R(v))\dd v\}
\Psi(t-j)$. Then $\int\rho|W(j)|^{r_a}\dd x\le Y_{r_a}(0)$, and, since
$\frac2{r_a}>1$, $Y=\int\rho|W|^2\dd x$ satisfies $Y'+d_2\Vert\nabla W\Vert^2\le0$.
The Gagliardo--Nirenberg inequality
$\Vert W\Vert^2\le C\Vert\nabla W\Vert^{2\vartheta}\Vert W\Vert_{L^{r_a}}^{2(1-\vartheta)}$,
$\vartheta=\frac{2\varsigma}{1+2\varsigma}$, gives, with $N=C\Vert\Psi(t)\Vert_{L^{r_a}}$,
\[
 Y'+c_*N^{-\frac1\varsigma}Y^{1+\frac1{2\varsigma}}\le0,\qquad
 Y(j)\le\min\bigl\{Y(0),CN^2j^{-2\varsigma}\bigr\},
\]
since $\frac1\vartheta=1+\frac1{2\varsigma}$ and
$\frac{2(1-\vartheta)}\vartheta=\frac1\varsigma$. As
$\exp\{\int_r^t\frac{\dd v}{64\tau_v}\}=(\frac{\tau_t}{\tau_r})^{\frac1{64}}$,
$\int_s^\infty b_R\le C_M$ and $\Vert\Psi(t)\Vert_{L^{r_a}}+\Vert\Psi(t)\Vert
\le C\Vert K_1\Vert_{L^{r_a}\cap L^2}$, this is the bound of $\Vert\Psi(r)\Vert$
in \eqref{e:adj}. For $a\ge2$ the same identity gives
$\frac{\dd Y_a}{\dd j}+\frac\kappa{c_v}\Vert\nabla|\Psi|^{\frac a2}\Vert^2\le C_ba^2Y_a$.
On a reverse unit interval, starting at $r+1$, take $a_k=2(\frac53)^k$ and
$b_k=\frac12(1-2^{-k})$, and multiply this inequality for $a=a_k$ by a time
cutoff rising from $0$ to $1$ between $b_k$ and $b_{k+1}$. With
$\int_{b_{k+1}}^1\Vert v\Vert_{L^{10/3}}^{10/3}\le C\bigl(\sup_{[b_{k+1},1]}
\Vert v\Vert^2\bigr)^{\frac23}\int_{b_{k+1}}^1\Vert\nabla v\Vert^2$ for
$v=|\Psi|^{a_k/2}$,
\[
 \Vert\Psi\Vert_{L^{a_{k+1}}([b_{k+1},1]\times\R^3)}
 \le\bigl(C_b(a_k^2+2^k)\bigr)^{\frac1{a_k}}
   \Vert\Psi\Vert_{L^{a_k}([b_k,1]\times\R^3)},
\]
and $\sum_k\frac{1+k+\log a_k}{a_k}<\infty$; hence
$\Vert\Psi(r)\Vert_{L^\infty}\le C_b\Vert\Psi(r+1)\Vert$ if $t-r\ge2$. If
$t-r\le2$, the maximum principle for the equation of $\Psi$, whose
coefficient $\frac{\kappa(g_T^2-\pa_1g_T)}{c_v\rho}$ is bounded by
\eqref{e:E6}, gives $\Vert\Psi(r)\Vert_{L^\infty}\le C_b\Vert\Psi(t)\Vert_{L^\infty}
\le C_b$. This proves \eqref{e:adj}.

By \eqref{eq:planar} and $\Vert(\psi,\Pi)(r)\Vert^2\le\tau_r^{q_a}N_o^2$,
\[
\begin{gathered}
 \Vert F_\nabla\Vert\le C_b\Vert\nabla(\Pi,\psi)\Vert,\qquad
 \Vert F_Q\Vert_{L^1}\le C\Vert\nabla(\Pi,\psi)\Vert^2,\qquad
 \Vert F_c\Vert\le C_M\delta\tau^{-\frac12}\Vert\psi\Vert
   \le C_M\delta\mathcal X\tau^{-\frac12+\frac{q_a}{2}},\\
 \Vert F_R\Vert\le C_bc_R\Vert(\psi,\Pi)\Vert^{\frac12}\Vert\nabla(\Pi,\psi)\Vert^{\frac12}
   \le C_bc_R\mathcal X^{\frac12}\tau^{\frac{q_a}{4}}\Vert\nabla(\Pi,\psi)\Vert^{\frac12},
\end{gathered}
\]
and $\int_s^t\tau^{-q_a}\Vert\nabla(\Pi,\psi)\Vert^2\dd r\le\tau_s^{-q_a}N_o^2\le\mathcal X^2$. Let $t\ge2s+2$. On $[s,t/2]$,
$K(t,r)\le C\tau_t^{1/64-\varsigma}\tau_r^{-1/64}$. The preceding source bounds,
Cauchy--Schwarz, and H\"older with exponents $(2,4,4)$ for $F_R$ give
\[
\begin{aligned}
 &\int_s^{t/2}K\bigl(\Vert F_\nabla\Vert+\Vert F_c\Vert+\Vert F_R\Vert
       +\Vert F_Q\Vert_{L^1}\bigr)\dd r\\
 &\qquad\le C_b\Bigl(\mathcal X\tau_t^{-(\varsigma-1/2-q_a/2)}
       +\mathcal X\tau_t^{-(\varsigma-1/4-q_a/2)}
       +\mathcal X^2\tau_t^{-(\varsigma-1/64)}\Bigr).
\end{aligned}
\]
On $[t/2,t]$, $K(t,r)\le C(1+t-r)^{-\varsigma}$ and
\[
 \sup_{t/2\le r\le t}\Vert\nabla(\Pi,\psi)(r)\Vert
 \le C\mathcal X\tau_t^{-1/8+q_a/2},\qquad
 \int_{t/2}^t\Vert\nabla(\Pi,\psi)\Vert^2\dd r
 \le\mathcal X^2\tau_t^{q_a}.
\]
Split the first convolution below according to $t-r\le\tau_t^{1/4}$.
Use the pointwise bound on this part and Cauchy--Schwarz on its complement:
\[
\begin{aligned}
 \int_{t/2}^tK\Vert\nabla(\Pi,\psi)\Vert\dd r
 &\le C\mathcal X\Bigl(\tau_t^{-1/8+q_a/2+(1-\varsigma)/4}
       +\tau_t^{q_a/2+(1/2-\varsigma)/4}\Bigr)
 \le C\mathcal X\tau_t^{-d_*},\\
 \int_{t/2}^tK\Vert\nabla(\Pi,\psi)\Vert^2\dd r
 &\le C\mathcal X\tau_t^{-1/8+q_a/2}
       \int_{t/2}^tK\Vert\nabla(\Pi,\psi)\Vert\dd r
 \le C\mathcal X^2\tau_t^{-(d_*+1/8-q_a/2)},\\
 \int_{t/2}^tK\Vert F_c\Vert\dd r
 &\le C_M\delta\mathcal X\tau_t^{-1/2+q_a/2}\int_{t/2}^tK\dd r
 \le C_M\delta\mathcal X\tau_t^{-(\varsigma-1/2-q_a/2)},\\
 \int_{t/2}^tK\Vert F_R\Vert\dd r
 &\le C_b\mathcal X^{1/2}\tau_t^{q_a/4}\Vert c_R\Vert_{L^2(s,t)}
       \Bigl(\int_{t/2}^tK^2\Vert\nabla(\Pi,\psi)\Vert\dd r\Bigr)^{1/2}
 \le C_b\mathcal X\tau_t^{-(1/16-q_a/2)}.
\end{aligned}
\]
Here $1/2<\varsigma<1$, so $\int_{t/2}^tK^2\dd r\le C$, and both powers
in the first line equal $-d_*$. The last line has strictly faster
decay, since $1/16-q_a/2=15/256>d_*$.
The initial term in
\eqref{e:dual} is at most $C_bI_s^{1/2}K(t,s)
\le C_{b,s}I_s^{1/2}\tau_t^{-(\varsigma-1/64)}$.
The earlier-time exponents are also larger than $d_*$, since
$0<q_a<1/64$ and $11/16<\varsigma<3/4$. If
$s\le t<2s+2$, then $K\le C_s$, and the same source bounds give
$C_{b,s}(I_s^{1/2}+\mathcal X+\mathcal X^2)$. This proves \eqref{e:thout}.
\end{proof}

\subsection{Decay and global existence}\label{s:glob}

\begin{lemma}\label{l:dec}
Let $S\ge S_C$, $0\le\eta_P<\tfrac12(d_*-q_a)$ and $C_G\ge1$, and suppose that every
solution of \eqref{eq:NSF} on $\R^3\times[0,T_e]$, $S<T_e<\infty$, with $q$ in
\eqref{eq:class}, $\rho,\theta$ bounded below by positive constants and
$\mu_*\le\mu_0$ on $[S,T_e]$, satisfies
\[
 N_o(S,t)^2\le C_GI_S\Bigl(\frac{\tau_t}{\tau_S}\Bigr)^{\eta_P}\qquad(S\le t<T_e).
\]
Then there are $\epsilon_0>0$ and $C$ such that for every
$q_0\in H^2(\R^3)$ with $\epsilon=\Vert q_0\Vert_{H^2}\le\epsilon_0$,
\eqref{eq:NSF} with $q(0)=q_0$ has a solution on $\R^3\times[0,\infty)$ with
$q$ in \eqref{eq:class} and $\rho,\theta$ bounded below by positive constants
on every $[0,T]$, and, for $t\ge0$,
\begin{equation}\label{e:decg}
 \ell_{\rm low}(t)+\Vert q(t)\Vert_{H^2_{\rm uloc}}+\Vert(q,\Pi)(t)\Vert_{L^\infty}
 \le C\epsilon\tau_t^{-(d_*-\eta_P)},\qquad
 \Vert(\Pi,\psi)(t)\Vert_{L^3}\le C\epsilon .
\end{equation}
\end{lemma}

\begin{proof}
Let $\epsilon\le\epsilon_{b,S+1}$ of Lemma~\ref{l:fp}. The solution on
$[0,S+1]$ given there satisfies \eqref{e:fp}, and
$\Vert\Pi\Vert_{H^2}\le C_{b,S}\Vert q\Vert_{H^2}$ by \eqref{e:pmul}. As
$E_*\le C_{b,S}\Vert q\Vert_{H^2}^2$,
$D_*\le C_{b,S}(\Vert q\Vert_{H^2}^2+\Vert\nabla(\psi,\zeta)\Vert_{H^2}^2)$ and
$w_{S_C}\le1$,
\begin{equation}\label{e:bridge}
 I_S+I_{S_C}+N_o(S_C,S)^2\le C_{b,S}\epsilon^2,\qquad
 \mu_*\le C_{b,S}\epsilon\ \text{ on }[0,S+1].
\end{equation}
Let $T^*$ be the supremum of all $T>S$ such that the solution extends to
$[0,T]$ with $q$ in \eqref{eq:class}, $\rho,\theta$ bounded below by positive
constants and $\mu_*\le\mu_0$ on $[S,T]$; by Lemma~\ref{l:uniq} the
extensions coincide where both are defined. By \eqref{e:bridge},
$T^*\ge S+1$ if $C_{b,S}\epsilon\le\mu_0$.

Let $S<T_e<T^*$ and $S\le t<T_e$. By the hypothesis and \eqref{e:bridge},
$N_o(S,t)^2\le C_GC_{b,S}\epsilon^2\tau_t^{\eta_P}$. Since
$w_{S_C}(r)=(\frac{\tau_{S_C}}{\tau_S})^{q_a}w_S(r)$ for $r\ge S$, splitting the
supremum and the integrals at $S$ gives
\begin{equation}\label{e:wcmp}
 N_o(S_C,t)^2\le N_o(S_C,S)^2+\Bigl(\frac{\tau_{S_C}}{\tau_S}\Bigr)^{q_a}N_o(S,t)^2,
\end{equation}
so that $I_{S_C}+N_o(S_C,t)^2\le C_b\epsilon^2\tau_t^{\eta_P}$. By \eqref{e:bridge},
$\mu_*\le\mu_0$ on $[S_C,T_e]$, and Lemma~\ref{l:age} gives
$\tau_r^{\frac14-q_a}\Vert\nabla(\Pi,\psi)(r)\Vert^2\le\tau_{S_C}^{-q_a}\omega(r)\Vert\nabla(\Pi,\psi)(r)\Vert^2
\le C_b\epsilon^2\tau_t^{\eta_P}$ for $S\le r\le t$. Hence, with a suitable
$C_b$,
\begin{equation}\label{e:ampb}
 \mathcal X^2=C_b\epsilon^2\tau_t^{\eta_P}:\qquad
 N_o(S,t)\le\mathcal X,\qquad
 \tau_r^{\frac14-q_a}\Vert\nabla(\Pi,\psi)(r)\Vert^2\le\mathcal X^2\quad(S\le r\le t),
\end{equation}
and Lemma~\ref{l:adj} with $s=S$ and $\epsilon\le1$ gives
$\sup_{\Vert J\Vert_{L^\infty}\le1}\Vert P_1(J\frac\zeta{\theta_b})(t)\Vert_{L^\infty}
\le C_b\epsilon\tau_t^{-(d_*-\eta_P)}$. Since $K_1\in L^{\frac65}$,
\begin{equation}\label{e:mechv}
 \Vert P_1\psi(t)\Vert_{L^\infty}+\Bigl\Vert P_1\frac{\Pi}{R\theta_b}(t)\Bigr\Vert_{L^\infty}
 \le C_b\Vert(\psi,\Pi)(t)\Vert_{L^6}\le C_b\Vert\nabla(\Pi,\psi)(t)\Vert
 \le C_b\epsilon\tau_t^{\frac{\eta_P}2-\frac18+\frac{q_a}2} .
\end{equation}
By \eqref{e:state}, $\zeta=\theta_b\frac\zeta{\theta_b}$ and
$\phi=\frac\Pi{R\theta_b}-\rho\frac\zeta{\theta_b}$; we apply the preceding bound
with $J=\frac{\theta_b}{\Vert\theta_b\Vert_{L^\infty}}$ and
$J=\frac\rho{\Vert\rho\Vert_{L^\infty}}$. With
$\frac18-\frac{q_a}2-\frac{\eta_P}2>d_*-\eta_P$, this gives
$\ell_{\rm low}(r)\le C_b\epsilon\tau_r^{-(d_*-\eta_P)}$ for $S\le r<T_e$.
Lemma~\ref{l:damp}(a) with $s=S$, $d=d_*-\eta_P$ and $K=C_b\epsilon$, and
$\Pi=R(\theta_b\phi+\rho_b\zeta+\phi\zeta)$, give the first bound in
\eqref{e:decg} for $S\le t<T_e$. For $F\in\{\Pi,\psi\}$,
$\Vert F(t)\Vert^2\le E_*(t)\le(\frac{\tau_t}{\tau_S})^{q_a}N_o(S,t)^2
\le C_b\epsilon^2\tau_t^{q_a+\eta_P}$, so
\[
 \Vert F(t)\Vert_{L^3}\le\Vert F(t)\Vert^{\frac23}\Vert F(t)\Vert_{L^\infty}^{\frac13}
 \le C_b\epsilon\tau_t^{-\frac13(d_*-q_a-2\eta_P)}\le C_b\epsilon,
\]
since $d_*-q_a-2\eta_P>0$ by hypothesis. By continuity of
$q$ in $H^2$ these bounds hold at $t=T_e$, and with \eqref{e:mu} they give
$\mu_*\le C_b\epsilon$ on $[S,T_e]$, with $C_b$ independent of $T_e$.

Let $\epsilon_0\le\min\{1,\epsilon_{b,S+1}\}$ with $C_b\epsilon_0\le\frac{\mu_0}2$
and $C_{b,S}\epsilon_0\le\frac{\mu_0}2$. If $T^*<\infty$, then
$\mu_*\le\frac{\mu_0}2$ on every $[S,T']$, $T'<T^*$, and \eqref{e:stb} gives
positive lower bounds of $\rho,\theta$ on $[S,T^*)$; by
Lemma~\ref{l:damp}(c) the solution extends to $[0,T^*+\ell]$, and by
continuity in $H^2$, $\mu_*\le\mu_0$ on $[S,T^*+\ell']$ for some
$\ell'\in(0,\ell]$, which contradicts the definition of $T^*$. Thus
$T^*=\infty$, and \eqref{e:decg} holds for $t\ge S$. For $0\le t\le S$,
\eqref{e:fp}, $H^2\subset L^\infty$, $H^1\subset L^3$ and
$\Vert\Pi\Vert_{H^2}\le C_{b,S}\Vert q\Vert_{H^2}$ give \eqref{e:decg} with
$C=C_{b,S}\tau_S^{d_*}$.
\end{proof}

\section{Proof of the main theorem}\label{sec:proof}

For each fixed datum choose an auxiliary $M$ as in
\eqref{e:dat-aux}. We verify that every wave-strength restriction
depends only on $A$ and the fixed data. Throughout this section
$C_A$ has only this dependence, whereas $C_b$ and the entrance times
may depend on the complete datum.

\begin{proof}[Proof of Theorem~\ref{thm:main}]
We use Lemma~\ref{l:background} and the uniformity of the imported
background estimates in Section~\ref{s:bg-interface}. Increase the
entrance time so that the bounds in \eqref{e:E6} are at most one.

For the three-dimensional estimates, the state-box constants of
Lemma~\ref{l:st} are now uniform. The acoustic gap \eqref{e:gap}
also has a uniform positive lower bound, since
$e^{-\alpha|s_\sigma m_\sigma|}\ge\beta_A>0$ on the bounded mass
range. Hence the acoustic recovery estimates in Lemmas~\ref{l:cor}
and~\ref{l:sp}, and \eqref{e:coer}, have constants depending only
on $A$ and the fixed test parameters. The source tests require one
distinction: the first line of \eqref{e:srct} may have constant
$K_\infty=C_{b,*}$, whereas the second has $K_E=C_{A,*}$.
The proof of Lemma~\ref{l:tw}(iii) places $K_\infty$ only with
$\mu_*$ or a vanishing tail; its coefficient of $\delta$ uses
$K_E$. For the contact-pressure term, integration by parts and
\eqref{e:eta}, \eqref{e:E1Z} give explicitly
\[
\begin{aligned}
 &\int_s^tw_s\left|\int_{\R^3}\gamma p_b\partial_1(b\eta)T_p\dd x\right|\dd r\\
 &\quad\le C_A\Vert b\eta\Vert_{L^2_w}
       \left(\Vert\nabla T_p\Vert_{L^2_w}
                       +\Vert(\partial_1p_b)T_p\Vert_{L^2_w}\right)\\
 &\quad\le\bigl(C_AK_E\delta+\omega_{b,K_E}(s)\bigr)\mathfrak B(s,t).
\end{aligned}
\]
This is also the separation used in
Lemmas~\ref{l:low} and~\ref{l:mix}.
Thus the error coefficients in the weighted closure have the form
\[
 C_{A,*}(\delta+\delta^2)+C_{b,*}\mu_*+\omega_{b,*}(s).
\]
Here the star denotes the frequency-test parameters chosen in the
order at the end of the proof of Lemma~\ref{l:cl}.

In particular, \eqref{e:coer} gives $C_R=C_R(A,\theta_a)$ and
$C_J=6C_R$. Choose $\epsilon_A=(4C_J)^{-1}$. The proof of
Lemma~\ref{l:cw} gives $c_1=c_1(A)$, $c_2=c_2(A)$, since its free
Young parameter is fixed from $\epsilon_A$ and the preceding
$C_A$ constants. Therefore
\[
 C_1=2C_Jc_1(A),\qquad C_2=2C_Jc_2(A)
\]
are uniform in the higher initial norms. The final three absorptions
in Lemma~\ref{l:cl} are precisely
\[
\begin{gathered}
 9C_R\epsilon_{\rm rec}\le\frac18,\qquad
 3C_RC_{A,*}(\delta+\delta^2)\le\frac18,\\
 3C_R\bigl(C_{b,*}\mu_*+\omega_{b,*}(s)\bigr)+\omega_b(s)\le\frac14.
\end{gathered}
\]
The first two determine low-data parameters and a wave threshold;
the last is obtained after the datum is fixed by increasing $s$
and decreasing $\mu_*$. Lemmas~\ref{l:age} and~\ref{l:damp} use
the same state-box coercivity and vanishing higher derivative norms.
The additional wave restriction in Lemma~\ref{l:adj} involves only
the contact Gaussian constant and the fixed data. The fixed state-box
restriction in Lemma~\ref{l:low-damp} is imposed here as well; it does
not involve the three-dimensional bounds. Hence all the
preceding estimates hold for a threshold $\delta_{\rm u}(A)>0$,
with background-dependent entrance times and finite constants.
Before fixing $\delta_{\rm u}(A)$, take $a_c=\alpha_G/2$ and
$a_s$ so that
\begin{equation}\label{e:gp}
 0<a_s\le\min\Bigl\{a_c,\frac{c_G}4\Bigr\},\qquad
 A^{a_c}_f\le A^{a_s}_f,\quad
 \mathscr N^{a_c}_{r,T}(f)\le\mathscr N^{a_s}_{r,T}(f),
\end{equation}
and then $0<b_0<a_s/64$.
Put $\eta_P=1/1024$ and choose
$\delta_0(A)\le\min\{\delta_{\rm u}(A),1,\theta_-/2\}$,
decreasing the threshold of Lemma~\ref{l:background} if necessary, so that
\begin{equation}\label{e:del3}
 4C_2(A)C_P\delta_0(A)^4\le\frac{\eta_P}{2}.
\end{equation}
Fix the datum and $0<\delta\le\delta_0(A)$. Take $r$ and
$P_b-p_*=g_1+g_2$ from Lemma~\ref{l:prs}, with
$\xi=\eta_P/(2C_2(A)\delta^2)$, then $S_C$ from Lemma~\ref{l:age},
and $S\ge\max\{r,S_C,S_W\}$ above the preceding entrance times with
\[
 C_1(A)\delta\,\mathscr N^{a_s}_{S,\infty}(g_1)\le\frac12.
\]
This is possible since $\mathscr N^{a_s}_{r,\infty}(g_1)<\infty$.
Choose $\mu_0\le\mu_b$ below the preceding smallness thresholds.

Let $S<T_e<\infty$ and let a solution on $[0,T_e]$ have
$q$ in \eqref{eq:class}, positive lower bounds for $\rho,\theta$ and
$\mu_*\le\mu_0$ on $[S,T_e]$. By \eqref{e:gp},
$C_1(A)\delta\mathscr N^{a_c}_{S,T_e}(g_1)\le1/2$.
Moreover, \eqref{e:prlog}, \eqref{e:del3} and
$(2+t)(1+S)\le(1+t)(2+S)$ give, for $S\le t\le T_e$,
\begin{equation}\label{e:pbapp}
\begin{aligned}
 \int_S^tC_2(A)\delta^2A^{a_c}_{g_2}\dd\ell
 &\le C_2(A)\delta^2\mathscr O^{a_s}_{r,t}(g_2)\\
 &\le C_2(A)\delta^2C_{b,r,\xi}
     +C_2(A)\delta^2(4C_P\delta^2+\xi)\log\frac{2+t}{2+r}\\
 &\le C_{P,b}+\eta_P\log\frac{\tau_t}{\tau_S},\\
 C_{P,b}&=C_2(A)\delta^2C_{b,r,\xi}
           +\eta_P\log\frac{2+S}{2+r}.
\end{aligned}
\end{equation}
In particular $A^{a_c}_{g_2}\in L^1(S,T_e)$, so \eqref{e:pc} and
\eqref{e:bpbd} hold with the chosen $\eta_P=1/1024$. Lemma~\ref{l:cl} yields
\[
 N_o(S,t)^2\le C_GI_S\Bigl(\frac{\tau_t}{\tau_S}\Bigr)^{\eta_P},
 \qquad C_G=\max\{1,C_{0,b}e^{C_{P,b}}\},
\]
uniformly in $T_e$. Since $d_*-q_a-2\eta_P=47/1024>0$, Lemmas~\ref{l:fp}
and~\ref{l:dec} first give the global solution and \eqref{e:r3}
with a background constant $C_b$ whenever
$\Vert q_0\Vert_{H^2}\le\epsilon_b^{(2)}$, for some
$\epsilon_b^{(2)}>0$. Here
$d_*-\eta_P=7/128>1/20$.

We now only assume \eqref{e:3small}. Lemma~\ref{l:thermal-layer}
gives a positive strong solution on $[0,1]$ with
\[
 \Vert q(1)\Vert_{H^2}\le K_1,\qquad
 \sup_{0\le t\le1}\Vert q(t)\Vert\le B_1\Vert q_0\Vert,
\]
where $K_1,B_1$ depend only on the background and the five prescribed
bounds, independently of $\Vert\nabla^2\zeta_0\Vert$.
Increase a fixed background time $S_0>\max\{S,1\}$ so that
Lemmas~\ref{l:low-damp} and~\ref{l:late-restart} apply.
The time-translated version of Lemma~\ref{l:low-finite} on $[1,S_0]$
gives a bound $K_0=K_{b,1,S_0}(K_1)$ when the original $L^2$ norm
is sufficiently small. Let $c_b(1+t)^{-m_b}$ be the admissible size in
Lemma~\ref{l:late-restart}. Choose $T>S_0$ so that
\[
 C_be^{-c(T-S_0)}K_0^2
 \le\frac14c_b^2(1+T)^{-2m_b}.
\]
This is possible because the left side decays exponentially.
Apply Lemma~\ref{l:low-finite} on $[1,T]$, obtaining
$K_T=K_{b,1,T}(K_1)$ and, with the initial-layer estimate,
\[
 \sup_{0\le t\le T}\Vert q(t)\Vert
 \le B_T\Vert q_0\Vert=:h.
\]
Decrease the initial $L^2$ threshold, after $T$ and $K_T$ have been
fixed, so that
\[
 h\le h_b(K_T),\qquad
 C_{b,K_T}h^{2/7}
 \le\frac14c_b^2(1+T)^{-2m_b}.
\]
The two finite-interval solutions coincide by Lemma~\ref{l:uniq},
so the bound $K_0$ at $S_0$ remains valid. Lemma~\ref{l:low-damp}
on $[S_0,T]$ now gives
\[
 \Vert q(T)\Vert_{H^2}^2
 \le\frac12c_b^2(1+T)^{-2m_b}.
\]
Lemma~\ref{l:late-restart} extends this positive strong solution
globally and gives the exponent $1/20$.

The constants in this construction are uniform for the class
\eqref{e:3small}, hence
\[
 \Vert(q,\Pi)(t)\Vert_{L^\infty}
 \le D(1+t)^{-1/20}\qquad(t\ge0).
\]
Here the state bounds and reduced energy in
Lemma~\ref{l:thermal-layer} cover $0\le t\le1$, the full $H^2$ bound
covers $1\le t\le T$, and the restarted estimate covers $t\ge T$.
Thus $D$ has the stated dependencies. If
$\epsilon=\Vert q_0\Vert_{H^2}\le\epsilon_b^{(2)}$, the earlier
small-$H^2$ argument gives the factor $C_b\epsilon$. Otherwise,
$D\le D\epsilon/\epsilon_b^{(2)}$.
This proves \eqref{e:r3} in both cases.
Lemma~\ref{l:uniq} gives uniqueness,
$\phi_t\in C([0,T_e];H^1)$ and \eqref{eq:far}.
Lemma~\ref{l:background} and the triangle inequality give \eqref{e:rt}.
Every choice involving the three-dimensional bounds is made after
$\delta_0(A)$; these choices affect $\epsilon_0$, the entrance times
and the finite constants, but not the wave-strength threshold.

\end{proof}

\appendix
\section{Local strong solutions}\label{app:local}

The local theory is proved under the background assumptions \eqref{e:locb}
stated in Section~\ref{sec:pre}.

\begin{lemma}\label{l:local}
Let $0<T_f<\infty$, and let $(\rho_b,u_b,\theta_b)$ solve
\eqref{eq:bg-eqs} on $\R\times[0,T_f]$ with \eqref{e:locb}.
\begin{enumerate}
\item[(i)] For $K_0,\underline c_0>0$ there are $h_0\in(0,\frac12]$ and
$C_{\rm loc}$, depending only on $K_0$, $\underline c_0$, $T_f$,
$(R,\gamma,\mu,\lambda,\kappa)$ and the background, such that for every
$t_0\in[0,T_f)$ and every
$q^{t_0}=(\phi^{t_0},\psi^{t_0},\zeta^{t_0})\in H^2(\R^3)$ with
\[
 \Vert q^{t_0}\Vert_{H^2}\le K_0,\qquad
 \rho_b(\cdot,t_0)+\phi^{t_0}\ge\underline c_0,\qquad
 \theta_b(\cdot,t_0)+\zeta^{t_0}\ge\underline c_0,
\]
\eqref{eq:NSF} has a unique strong solution on $\R^3\times I$,
$I=[t_0,\min\{t_0+h_0,T_f\}]$, with $q(t_0)=q^{t_0}$, and
\begin{equation}\label{e:loc}
 \phi_t\in C(I;H^1),\qquad \rho,\theta\ge\frac{\underline c_0}{2},\qquad
 \sup_{r\in I}\Vert q(r)\Vert_{H^2}^2
    +\int_I\bigl(\Vert(\psi,\zeta)\Vert_{H^3}^2
                         +\Vert q_t\Vert_{H^1}^2\bigr)\dd r\le C_{\rm loc}.
\end{equation}
\item[(ii)] Let $0<T_{\rm fin}<T_f$, and let $(\rho,u,\theta)$ be a
strong solution of \eqref{eq:NSF} on $\R^3\times[0,T]$ for every
$T<T_{\rm fin}$. If
$\sup_{0\le t<T_{\rm fin}}\Vert q(t)\Vert_{H^2}<\infty$ and
$\inf_{\R^3\times[0,T_{\rm fin})}\min\{\rho,\theta\}>0$, then for some
$\ell>0$ the solution extends to a strong solution on
$\R^3\times[0,T_{\rm fin}+\ell]$.
\end{enumerate}
\end{lemma}

Part (i) is proved below; uniqueness and part (ii) follow from
Lemmas~\ref{l:uniq} and~\ref{l:cont}, respectively.

Let the background satisfy \eqref{e:locb} on $[0,T_f]$. The derivation of
\eqref{eq:ph-phi}, \eqref{eq:ph-psi} and \eqref{e:peq3} uses only
\eqref{eq:bg-eqs}, in the form \eqref{e:bal}, for the background. Hence,
for $q$ with $\rho,\theta>0$, \eqref{eq:NSF} is equivalent to
\begin{equation}\label{e:locq}
 D_u\phi=-\rho\Div\psi-h_b\phi-(\partial_1\rho_b)\psi_1,\qquad
 \psi_t-\frac1\rho\mathcal L\psi=F_\psi,\qquad
 \zeta_t-\frac\kappa{c_v\rho}\Delta\zeta=F_\zeta,
\end{equation}
where $\rho=\rho_b+\phi$, $\theta=\theta_b+\zeta$, $u=u_be_1+\psi$,
$\Pi=R(\theta\phi+\rho_b\zeta)$ and, with $f_m,f_z$ of \eqref{e:ps3},
\[
 F_\psi=\frac{f_m-\nabla\Pi}{\rho}-u\cdot\nabla\psi,\qquad
 F_\zeta=f_z-(\gamma-1)\theta\Div\psi-u\cdot\nabla\zeta .
\]
By \eqref{e:bal},
$A_b=\frac{\nu\partial_1h_b-\partial_1p_b}{\rho_b}$ and
$B_b=\frac{\kappa\partial_1^2\theta_b+\nu h_b^2-p_bh_b}{c_v\rho_b}$, so
\eqref{e:locb} gives $(A_b,B_b)\in L^2(0,T_f;H^1(\R))$ and
\[
 k_{\rm loc}=1+\Vert\partial_1(u_b,\theta_b)\Vert_{H^2(\R)}^2
       +\Vert(A_b,B_b)\Vert_{H^1(\R)}^2\in L^1(0,T_f).
\]
For $|\alpha|\le2$, the first equation of \eqref{e:locq} gives
$D_u\partial^\alpha\phi=G_\alpha$ with
\[
 G_\alpha=-\rho\partial^\alpha\Div\psi-[\partial^\alpha,\rho]\Div\psi
  -[\partial^\alpha,u\cdot\nabla]\phi
  -\partial^\alpha\bigl(h_b\phi+(\partial_1\rho_b)\psi_1\bigr).
\]
Let $[t_1,t_2]\subset[0,T_f]$, $q\in L^\infty(t_1,t_2;H^2)$,
$(\psi,\zeta)\in L^2(t_1,t_2;H^3)$, $q_t\in L^2(t_1,t_2;H^1)$ and
$\inf\rho>0$. Then $q$ is weakly continuous into $H^2$,
$\nabla u=h_be_1\otimes e_1+\nabla\psi\in L^2(t_1,t_2;L^\infty)$, and
\eqref{e:pmul}, \eqref{eq:planar} for $(\partial_1^3\rho_b)\psi_1$ and
$(\partial_1^2h_b)\phi$, and \eqref{e:locb} give
$G_\alpha\in L^2(t_1,t_2;L^2)$. For a mollifier $j_\varepsilon$ in $x$,
$f^\varepsilon=j_\varepsilon*\partial^\alpha\phi$ satisfies
$D_uf^\varepsilon=j_\varepsilon*G_\alpha+r_\varepsilon$ with
$r_\varepsilon=u\cdot\nabla f^\varepsilon
-j_\varepsilon*(u\cdot\nabla\partial^\alpha\phi)$,
$\Vert r_\varepsilon\Vert_{L^2}\le C\Vert\nabla u\Vert_{L^\infty}
\Vert\partial^\alpha\phi\Vert_{L^2}$ and $r_\varepsilon\to0$ in
$L^1(t_1,t_2;L^2)$. The energy identity for $f^\varepsilon$ and
$\varepsilon\to0$ give, for $t_1\le r\le t\le t_2$,
\begin{equation}\label{e:phid}
 \Vert\partial^\alpha\phi(t)\Vert_{L^2}^2-\Vert\partial^\alpha\phi(r)\Vert_{L^2}^2
 =\int_r^t\!\!\int_{\R^3}\bigl(\Div u\,|\partial^\alpha\phi|^2
       +2G_\alpha\partial^\alpha\phi\bigr)\dd x\dd v .
\end{equation}

\begin{proof}[Proof of Lemma~\ref{l:uniq}]
Let $q$ be one of the solutions. By the first equation of \eqref{e:locq},
\[
 \phi_t=-u\cdot\nabla\phi-\rho\Div\psi-h_b\phi-(\partial_1\rho_b)\psi_1 ;
\]
\eqref{e:pmul} with $j=1$,
$\Vert FG\Vert_{H^1}\le C\Vert F\Vert_{H^2}\Vert G\Vert_{H^1}$,
\eqref{e:locb} and $q\in C([t_1,t_2];H^2)$ make every term continuous
into $H^1$, so $\phi_t\in C([t_1,t_2];H^1)$. Since
$H^2(\R^3)\subset C_0(\R^3)$, $\sup_{x'}|q(x_1,x',t)|\to0$ as
$x_1\to\pm\infty$; adding the end states of the background gives
\eqref{eq:far}.

Let $q^1,q^2$ be two solutions with $q^1(t_1)=q^2(t_1)$, and
$d=(d_\phi,d_\psi,d_\zeta)=q^1-q^2$. In this proof $C$ depends on
$\sup_{[t_1,t_2]}\Vert q^i\Vert_{H^2}$, the lower bounds of
$\rho^i,\theta^i$ and the background. By \eqref{e:locq},
\[
\begin{aligned}
 \partial_td_\phi+u^1\cdot\nabla d_\phi
 &=-d_\psi\cdot\nabla\phi^2-\rho^1\Div d_\psi-d_\phi\Div\psi^2
   -h_bd_\phi-(\partial_1\rho_b)d_{\psi,1},\\
 \partial_td_\psi-\frac1{\rho^1}\mathcal Ld_\psi
 &=\Bigl(\frac1{\rho^1}-\frac1{\rho^2}\Bigr)\mathcal L\psi^2
   +F_\psi^1-F_\psi^2,\\
 \partial_td_\zeta-\frac\kappa{c_v\rho^1}\Delta d_\zeta
 &=\frac\kappa{c_v}\Bigl(\frac1{\rho^1}-\frac1{\rho^2}\Bigr)\Delta\zeta^2
   +F_\zeta^1-F_\zeta^2 .
\end{aligned}
\]
Test the three equations with $d_\phi$, $d_\psi$, $d_\zeta$; this is
justified by $d_\phi,\partial_td_\phi\in C([t_1,t_2];H^1)$ and
$(d_\psi,d_\zeta)\in L^2(t_1,t_2;H^1)$,
$\partial_t(d_\psi,d_\zeta)\in L^2(t_1,t_2;H^{-1})$. The diffusion terms
give $C^{-1}\Vert\nabla(d_\psi,d_\zeta)\Vert_{L^2}^2$ up to
$\int_{\R^3}|\nabla(\rho^1)^{-1}||\nabla(d_\psi,d_\zeta)||(d_\psi,d_\zeta)|\dd x$,
and $\nabla(\rho^1)^{-1}=-(\rho^1)^{-2}((\partial_1\rho_b)e_1+\nabla\phi^1)$
with $\Vert\nabla\phi^1\Vert_{L^6}\Vert(d_\psi,d_\zeta)\Vert_{L^3}
\le C\Vert(d_\psi,d_\zeta)\Vert_{L^2}^{\frac12}
\Vert\nabla(d_\psi,d_\zeta)\Vert_{L^2}^{\frac12}$. The transport terms
are bounded by
$C(\Vert h_b\Vert_{L^\infty}+\Vert\nabla\psi^1\Vert_{L^\infty})\Vert d\Vert_{L^2}^2$,
and $d_\psi\cdot\nabla\phi^2$ by
$\Vert d_\psi\Vert_{L^6}\Vert\nabla\phi^2\Vert_{L^3}$. The reciprocal
differences satisfy $|\frac1{\rho^1}-\frac1{\rho^2}|\le C|d_\phi|$, and
their products with $\nabla^2\psi^2$, $\nabla^2\zeta^2$ and $\nabla\Pi^2$
are bounded in $L^2\cdot L^3\cdot L^6$ by
$C\Vert d_\phi\Vert_{L^2}\Vert(\nabla^2\psi^2,\nabla^2\zeta^2,\nabla\Pi^2)
\Vert_{L^3}\Vert(d_\psi,d_\zeta)\Vert_{L^6}$. In the pressure term,
$d_\Pi=R(\theta^1d_\phi+\rho^2d_\zeta)$ and
$\int_{\R^3}\frac{\nabla d_\Pi}{\rho^1}\cdot d_\psi\dd x
=-\int_{\R^3}d_\Pi\Div\frac{d_\psi}{\rho^1}\dd x$. In the acceleration
terms, $\frac{\phi^1}{\rho^1}-\frac{\phi^2}{\rho^2}
=\frac{\rho_bd_\phi}{\rho^1\rho^2}$ and
$\Vert(A_b,B_b)\Vert_{L^\infty}\le C\Vert(A_b,B_b)\Vert_{H^1(\R)}$. In the
heating term,
$|Q_\mu(\nabla\psi^1)-Q_\mu(\nabla\psi^2)|
\le C(|\nabla\psi^1|+|\nabla\psi^2|)|\nabla d_\psi|$. The other terms
are linear in $d$ with coefficients bounded by $C$. Young's inequality
absorbs the gradients and gives
\[
\begin{gathered}
 \frac{\dd}{\dd t}\Vert d\Vert_{L^2}^2
   +C^{-1}\Vert\nabla(d_\psi,d_\zeta)\Vert_{L^2}^2
 \le k\Vert d\Vert_{L^2}^2,\\
 k=C\Bigl[1+\Vert(A_b,B_b)\Vert_{H^1(\R)}
   +\sum_{i=1}^2\bigl(\Vert\nabla(\psi^i,\zeta^i)\Vert_{L^\infty}^2
                  +\Vert\nabla^2(\psi^i,\zeta^i)\Vert_{L^3}^2\bigr)\Bigr],
\end{gathered}
\]
with $k\in L^1(t_1,t_2)$ by \eqref{eq:class}. Gronwall's inequality and
$d(t_1)=0$ give $d=0$.
\end{proof}

Let $K\ge1$ and $m>0$. Suppose $\Vert q\Vert_{H^2}\le K$ and
$\rho,\theta\ge\frac m4$, and let
$C$ depend on $K$, $m$, $(R,\gamma,\mu,\lambda,\kappa)$ and the
background. By \eqref{e:pmul}, \eqref{e:locb}, the algebra property of
$H^2$ and
$\Vert\nabla f\Vert_{L^\infty}+\Vert\nabla^2f\Vert_{L^3}
\le C\Vert\nabla f\Vert_{H^1}^{\frac12}\Vert\nabla f\Vert_{H^2}^{\frac12}$,
\begin{equation}\label{e:locF}
\begin{gathered}
 \Vert\rho^{-1}-\rho_b^{-1}\Vert_{H^2}+\Vert\Pi\Vert_{H^2}\le C\Vert q\Vert_{H^2},\qquad
 \Vert\nabla(\psi,\zeta)\Vert_{L^\infty}+\Vert\nabla^2(\psi,\zeta)\Vert_{L^3}
 \le C\Vert q\Vert_{H^2}^{\frac12}\Vert(\psi,\zeta)\Vert_{H^3}^{\frac12},\\
 \Vert F_\psi\Vert_{H^1}+\Vert F_\zeta\Vert_{H^1}
 \le C\bigl(k_{\rm loc}^{\frac12}\Vert q\Vert_{H^2}+\Vert q\Vert_{H^2}^{\frac32}\Vert(\psi,\zeta)\Vert_{H^3}^{\frac12}\bigr);
\end{gathered}
\end{equation}
here $\frac\phi\rho A_b$ and $\frac\phi\rho B_b$ are bounded in $H^1$ only,
by \eqref{e:pmul} with $j=1$, and
$\Vert Q_\mu(\nabla\psi)\Vert_{H^1}\le C\Vert\nabla\psi\Vert_{L^\infty}
\Vert\nabla\psi\Vert_{H^1}$. Test the second and third equations of
\eqref{e:locq} with $\psi+\mathcal L(\mathcal L\psi)$ and
$\zeta+\Delta^2\zeta$, and $D_u\partial^\alpha\phi=G_\alpha$,
$|\alpha|\le2$, with $\partial^\alpha\phi$. With
\[
 E_{\rm loc}=\Vert\phi\Vert_{H^2}^2+\Vert\psi\Vert_{L^2}^2
   +\Vert\mathcal L\psi\Vert_{L^2}^2+\Vert\zeta\Vert_{L^2}^2
   +\Vert\Delta\zeta\Vert_{L^2}^2,
\]
$C_*^{-1}\Vert q\Vert_{H^2}^2\le E_{\rm loc}\le C_*\Vert q\Vert_{H^2}^2$ and
$\Vert(\psi,\zeta)\Vert_{H^3}^2\le C_*(\Vert q\Vert_{H^2}^2+\Vert\nabla\mathcal L\psi\Vert_{L^2}^2
+\Vert\nabla\Delta\zeta\Vert_{L^2}^2)$ with $C_*$ depending only on
$\mu,\lambda$. The highest velocity pairing is
\[
\begin{aligned}
 \int_{\R^3}\frac{\mathcal L\psi}{\rho}\cdot\mathcal L(\mathcal L\psi)\dd x
 ={}&-\int_{\R^3}\frac{\mu|\nabla\mathcal L\psi|^2
                 +(\mu+\lambda)(\Div\mathcal L\psi)^2}{\rho}\dd x\\
 &-\int_{\R^3}\Bigl[\mu\sum_{k,l}(\mathcal L\psi)_k\,\partial_l\rho^{-1}\,
              \partial_l(\mathcal L\psi)_k
   +(\mu+\lambda)(\nabla\rho^{-1}\cdot\mathcal L\psi)
                       \Div\mathcal L\psi\Bigr]\dd x,
\end{aligned}
\]
and
\[
 \int_{\R^3}\frac{\Delta\zeta}{\rho}\,\Delta^2\zeta\dd x
 =-\int_{\R^3}\frac{|\nabla\Delta\zeta|^2}{\rho}\dd x
  -\int_{\R^3}\Delta\zeta\,\nabla\rho^{-1}\cdot\nabla\Delta\zeta\dd x .
\]
Since
$\nabla\rho^{-1}=-\rho^{-2}((\partial_1\rho_b)e_1+\nabla\phi)$ and
$\Vert\nabla\phi\Vert_{L^6}\Vert\mathcal L\psi\Vert_{L^3}
\le C\Vert q\Vert_{H^2}^{\frac32}\Vert(\psi,\zeta)\Vert_{H^3}^{\frac12}$, the diffusion terms other than the
dissipation are at most $C(\Vert q\Vert_{H^2}\Vert(\psi,\zeta)\Vert_{H^3}+\Vert q\Vert_{H^2}^{\frac32}\Vert(\psi,\zeta)\Vert_{H^3}^{\frac32})$. Moving one
derivative to $F_\psi,F_\zeta$ in the highest forcing pairings, the
forcing terms are at most $C(k_{\rm loc}^{\frac12}\Vert q\Vert_{H^2}\Vert(\psi,\zeta)\Vert_{H^3}+\Vert q\Vert_{H^2}^{\frac32}\Vert(\psi,\zeta)\Vert_{H^3}^{\frac32})$
by \eqref{e:locF}. In the density terms,
$\int_{\R^3}(\frac12\Div u|\partial^\alpha\phi|^2
+G_\alpha\partial^\alpha\phi)\dd x$, the bounds
$\Vert\Div u\Vert_{L^\infty}\le C(k_{\rm loc}^{\frac12}+\Vert q\Vert_{H^2}^{\frac12}\Vert(\psi,\zeta)\Vert_{H^3}^{\frac12})$,
\eqref{e:pmul} and \eqref{eq:planar} for $(\partial_1^3\rho_b)\psi_1$ and
$(\partial_1^2h_b)\phi$ give
$\Vert G_\alpha\Vert_{L^2}\le C(\Vert(\psi,\zeta)\Vert_{H^3}+k_{\rm loc}^{\frac12}\Vert q\Vert_{H^2}+\Vert q\Vert_{H^2}^{\frac32}\Vert(\psi,\zeta)\Vert_{H^3}^{\frac12})$
and the bound $C(k_{\rm loc}^{\frac12}\Vert q\Vert_{H^2}^2+\Vert q\Vert_{H^2}\Vert(\psi,\zeta)\Vert_{H^3}+\Vert q\Vert_{H^2}^{\frac52}\Vert(\psi,\zeta)\Vert_{H^3}^{\frac12})$. By
Young's inequality,
\begin{equation}\label{e:locen}
 E_{\rm loc}'+C^{-1}\Vert(\psi,\zeta)\Vert_{H^3}^2\le Ck_{\rm loc},\qquad
 \Vert q_t\Vert_{H^1}^2\le C(k_{\rm loc}+\Vert(\psi,\zeta)\Vert_{H^3}^2),
\end{equation}
the second by \eqref{e:locq} and \eqref{e:locF}. If $\Vert q\Vert_{H^2}\le1$ and
$\rho,\theta\ge\frac12\inf_{\R\times[0,T_f]}\min\{\rho_b,\theta_b\}$,
Young's inequality bounds each term above by
$\epsilon \Vert(\psi,\zeta)\Vert_{H^3}^2+C_\epsilon k_{\rm loc}\Vert q\Vert_{H^2}^2$, which gives
\begin{equation}\label{e:locsm}
 E_{\rm loc}'+C_b^{-1}\Vert(\psi,\zeta)\Vert_{H^3}^2\le C_bk_{\rm loc}E_{\rm loc},\qquad
 \Vert q_t\Vert_{H^1}^2\le C_b(k_{\rm loc}E_{\rm loc}+\Vert(\psi,\zeta)\Vert_{H^3}^2),
\end{equation}
with $C_b$ depending only on $(R,\gamma,\mu,\lambda,\kappa)$ and the
background.

\begin{proof}[Proof of Lemma~\ref{l:local}(i)]
Let $J_n$ be the Fourier projection onto $\{|\xi|\le n\}$, and let
$\mathcal F(t,q)$ be the right-hand side of \eqref{e:locq} solved for
$q_t$. On each set
$\{q\in J_nH^2:\Vert q\Vert_{H^2}\le K,\ \rho,\theta\ge\frac m4\}$, the map
$J_n\mathcal F(t,\cdot)$ is Lipschitz in $H^2$ with constant
$C_nk_{\rm loc}(t)^{\frac12}$; successive approximation on intervals on
which $\int k_{\rm loc}^{\frac12}$ is small gives an absolutely
continuous $q_n$ with
\[
 q_n(t)=J_nq^{t_0}+\int_{t_0}^tJ_n\mathcal F(r,q_n(r))\dd r,
\]
as long as $q_n$ stays in this set. Since $J_n$ is self-adjoint and
commutes with $\partial^\alpha$, $\mathcal L$ and $\Delta$, the pairings
above are unchanged by $J_n$, and \eqref{e:locen} holds for $q_n$ with
$C$ independent of $n$.

For $f\in H^2(\R^3)$,
$\Vert(1-J_n)f\Vert_{L^\infty}\le\int_{|\xi|>n}|\hat f|\dd\xi
\le Cn^{-\frac12}\Vert f\Vert_{H^2}$. Take $m=\underline c_0$ and
$K=2C_*(1+K_0)$. As long as $\Vert q_n\Vert_{H^2}\le K$ and
$\rho_n,\theta_n\ge\frac m4$, \eqref{e:locen}, the interpolation
$\Vert f\Vert_{L^\infty}\le C\Vert f\Vert_{H^1}^{\frac14}
\Vert f\Vert_{H^2}^{\frac34}$ and $E_{\rm loc}(t_0)\le C_*K_0^2$ give
\[
\begin{aligned}
 E_{\rm loc}(t)+C^{-1}\int_{t_0}^t\Vert(\psi,\zeta)\Vert_{H^3}^2\dd r
 &\le C_*K_0^2+C\int_{t_0}^tk_{\rm loc}\dd r,\\
 \Vert q_n(t)-q_n(t_0)\Vert_{H^1}
 &\le C|t-t_0|^{\frac12}\Bigl(1+\int_{t_0}^tk_{\rm loc}\dd r\Bigr)^{\frac12},\\
 \Vert q_n(t)-q_n(t_0)\Vert_{L^\infty}
 &\le C|t-t_0|^{\frac18}\Bigl(1+\int_{t_0}^tk_{\rm loc}\dd r\Bigr)^{\frac18}.
\end{aligned}
\]
Since $k_{\rm loc}\in L^1(0,T_f)$,
$\sup\{\int_{I'}k_{\rm loc}\dd r:I'\subset[0,T_f],\ |I'|\le h\}\to0$ as
$h\to0$, and $(\rho_b,\theta_b)$ is uniformly continuous from $[0,T_f]$
into $L^\infty(\R)$. Hence there is $h_0\in(0,\frac12]$, depending only
on $K_0,\underline c_0,T_f$, the physical parameters and the background,
such that, whenever $0\le t'<t\le T_f$ and $t-t'\le h_0$,
\[
 C\int_{t'}^tk_{\rm loc}\dd r\le1,\qquad
 2Ch_0^{\frac18}\le\frac{\underline c_0}8,\qquad
 \Vert(\rho_b,\theta_b)(t)-(\rho_b,\theta_b)(t')\Vert_{L^\infty}
 \le\frac{\underline c_0}8 .
\]
With $Cn^{-\frac12}K_0\le\frac{\underline c_0}8$, the preceding bounds give
$\Vert q_n\Vert_{H^2}^2\le C_*(C_*K_0^2+1)<K^2$ and
$\rho_n,\theta_n\ge\frac58\underline c_0$, so $q_n$ exists on
$I=[t_0,\min\{t_0+h_0,T_f\}]$ and satisfies \eqref{e:locen} there,
uniformly in $n$ and $t_0$.

By \cite[Corollary~4]{Simon1987} and interpolation, along a subsequence
$q_n\to q$ in $C(I;H^a(B))$ for $a<2$ and every
ball $B$, $q_n\rightharpoonup q$ weakly-$*$ in $L^\infty(I;H^2)$,
$(\psi_n,\zeta_n)\rightharpoonup(\psi,\zeta)$ in $L^2(I;H^3)$ and
$\partial_tq_n\rightharpoonup q_t$ in $L^2(I;H^1)$. With $a>\frac32$ the
coefficients converge locally uniformly, the products converge in
$L^1_{\rm loc}$, and the limit solves \eqref{e:locq}, hence
\eqref{eq:NSF}, with $q(t_0)=q^{t_0}$ and $\rho,\theta\ge\frac58\underline c_0$.
Since $(\psi,\zeta)\in L^2(I;H^3)$ and $(\psi_t,\zeta_t)\in L^2(I;H^1)$,
$(\psi,\zeta)\in C(I;H^2)$. By \eqref{e:phid}, $\Vert\partial^\alpha\phi\Vert_{L^2}$
is continuous for $|\alpha|\le2$; with the weak continuity of $\phi$ in
$H^2$ this gives $\phi\in C(I;H^2)$. Lemma~\ref{l:uniq} gives
$\phi_t\in C(I;H^1)$, and lower semicontinuity gives the bound in
\eqref{e:loc}.
\end{proof}

\begin{proof}[Proof of Lemma~\ref{l:cont}]
Choose $K_0$ above the supremum and $\underline c_0$ below the infimum in
\eqref{e:cont}, let $h_0$ be given by Lemma~\ref{l:local}(i), and put
$h=\min\{h_0,\frac{T_f-T_{\rm fin}}2\}$. For
$r\in(T_{\rm fin}-\frac h2,T_{\rm fin})$, Lemma~\ref{l:local}(i) with
$t_0=r$ and $q^{t_0}=q(r)$ gives a solution on $[r,r+h]$, and
$r+h>T_{\rm fin}+\frac h2$. By Lemma~\ref{l:uniq} it coincides with
$q$ on $[r,T_{\rm fin})$; hence $q(t)$ converges in $H^2$ as
$t\to T_{\rm fin}$, and the two solutions together give the extension
with $\ell=\frac h2$.
\end{proof}

\begin{proof}[Proof of Lemma~\ref{l:fp}]
Put $b_0=\inf_{\R\times[0,T]}\min\{\rho_b,\theta_b\}$, let $\eta\in(0,1]$
be such that $\Vert q\Vert_{H^2}\le\eta$ implies
$\Vert(\phi,\zeta)\Vert_{L^\infty}\le\frac{b_0}4$, and choose
$\epsilon_{b,T}$ with
\[
 C_*^2\epsilon_{b,T}^2\exp\Bigl(C_b\int_0^Tk_{\rm loc}\dd r\Bigr)
 \le\frac{\eta^2}4 .
\]
Apply Lemma~\ref{l:local}(i) with $K_0=\eta$ and
$\underline c_0=\frac34b_0$. While $\Vert q\Vert_{H^2}\le\eta$, we have
$\rho,\theta\ge\frac34b_0$, and Gronwall's inequality in \eqref{e:locsm} gives
\[
 E_{\rm loc}(q(t))\le E_{\rm loc}(q_0)
       \exp\Bigl(C_b\int_0^tk_{\rm loc}\dd r\Bigr)
 \le\frac{\eta^2}{4C_*},\qquad
 \Vert q(t)\Vert_{H^2}\le\frac\eta2 .
\]
This strict improvement, continuity in $H^2$, and the uniform restart time
in Lemma~\ref{l:local}(i) extend the solution throughout $[0,T]$;
Lemma~\ref{l:uniq} identifies the solutions on overlapping intervals.
In particular, $\rho\ge\frac34\inf\rho_b$ and $\theta\ge\frac34\inf\theta_b$. Integrating
\eqref{e:locsm} once more,
\[
 \int_0^T\bigl(\Vert(\psi,\zeta)\Vert_{H^3}^2+\Vert q_t\Vert_{H^1}^2\bigr)\dd r
 \le C_bE_{\rm loc}(q_0)\Bigl(1+\int_0^Tk_{\rm loc}\dd r\Bigr)
       \exp\Bigl(C_b\int_0^Tk_{\rm loc}\dd r\Bigr),
\]
and $\Vert\nabla(\psi,\zeta)\Vert_{H^2}\le \Vert(\psi,\zeta)\Vert_{H^3}$ gives \eqref{e:fp}. The
constants depend only on $T$, the physical parameters, the lower bounds
in \eqref{e:locb} and $\int_0^Tk_{\rm loc}\dd r$, hence on the norms in
\eqref{e:locb} on $[0,T]$.
\end{proof}

\subsection{Small energy with a prescribed higher-norm bound}
\label{app:low3}

We use small relative energy, as in the constant-state analysis
of Huang and Li~\cite{HuangLi2018}. The relative-entropy framework
for the full system is developed in~\cite{FeireislNovotny2012}; the
identity for our planar background is derived below. For positive
$\rho,\theta$, define
\begin{equation}\label{e:rel3}
\begin{aligned}
 \mathcal E_{\rm rel}(t)=\int_{\R^3}\Bigl[&\frac12\rho|\psi|^2
 +R\theta_b\left(\rho\log\frac\rho{\rho_b}-\rho+\rho_b\right)\\
 &+c_v\rho\left(\theta-\theta_b
                  -\theta_b\log\frac\theta{\theta_b}\right)\Bigr]\dd x.
\end{aligned}
\end{equation}

\begin{lemma}\label{l:low-finite}
Let the background solve \eqref{eq:bg-eqs} with \eqref{e:locb}, and
let $0<T\le T_f$ and $M_3\ge1$. There are
$\varepsilon_{b,T}(M_3)>0$ and $K_{b,T}(M_3)<\infty$ such that
\[
 q_0\in H^2(\R^3),\qquad
 \Vert q_0\Vert_{H^2}\le M_3,\qquad
 \Vert q_0\Vert_{L^2}\le\varepsilon_{b,T}(M_3)
\]
give a unique strong solution on $[0,T]$ in the class
\eqref{eq:class}, with
\[
 \rho\ge\frac12\rho_b,\qquad \theta\ge\frac12\theta_b,\qquad
 \sup_{0\le t\le T}\Vert q(t)\Vert_{H^2}\le K_{b,T}(M_3),
\]
and
\begin{equation}\label{e:low-finite}
 \sup_{0\le t\le T}\Vert q(t)\Vert_{L^2}^2
 +\int_0^T\Vert(\nabla\psi,\nabla\zeta)\Vert_{L^2}^2\dd t
 \le C_{b,T}\Vert q_0\Vert_{L^2}^2.
\end{equation}
Here $C_{b,T}$ is independent of $M_3$ and $q_0$; the permitted
$\varepsilon_{b,T}(M_3)$ may depend on $M_3$. No higher initial
regularity or compatibility condition is imposed. The same statement
holds on every fixed interval $[t_0,T]\subset[0,T_f]$, with datum
$q(t_0)$ and constants depending also on $t_0$.
\end{lemma}

\begin{proof}
It suffices to treat $t_0=0$; time translation gives the stated
fixed-interval version. Subtract the common velocity $u_*$ by a
Galilean transformation.
Put $m=\inf_{\R\times[0,T]}\min\{\rho_b,\theta_b\}>0$.
Constants $C_b,c_b>0$ below depend on the background on $[0,T]$,
its positive lower bounds and the physical coefficients, but not on
$M_3$ or a bound for $\Vert q\Vert_{H^2}$. We first work on an interval
where
\[
 \Vert(\phi,\zeta)\Vert_{L^\infty}\le\min\{m/2,1\}.
\]
The densities and temperatures then lie in a fixed positive interval,
and
\begin{equation}\label{e:rel3-coerc}
 c_b\Vert q\Vert_{L^2}^2\le\mathcal E_{\rm rel}
 \le C_b\Vert q\Vert_{L^2}^2.
\end{equation}

To obtain an estimate independent of the perturbation derivatives, put
\[
 f=\rho\log\frac\rho{\rho_b}-\rho+\rho_b,\qquad
 j=\theta-\theta_b-\theta_b\log\frac\theta{\theta_b},\qquad
 \ell=\log\frac\theta{\theta_b},\qquad
 \alpha=1-\frac{\theta_b}\theta.
\]
Using the continuity equations, test \eqref{eq:ph-psi} by $\psi$
and the temperature equation by $\alpha$. The pressure terms cancel
because $\Pi=R(\theta_b\phi+\rho\zeta)$. The background temperature
equation gives
\[
 \kappa\partial_1^2\theta_b
 =c_v\rho_bB_b-\nu h_b^2+R\rho_b\theta_bh_b.
\]
Consequently the resulting local identity, integrated in space, is
\begin{equation}\label{e:rel3-id}
\begin{aligned}
 \mathcal E_{\rm rel}'&+\int_{\R^3}\left[
   \frac{\theta_b}\theta Q_\mu(\nabla\psi)
       +\kappa\theta_b|\nabla\ell|^2\right]\dd x\\
 =\int_{\R^3}\Bigl[&-\rho h_b\psi_1^2-\phi A_b\psi_1
       +(Rf-c_v\phi\ell)B_b\\
 &+(Rf-c_v\rho\ell)\psi_1\partial_1\theta_b
       -R\theta_b\phi\psi_1\partial_1\log\rho_b\\
 &-Rh_b(\rho_bj+\phi\zeta)
       +\nu h_b^2(\alpha-\ell)+\alpha H[\psi]\Bigr]\dd x.
\end{aligned}
\end{equation}
The cancellations are performed before integration: the relative
conductive flux is
\[
 \kappa\left[\alpha\nabla\zeta
       +(\alpha-\ell)(\partial_1\theta_b)e_1\right]\in L^1(\R^3),
\]
since $|\alpha|+|\ell|\le C_b|\zeta|$ and
$|\alpha-\ell|\le C_b\zeta^2$. The relative pressure and viscous
fluxes are also in $L^1$, by $\Pi,\psi,\nabla\psi\in L^2$.
Thus spatial cutoffs and smoothing justify \eqref{e:rel3-id} in the
stated strong class, without assuming $q\in L^1$.

On the positive interval, $f\le C_b\phi^2$, $j\le C_b\zeta^2$,
and $\alpha-\ell=1-e^{-\ell}-\ell\le0$.
All terms on the right of \eqref{e:rel3-id}, except $\alpha H[\psi]$,
are bounded above by $C_bk_{\rm loc}\mathcal E_{\rm rel}$; in the term
$f\psi_1$, use $\phi^4\le C_b\phi^2$.
Moreover,
\[
 \left|\int\alpha H[\psi]\dd x\right|
 \le\eta\Vert\nabla\psi\Vert_{L^2}^2
       +C_{b,\eta}\Vert h_b\Vert_{L^\infty}^2\Vert\zeta\Vert_{L^2}^2.
\]
Here $Q_\mu\ge0$ and
$\int Q_\mu(\nabla\psi)
 =\mu\Vert\nabla\psi\Vert_{L^2}^2
   +(\mu+\lambda)\Vert\Div\psi\Vert_{L^2}^2$.
Also
\[
 \Vert\nabla\zeta\Vert_{L^2}^2
 \le C_b\Vert\nabla\ell\Vert_{L^2}^2
   +C_b\Vert\partial_1\theta_b\Vert_{L^\infty}^2
                    \Vert\zeta\Vert_{L^2}^2.
\]
Choosing $\eta$ small in \eqref{e:rel3-id} yields
\begin{equation}\label{e:rel3-low}
 \mathcal E_{\rm rel}'+c_b\Vert(\nabla\psi,\nabla\zeta)\Vert_{L^2}^2
 \le C_bk_{\rm loc}\mathcal E_{\rm rel}.
\end{equation}
Since $k_{\rm loc}\in L^1(0,T)$, this proves
\eqref{e:low-finite} as long as the positive interval is retained.

We next retain a possibly large $H^2$ bound. For this calculation write
\[
 X=\Vert q\Vert_{H^2},\qquad
 Y=\Vert(\psi,\zeta)\Vert_{H^3},\qquad
 a=\Vert q\Vert_{L^\infty}+\Vert\nabla q\Vert_{L^3}.
\]
The Gagliardo--Nirenberg inequalities give
\begin{equation}\label{e:rel3-gn}
 a\le C\Vert q\Vert_{L^2}^{1/4}X^{3/4},\qquad
 \Vert\nabla\psi\Vert_{L^\infty}\le Ca^{1/3}Y^{2/3},\qquad
 \Vert\nabla\psi\Vert_{L^{12}}^2\le CaY.
\end{equation}
Reconsider the tests defining $E_{\rm loc}$ above, assuming $a\le1$.
In the diffusion commutators, the perturbation part satisfies
\[
 \Vert\nabla\phi\Vert_{L^3}
 \bigl(\Vert\mathcal L\psi\Vert_{L^6}
               \Vert\nabla\mathcal L\psi\Vert_{L^2}
       +\Vert\Delta\zeta\Vert_{L^6}
               \Vert\nabla\Delta\zeta\Vert_{L^2}\bigr)
 \le CaY^2.
\]
The part containing $\partial_1\rho_b$ is bounded by
$\eta Y^2+C_{b,\eta}X^2$.
The forcing terms in \eqref{e:locq} satisfy the sharper bound
\begin{equation}\label{e:rel3-force}
 \Vert F_\psi\Vert_{H^1}+\Vert F_\zeta\Vert_{H^1}
 \le C_bk_{\rm loc}^{1/2}X+C_baY.
\end{equation}
Indeed, the pressure, transport and background factors are bounded
using \eqref{e:pmul} and $\Vert\nabla q\Vert_{L^3}\le a$.
For the only quadratic derivative source, use
\[
\begin{aligned}
 \Vert\nabla Q_\mu(\nabla\psi)\Vert_{L^2}
 &\le C\Vert\nabla\psi\Vert_{L^3}
                   \Vert\nabla^2\psi\Vert_{L^6}\le CaY,\\
 \Vert\nabla\phi\,Q_\mu(\nabla\psi)\Vert_{L^2}
 &\le C\Vert\nabla\phi\Vert_{L^3}
                   \Vert\nabla\psi\Vert_{L^{12}}^2\le Ca^2Y;
\end{aligned}
\]
the undifferentiated heating term is at most $CaX$.

In the density equation, the terms containing three derivatives of
$\psi$ are bounded by $C_bXY$. The two top perturbation commutators
have the bounds
\[
\begin{aligned}
 \int|\nabla^2\psi||\nabla\phi||\nabla^2\phi|\dd x
 &\le CaXY,\\
 \Vert\nabla\psi\Vert_{L^\infty}X^2
 &\le Ca^{1/3}Y^{2/3}X^2
 \le\eta Y^2+C_\eta a^{1/2}X^3.
\end{aligned}
\]
All background terms are bounded by
$C_bk_{\rm loc}^{1/2}X^2$; for
$(\partial_1^3\rho_b)\psi_1$ and $(\partial_1^2h_b)\phi$, use
\eqref{eq:planar}. Thus the same tests, with
\eqref{e:rel3-force}, give, for a sufficiently small fixed $a_b>0$,
\begin{equation}\label{e:rel3-high}
 E_{\rm loc}'+c_bY^2
 \le C_bk_{\rm loc}E_{\rm loc}+C_ba^{1/2}E_{\rm loc}^{3/2}
 \qquad(a\le a_b).
\end{equation}
This calculation uses no small bound for $X$.

Choose $a_b\le\min\{1,m/4\}$. If $X\le K$ and
$a\le\min\{a_b,c_bK^{-2}\}$, decreasing $c_b$ if necessary,
\eqref{e:rel3-high} implies
\[
 E_{\rm loc}'+c_bY^2\le C_bk_{\rm loc}E_{\rm loc}.
\]
Choose $K$ larger than four times the resulting Gronwall bound for
$X$ when $X(0)\le M_3$, and then choose
$\varepsilon_{b,T}(M_3)$ so small that
\[
 C_{b,T}\varepsilon_{b,T}(M_3)^{1/4}K^{3/4}
 \le\frac12\min\{a_b,c_bK^{-2}\}.
\]
The initial bounds hold strictly by \eqref{e:rel3-gn}; in particular
$\rho_0\ge3\rho_b(\cdot,0)/4$ and
$\theta_0\ge3\theta_b(\cdot,0)/4$, so Lemma~\ref{l:local} applies.
On their maximal interval, \eqref{e:rel3-coerc},
\eqref{e:rel3-low} and \eqref{e:rel3-gn} improve the bound for $a$
by a factor of two, while \eqref{e:rel3-high} improves $X\le K$
to $X\le K/2$. In particular $\rho\ge3\rho_b/4$ and
$\theta\ge3\theta_b/4$. The uniform restart time in
Lemma~\ref{l:local} and uniqueness therefore extend the solution
through $[0,T]$.
Integrating \eqref{e:rel3-high} gives $Y\in L^2(0,T)$;
\eqref{e:locq} and \eqref{e:rel3-force} give $q_t\in L^2(0,T;H^1)$.
The local $H^2$ continuity is preserved on overlapping intervals.
This proves the stated class, positivity and bounds.
\end{proof}

\begin{lemma}\label{l:thermal-layer}
Let the background satisfy \eqref{e:locb} on $[0,T_f]$, where
$T_f>1$. Fix $M_\rho,M_u,M_\theta\ge1$ and
$0<m_\theta\le\Theta_\theta<\infty$. There exist
$\varepsilon>0$, $K<\infty$, and $0<c_\theta<C_\theta<\infty$,
depending only on these parameters, the background on $[0,1]$, and
the physical coefficients, such that
\[
 \begin{gathered}
 q_0\in H^2(\R^3),\qquad
 \Vert\nabla^2\phi_0\Vert_{L^2}\le M_\rho,\qquad
 \Vert\nabla^2\psi_0\Vert_{L^2}\le M_u,\qquad
 \Vert\nabla\zeta_0\Vert_{L^2}\le M_\theta,\\
 m_\theta\le\theta_b(\cdot,0)+\zeta_0\le\Theta_\theta,
 \qquad \Vert q_0\Vert_{L^2}\le\varepsilon
 \end{gathered}
\]
give a unique strong solution on $[0,1]$ in \eqref{eq:class}, with
\begin{equation}\label{e:thermal-layer}
\begin{gathered}
 \rho\ge\tfrac12\rho_b,\qquad
 c_\theta\le\theta\le C_\theta,\\
 \sup_{0\le t\le1}
  \bigl(\Vert(\phi,\psi)(t)\Vert_{H^2}^2
              +\Vert\zeta(t)\Vert_{H^1}^2\bigr)
 +\int_0^1\bigl(\Vert\psi\Vert_{H^3}^2
              +\Vert\zeta\Vert_{H^2}^2\bigr)\dd t\le K^2,\\
 \sup_{0<t\le1}t\Vert\zeta(t)\Vert_{H^2}^2
       +\int_0^1t\Vert\zeta\Vert_{H^3}^2\dd t\le K^2.
\end{gathered}
\end{equation}
Moreover,
\[
 \sup_{0\le t\le1}\Vert q(t)\Vert_{L^2}^2
 +\int_0^1\Vert(\nabla\psi,\nabla\zeta)\Vert_{L^2}^2\dd t
 \le C\Vert q_0\Vert_{L^2}^2.
\]
The constants are independent of $\Vert\nabla^2\zeta_0\Vert_{L^2}$.
The last constant $C$ and the temperature bounds may be chosen
independent of $M_\rho,M_u,M_\theta$.
\end{lemma}

\begin{proof}
Work in the Galilean frame with $u_*=0$. We use \eqref{e:locq}
and the relative energy \eqref{e:rel3}.
Write, within this proof,
\[
 \begin{aligned}
 X^2&=\Vert(\phi,\psi)\Vert_{H^2}^2
                    +\Vert\zeta\Vert_{H^1}^2,&
 Y^2&=\Vert\psi\Vert_{H^3}^2+\Vert\zeta\Vert_{H^2}^2,\\
 a&=\Vert(\phi,\psi)\Vert_{L^\infty}
                      +\Vert\nabla(\phi,\psi)\Vert_{L^3},&
 h&=\sup_{0\le r\le t}\Vert q(r)\Vert_{L^2}.
 \end{aligned}
\]
Initially $X\le C(1+M_\rho+M_u+M_\theta)$ if
$\Vert q_0\Vert_{L^2}\le1$, by interpolation. We choose a fixed
density interval containing $[\frac12\inf\rho_b,2\sup\rho_b]$.
On any interval where $a$ is small, the density stays in this
interval. All constants below are uniform for density in this
interval and temperature in a fixed positive interval.

We specify the temperature interval before choosing the bound for
$X$. Set
\[
 B_j=\int_0^1\Vert h_b(t)\Vert_{L^\infty}^j\dd t\quad(j=1,2),
 \qquad
 c_\theta=\tfrac12m_\theta e^{-C(B_1+1)},\qquad
 C_\theta=2\{\Theta_\theta+C(B_2+1)\}e^{C(B_1+1)},
\]
where $C$ is large enough depending only on the density interval and
the physical coefficients. Suppose temporarily that
\[
 \int_0^t\Vert\nabla\psi\Vert_{L^\infty}\dd r\le1,
 \qquad
 \int_0^t\Vert\nabla\psi\Vert_{L^\infty}^2\dd r\le1.
\]
Since $Q_\mu(\nabla u)\ge0$, the full temperature equation and
parabolic comparison give
\[
 \begin{aligned}
 \inf_x\theta(x,t)
 &\ge m_\theta\exp\left(-C\int_0^t
            \{\Vert h_b\Vert_{L^\infty}
                    +\Vert\nabla\psi\Vert_{L^\infty}\}\dd r\right)
       \ge2c_\theta,\\
 \Vert\theta(t)\Vert_{L^\infty}
 &\le\left[\Theta_\theta+C\int_0^t
            \{\Vert h_b\Vert_{L^\infty}^2
                    +\Vert\nabla\psi\Vert_{L^\infty}^2\}\dd r\right]
       \exp\left(C\int_0^t
            \{\Vert h_b\Vert_{L^\infty}
                    +\Vert\nabla\psi\Vert_{L^\infty}\}\dd r\right)
       \le\tfrac12 C_\theta.
 \end{aligned}
\]
These inequalities follow by comparison with scalar time-dependent
sub- and supersolutions; spatial cutoffs and smoothing justify
comparison in the local strong class. They require only
$Q_\mu\ge0$, and do not require $2\mu+3\lambda>0$.

On this state interval, the proof of \eqref{e:rel3-low} applies
unchanged, with constants depending on $c_\theta,C_\theta$:
\[
 \mathcal E_{\rm rel}\asymp\Vert q\Vert_{L^2}^2,\qquad
 \mathcal E_{\rm rel}'+c\Vert(\nabla\psi,\nabla\zeta)\Vert_{L^2}^2
       \le Ck_{\rm loc}\mathcal E_{\rm rel}.
\]
In particular $h\le C\Vert q_0\Vert_{L^2}$. No perturbation derivative
or initial temperature Hessian occurs in $C$.

We next estimate the reduced energy. The temperature equation is
tested only with $-\Delta\zeta$. Its principal term is exactly
$\int\kappa(c_v\rho)^{-1}|\Delta\zeta|^2$, with no derivative on
$\rho^{-1}$. The following source estimates hold when $a\le1$:
\begin{align*}
 \Vert F_\zeta\Vert_{L^2}
  &\le Ck_{\rm loc}^{1/2}X,\\
 \Vert F_\psi\Vert_{H^1}
  &\le Ck_{\rm loc}^{1/2}X+C\Vert\zeta\Vert_{H^2}.
\end{align*}
For clarity, all source classes can be checked as follows. In
$F_\zeta$, the terms with $\partial_1\theta_b,B_b,h_b$ and
$H[\psi]$ are bounded by $Ck_{\rm loc}^{1/2}X$, using
\eqref{e:pmul}, the fixed state interval and $a\le1$.
The terms $\theta\Div\psi$ and $u\cdot\nabla\zeta$ have the same
bound. For the remaining quadratic source,
\[
 \Vert Q_\mu(\nabla\psi)\Vert_{L^2}
 \le C\Vert\nabla\psi\Vert_{L^4}^2
 \le C\Vert\psi\Vert_{L^\infty}
                          \Vert\nabla^2\psi\Vert_{L^2}
 \le CaX.
\]
For the momentum pressure term, the useful exact expansion is
\[
 \frac{\nabla\Pi}{\rho}
 =R\frac\theta\rho\nabla\phi+R\nabla\zeta
       +R\frac{\phi\partial_1\theta_b
                     +\zeta\partial_1\rho_b}{\rho}e_1.
\]
The last term is bounded in $H^1$ by $Ck_{\rm loc}^{1/2}X$;
terms such as $\zeta\nabla\phi$ use
$\Vert\zeta\Vert_{L^6}\Vert\nabla\phi\Vert_{L^3}\le CaX$.
Differentiating the first term gives bounded multiples of
$\nabla^2\phi$, background terms, $\nabla\phi\otimes\nabla\phi$,
and $\nabla\phi\otimes\nabla\zeta$. Here
\[
 \Vert\nabla\phi\otimes\nabla\phi\Vert_{L^2}
   \le Ca\Vert\phi\Vert_{H^2},\qquad
 \Vert\nabla\phi\otimes\nabla\zeta\Vert_{L^2}
   \le\Vert\nabla\phi\Vert_{L^3}
                     \Vert\nabla\zeta\Vert_{L^6}
   \le Ca\Vert\zeta\Vert_{H^2}.
\]
The momentum background terms $f_m/\rho$ are bounded by
$Ck_{\rm loc}^{1/2}X$; the transport term uses
$\Vert\nabla\psi\Vert_{L^3}\Vert\nabla\psi\Vert_{L^6}\le CaX$.
This proves both source estimates without a temperature $H^2$
bound at the initial time.

Test the momentum equation with $\mathcal L(\mathcal L\psi)$ and
the differentiated density equations as in
\eqref{e:rel3-high}. The diffusion commutator with $\nabla\phi$ is
bounded by $Ca\Vert\psi\Vert_{H^3}^2$; its background part is
bounded by $\eta\Vert\psi\Vert_{H^3}^2+C_\eta X^2$.
The density tests have exactly the previous bounds
\[
 CXY+CaXY
   +C\Vert\nabla\psi\Vert_{L^\infty}X^2
   +Ck_{\rm loc}^{1/2}X^2,
 \qquad
 \Vert\nabla\psi\Vert_{L^\infty}
      \le Ca^{1/3}\Vert\psi\Vert_{H^3}^{2/3}.
\]
Consequently their last perturbation contribution is at most
$\eta Y^2+C_\eta a^{1/2}X^3$.
The new pressure product is at most $CaY^2$ and is absorbed for
small $a$.
Choose a sufficiently large fixed weight $L_\theta$ for the temperature
$H^1$ test. This absorbs the linear
$C\Vert\zeta\Vert_{H^2}\Vert\psi\Vert_{H^3}$ pressure pairing.
Adding a sufficiently large fixed multiple $L_0$ of the relative energy gives
\[
 E=\Vert\phi\Vert_{H^2}^2+\Vert\mathcal L\psi\Vert_{L^2}^2
              +L_\theta\Vert\nabla\zeta\Vert_{L^2}^2
              +L_0\mathcal E_{\rm rel}\asymp X^2
\]
and, for sufficiently small fixed $a_*>0$,
\[
 E'+cY^2
 \le Ck_{\rm loc}E+C a^{1/2}E^{3/2}
 \qquad(a\le a_*).
\]
The lower-order parts of $Y$ are included by adding $CE$ to the
right, since $k_{\rm loc}\ge1$.

Fix $K$ greater than four times the resulting Gronwall bound for
$X$ and $(\int Y^2)^{1/2}$ with the stated bound for $X(0)$.
If $X\le K$, then
\[
 a\le Ch^{1/4}K^{3/4}.
\]
Choose $h$ sufficiently small, after $K$, so that the last
term of the reduced energy inequality is bounded by $E$ and
$a\le a_*$. This yields the improved bound $X\le K/2$ and
$\int_0^tY^2\dd r\le K^2/4$, after enlarging $K$ once by a fixed
factor. Moreover,
\[
 \begin{aligned}
 \int_0^t\Vert\nabla\psi\Vert_{L^\infty}\dd r
 &\le Ch^{1/6}\Bigl(\int_0^t\Vert\psi\Vert_{H^3}^2\dd r\Bigr)^{5/12},\\
 \int_0^t\Vert\nabla\psi\Vert_{L^\infty}^2\dd r
 &\le Ch^{1/3}\Bigl(\int_0^t\Vert\psi\Vert_{H^3}^2\dd r\Bigr)^{5/6},
 \end{aligned}
\]
by $t\le1$ and
$\Vert\nabla\psi\Vert_{L^\infty}
\le C\Vert\psi\Vert_{L^2}^{1/6}
                   \Vert\nabla^3\psi\Vert_{L^2}^{5/6}$.
Decreasing $h$ further makes both integrals less than $1/2$.
The maximum-principle bounds therefore improve the fixed
temperature interval, and $a$ improves the density interval.
As $h\le C\Vert q_0\Vert_{L^2}$, all these conditions follow from
a single choice of $\varepsilon$ with the stated dependencies.

It remains to ensure that a large but finite temperature Hessian
cannot shorten the existence interval. We claim the source bound
\[
 \Vert F_\zeta\Vert_{H^1}
 \le C\{k_{\rm loc}^{1/2}X
                 +\Vert\zeta\Vert_{H^2}
                 +a\Vert\psi\Vert_{H^3}\}.
\]
Indeed, differentiating $\theta\Div\psi$ or
$u\cdot\nabla\zeta$ produces the new product
$\nabla\psi\,\nabla\zeta$, bounded by
$\Vert\nabla\psi\Vert_{L^3}\Vert\nabla\zeta\Vert_{L^6}
\le Ca\Vert\zeta\Vert_{H^2}$; the other differentiated terms are
bounded by the first two terms on the right. The background
sources are bounded with \eqref{e:pmul}; in particular only
$B_b\in H^1(\R)$ is used. The differentiated heating terms satisfy
\[
 \begin{aligned}
 \Vert\nabla Q_\mu(\nabla\psi)\Vert_{L^2}
 &\le C\Vert\nabla\psi\Vert_{L^3}
                    \Vert\nabla^2\psi\Vert_{L^6}
       \le Ca\Vert\psi\Vert_{H^3},\\
 \Vert\nabla\phi\,Q_\mu(\nabla\psi)\Vert_{L^2}
 &\le C\Vert\nabla\phi\Vert_{L^3}
                    \Vert\nabla\psi\Vert_{L^{12}}^2
       \le Ca^2\Vert\psi\Vert_{H^3}.
 \end{aligned}
\]
Here $\Vert\nabla\psi\Vert_{L^{12}}^2
\le C\Vert\psi\Vert_{L^\infty}\Vert\nabla^3\psi\Vert_{L^2}$.
The remaining density coefficient derivative in the heating term
is a bounded background factor times $Q_\mu$, already estimated.

Test the temperature equation with $\Delta^2\zeta$. The background
part of the diffusion commutator is bounded by
$\eta\Vert\nabla\Delta\zeta\Vert_{L^2}^2
+C_\eta\Vert\Delta\zeta\Vert_{L^2}^2$, while its perturbation part
is at most
\[
 C\Vert\nabla\phi\Vert_{L^3}
     \Vert\Delta\zeta\Vert_{L^6}
     \Vert\nabla\Delta\zeta\Vert_{L^2}
 \le Ca\Vert\nabla\Delta\zeta\Vert_{L^2}^2.
\]
The preceding $H^1$ source estimate therefore gives
\[
 \frac{\dd}{\dd t}\Vert\Delta\zeta\Vert_{L^2}^2
       +c\Vert\nabla\Delta\zeta\Vert_{L^2}^2
 \le C\{k_{\rm loc}X^2+\Vert\zeta\Vert_{H^2}^2
                           +a^2\Vert\psi\Vert_{H^3}^2\}.
\]
Let $T_*$ be the maximal strong existence time. For
$0<T<\min\{1,T_*\}$, multiply by $t$ and integrate over $[0,T]$. Since
$t\Vert\Delta\zeta(t)\Vert_{L^2}^2\to0$ as $t\downarrow0$
for each of the local $H^2$ solutions, and every term on the right,
including the extra $\Vert\Delta\zeta\Vert_{L^2}^2$, is integrable
by the reduced energy estimate, this proves
\[
 \sup_{0<t\le T}t\Vert\zeta(t)\Vert_{H^2}^2
       +\int_0^T t\Vert\zeta\Vert_{H^3}^2\dd t\le C K^2,
\]
with $C$ depending only on the already fixed parameters. This bound is
independent of the initial temperature Hessian.

Lemma~\ref{l:local}(i) applies initially because $q_0\in H^2$,
the temperature has the stated positive lower bound, and small
$\Vert q_0\Vert_{L^2}$ gives $\rho_0\ge\frac34\rho_b(\cdot,0)$.
If $T_*\le1$, the weighted estimate would give a uniform full
$H^2$ bound on
$[T_*/2,T_*)$. On the earlier compact interval the local $H^2$
continuity supplies a finite bound, possibly depending on the
finite initial Hessian. The state bounds are uniform and positive.
Lemma~\ref{l:local}(ii) then extends the solution beyond $T_*$,
a contradiction. Thus the solution exists on $[0,1]$.
The same scalar second-derivative inequality, integrated without
the time weight, also gives
\[
 \sup_{0\le t\le1}\Vert\zeta(t)\Vert_{H^2}^2
       +\int_0^1\Vert\zeta\Vert_{H^3}^2\dd t
 \le C\bigl(\Vert\zeta_0\Vert_{H^2}^2+K^2\bigr).
\]
Together with \eqref{e:locq}, this preserves the unweighted class
\eqref{eq:class}; this regularity bound may depend on the finite
initial temperature Hessian and is not used to choose $\varepsilon$.
Enlarging $K$ gives both uniform energy bounds; uniqueness follows from
Lemma~\ref{l:uniq}.
\end{proof}

\begin{lemma}\label{l:low-damp}
There are $S_0\ge1$, $c>0$ and $C_b$ such that, for every $K\ge1$,
there are $h_b(K)>0$ and $C_{b,K}$ with the following property.
Let $S_0\le s<T$ and let a positive strong solution satisfy
\[
 \sup_{s\le r\le T}\Vert q(r)\Vert_{H^2}\le K,
 \qquad \sup_{s\le r\le T}\Vert q(r)\Vert\le h\le h_b(K).
\]
Then, for $s\le t\le T$,
\begin{equation}\label{e:low-damp}
 \Vert q(t)\Vert_{H^2}^2
 \le C_be^{-c(t-s)}\Vert q(s)\Vert_{H^2}^2+C_{b,K}h^{2/7}.
\end{equation}
\end{lemma}
\begin{proof}
As in the proof of Lemma~\ref{l:damp}, choose the wave-strength
threshold and $S_0$ so that, for $r\ge S_0$,
\[
 \Vert(\rho_b-\rho_-,u_b-u_*,\theta_b-\theta_-)\Vert_{L^\infty}
 +\Vert\partial_1(\rho_b,u_b,\theta_b)\Vert_{H^2(\R)}
 +\Vert(A_b,B_b)\Vert_{H^1(\R)}\le\beta,
\]
where $\beta>0$ is a fixed small number.
Interpolation gives
\begin{equation}\label{e:low-coeff}
 \Vert q\Vert_{L^\infty}+\Vert\nabla q\Vert_{L^3}
 \le Ch^{1/4}K^{3/4}=:a.
\end{equation}
We choose $h_b(K)\le1$ so that $a$ is small. The density and temperature
then lie in a fixed positive state box.

Take $P_0$ from the proof of Lemma~\ref{l:damp} and put
$q^h=(1-P_0)q$. Use its quadratic energy and interactive term with
the spatial weight replaced by $1$:
\[
\begin{aligned}
 \mathcal H={}&\frac12\sum_{k=0}^2\binom2k\int
 \left(\frac{R\theta_-}{\rho_-}|\nabla^k\phi^h|^2
       +\rho_-|\nabla^k\psi^h|^2
       +\frac{c_v\rho_-}{\theta_-}|\nabla^k\zeta^h|^2\right)\dd x\\
 &+\epsilon_c\int\left(\psi^h\cdot\nabla\phi^h
             +\sum_i\partial_i\psi^h\cdot\nabla\partial_i\phi^h\right)\dd x,\\
 \mathcal D={}&\Vert\nabla\phi^h\Vert_{H^1}^2
              +\Vert\nabla(\psi^h,\zeta^h)\Vert_{H^2}^2.
\end{aligned}
\]
For a sufficiently small fixed $\epsilon_c>0$, the spectral gap gives
\[
 C^{-1}\Vert q^h\Vert_{H^2}^2\le\mathcal H
 \le C\Vert q^h\Vert_{H^2}^2\le C\mathcal D,
 \qquad
 \Vert P_0q\Vert_{H^k}+\Vert\nabla^kP_0q\Vert_{L^\infty}\le C_kh.
\]

We repeat the tests in the proof of Lemma~\ref{l:damp}, using
\eqref{e:low-coeff} for the coefficient and parabolic commutator terms.
The top density terms are estimated directly, without assuming small
$H^2_{\rm uloc}$ norm. Indeed,
\[
 \Vert\nabla\psi\Vert_{L^\infty}
       +\Vert\nabla^2\psi\Vert_{L^3}
 \le C\Vert\psi\Vert^{1/6}\Vert\nabla^3\psi\Vert^{5/6},
\]
so that, for every $\eta>0$,
\begin{equation}\label{e:low-density-comm}
\begin{aligned}
 &\int\left(|\nabla\psi||\nabla^2\phi^h|^2
     +|\nabla^2\psi||\nabla\phi^h||\nabla^2\phi^h|\right)\dd x\\
 &\qquad\le Ch^{1/6}K^2\Vert\nabla^3\psi\Vert^{5/6}
 \le\eta\Vert\nabla^3\psi\Vert^2
             +C_\eta h^{2/7}K^{24/7}.
\end{aligned}
\end{equation}
All remaining products in \eqref{e:nphys} have small coefficients
$C_b(\beta+a)$ in the required $H^1$ source norms. For example,
\[
\begin{aligned}
 \Vert\nabla\phi\,\nabla^2(\psi,\zeta)\Vert
 &\le\Vert\nabla\phi\Vert_{L^3}
             \Vert\nabla^2(\psi,\zeta)\Vert_{L^6},\\
 \Vert\nabla\phi\,|\nabla\psi|^2\Vert
 &\le\Vert\nabla\phi\Vert_{L^3}\Vert\nabla\psi\Vert_{L^{12}}^2
 \le C\Vert\nabla\phi\Vert_{L^3}\Vert\psi\Vert_{L^\infty}
                  \Vert\nabla^3\psi\Vert .
\end{aligned}
\]
Products of two first derivatives are bounded by the first inequality
of \eqref{e:gn4}. The background terms use
$\Vert b(x_1)f\Vert\le C\Vert b\Vert_{L^2(\R)}\Vert f\Vert_{H^1}$;
this includes $(\partial_1^3\rho_b)\psi_1$.
The terms containing $P_0$ are smoothing and satisfy the same bounds;
in the density row put $w_*=(u_b-u_*)e_1+\psi$ and use
\[
 (1-P_0)(w_*\cdot\nabla\phi)
 =w_*\cdot\nabla\phi^h-P_0(w_*\cdot\nabla\phi)
                       +w_*\cdot\nabla P_0\phi.
\]
Thus the constant-coefficient cancellations and the differentiated
interactive term give
\[
 \mathcal H'+c_1\mathcal D
 \le C_b(\beta+a)\mathcal D+\eta\mathcal D
                         +C_{b,K,\eta}h^{2/7}.
\]
Here the low-frequency parts of $\nabla^3(\psi,\zeta)$ contribute
only $C_{b,K}h^2$. Choose $\eta$, then $\beta$, and finally $h_b(K)$
to absorb the first two terms. The spectral gap yields
\[
 \mathcal H'+c\mathcal H+c\mathcal D\le C_{b,K}h^{2/7}.
\]
Integration and $\Vert P_0q\Vert_{H^2}^2\le Ch^2\le Ch^{2/7}$
prove \eqref{e:low-damp}.
\end{proof}

\begin{lemma}\label{l:late-restart}
Fix a background for which Lemmas~\ref{l:prs}, \ref{l:cl},
\ref{l:age}, \ref{l:damp} and~\ref{l:adj} hold. Assume
$4C_2C_P\delta^4\le1/2048$, where $C_2$ and $C_P$ are the constants
in Lemmas~\ref{l:cl} and~\ref{l:prs}. There are $S_0\ge1$, $c_b>0$ and a finite $m_b>0$ such that, for every
$s\ge S_0$, a datum $q^s\in H^2$ with
\[
 \Vert q^s\Vert_{H^2}\le c_b(1+s)^{-m_b}
\]
generates a unique global positive strong solution on $[s,\infty)$.
For $t\ge s$, it satisfies
\[
 \Vert(q,\Pi)(t)\Vert_{L^\infty}
 \le C_{b,s}\Vert q^s\Vert_{H^2}(1+t)^{-1/20}.
\]
\end{lemma}
\begin{proof}
Restart the proof of Lemma~\ref{l:dec} with $S_C=S=s$.
All background entrance conditions hold once $s\ge S_0$.
The local estimates on $[s,s+1]$ have constants bounded uniformly in
$s$, since the late-time background state, coefficient norms and
unit-interval integrals are uniformly bounded. No estimate on $[0,s]$
is used.

The remaining dependence on $s$ is polynomial. The constant $C_{0,b}$
in Lemma~\ref{l:cl} is uniform in the starting time, and the pressure
budget \eqref{e:prlog}, with the stated restriction on $\delta$, gives
\[
 C_{P,b}(s)\le C_b+\eta_P\log(1+s),\qquad
 C_G(s)\le C_b(1+s)^{\eta_P},\qquad \eta_P=1/1024.
\]
In Lemma~\ref{l:age}, the coefficient-tail integrals are uniformly
bounded and the starting-time factor is
$\omega(s)=(1+s)^{1/4}$. Its constant is therefore at most
$C_b(1+s)^{1/4}$. The adjoint kernel in \eqref{e:adj} is uniform
in $s$. On $s\le t<2s+2$, its time $L^2$ norm is uniformly bounded;
the source estimates give at worst $(1+s)^{q_a}$, and multiplication
by $(1+t)^{d_*}$ gives a fixed further power. Thus the constant in
\eqref{e:thout} is at most $C_b(1+s)^{m}$ for a fixed finite $m$.
Finally, the explicit exponential-convolution formula in
Lemma~\ref{l:damp} gives at most the factor $(1+s)^d$.

Consequently the same reconstruction as in Lemma~\ref{l:dec} yields
\[
 \mu_*\le C_b(1+s)^{m_b}\Vert q^s\Vert_{H^2}
\]
for some fixed finite $m_b$, on every bootstrap interval.
Choose $c_b$ so that this is at most $\mu_0/2$ and the uniform
unit-interval startup conditions hold. The continuation and decay
argument of that lemma then applies.
\end{proof}

\section*{Acknowledgments}
The authors acknowledge the use of generative AI tools in developing and
checking mathematical arguments and in drafting and revising the manuscript.

\begingroup
\small
\urlstyle{same}

\endgroup
\end{document}